\def\real{{\tt I\kern-.2em{R}}}
\def\nat{{\tt I\kern-.2em{N}}}
\def\snat{{\rm I\kern-.2em{N}}}
\def\eps{\epsilon}
\def\realp#1{{\tt I\kern-.2em{R}}^#1}
\def\natp#1{{\tt I\kern-.2em{N}}^#1}
\def\hyper#1{\,^*\kern-.2em{#1}}
\def\hy#1{\,^*\kern-.2em{#1}}
\def\monad#1{\mu (#1)}

\def\St#1{{\tt st}#1}
\def\st#1{{\tt st}(#1)}
\def\hyperreal{{^*{\real}}}
\def\hyperrealp#1{{\tt ^*{I\kern-.2em{R}}}^#1} 
\def\hypernat{{^*{\nat }}}
\def\hypernatp#1{{{^*{{\tt I\kern-.2em{N}}}}}^#1} 
\def\eskip{\hskip.25em\relax}

\def\Hyper#1{\hyper {\eskip #1}}

\def\leaderfill{\leaders\hbox to 1em{\hss.\hss}\hfill}
\def\srealp#1{{\rm I\kern-.2em{R}}^#1}

\def\power#1{{{\cal P}(#1)}}
\def\iff{\leftrightarrow}
\def\qed{{\vrule height6pt width3pt depth2pt}\par\medskip}
\def\pars{\par\smallskip}
\def\parm{\par\medskip}
\def\r#1{{\rm #1}}
\def\b#1{{\bf #1}}
\def\ref#1{$^{#1}$}

\def\m@th{\mathsurround=0pt}
\def\rightarrowfill{$\m@th \mathord- \mkern-6mu \cleaders\hbox{$\mkern-2mu 
\mathord- \mkern-2mu$}\hfil \mkern-6mu \mathord\rightarrow$}
\def\leftarrowfill{$\mathord\leftarrow
\mkern -6mu \m@th \mathord- \mkern-6mu \cleaders\hbox{$\mkern-2mu 
\mathord- \mkern-2mu$}\hfil $}
\def\noarrowfill{$\m@th \mathord- \mkern-6mu \cleaders\hbox{$\mkern-2mu 
\mathord- \mkern-2mu$}\hfil$}
\def\orgate{$\bigcirc \kern-.80em \lor$}
\def\andgate{$\bigcirc \kern-.80em \land$}
\def\inverter{$\bigcirc \kern-.80em \neg$}
\def\vd{\vdash}
\def\mod{\models}
\def\l{\land}
\def\lint{\lceil}
\def\rint{\rceil}
\def\ass{\underline{a}}
\def\em#1{\hskip #1em\relax}
\def\hs{\hrule\smallskip}
\def\hm{\hrule\medskip}
\def\s{\smallskip}
\def\c{\cal}
\def\p{^{\prime}}
\def\str{{\c M}_I}
\def\vd{\vdash}
\def\mod{\models}
\def\ol{\overline{\Gamma}}
\def\ss{\s\s\s\hrule\s\hrule\s\hrule\s\s\s}
\settabs 2 \columns 
\magnification=1000
\tolerance 10000
\hsize 6.50 true in
\vsize 9.00 true in
\baselineskip=12pt

\font\bigbf=cmb10 scaled 2000
\font\bigcr=cmb10 scaled 1600
\centerline{}
\vskip 1.75in
$$\vbox{\offinterlineskip
\hrule
\halign{&\vrule#&
\strut\quad\hfil#\hfil\quad\cr
height2pt&\omit&\cr
&\quad&\cr
&{{\bigbf LOGIC FOR EVERYONE}}&\cr
&\quad&\cr
height2pt&\omit&\cr}
\hrule}$$
\vskip 0.5in
\centerline{{\bigcr Robert A. Herrmann}}
\vskip 3in
\vfil\eject
{\quad}
\vskip 4.0in
\centerline{Previous titled ``Logic For Midshipmen''}\bigskip
\centerline{Mathematics Department}
\centerline{U. S. Naval Academy}
\centerline{572C Holloway Rd.}
\centerline{Annapolis$,$ MD 21402-5002}
{\quad}
\medskip\vfil
\eject

\centerline{}
\centerline{\bf CONTENTS}
\bigskip
\indent Chapter 1\pars
{\bf Introduction}\pars
\line{\indent 1.1 \ \ Introduction\leaderfill 5}
\bigskip
\indent Chapter 2\pars
{\bf The Propositional Calculus}\pars
\line{\indent 2.1 \ \ Constructing a Language by Computer\leaderfill 9}
\line{\indent 2.2 \ \ The Propositional Language\leaderfill 9}
\line{\indent 2.3 \ \ Slight Simplification$,$ Size$,$ Common Pairs\leaderfill 13}
\line{\indent 2.4 \ \ Model Theory --- Basic Semantics\leaderfill 15}
\line{\indent 2.5 \ \ Valid Formula\leaderfill 19}
\line{\indent 2.6 \ \ Equivalent Formula\leaderfill 22}
\line{\indent 2.7 \ \ The Denial$,$ Normal Form$,$ Logic Circuits\leaderfill 26}
\line{\indent 2.8 \ \ The Princeton Project$,$ Valid Consequences\leaderfill 32}
\line{\indent 2.9 \ \ Valid Consequences\leaderfill 35}
\line{\indent 2.10 \ \ Satisfaction and Consistency \leaderfill 38}
\line{\indent 2.11 \ \ Proof Theory \leaderfill 42}
\line{\indent 2.12 \ \ Demonstrations$,$ Deduction from Premises\leaderfill 45}
\line{\indent 2.13 \ \ The Deduction Theorem\leaderfill 47}
\line{\indent 2.14 \ \ Deducibility Relations\leaderfill 50}
\line{\indent 2.15 \ \ The Completeness Theorem\leaderfill 52}
\line{\indent 2.16 \ \ Consequence Operators\leaderfill 54}
\line{\indent 2.17 \ \ The Compactness Theorem\leaderfill 57}
\bigskip                                                      
\indent Chapter 3\pars
{\bf Predicate Calculus}\pars
\line{\indent 3.1 \ \ First-Order Language\leaderfill 63}
\line{\indent 3.2 \ \ Free and Bound Variable Occurrences\leaderfill 67} 
\line{\indent 3.3 \ \ Structures\leaderfill 70}
\line{\indent 3.4 \ \ Valid Formula in $Pd.$\leaderfill 76}
\line{\indent 3.5 \ \ Valid Consequences and Models\leaderfill 82}
\line{\indent 3.6 \ \ Formal Proof Theory\leaderfill 86}
\line{\indent 3.7 \ \ Soundness and Deduction Theorem for $Pd^\prime$\leaderfill 87}
\line{\indent 3.8 \ \ Consistency$,$ Negation Completeness$,$ Compactness$,$ Infinitesimals \leaderfill 91}
\line{\indent 3.9 \ \ Ultralogics and Natural Systems\leaderfill 97}
\bigskip
\indent {\bf Appendix}\pars
\line{\indent\indent \ \ Chapter 2\leaderfill 102}
\line{\indent\indent \ \ Chapter 3 \leaderfill 105}
\bigskip
\line{\indent{\bf Answers to Some Exercises}\leaderfill 109}
\bigskip
\line{\indent{\bf Index}\leaderfill 123}
\vfil\eject
{\quad}\vfil\eject
\centerline{\bf Chapter 1 - INTRODUCTION}\par
\medskip
\noindent {\bf 1.1 Introduction.}\pars
The discipline known as {{\it Mathematical Logic}} will not specifically be 
defined within this text. Instead$,$ you will study some of the concepts in this 
significant discipline by actually doing mathematical logic. Thus$,$ you will
be  able to surmise for yourself what the mathematical logician is attempting 
to accomplish. \pars
Consider the following three arguments taken from the disciplines of military 
science$,$ biology$,$ and set-theory$,$ where the symbols (a)$,$ (b)$,$ (c)$,$ (d)$,$ (e) are 
used only to locate specific sentences. \parm
{\leftskip 0.25in \rightskip 0.25 in (1) (a) If armored vehicles are used$,$ then 
the battle will be won. (b) If the infantry walks to the battle field$,$ then the 
enemy is warned of our presence. (c) If the enemy is warned of our presence and 
armored vehicles are used$,$ then we  sustain many casualties. (d) If the battle 
is won and we sustain many casualties$,$ then we will not be able to advance to 
the next objective. (e) Consequently$,$ if the infantry walks to the battle field 
and armored vehicles are used$,$ then we will not be able to advance to the next 
objective.\pars                          
(2) (a) If bacteria grow in culture A$,$ then the bacteria growth is normal. 
(b) If an 
antibiotic is added to culture A$,$ then mutations are formed. (c) If mutations are 
formed and bacteria grow in culture A$,$ then the growth medium is enriched. 
(d) If 
the bacteria growth is normal and the growth medium is enriched$,$ then there is 
an increase in the growth rate. (e) Thus$,$ if an antibiotic is added to culture 
A and bacteria grow in culture A$,$ then there is an increase in the growth 
rate.\pars
(3) (a) If $b \in B$$,$ then $(a,b) \in A \times B$. (b) If $c \in C,$ 
then $s \in S$. (c) If $s \in S$ and $b \in B$$,$ then $a \in A$. (d) 
If   $(a,b) \in A 
\times B$        
and $a \in A,$ then $(a,b,c,s) \in A \times B\times C\times S$. (e) Therefore$,$ 
if $c \in C$ and $b \in B$$,$ then  $(a,b,c,s) \in A \times B\times C\times 
S$.
\par}\parm
 With respect to the three cases above$,$ the statements that 
appear before the words ``Consequently$,$ Thus$,$ Therefore'' need not be assumed 
to be  ``true in reality.''  The actual logical pattern being presented is 
not$,$ as yet$,$ relative to the concept of what is ``true in reality.'' How can 
we analyze the logic behind each of these arguments? First$,$ notice that each 
of the above arguments employs a technical language peculiar to the specific 
subject under discussion.  This technical language should not affect the logic 
of each argument. The logic is something  ``pure'' in character which should 
be independent of such phrases as $a \in A.$ Consequently$,$ we could substitute 
abstract symbols -- symbols that carry no meaning and have no internal 
structure -- for each of the phrases such as the one ``we will not be able to 
advance to the next objective.'' Let us utilize the symbols P$,$ Q$,$ R$,$ S$,$ T$,$ H 
as replacements for these phrases with their technical terms. \pars
Let P = {\it armored vehicles are used}$,$ Q =  {\it the battle will be won}$,$ R =
{\it the infantry walks to the battle field}$,$ S =  {\it the 
enemy is warned of our presence}$,$ H =  {\it we  sustain many casualties}$,$ 
T= {\it we will not be able to advance to the next objective}. Now the words
{\it Consequently$,$ Thus$,$ Therefore} are replaced by the symbol $\vd,$  where 
the $\vd$ represents the processes the human mind (brain) goes through to 
``logically arrive at the statement'' that follows these words.\pars
Mathematics$,$ in its most fundamental form$,$ is based upon human experience and 
what we do next is related totally to such an experience. You must 
intuitively know your left from your right$,$ you must intuitively know what is 
means to  ``move from the left to the right$,$'' you must know what it means 
to ``substitute'' one thing for another$,$ and you must intuitively know one 
alphabet letter from another although different individuals may write them 
in slightly different forms. Thus P is the same as {\it P}$,$ etc. 
Now each of the above sentences contains the words  {\it If} and {\it then.} 
These two words are not used when we analyze the above three logical arguments
they will intuitively be understood. They will be part of the  symbol $\to.$ 
Any time you have a statement such as ``If P$,$ then Q'' this will be symbolized
as  $P \to  Q.$ There is one other important word in these statements. 
This word is  {\it and}. We symbolize this word {\it and} by the symbol $\land.$
What do these three arguments look like when we translate them into these 
defined 
symbols? Well$,$ in the next display$,$ I've used the  ``comma'' to separated the 
sentences and parentheses to remove any possible misunderstandings that might 
occur. When the substitutions are made in argument (1) and we write the sentences 
(a)$,$ (b)$,$ (c)$,$ (d)$,$ (e) from left to right$,$ the logical argument looks like
$$P \to  Q,\  R \to  S,\ ( S \land  P) \to H,\ ( Q \land  H) 
\to T \vd ( R \land P) \to  T. \eqno (1)^\prime$$
Now suppose that you use the same symbols P$,$ Q$,$ R$,$ S$,$ H$,$ T for the phrases in 
the sentence (a)$,$ (b)$,$ (c)$,$ (d)$,$ (e) (taken in the same order from left to 
right) for arguments (2)$,$ (3). Then these next two arguments would look 
like
$$P \to  Q,\  R \to  S,\ ( S \land  P) \to H,\ ( Q \land  H) 
\to T \vd ( R \land P) \to  T. \eqno (2)^\prime$$
$$P \to  Q,\  R \to  S,\ ( S \land  P) \to H,\ ( Q \land  H) 
\to T \vd ( R \land P) \to  T. \eqno (3)^\prime$$
\noindent Now$,$ from human experience$,$ compare these three patterns 
(i.e. compare them as if they are geometric configurations written left to right).
It is obvious$,$ is it not$,$ that they are the ``same.'' What this means for us 
is that the  logic behind the three arguments (1)$,$ (2) $,$ (3) appears to be the 
same logic. All we need to do is to analyze one of the patterns such as 
(1)$^\prime$ in order to understand the process more fully. For example$,$ is 
the logical argument represented by (1)$^\prime$ correct?\pars 
    One of the most important basic questions is how can we mathematically 
analyze such a logical pattern when we must use a language for the 
mathematical discussion as well as some type of logic for the analysis?
Doesn't this yield a certain type of  {\it double think} or an obvious 
paradox? This will certainly be the case if we don't proceed very carefully.
In 1904$,$ David Hilbert gave the following solution to this problem which we 
re-phrase in terms of the modern computer. A part of Hilbert's method can be 
put into the following form.\pars
The abstract language involving the symbols {\it P$,$ Q$,$ R$,$ S$,$ T$,$ H$,$} $\vd,\ 
\land, \to$ are part of the computer language for a  ``logic computer.'' The 
manner in which these symbols are combined together to form correct logical 
arguments can be checked or verified by a fixed computer program. However$,$ 
outside of the computer we use a language to write$,$ discuss and use mathematics to construct$,$ 
study and analyze the computer programs before they are entered into various 
files. Also$,$ we analyze the actual computer operations and construction using 
the same outside language. Further$,$ we don't specifically explain the human logic that 
is used to do all of this analysis and construction. Of course$,$ the symbols  
{\it P$,$ Q$,$ R$,$ S$,$ T$,$ H$,$} $\vd,\ 
\land, \to$ are a small part of the language we use. What we have is two 
languages. The language the computer understands and the much more 
complex and very large language --- in this case English --- that is employed 
to analyze and discuss the computer$,$ its programs$,$ its operations and the 
like. Thus$,$ we do our mathematical analysis of the logic computer in what is 
called a  {{\it metalanguage}} (in this case English) and we use the simplest 
possible human logic called the {{\it metalogic}} which we don't formally state. 
Moreover$,$ we use the simplest and most convincing mathematical procedures ---
procedures that we call {{\it metamathematics.}}  These procedures are those that 
have the largest amount of empirical evidence that they are consistent. In the 
literature the term {{\it meta}} is sometimes replaced by the term {{\it 
observer.}}  Using this compartmentizing procedure for the languages$,$ one 
compartment the computer language and another compartment a larger metalanguage outside of the 
computer$,$ is what prevents the mathematical study of logic from being 
``circular'' or a  ``double think'' in character. I mention that the metalogic 
is composed of a set of logical procedures that are so basic in character that 
they are universally held as correct. We use simple principles to investigate 
some highly complex logical concepts in a step-b-step effective manner.\pars
It's clear that in order to analyze mathematically human deductive procedures a 
certain philosophical stance must be taken. We must believe that the 
mathematics employed is itself correct logically and$,$ indeed$,$ that it is 
powerful enough to analyze \underbar{all} significant concepts associated with 
the discipline known as   ``Logic.'' The major reason we accept this 
philosophical stance is that the mathematical methods employed have 
applications to thousands of areas completely different from one another. If 
the mathematical methods utilized are somehow in error$,$ then these errors 
would have appeared in all of the thousands of other areas of application. 
Fortunately$,$ mathematicians attempt$,$ as best  as they can$,$ to remove all 
possible error from their work since they are aware of the fact that their 
research findings will be used by many thousands of individuals who accept 
these finding as absolutely correct logically.\pars
It's the facts expressed above that leads one to believe that the carefully 
selected mathematical procedures used by the mathematical logician are as 
absolutely correct as can be rendered by the human mind. 
Relative to the 
above arguments$,$ is it important that they be logically correct? The 
argument as stated in biological terms is an actual experimental scenario 
conducted at the University of Maryland Medical School$,$ from 1950 -- 51$,$ 
by Dr. Ernest C. 
Herrmann$,$ this author's brother. I actually aided$,$ as a teenager$,$ with the 
basic mathematical aspects for this experiment. It was shown that the 
continued use of an antibiotic not only produced resistant mutations but the 
antibiotic was also an enriched growth medium for such mutations. Their rate of 
growth increased with continued use of the same antibiotic. This led to a 
change in medical procedures$,$ at that time$,$ where combinations of antibiotics 
were used to 
counter this fact and the saving of many more lives. But$,$ the successful 
conclusion of this experiment actually 
led to a much more significant result some years later when my brother 
discovered the first useful anti-viral agent. The significance of this 
discovery  
is obvious and$,$ moreover$,$ with this discovery began the entire scientific 
discipline that studies and produces  anti-viral drugs and agents. \pars
From 1979 through 1994$,$ your author worked on one  problem and two questions
 as they were 
presented to him by John Wheeler$,$ the Joseph Henry Professor of Theoretical 
Physics at Princeton University. These are suppose to be the ``greatest 
problem and questions on the books of physics.'' 
The first problem is called the {\it General 
Grand Unification Problem.} This means to develop some sort of theory that 
will unify$,$ under a few theoretical properties$,$ all of the scientific 
theories for the 
behavior of all of the Natural systems that comprise our universe. Then  
the two other questions are ``How did our universe come into being?'' and
``Of what is empty space composed?'' As research progressed$,$  
findings were announced in various scientific journals. The first announcement 
appeared in 1981 in the {\it Abstracts of papers presented before the American 
Mathematical Society}$,$ 2(6)$,$ \#83T-26-280$,$ p. 527. Six more announcements were 
made in this journal$,$ the last one being in 1986$,$ 7(2)$,$\# 86T-85-41$,$ p. 238$,$ entitled
``A solution of the grand unification problem.'' Other important papers were 
published discussing the methods and results obtained. One of these was 
published in 1983$,$ ``Mathematical philosophy and developmental processes$,$''
{\it Nature and System$,$} 5(1/2)$,$ pp. 17-36. Another one was the 1988 paper$,$
``Physics is legislated by a cosmogony$,$'' {\it Speculations in Science and 
Technology}$,$ 11(1)$,$ pp. 17-24. There have been other publications using 
some of the procedures that were developed to solve this problem and answer the two 
questions. The last paper$,$ which contained the entire solution and almost all 
of the actual mathematics$,$ was presented before the Mathematical Association 
of America$,$ on 12 Nov.$,$ 1994$,$ at Western Maryland College. \pars
Although there are numerous applications of the methods presented within this 
text to the sciences$,$ it is  shown in section 3.9 
that there exists an elementary 
{\it ultralogic} as well as an {\it ultraword}. The properties associated with 
these two entities 
should give you a 
strong indication as to how the above discussed theoretical problem has been 
solved and 
how the
two physical questions have been answered. \vfil\eject 
\centerline{NOTES}\vfil\eject
\centerline{\bf Chapter 2 - THE PROPOSITIONAL CALCULUS}\par
\medskip
\noindent {\bf 2.1 Constructing a Language By Computer.}\pars
Suppose that you are given the symbols P$,$ Q$,$ $\land$$,$ and left
parenthesis ($,$ right parenthesis ). You want to start with the set
$L_0 = \{P,Q\}$ and construct the complete set of different (i.e. not 
geometrically congruent in the plane) strings of symbols
$L_1$ that can be formed by putting the $\land$ between two of the symbols 
from the set $L_0$$,$ {\bf with repetitions allowed}$,$ and putting the ( on the 
left and the ) on the right of the construction. {\it Also you must include the 
previous set $L_0$as a subset of $L_1.$} I hope you see easily that 
the complete set formed from these (metalanguage) rules would be \par
$$L_1 = \{P, Q, (P\land P), (Q \land  Q), (P\land Q), (Q \land P)\} \eqno (2.1.1)$$ 
\par
Now suppose that you start with $L_1$ and follow the same set of rules and 
construct the complete set of symbol strings $L_2.$  This would give
$$L_2 = \{P, Q, (P\land P), (P\land Q), (P \land (P\land P)), (P\l ( P \l Q)),
(P \l (Q \l P)),$$ $$ (P\l (Q\l Q)), (Q \l P), (Q \l Q), (Q \l (P \l P)), 
(Q\l (P\l 
Q)), (Q \l (Q \l P)),$$ 
$$ (Q\l (Q \l Q)), ((P\l P) \l P), ((P\l P) \l Q), 
((P\l P) 
\l (P \l P)), ((P\l P) \l (P \l Q)), $$ 
$$ ((P\l P) \l (Q \l P)),((P\l P) \l (Q \l Q)), ((P\l Q) \l P), ((P\l Q) \l Q), $$
$$((P \l Q) \l (P \l P)),((P \l Q)\l (P \l Q)),((P\l Q) \l (Q \l P)),$$ 
$$((P\l Q) \l (Q \l Q)), ((Q \l P ) \l P), ((Q\l P) \l Q),$$
$$((Q \l P )\l (P\l P)),((Q\l P) \l (P\l Q)),((Q \l P)\l (Q \l P)),$$
$$((Q \l P) \l (Q \l Q)), ((Q \l Q)\l P), ((Q\l Q) \l Q), $$ 
$$((Q\l Q) \l (P \l P)), ((Q\l Q) \l (P\l Q)),$$ 
$$((Q\l Q) \l (Q\l P)),((Q\l Q)\l (Q \l Q)) \}.\eqno (2.1.2)$$ \par
Now I did not form the$,$ {{\it level two}}$,$ $L_2$ by guess. 
I wrote a simple computer program that displayed this result. If I now follow 
the same instructions and form level three$,$ $L_3,$ I would print out a set 
that takes four pages of small print to express. But you have the intuitive 
idea$,$ the metalanguage rules$,$ as to what you would do if you had the 
previous level$,$ say $L_3,$ and wanted to find the strings of symbols that 
appear in $L_4.$ But$,$ the computer would have a little difficulty in 
printing out the set of all different strings of symbols or what are called
{{\it formulas$,$}} (these are also called well-defined formula by many 
authors and$,$ in that case$,$ the name is abbreviated by the symbol {{\it 
wffs}}). 
Why? Since there are 2$,$090$,$918 different formula in $L_4.$
Indeed$,$ the computer could not produce even internally all of the formulas in 
level nine$,$ $L_9,$ since there are more than $2.56 \times 10^{78}$ different 
symbol strings in this set. This number is greater than the estimated number
of atoms in the observable universe.   But you will soon able to show that  
$(((((((((P\l Q)\l (Q\land Q)))))))))\in L_9$ ($\in$ means {\it member of}) 
and this formula is not a 
member of any other level that comes before $L_9.$ You'll also be able to show 
that $(((P \l Q)\l (P\l Q))$ is not a formula at all. But all that is still to 
come. \par
In the next section$,$ we begin a serious study of formula$,$ where we can 
investigate properties associated with these symbol strings on any level of 
construction and strings that contain many more  {{\it atoms}}$,$ these are the 
symbols in $L_0$$,$ and many more {{\it connectives}}$,$ these are symbols like
$\land, \ \to$ and more to come.\pars
\ss
\noindent {\bf 2.2 The Propositional Language.}
\medskip
The are many things done in mathematical logic that are a mathematical 
formalization of obvious and intuitive things such as the above construction of 
new symbol strings from old symbol strings. {\bf The intuitive concept comes 
first and then the formalization comes after this.} In many cases$,$ I am going 
to put the actual accepted mathematical formalization in the appendix. If you 
have a background in mathematics$,$ then you can consult the appendix for 
the formal mathematical definition. As I define things$,$ I will indicate that 
the deeper stuff appears in the appendix by writing (see appendix). \pars
We need a way to talk about formula in general. That is we need symbols that 
act like {{\it formula variables}}. This means that these symbols represent 
\underbar{any} formula in our formal language$,$ with or without additional 
restrictions such as the level $L_n$ in which they are members. \pars
\vfil\eject
\hrule\smallskip
{\bf Definition 2.2.1.} Throughout this text$,$ the symbols {\it A$,$ B$,$ C$,$ D$,$ E$,$ 
F} (letters at the front of the alphabet) will denote formula variables. \parm
\hrule\medskip
In all that follows$,$ we use the following {{\it interpretation metasymbol}}$,$
``$\lint\ \rint$:'' I'll show you the meaning of this by example. The symbol 
will be presented in the following manner. \pars
\centerline{$\lint A \rint$: . . . . . . . . . . . . .}
\noindent   There will be stuff written 
where the dots . . . . . . . . . . . . . . . are placed. 
Now what you do is the 
substitute for the formula {\it A}$,$ in ever place that it appears$,$ 
the stuff that appears where the . . . . . . . 
. . . . . . .  are located. For example$,$ suppose that \par
\centerline{$\lint A\rint$: it 
rained all day$,$ $\lint \land \rint$: and}
\noindent  Then for formula $A \land A$$,$ the interpretation 
$\lint A \land A\rint$: would read\par
\centerline{it rained all day and it rained all day}
\noindent You could then adjust this so that it  
corresponds to the correct English format.  This gives\par
\centerline{It rained all day and it rained 
all day.}
\pars
Although it is not necessary that we use all of the following {{\it logical 
connectives}}$,$ using them makes it much easier to deal with ordinary everyday 
logical arguments. \parm
\hrule\smallskip  
{\bf Definition 2.2.2.} The following is the list of basic {{\it logical 
connectives}} with their technical names. \par
{\+ \indent (i) $\neg$ (Negation)& (iv) $\to$ (The conditional)\cr
\+ \indent (ii) $\l$ (Conjunction)& (v) $\iff$ (Biconditional)\cr
\+ \indent (iii) $\lor$ (Disjunction)& \cr} \pars
\hrule
\smallskip
REMARK: Many of the symbols in Definition 2.2.2 carry other names throughout 
the literature and even other symbols are used. \parm
To construct a formal language from the above logical connectives$,$ you 
consider (ii)$,$ (iii)$,$ (iv)$,$ (v) as {{\it binary}} connectives$,$ where this 
means that some formula is placed immediately to the left of each of 
them and some formula is placed immediately to the right. BUT$,$ the symbol 
$\neg$ is special. It is called an {{\it unary}} connective and formulas are 
formed as follows: your write down $\neg$ and place a formula immediately to 
the right and only the right of $\neg.$ Hence if $A$ is a formula$,$ then 
$\neg A$ is also a formula. \parm
\hrule\smallskip
{\bf Definition 2.2.3.} The construction of the propositional language $L$ (see 
appendix).\pars 
{\leftskip 0.25in \rightskip 0.25in (1) Let $P,\ Q,\ R,\ S,\ P_1, \ Q_1, \ R_1,\ S_1,\ P_1
, \ Q_2, \ R_2,\ S_2, \ldots $ be an infinite set of starting formula
called the set of {{\it atoms}}. \pars
(2) Now$,$ as our starting level$,$ take any nonempty subset of these atoms$,$ and 
call it $L_0.$\pars 
(3) You construct$,$ in a step-by-step manner$,$ the next level $L_1.$ You first 
consider as members of $L_1$ all the elements of $L_0.$ Then for each  
 and every member $A$ in $L_0$ (i.e. $A \in L_0$) you add $(\neg A)$ to $L_1.$ 
Next you 
take each and every pair of members $A,\  B$ from $L_0$ 
{\bf where repetition is allowed} 
(this means that $B$ could be the same as $A$)$,$ and add the new 
formulas $(A \land B),\ (A \lor B),\ (A \to B),\ (A \iff B).$ The result of 
this construction  is the set of formula $L_1.$ Notice that in $L_1$ every 
formula except for an atom has a left parenthesis ( and a right parenthesis ) 
attached to it. These parentheses are called {{ \it extralogical symbols}}. \pars
(4) Now repeat the construction using $L_1$ in place of $L_0$ and you get $L_2$. 
\pars
(5) This construction now continues step-by-step so that for any natural 
number $n$ you have a level $L_n$ constructed from the previous level
and level $L_n$ contains the previous levels. \par
(6) Finally$,$ a formula $F$ is a member of the propositional language $L$ if 
and only if there is some natural number $n \geq 0$ such that $F \in 
L_n.$\par}\pars
\hrule\parm
{\bf Example 2.2.1} The following are examples of formula and the particular 
level $L_i$ indicated is the first level in which they appear. Remember that 
$\in$ means  ``a member or element of''.\pars
$P \in L_0;\ (\neg P)\in L_1;\ (P\land (Q \to R))\in L_2;\ ((P\land Q)\land R) 
\in L_2;\ (P\land (Q\land R)) \in L_2;\ ((P\to Q)\lor(Q \to S))\in L_2;\ 
(P \to (Q \to (R\to S_2)))\in L_3.$\pars
\vfil\eject
{\bf Example 2.2.2} The following are examples of strings of symbols that are 
NOT in $L.$\pars
$(P);\ ((P\to Q);\ \neg(P);\ ()Q;\ (P\to (Q));\ (P=(Q \to S)).$\pars
Unfortunately$,$ some more terms must be defined so that we can communicate 
successfully. Let $A \in L.$ The {${\rm size} (A)$} is the smallest $n\geq 0$ 
such that $A \in L_n.$ Note that if ${\rm size}(A)$ = n$,$ then $A \in L_m$ for 
each level $m$ such that $m\geq n.$ And$,$ of course$,$ $A\not\in L_k$ for all 
$k$$,$ if any$,$ such that $0 \leq k < n.$ ($\not\in$ is read  ``not a member 
of''). Please note what symbols are metasymbols and that they are not symbols within 
the formal language $L.$ \pars
There does not necessary exist a unique interpretation of the above formula in 
terms of English language expressions. There is a very basic interpretation$,$ 
but there are others that experience indicates are logically equivalent to the 
basic interpretations. The symbol $\nat$ means the set $\{0,1,2,3,4,5,\ldots 
\} $ of natural numbers including zero.\parm
\hrule\smallskip
{\bf Definition 2.2.4} The basic English language interpretations.\pars
{\leftskip 0.25in \rightskip 0.25in (i) $\lint \neg \rint$: not$,$ (it is not the  case that).\pars
(ii) $\lint \land \rint$: and\pars
(iii) $\lint \lor \rint$: or\pars
(iv) For any $A \in L_0,$ $\lint A \rint$: a simple declarative sentence$,$ a 
sentence which contains no interpreted logical connectives OR a set of English 
language symbols that is NOT considered as decomposed into distinct 
parts.\pars
(v) For any $n \in \nat, \ A,\ B \in L_n$; $\lint A\lor B \rint$: $\lint A 
\rint$ or $\lint B \rint$.\pars       
(vi) For any $n \in \nat, \ A,\ B \in L_n$; $\lint A\land B \rint$: $\lint A 
\rint$ and $\lint B \rint$.\pars        
(vii) For any $n \in \nat, \ A,\ B \in L_n$; $\lint A\to B \rint$: if $\lint A 
\rint$$,$ then $\lint B \rint$.\pars 
(viii) For any $n \in \nat, \ A,\ B \in L_n$; $\lint A\iff B \rint$: $\lint A 
\rint$ if and only if $\lint B \rint$.\pars
(ix) The above interpretations are then continued  ``down'' the levels $L_n$ 
until they stop at level $L_0.$}\pars
\hrule\medskip
Please note that the above is not the only translations that can be applied to 
these formulas. Indeed$,$ the electronic hardware known as switching circuits or 
gates can also be used to interpret these formulas. This hardware 
interpretation is what has produced the modern electronic computer.\pars
Unfortunately$,$ when translating from English or conversely the members of $L,$ 
the 
above basic interpretations must be greatly expanded. The following is a list 
for reference purposes of the usual English constructions that can be properly 
interpreted by members of $L.$\parm
\hrule\smallskip
\+ (x) For any $n \in \nat, \ A,\ B \in L_n$; $\lint A\iff B \rint$: &\cr
\+ (a) $\lint A \rint$ if $\lint B \rint,$ and 
$\lint B \rint$ if $\lint A \rint$. & (g)  $\lint A \rint$ exactly if
  $\lint B \rint$.\cr
\+ (b) If $\lint A \rint$$,$ then  $\lint B \rint,$ and 
conversely.&(h) $\lint A \rint$ is material equivalent to $\lint B \rint.$\cr
\+ (c) $\lint A \rint$ is (a) necessary and sufficient (condition) for 
       $\lint B \rint$ & \cr
\+ (d) $\lint A \rint$ is equivalent to  $\lint B \rint$. (sometimes used in 
this manner)&\cr
\+ (e) $\lint A \rint$ exactly when $\lint B \rint.$ & (i) 
$\lint A \rint$ just in case  $\lint B \rint.$ \cr
\+ (f) If and only if $\lint A \rint$$,$ (then) $\lint B \rint.$ &\cr
\pars
\hrule\smallskip
\+ (xi) For any $n \in \nat, \ A,\ B \in L_n$; $\lint A\to B \rint$: &\cr
\+ (a) $\lint B \rint$ if $\lint A \rint.$ & (h) $\lint A \rint$ only if $\lint B 
\rint.$\cr
\+(b) When $\lint A \rint,$ then  $\lint B \rint.$ & (i) $\lint B \rint$ when 
$\lint A \rint.$\cr
\+(c) $\lint A \rint$ only when  $\lint B \rint.$ & (j) In  case $\lint A \rint,$
$\lint B \rint.$\cr
\+(d) $\lint B \rint$ in case $\lint A \rint.$& (k) $\lint A \rint$ only in case
$\lint B \rint.$\cr
\+(e)  $\lint A \rint$ is a sufficient condition for  $\lint B \rint.$&\cr
\+(f) $\lint B \rint$ is a necessary condition for $\lint A \rint.$&\cr
\+(g) $\lint A \rint$ materially implies $\lint B \rint.$&(l) $\lint A 
\rint$ implies $\lint B \rint.$\cr\pars
\hrule
\smallskip
\+(xii) For any $n \in \nat, \ A,\ B \in L_n$; $\lint A\land B \rint$: &\cr
\+(a) Both $\lint A \rint$ and $\lint B \rint.$&(e) Not only $\lint A \rint$ but
$\lint B \rint.$ \cr
\+(b) $\lint A \rint$ but $\lint B \rint.$&(f) $\lint A \rint$ while $\lint B 
\rint.$ \cr
\+(c) $\lint A \rint$ although $\lint B \rint.$&(g) $\lint A \rint$ despite 
$\lint B \rint.$ \cr
\+(d) $\lint A \rint$ yet $\lint B \rint.$&\cr\pars
\hrule
\smallskip
\+(xiii) For any $n \in \nat, \ A,\ B \in L_n$; $\lint A\lor B \rint$: &\cr
\+(a) $\lint A \rint$ or $\lint B \rint$ or both.&(d) $\lint A \rint$ and/or
$\lint B\rint$\cr
\+(b) $\lint A \rint$ unless $\lint B \rint.$& (e) Either $\lint A \rint$ or 
$\lint B \rint.$ (usually) \cr                        
\+(c) $\lint A \rint$  except when $\lint B \rint.$ (usually)&\cr\pars
\hrule\smallskip
\+(xiv) For any $n \in \nat, \ A,\ B \in L_n$; $\lint ( A\lor B 
)\land(\neg(A\land B))\rint$: &\cr       
\+(a) $\lint A \rint$ or $\lint B \rint$ not both.&(c) $\lint A \rint$ or else 
$\lint B \rint.$ (usually)\cr
\+(b) $\lint A \rint$ or $\lint B \rint.$ (sometimes)&(d) Either $\lint A \rint$ or 
$\lint B \rint.$ (sometimes)\cr\pars
\hrule\smallskip
\noindent (xv) For any $n \in \nat, \ A,\ B \in L_n$; $\lint (\neg(A \iff 
B))\rint$: $\lint ((\neg A) \iff B)))\rint$:\par   
\+(a) $\lint A \rint$ unless $\lint B \rint.$ (sometimes)&\cr\pars
\hrule\smallskip
\+(xvi) For any $n \in \nat, \ A,\ B \in L_n$; $\lint (A \iff (\neg B))\rint$: 
&\cr
\+(a) $\lint A \rint$ except when  $\lint B \rint.$ (sometimes)&\cr\pars
\hrule\smallskip
\+(xvii) For any $n \in \nat, \ A,\ B \in L_n$; $\lint(\neg (A \lor B))\rint$: 
&\cr
\+(a) Neither $\lint A \rint$ nor $\lint B \rint.$ &\cr\pars
\hrule\smallskip
\+(xviii) For any $n \in \nat,\ A \in L_n$; $\lint(\neg A )\rint$: 
&\cr\pars
Not $\lint A \rint$ (or the result of transforming $\lint A \rint$ to give the 
intent of  ``not'' such as  ``$\lint A \rint$ doesn't hold'' or ``$\lint A \rint$ 
isn't so.''\pars
\hrule\smallskip\hrule\bigskip
\centerline{\bf EXERCISES 2.2}
\medskip
\centerline{In what follows assume that $P,\ Q,\ R,\ S \in L_0.$}
\noindent 1. Let $A$ represent each of the following strings of symbols.
Determine if $A \in L$ or $A \not\in L.$ State your conclusions.
\+\indent (a) $A = (P \lor(Q\to (\neg S))$& (f) $A = )P) \lor((\neg S)))$\cr
\+\indent (b) $A = (P\iff(Q\lor S))$& (g) $A = (P\iff (\neg(R\iff S)))$\cr
\+\indent (c) $A = (P \to (S \land R))$& (h) $A = (R \land (\neg(R\lor S)) 
\to P)$\cr
\+\indent (d) $A = ((P)\to (R\land S))$& (i) $A = (P\land(P\land P)\to Q)$\cr
\+\indent (e) $A = (\neg P) \to (\neg(R\lor S))$& (j) $A =((P\land P)\to P\to 
P)$\cr\parm
\noindent 2. Each of the following formula $A$ are members of $L.$ Find the
${\rm size}(A)$ of each. \pars
\+ (a) $A = ((P\lor Q) \to (S\to R))$& (c) $A = (P\lor(Q \land(R\land S)))$\cr
\+ (b) $A = (((P\lor Q)\to R)\iff S)$& (d) $ A = (((P\lor Q)\iff (P\land Q)) 
\to S)$\cr\parm
\noindent 3. Use the indicated atomic symbol to translate each of the 
following into a member of $L.$\pars
(a) Either (P) the port is open or (Q) someone left the shower on.\pars
(b) If (P) it is foggy tonight$,$ then either (Q) the Captain will stay in his 
cabin or (R) he will call me to extra duty.\pars
(c) (P) Midshipman Jones will sit$,$ and (Q) wait or (R) Midshipman George will 
wait.\pars
(d) Either (Q) I will go by bus or (R) (I will go) by airplane. \pars
(e) (P) Midshipman Jones will sit and (Q) wait$,$ or (R) Midshipman George will 
wait.\pars
(f) Neither (P) Army nor (Q) Navy won the game. \pars
(g) If and only if the (P) sea-cocks are open$,$ (Q) will the ship sink; (and) should 
the ship sink$,$ then (R) we will go on the trip and (S) miss the dance. \pars
(h) If I am either (P) tired or (Q) hungry$,$ then (R) I cannot study.\pars
(i) If (P) Midshipman Jones gets up and (Q) goes to class$,$ (R) she will pass 
the quiz; and if she does not get up$,$ then she will fail the quiz. \parm
\noindent 4. Let $\lint P \rint$: it is nice; $\lint Q \rint$: it is hot;
  $\lint R \rint$: it is cold; $\lint S \rint$:  it is small. Translate 
(interpret)  the following formula into acceptable non-ambiguous English
sentences.\pars
\+\indent (a) $(P \to (\neg(Q \land R)))$& (d) $((S \to Q)\lor P)$\cr
\+\indent (b) $(S \iff P)$& (e) $(P \iff ((Q \land (\neg R))\lor S))$\cr
\+\indent (c) $(S \land(P\lor Q)) $& (f) $((S \to Q)\lor P)$\cr\pars\ss
\noindent {\bf 2.3 Slight Simplification$,$ Size$,$ Common Pairs and Computers.}
\medskip
Each formula has a unique size $n$$,$ where $n$ is a natural number$,$ $\nat,$ 
greater than  or equal to zero. Now if ${\rm size}(A) = n,$ then 
$A \in L_m$ for all $m \geq n,$ and $A \not\in L_m$ for all $m < n.$  For each 
formula that is not an atom$,$ there appears a certain number of left ``('' and 
right ``)'' parentheses. These parentheses occur in what is called {{\it common 
pairs}}. Prior to the one small simplification we may make to a formula$,$ we'll
learn how to calculate which parentheses are common pairs. The common pairs 
are the parentheses that are included in a specific construction step for a 
specific level $L_n.$ The method we'll use can be mathematically established; 
however$,$ its demonstration is somewhat long and tedious. Thus the  ``proof''
will be omitted. The following is the {{\it common pair rule}}. \parm
\hrule
\smallskip
{\bf Rule 2.3.1.} This is the common pair rule (CPR). Suppose that we are 
given an expression that is thought to be a member of $L.$\pars
(1)  Select any
left parenthesis ``(.'' Denote this parenthesis by the number +1. \pars
(2) Now moving towards the right$,$ each time you arrive at another left
parenthesis  ``('' add the number 1 to the previous number.\pars
(3) Now moving towards the right$,$ each time you arrive at a 
right parenthesis  ``)'' subtract the number 1 from the previous number.\pars
(4) {\bf The first time} you come to a parenthesis that yields a ZERO by the 
above cumulative algebraic summation process$,$ then that particular parenthesis
is the companion parenthesis with which the first parenthesis you started with 
forms a common pair.\pars
\hrule
\medskip
The common pair rule will allow us to find out what expressions within a 
formula 
are also formula. This rule will also allow us to determine the size of a 
formula. A formula is written in {{\it atomic form}} if only atoms$,$ connectives$,$ 
and parentheses appear in the formula.
\parm
\hrule \smallskip
{\bf Definition 2.3.1} Non-atomic subformula.\pars
Given an $A \in L$ (written in atomic form). A subformula is any expression 
that appears between and includes a common pair of parentheses. \pars
\hrule
\medskip
Note that according to Definition 2.3.1$,$ the formula $A$ is a subformula. I 
now state$,$ without proof$,$ the theorem that allows us to determine the size
of a formula.\parm
{\bf Theorem 2.3.1} {\sl Let $A \in L$ and $A$ is written in atomic form.
If there does not exist a parenthesis in $A,$ then  $A \in L_0$  and $A$ has 
size zero. If there exists a left most parenthesis  ``('' {\rm [}i.e. no more parentheses 
appear on the left in the expression{\rm ]}$,$ then beginning  with this parenthesis
the common pair rule will yield a right most parenthesis for the common pair.
During this common pair procedure$,$ the  largest natural number obtained will 
be the size of $A.$ }\parm
{\bf Example 2.3.1.}  The numbers in the next display are obtained by 
starting with parenthesis $a$ and show that the ${\rm size}(A) = 3$.\pars
$$ \matrix{A =&(&(&(&P\land Q&)&\lor R&)&\to &(&S\lor P&)&)\cr  
\quad&a&b&c&\quad&d&\quad&e&\quad&f&\quad&g&h\cr
\quad&1&2&3&\quad&2&\quad&1&\quad&2&\quad&1&0\cr}$$
Although the common pairs are rather obvious$,$ the common pair rule can be 
used. This gives (c$,$d)$,$ (b$,$e)$,$ (f$,$g)$,$ and (a$,$h) as common pairs. Hence$,$
the subformula are $A,\ (P\land Q),\ ((P\land Q)\lor R),\ (S \lor P)$ and you 
can apply the common pair rule to each of these to find their sizes.
Of course$,$ the rule is most useful when the formula are much more 
complex.\parm
The are various simplification processes that allow for the removal of many of 
the parenthesis. One might think that logicians like to do this since those 
that do not know the  
simplifications would not have any knowledge as to what the formula actually 
looks like.  The real reason is to simply write less. These simplification 
rules are relative to a listing of the strengths of the connectives. However$,$
for this beginning course$,$ such simplification rules are not important with 
one exception. \parm
\hrule
\smallskip
{\bf Definition 2.3.2} The one simplification that may be applied is the 
removal of the outermost left parenthesis ``('' and the outermost
right  parenthesis ``).'' It should be obvious when such parentheses have been 
removed. BUT$,$ they must be inserted prior to application of the common pair rule 
and Theorem 2.3.1.\pars
\hrule
\medskip
One of the major applications of the propositional calculus is in the design 
of the modern electronic computer. This design is based upon the basic and 
simplest behavior of the logical network which itself is based upon the simple 
idea of switching circuits. Each switching device is used to produce the 
various types of  ``gates.'' These gates will not specifically be identified 
but the basic switches will be identified.  Such switches are conceived of as 
simple single pole relays. A switch may be {{\it normally open}} when no current flows 
through the coil. One the other hand$,$ the switch could be {{normally closed}} 
when no current flows. When current flows through the relay coil$,$ 
the switch takes the opposite state. The coil is not shown only the circuit 
that is formed or broken when the coil is energized for the $P$ or $Q$ relay. 
The action is ``current through a coil'' and leads to or prevents current 
flowing through the indicated line. \pars
For what follows the atoms of our propositional calculus represent relays in 
the normally open position. Diagrammatically$,$ a normally open relay (switch) is 
symbolized as follows: \parm
\line{(a) $\lint P \rint: \ \ \ \ \underline{\em {7}}{\kern -.1em|} \ \ P 
\  \backslash\kern -.1em\underline{\em {7}}$\hfil}\parm
\noindent Now for each atom $P,$  let $\neg P$ represent a normally closed 
relay. Diagrammatically$,$ a normally closed relay is symbolized by:\parm
\line{(b) $\lint \neg P \rint: \ \ \ \ \underline{\em {6.5}}\kern -.1em|\kern -.1em
\underline{\ \neg P\ }\kern -.1em|\kern -.1em\underline{\em {6.5}}$\hfil}\parm
\indent Now any set of normally open or closed relays  can be wired together 
in various ways. For the normally open ones$,$ the following will model the 
binary connectives $\to,\ \iff, \land,\ \lor.$ This gives a non-linguistic 
model. In the following$,$ $P, \ Q$ are atoms and the relays are normally 
open.\parm
\line{\hfil$\kern -.9em\underline{\em {4}}\kern -.1em|\ \ P\ \ \ 
\backslash\kern -.1em\underline{\em {4.25}}$\hfil}
\baselineskip=10pt
\line{${\rm (i)}\ \lint P \lor Q\rint: \hfil \ \ \ \kern -.3em\underline{\em {4.0}}\kern -.1em|$
\hskip 3.9cm$\kern -.6em |\kern -.1em\underline{\em {4}}$\hskip 2.63cm\hfil}
\line{\hfil$\kern -1.00em|\kern -.1em\underline{\em {4}}\kern -.1em|\ \ Q\ \ \ 
\backslash\kern -.1em\underline{\em {4}}\kern -.1em|$\hfil}
\baselineskip=14pt
\bigskip           
\line{(ii) $\lint P \land Q\rint: \hfil\underline{\em {5.5}}
{\kern -.1em|} \ \ P 
\  \backslash\kern -.1em\underline{\em {3}}{\kern -.1em|} \ \ Q 
\  \backslash\kern -.1em\underline{\em {5.5}} $\hfil}
\medskip
\line{\hfil$\kern -.9em\underline{\em {4}}\kern -.1em|\kern -.1em\underline{\ 
\ \neg P \ \ }\kern -.1em|\kern -.1em\underline{\em {4.25}}$\hfil}
\baselineskip=10pt
\line{${\rm (iii)}\ \lint P \to Q\rint:\hfil \ \kern -.58em\underline{\em {4.0}}\kern -.1em|$
\hskip 3.9cm$\kern -.6em |\kern -.1em\underline{\em {4}}$\hskip 2.63cm\hfil}
\line{\hfil$\kern -1.00em|\kern -.1em\underline{\em {4}}\kern -.1em|\ \ Q\ \ \ 
\backslash\kern -.1em\underline{\em {4}}\kern -.1em|$\hfil}
\baselineskip=12pt\medskip
\noindent (iv) In order to model the expression $P\iff Q$ we need two coils. The 
$P$ coil has a switch at both ends$,$ one normally open the other normally 
closed. The 
$Q$ coil has a switch at both ends$,$ one normally open the other normally 
closed. But$,$ the behavior of the two coils is opposite from one another as 
shown in the following diagram$,$ where (iii) denotes the previous diagram. 
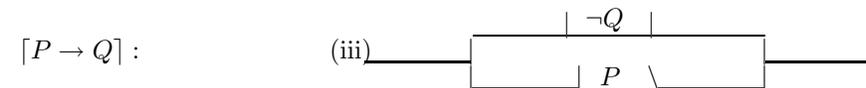
\parm                                                                                                                                      
\line{\hfil$\kern -.9em\underline{\em {3.5}}\kern -.1em|\kern -.1em$\underbar{\ 
\ \raise 2pt\hbox{$\neg Q$} \ \ }$\kern -.1em|\kern -.1em\underline{\em {4.25}}$\hfil}
\baselineskip=10pt
\line{$\ \lint P \to Q\rint:\hfil {\rm (iii)}\kern -.3em\underline{\em {4.0}}\kern -.1em|$
\hskip 3.9cm$\kern -.6em |\kern -.1em\underline{\em {4}}$\hskip 2.63cm\hfil}
\line{\hfil$\kern -1.00em|\kern -.1em\underline{\em {4}}\kern -.1em|\ \ P\ \ \ 
\backslash\kern -.1em\underline{\em {4}}\kern -.1em|$\hfil}
\baselineskip=14pt
\bigskip\bigskip
\hrule\smallskip\hrule\bigskip
\centerline{\bf EXERCISES 2.3}
\bigskip 
\noindent 1. When a formula is written in atomic form$,$ the (i) atoms$,$ (ii) 
(not necessary distinct) connectives$,$ and (iii) the parentheses are displayed. 
Of the three collections (i)$,$ (ii)$,$ (iii)$,$ find the one collection that can be 
used to determine immediately by counting the (a) number of common pairs and 
(b) the number of subformula. What is it that you count?\parm 
\noindent 2. For each of the following formula$,$ use the indicated letter and 
list as order pairs$,$ as I have done for Example 2.3.1 on page 16$,$ the letters 
that identify the common pairs of parentheses.\pars
\indent {$(A)=((P\to(Q\lor R))\iff(\neg S))$}\par         
\baselineskip=10pt
 {${\rm \ \ \  \ \ \ \ \ \ \ ab\ \ \ \ \ \ c \ \ \ \ \ \ de\ \ \ \ \ f\ \ \ \ gh}$} 
\baselineskip=14pt        
\parm
\indent $(B)=((P\lor(Q\lor(S\lor Q))) \to (\neg(\neg R)))$\par         
\baselineskip=10pt
${\rm \ \ \  \ \ \ \ \ \ \ ab\ \ \ \ c\ \  \ \ \  d\ \ \ \ \ \ \ efg\ \ \ \ 
 \ h\ \ i \ \ \ \ jkm}$ 
\baselineskip=14pt      
\parm  
\indent $(C)=(((\neg (P\iff (\neg R)))\to ((\neg Q)\iff (R\lor P)))\iff S) $\par         
\baselineskip=10pt
${\rm\  \ \ \ \ \ \ \ abc \ d\ \  \ \ \ \ \  e\ \ \ \ fgh\ \ \ \ 
 \ ij\ \ \ \ k \ \ \ \ m \ \ \ \ \ \ \ nop\ \ \ \ \ q}$ 
\baselineskip=14pt   
\parm
\noindent 3. Find the size of each of the formula in problem 3 above.\parm 
\noindent 4. Although it would not be the most efficient$,$ (we will learn how 
to find logically equivalent formula so that we can make them more efficient)$,$
use the basic relay (switching) circuits described in this section (i.e. 
combine them together) so that the circuits will model the following 
formula.\pars
\+\indent (a) $((P\lor Q)\land (\neg R))$ & (c) 
$(((\neg P)\land Q) \lor ((\neg Q)\land P))$\cr
\+\indent  (b) $((P\to Q)\lor (Q \to P))$ & (d) $((P\land Q )\land (R \lor 
S))$\cr \s\ss
\noindent {\bf 2.4 Model Theory -- Basic Semantics.}
\medskip
Prior to 1921 this section could not have been rigorously written. It was not 
until that time when Emil Post convincingly established that the  {{\it 
semantics}} for the language $L$ and the seemingly more complex formal 
approach to Logic as practiced in the years prior to 1921 are equivalent. The 
semantical ideas had briefly been considered for some years prior to 1921. 
However$,$ Post was apparently the first to investigate rigorously such 
concepts. It has been said that much of modern mathematics and various 
simplifications came about since ``we are standing on the shoulders of 
giants.'' This is$,$ especially$,$ true when we consider today's simplified
semantics for $L.$ \pars
 Now the term {{\it semantics}} means that we are going to give a special 
meaning to each member of $L$ and supply rules to obtain these
``meanings'' from the atoms and connectives whenever a formula is written in
atomic form (i.e. only atoms$,$ connectives and parentheses appear). These 
meanings will  ``mirror'' or ``model'' the behavior of the classical concepts 
of  ``truth'' and ``falsity.'' HOWEVER$,$ so as not to impart any philosophical 
meanings to our semantics until letter$,$ we replace   ``truth'' by the letter 
$T$
and ``falsity'' by the $F.$ \parm
\hrule
\smallskip
In the applications of the following semantical rules to the real world$,$ 
it 
is often better to model the $T$ by the concept  ``occurs in reality''
and the $F$ by the concept ``does not occur in reality.'' Further$,$ in many 
cases$,$ the words  ``in reality'' many not be justified. \pars
\hrule
\medskip
\smallskip
\hrule
\smallskip
{\bf Definition 2.4.1} The following is the idea of an {{\it assignment}}.
Let $A \in L$ and assume that $A$ is written in atomic form. Then there exists some natural number $n$ such that $A \in L_n$ 
and ${\rm size} (A)=n.$ Now there is in $A$ a finite list of distinct
atoms$,$ say $(P_1, P_2, \ldots, P_m)$ reading left to right. We will assign to 
each $P_i$ in the list the symbol T or the symbol F in as many different ways 
as possible. If there are $n$ different atoms$,$ there will be $2^n$ different 
arrangements of such Ts and Fs. These are the  {{\it values}} of the assignment. This can be 
diagrammed as follows: \par
\centerline{$(P_1,P_2,\ldots,P_m)$}\par
\centerline{$(\updownarrow\ \ \ \updownarrow\ \  \cdots\ \updownarrow)$}\par
\centerline{$(\ T\ , \ F\ ,\ldots,\ T)$}\parm
\hrule
\smallskip
{\bf Example 2.4.1} This is the example that shows how to give a standard fixed 
method to display and find all of the different arrangements of the Ts and Fs
for$,$ say three atoms$,$ $P, \ Q,\ R.$ There would be a total of 8 different 
arrangements. Please note how I've generated the first$,$ second and third 
columns of the following  ``assignment'' table. This same idea can be used to 
 generate quickly an assignment table for any finite number of atoms.\vfil\eject
\newbox\medstrutbox
\setbox\medstrutbox=\hbox{\vrule height14.5pt depth9.5pt width0pt}
\def\medstrut{\relax\ifmmode\copy\medstrutbox\else\unhcopy\medstrutbox\fi}

\bigskip

\newdimen\leftmargin
\leftmargin=0.0truein
\newdimen\widesize
\widesize=2.0truein
\advance\widesize by \leftmargin       
\hfil\vbox{\tabskip=0pt\offinterlineskip
\def\tablerule{\noalign{\hrule}}

\halign to \widesize{\medstrut\vrule#\tabskip=0pt plus2truein

&\hfil\quad#\quad\hfil&\vrule#
&\hfil\quad#\quad\hfil&\vrule#
&\hfil\quad#\quad\hfil&\vrule#

\tabskip=0pt\cr\tablerule


&$P$&&$Q$&&$R$&\cr\tablerule
&$T$&&$T$&&$T$&\cr\tablerule
&$T$&&$T$&&$F$&\cr\tablerule
&$T$&&$F$&&$T$&\cr\tablerule
&$T$&&$F$&&$F$&\cr\tablerule
&$F$&&$T$&&$T$&\cr\tablerule
&$F$&&$T$&&$F$&\cr\tablerule
&$F$&&$F$&&$T$&\cr\tablerule
&$F$&&$F$&&$F$&\cr\tablerule
}}
\parm
The rows in the above table represent one assignment (these are also called 
truth-value assignments) and such an assignment 
will be denoted by symbols such as $\underline{a}=(a_1,a_2,a_3),$ where the 
$a_i$ take the value $T$ or $F$. In practice$,$ one could re-write these 
assignments in terms of the numbers 1 and 0 if one wanted to remove any 
association with the philosophical concept of  ``truth'' or ``falsity.''\pars
For this fundamental discussion$,$ we will assume that the list of atoms in a 
formula $A \in L$ is known and we wish to define in a appropriate manner the 
intuitive concept of the truth-value for the formula $A$ for a specific 
assignment $\ass.$ Hence you are given $\ass = (a_1,a_2,\ldots, a_n)$ that 
corresponds to the atoms $(P_1,P_2,\ldots, P_n)$  in $A$ and we want a 
definition that will allow us to determine inductively a unique 
``truth-value'' from $\{T,F\}$ that corresponds to $A$ for the assignment 
$\ass.$ \parm
\hrule\smallskip
{\bf Definition 2.4.2} The truth-value for a given formula $A$ and a given
assignment $\ass$ is denoted by $v(A,\ass).$ The procedure that we'll use is 
called a {{\it valuation procedure}}.\pars\hrule\parm
Prior to presenting the valuation procedure$,$ let's make a few observations. 
If you take a formula with $m$ distinct atoms$,$ then the assignments $\ass$ are
exactly the same for {\it any} formula with $m$ distinct atoms no matter 
what they are. \parm
(1) An assignment $\ass$ only depends upon the \underbar{number} of 
distinct atoms and not the actual atoms themselves. \parm
(2) Any rule that assigns a truth-value $T$ or $F$ to a formula $A,$ 
where $A$ is not an atom must depend only upon the connectives contained in 
the formula.\pars
(3) For any formula with $m$ distinct atoms$,$ changing the names of the $m$ 
distinct atoms to different atoms that are distinct will not change the 
assignments.\parm 
Now$,$ we have another observation based upon the table of assignments that 
appears on page 20. This assignment table is for three distinct atoms.
Investigation of just two of the columns yields the following:\parm
(4) For any assignment table for $m$ atoms$,$ any $n$ columns$,$ where $1\leq n 
\leq m$ can be used to obtain (with possible repetition) all of the 
assignments that correspond to a set of $n$ atoms.\parm
The actual formal inductively defined valuation procedure
is given in the appendix and is based upon the size of a formula. This formal 
procedure is not the actual way must mathematicians obtain $v(A,\ass)$ for a 
given $A,$ however. What's presented next is the usual informal (intuitive) 
procedure that's used. It's called the {{truth-table}} procedure and is based 
upon the ability of the human mind to take a general statement and to apply 
that statement in a step-by-step manner to specific cases. There are five
basic truth-tables.\parm
The $A,B$ are {\it any} two formula in $L.$ As indicated 
above we need only to define the truth-value for the five connectives. \parm
\hskip 0.1in {\vbox{\offinterlineskip
\hrule
\halign{&\vrule#&
\strut\hfil#\hfil\cr
height2pt&\omit&\cr
&$\underline{\em {.40} \ \ \ {\rm (i)}\ \ \ \hskip .40em\relax}$&\cr
&$\underline{\ A\hskip .75em\relax|\ \neg A\ }$&\cr
&$\underline{\hskip .2em\relax\ T\em {.1}\ |\ \ F\ \hskip .30em\relax}$&\cr
&$\ \em {-.1}F\em {.1}\ |\ \ T\ $&\cr
height2pt&\omit&\cr}
\hrule}}    
\hskip 0.75in{\vbox{\offinterlineskip
\hrule
\halign{&\vrule#&
\strut\hfil#\hfil\cr
height2pt&\omit&\cr
&$\underline{\em {2.6}{\rm (ii)}\em {2.6}}$&\cr
&$\underline{\em {.2}A\ |\ B\ |\ A \lor B\em {.1}}$&\cr
&$\underline{\em {.3}T\ \em {-.1}|\ T\ \em {.1}| \ \ \ T \ \ \ \em {.4}}$&\cr
&$\underline{\em {.3}T\ \em {-.1}|\ F\ \em {.1}| \ \ \ T \ \ \ \em {.4}}$&\cr
&$\underline{\em {.3}F\ \em {-.1}|\ T\ \em {.1}| \ \ \ T \ \ \ \em {.4}}$&\cr
&${\em {.00}F\em {.2}| \ F\ | \ \ F\ \ \ \ }$&\cr 
height2pt&\omit&\cr}
\hrule}}    
\hskip 0.75in{\vbox{\offinterlineskip
\hrule
\halign{&\vrule#&
\strut\hfil#\hfil\cr
height2pt&\omit&\cr
&$\underline{\em {2.5}{\rm (iii)}\em {2.5}}$&\cr 
&$\underline{\em {.2}A\ |\ B\ |\ A \land B\em {.1}}$&\cr
&$\underline{\em {.35}T\ \em {-.1}|\ T\ \em {.2}| \ \ \ T\ \ \ \em {.4}}$&\cr
&$\underline{\em {.35}T\ \em {-.1}|\ F\ \em {.2}| \ \ \ F\ \ \ \em {.4}}$&\cr
&$\underline{\em {.2}F\ |\ T\ \em {.2}| \ \ \ F\ \ \ \em {.4}}$&\cr
&$\em {-.1}F\em {.3}| \ F\ | \ \ \em {.2}F\em {-.2} \ \ \ \ $&\cr 
height2pt&\omit&\cr}
\hrule}}    
\hskip 0.75in{\vbox{\offinterlineskip
\hrule
\halign{&\vrule#&
\strut\hfil#\hfil\cr
height2pt&\omit&\cr
&$\underline{\em {3.0}{\rm (iv)}\em {3.0}}$&\cr 
&$\underline{\em {.4}A\em {.7}|\ B\ |\ A \to B\em {.2}}$&\cr
&$\underline{\em {.7}T\ |\ T\ \em {.2}|\em {.1} \ \ \ T\ \ \em {1.1}}$&\cr
&$\underline{\em {.7}T\ |\ F\ \em {.2}|\em {.1} \ \ \ F\ \ \em {1.1}}$&\cr
&$\underline{\em {.7}F\ |\ T\ \em {.3}| \ \ \ T\ \ \em {1.1}} $&\cr
&$\em {-.7}F\em {.4}| \ F\ \em {.2}| \ \ \em {.2}T\ \ \ \em {-.4}$&\cr 
height2pt&\omit&\cr}
\hrule}} \par\bigskip\bigskip   
\hskip 2.5in{\vbox{\offinterlineskip
\hrule
\halign{&\vrule#&
\strut\hfil#\hfil\cr
height2pt&\omit&\cr
&$\underline{\em {3.1}{\rm (v)}\em {3.1}}$&\cr 
&$\underline{\em {.4}A\em {.7}|\ B\ |\ A \iff B\em {.2}}$&\cr
&$\underline{\em {.7}T\ |\ T\ \em {.2}|\em {.1}\ \ \ T\ \ \em {1.1}}$&\cr
&$\underline{\em {.7}T\ |\ F\ \em {.2}|\em {.1}\ \ \ F\ \ \em {1.1}}$&\cr
&$\underline{\em {.7}F\ |\ T\ \em {.2}|\em {.1}\ \ \ F\ \ \em {1.1}}$&\cr
&$\em {-.9}F\em {.25}| \ F\ \em {.1}| \ \ \em {.2}T\ \ \em {-.4}$&\cr 
height2pt&\omit&\cr}
\hrule}} \par
\bigskip
Observe that the actual truth-value for a connective does not depend upon the 
symbols $A$ or $B$ but only upon the values $T$ or $F$. For this reason the
above general truth-tables can be replaced with a simple collection of 
statements relating the $T$ and $F$ and the connectives only. This is the 
quickest way to find the truth-values$,$ simply concentrate upon the connectives 
and use the following:\pars
 (i) ${\buildrel \neg F \over T},\  {\buildrel \neg T \over F}.$\pars
(ii) ${\buildrel T \lor T \over T},\ {\buildrel T \lor F \over T},
\ {\buildrel F \lor T \over T},\ {\buildrel F \lor F \over F}$\pars 
(iii) ${\buildrel T \land T \over T},\ {\buildrel T \land F \over F},
\ {\buildrel F \land T \over F},\ {\buildrel F \land F \over F}$\pars 
(iv) ${\buildrel T \to T \over T},\ {\buildrel T \to F \over F},
\ {\buildrel F \to T \over T},\ {\buildrel F \to F \over T}$\pars 
(v) ${\buildrel T \iff T \over T},\ {\buildrel T \iff F \over F},
\ {\buildrel F \iff T \over F},\ {\buildrel F \iff F \over T}$\parm 

The procedures will now be applied in a step-by-step manner to 
find the truth-values for a specific formula. This will be displayed as a 
truth-table with the values in the last column. Remember that the actual 
truth-table can contain many more atoms than those that appear in a given 
formula. By using all the distinct atoms contained in all the formulas$,$ one truth-table can be used to find the truth values for more than one formula.\pars 
The construction of a truth-table is best understood by example. In the 
following example$,$ the numbers 1$,$ 2$,$ 3$,$ 4 for the rows and the letters
a$,$ b$,$ c$,$ d$,$ e$,$ f that identify the columns are used here for reference only  
and are not used in the actual construction.\parm\vfil\eject
{\bf Example 2.4.2} Truth-values for the formula $A = (((\neg P)\lor R)\to 
(P\iff R)).$\par  
\newbox\medstrutbox
\setbox\medstrutbox=\hbox{\vrule height14.5pt depth9.5pt width0pt}
\def\medstrut{\relax\ifmmode\copy\medstrutbox\else\unhcopy\medstrutbox\fi}

\bigskip

\newdimen\leftmargin
\leftmargin=0.0truein
\newdimen\widesize
\widesize=5.0truein
\advance\widesize by \leftmargin       
\hfil\vbox{\tabskip=0pt\offinterlineskip
\def\tablerule{\noalign{\hrule}}

\halign to \widesize{\medstrut\vrule#\tabskip=0pt plus2truein

&\hfil\quad#\quad\hfil&\vrule#
&\hfil\quad#\quad\hfil&\vrule#
&\hfil\quad#\quad\hfil&\vrule#
&\hfil\quad#\quad\hfil&\vrule#
&\hfil\quad#\quad\hfil&\vrule#
&\hfil\quad#\quad\hfil&\vrule#
&\hfil\quad#\quad\hfil&\vrule#

\tabskip=0pt\cr\tablerule


&\quad&&$P$&&$R$&&$\neg P$&&$(\neg P)\lor R$&&$P\iff R$&&$v(A,\ass)$&\cr\tablerule
&(1)&&$T$&&$T$&&$F$&&$T$&&$T$&&$T$&\cr\tablerule
&(2)&&$T$&&$F$&&$F$&&$F$&&$F$&&$T$&\cr\tablerule
&(3)&&$F$&&$T$&&$T$&&$T$&&$F$&&$F$&\cr\tablerule
&(4)&&$F$&&$F$&&$T$&&$T$&&$T$&&$T$&\cr\tablerule
&\quad&&a&&b&&c&&d&&e&&f&\cr\tablerule}}
\parm
Now I'll go step-by-step through the above process using the general 
truth-tables on the previous page.\pars
(i) Columns (a) and (b) are written down with their usual patterns. \pars
(ii) Now go to column (a) only to calculate the truth-values in column (c). 
Note as mentioned previously there will be repetitions.\pars
(iii) Now calculate (d) for the $\lor$ connective using the truth-values
in columns (b) and (c).\pars
(iv) Next calculate (e) for the $\iff$ connective using columns (a) and (b).\pars
(v) Finally$,$ calculate (f)$,$ the value we want$,$ for the $\to$ connective using
columns (d) and (e).\parm

\hrule\smallskip\hrule\bigskip
\centerline{\bf EXERCISES 2.4}
\bigskip 
\noindent 1. First$,$ we assign the indicated truth values for the indicated atoms
$v(P) = T,\ v(Q) = F,\ v(R)=F$ and $v(S) = T.$ These values will yield one row 
of a truth-table$,$ one assignment $\ass.$ For this assignment$,$ find the truth-value 
for the indicated formula.  (Recall that $v(A,\ass)$ means the unique truth-
value for the formula $A$.)\pars
\+(a) $v((R\to(S\lor P)),\ass )$&(d) $v((((\neg S)\lor Q)\to(P\iff 
S)),\ass )$\cr  
\+(b) $v(((P\lor R)\iff (R \land(\neg S))),\ass )$&(e) $v((((P\lor (\neg Q))\lor 
R)\to ((\neg S)\land S)),\ass )$\cr
\+ (c) $v((S \iff (P\to ((\neg P)\lor S))),\ass )$&\cr \parm
\noindent 2. Construct complete truth tables for each of the following 
formula.\pars
\+(a) $(P\to(Q \to P))$&(c) $((P\to Q) \iff (P\lor (\neg Q)))$\cr
\+(b) $((P\lor Q)\iff (Q \lor P))$&(d) $((Q\land P)\to((Q \lor (\neg Q))\to 
(R\lor Q)))$\cr \parm 
\noindent 3. For each of the following determine whether or not the 
truth-value information given will yield a unique truth-value for the formula. 
State your conclusions. If the information is sufficient$,$ then give the
unique truth-value for the formula.\parm
\+(a) $(P \to Q)\to R,\ v(R) = T$&(d) $(R\to Q)\iff Q,\ v(R) = T$\cr
\+(b) $P\land (Q \to R),\ v(Q \to R) =F$&(e) $(P \to Q) \to R,\ v(Q) = F$\cr
\+(c) $(P \to Q)\to ((\neg Q) \to (\neg P))$&(f) $(P\lor (\neg P))\to R,\ 
v(R) = F$\cr 
\+ For (c)$,$ $v(Q) = T$&\cr\s\ss
\vfil\eject

\noindent {\bf 2.5 Valid Formula.}\parm
There may be something special about those formula that take the value $T$ for 
{\it any} assignment.\parm
\hrule\smallskip
{\bf Definition 2.5.1} (Valid formulas and contradictions). Let $A\in L.$ If 
for every assignment $\ass$ to the atoms in $A, \ v(A,\ass) = T,$  then $A$ is 
called a {{\it valid}} formula. If to every assignment $\ass$ to the atoms of 
$A,\ v(A,\ass) = F,$ then $A$ is called a {{\it (semantical)}} 
contradiction. If a formula $A$ is valid$,$ we use the notation $\mod A$ to 
indicate this fact.\pars
\hrule
\medskip
If we are given a formula in atomic form$,$ then a simple truth-table 
construction will determine whether or not it is a valid formula or a 
contradiction. Indeed$,$ $A$ is valid if and only if the column under the $A$ 
in its truth-table contains only $T$ in each position. A formula $A$ is a 
contradiction if and only if the column contains only $F$ in every position. 
From our definition$,$ we read the expression $\mod A$ ``$A$ is a valid 
formula.''  We read the notation $\not\mod A$  ``$A$ is not a valid formula.
Although a contradiction is not a valid formula$,$ there are infinitely many 
formula that are not valid AND not a contradiction.\parm
{\bf Example 2.5.1} Let $P,\ Q \in L_0.$\parm
\noindent (i) \hskip 0.25in $\mod P \to P$ \hskip 1.0in{\vbox{\offinterlineskip
\hrule
\halign{&\vrule#&
\strut\hfil#\hfil\cr
height2pt&\omit&\cr
&$\underline{\ P\hskip .75em\relax|\ P \to P\ }$&\cr
&$\underline{\hskip .2em\relax\ T\em {.1}\ |\ \ \ T \hskip 2.0em\relax}$&\cr
&$\ \em {-.1}F\em {.15}\ |\ \ \ T \em {2.0} $&\cr
height2pt&\omit&\cr}
\hrule}}
\par                                      
\bigskip
\noindent (ii) \hskip 0.25in $\mod P\to (Q \to P)$  \hskip 0.5in{\vbox{\offinterlineskip
\hrule
\halign{&\vrule#&
\strut\hfil#\hfil\cr
height2pt&\omit&\cr
&$\underline{\em {.4}P\em {.1}|\em {.1}\ Q\ |Q\to P|P\to (Q\to P)\em {.2}}$&\cr
&$\underline{\em {.4}T\ \em {-.1}|\em {.2}\ T\ \em {.00}|\em {.4}\ \ \ T\ \em {.7} |\em {.4} \ \ \ \ \ \ \ \ T\ \ \ \ \ \ \ \ }$&\cr
&$\underline{\em {.4}T\ \em {-.1}|\em {.2}\ F\ \em {.00}|\em {.4}\ \ T\ \ \em {.7}|\em {.4}\ \ \ \ \ \ \ \ T\ \ \ \ \ \ \ \ }$&\cr
&$\underline{\em {.4}F\ \em {-.1}|\em {.2}\ T\ \em {.00}|\em {.4}\ \  F\ \ \em {.7}|\em {.4}\ \ \ \ \ \ \ \ T\ \ \ \ \ \ \ \ }$&\cr
&$\em {-1.3}F\em {.3}|\ F\ \em {.1} | \ \ \em {.5}T\ \ \ \ \em {.1}|\em {.1}\ \ \ \ \ \ \ \ \em {.2}T \em {-.2}\ \ \ \ $&\cr 
height2pt&\omit&\cr}
\hrule}}\par                                        
\bigskip
\noindent (iii) \hskip 0.25in $\not\mod P \to Q$ \hskip 1.0in 
{\vbox{\offinterlineskip
\hrule
\halign{&\vrule#&
\strut\hfil#\hfil\cr
height2pt&\omit&\cr
&$\underline{\em {.4}P\em {.1}|\em {.2}\ Q\ |P\to Q\em {.7} }$&\cr
&$\underline{\em {.4}T\ \em {-.1}|\em {.2}\ T\ \em {.00}|\em {.4} \ \ \ \ T\  \ \ \ }$&\cr
&$\underline{\em {.4}T\ \em {-.1}|\em {.2}\ F\ \em {.00}|\em {.4} \ \ \ \ F\  \ \ \ }$&\cr
&$\underline{\em {.4}F\ \em {-.1}|\em {.2}\ T\ \em {.00}|\em {.4} \ \ \ \ T\ \ \ \ }$&\cr
&$\em {-1.3}F\em {.3}|\ F\ \em {.1} |\em {.1} \ \ \ \ \em {.2}T \em {-.2} $&\cr 
height2pt&\omit&\cr}
\hrule}}\par
\bigskip
\noindent (iv) \hskip 0.25in ${\rm contradiction}\ P \land (\neg P)$ \hskip 
0.20in {\vbox{\offinterlineskip
\hrule
\halign{&\vrule#&
\strut\hfil#\hfil\cr
height2pt&\omit&\cr
&$\underline{\em {.4}P\em {.4}|\ \neg P\em {.5} |P\land (\neg P)\em {.2} }$&\cr
&$\underline{\em {.5}T\ \em {-.1}|\ \ \ F\ \ | \em {.5} \ \ 
F\  \ \ \ \ \ \ }$&\cr
&${\em {.4}F\ \em {-.1}|\ \ \ T\ \ |\em {.4} \ \ 
F\  \ \ \ \ \ \ }$&\cr
height2pt&\omit&\cr}
\hrule}}    
\bigskip
We now begin the mathematical study of the validity concept.\parm
\hrule
\smallskip
{\it Throughout this 
text$,$ I will  ``prove'' various theorems in a way that is acceptable to the 
mathematical community.  Since the major purpose for this book is NOT
to produce trained mathematical logicians$,$ but$,$ rather$,$ to give the tools 
necessary to apply 
certain results from the  discipline to other areas$,$ I$,$ usually$,$ don't require 
a student to learn either these ``proofs'' or the methods used to obtain the 
proofs.}\pars
\hrule\parm
 {\bf Theorem 2.5.1} {\sl A formula $A\in L$ is valid if and only if 
$\neg A$ is a contradiction.}\pars
Proof. First$,$ notice that $\ass$ is an assignment for $A$ if and only if 
$\ass$ is a assignment  for $\neg A.$ Assume that $\mod A.$ Then for any 
$\ass$ for $A, \ v(A,\ass)=T.$ Consequently$,$ from the definition of $v,\ v(\neg 
A,\ass) = F.$ Since $\ass$ is any arbitrary assignment$,$ then
$v(\neg A, \ass) =F$ for all assignment $\ass.$\pars
Conversely$,$ let $\ass$ be any assignment to the atoms in $\neg A.$  Then 
$\ass$ is an assignment to the atoms in $A.$ Since $\neg A$ is a contradiction,
$v(\neg A,\ass) = F.$ From the truth-table (or the formal result in the 
appendix)$,$ it follows that $v(A,\ass) = T.$ Once again$,$ since $\ass$ is an 
arbitrary assignment$,$ this yields that $v(A,\ass) = T$ for all assignments 
and$,$ thus$,$ $\mod A.$\qed
Valid formula are important elements in our investigation of the propositional 
logic. It's natural to ask whether or not the validity of a formula is 
completely dependent upon its atomic components or its connectives? To answer 
this question$,$ we need to introduce the following {{\it substitution 
process}}.\parm
\hrule
\smallskip
{\bf Definition 2.5.2} (Atomic substitution process.) Let $A \in L$ be 
written in atomic form. Let $P_1,\ldots,P_n$ denote the atoms in $A.$
Now let $A_1,\ldots A_n$ be ANY (not necessarily distinct) members of $L.$ 
Define $A^*$ to be the result of substituting for each and every appearance
of an atom $P_i$ the corresponding formula $A_i$.\pars
\hrule\medskip
{\bf Theorem 2.5.2} {\sl Let $A\in L.$ If $\mod A,$ then $\mod A^*.$}\parm
Proof. Let $\ass$ be any assignment to the atoms in $A^*.$ In the step-by-step 
valuation process there is a level $L_m$ where the formula $A^*$ first 
appears. In the valuation process$,$ at level $L_m$ each constituent of $A^*$
takes on the value $T$ or $F$. Since the truth-value of $A^*$ only depends upon 
the connectives (they are independent of the symbols used for the formulas) 
and the truth-values of the $v(A_i, \ass)$ are but an 
assignment $\underline{b}$ that can be applied to the original atoms 
$P_1,\ldots,P_n,$ it follows that $v(A^*,\ass) = v(A,\underline{b}) = T.$ But$,$ 
$\ass$ is an arbitrary assignment for $A^*.$  Hence$,$ $\mod A^*$.\qed
{\bf Example 2.5.2} Assume that $P,\ Q \in L_0.$ Then we know that 
$\mod P\to (Q \to P).$ Now let $A,\ B\in L$ be any formula. Then 
$\mod A \to (B \to A).$ In particular$,$ letting $A = (P\to Q),\ B = (R\to S),$ 
where $P,\ Q,\ R,\ S \in L_0,$ then $\mod (P\to Q) \to ((R \to S)\to (P\to 
Q)).$\parm
Theorem 2.5.2 yields a simplification to the determination of a valid formula 
written with {\it some} connectives displayed. If you show that 
$v(A,\underline{c}) = T$
where you have created all of the possible assignments $\underline{c}$ not to the atoms of 
$A$ but only for the displayed components$,$ then $\mod A.$ (You think of the 
components as atoms.) Now the reason that this non-atomic method can be 
utilized follows from our previous results. Suppose that we have a list of 
components $A_1,\cdots, A_n$ and we substitute for 
each distinct component of $A$ a distinct atom in place of the components. 
Then any truth-value we give to the original components$,$ becomes an 
assignment $\ass$ for this newly created formula $A^\prime.$ Observe that 
using $A_1,\cdots, A_n$ it follows that $(A^\prime)^*= A.$ Now application of 
theorem 2.5.2 yields if $\mod A^\prime,$ then $\mod A.$ What this means is 
that whenever we wish to establish validity for a formula we may consider it 
written in component variables and make assignments only to these variables; if 
the last column is all Ts$,$ then the original formula is valid.\parm

\hrule\smallskip
{\bf WARNING}: We cannot use the simplified version to show that a formula is 
NOT valid. As a counter example$,$ let $A,\ B \in L.$ Then if we assume that 
$A,\ B$ behave like atoms and want to show that the composite formula
$A \to B$ is not valid and follow that procedure thinking it will show 
non-validity$,$ we would$,$ indeed$,$ have an $F$ at one row of the truth-table.
But if $A = B = P,$ which could be the case since $A,\ B$ are propositional 
language variables$,$ then we have a contradiction  since $\mod P \to P.$
Hence$,$ the formula can be considered as written in non-atomic form only if
it tests to be valid.\pars
\hrule
\medskip
It's interesting to note the close relation which exists between set-theory 
and logic. Assume that we interpret propositional symbols as names for sets 
which are subsets of an infinite set $X.$ Then interpret the conjunction as 
set-theoretic intersection (i.e. $\lint \land \rint:\  \cap$)$,$ the 
disjunction as set-theoretic union (i.e. $\lint \lor \rint:\  \cup$),
the negation as set-theoretic complementation with respect to $X$ 
(i.e. $\lint \neg \rint:\ X - 
{\rm or} \ X \setminus$)$,$ and the combination of validity with the 
biconditional to be set-theoretical equality (i.e $\lint \mod A \iff B\rint:\ 
A = B$). Now the valid formula $(P\land (( Q \lor R))) \iff ((P\land Q)\lor (P 
\land R))$ translates into the correct set-theoretic expression
$(P\cap (( Q \cup R))) = ((P\cap Q)\cup (P \cap R)).$ Now in this text we will 
NOT use the known set-theoretic facts to establish a valid formula even though some 
authors do so within the setting of the theory known as a Boolean algebra.
This idea would not be a circular approach since the 
logic used to determine these set-theoretic expressions is the metalogic of 
mathematics.\pars
In the next theorem$,$ we give$,$ FOR REFERENCE PURPOSES$,$ an important list of 
 formula each of which can be establish as valid by the simplified procedure
of using only language variables.\parm
{\bf Theorem 2.5.3} {\sl Let $A,\ B,\ C$ be any members of $L.$ Then 
the symbol $\mod$ can be place before each of the following formula.}
\+(1) $A \to (B \to A)$ &(8) $ B \to (A \lor B)$\cr 
\smallskip    
\+(2) $(A \to (B \to C))\to$& (9) $(A\to C) \to ((B\to C)\to $\cr 
\+ $\ \ \ \ ((A\to B)\to (A \to C)) $&$\ \ \ \ ((A\lor B)\to C))$\cr 
\smallskip    
\+(3) $(A\to B)\to$& (10) $(A \to B)\to ((A\to (\neg B)) \to (\neg A))$\cr
\smallskip 
\+$\ \ \ ((A\to (B\to C))\to (A \to C))$&\cr
\smallskip 
\+(4) $A \to (B \to (A\land B))$& (11) $(A \to B)\to$\cr 
\+& $\ \ \ \ ((B\to A) \to (A \iff B))$\cr\smallskip 
\+(5) $(A \land B) \to A$& (12) $(\neg(\neg A)) \to A$\cr\smallskip 
\+(6) $(A \land B)\to B $& (13) $(A \iff B)\to (A \to B)$\cr\smallskip 
\+(7) $A \to (A \lor B)$& (14) $(A \iff B)\to (B \to A)$\cr 
\line{\leaderfill}
\+(15) $A \to A   $& (17) $(A \to B) \to ((B\to C) \to$\cr 
\+& $\ \ \ \ (A \to C))$\cr\smallskip 
\+(16) $(A\to (B\to C))\iff $& (18) $(A \to(B\to C))\iff ((A\land B)\to C))$\cr 
\+ $\ \ \ \ (B \to (A \to C))$&\cr
\line{\leaderfill} 
\+ (19) $(\neg A)\to (A \to B)$& (20) $((\neg A) \to (\neg B))\iff (B \to A)$\cr 
\smallskip
\centerline{(21) $((\neg A) \to (\neg B))\to (B \to A)$}
\line{\leaderfill}
\+(22) $ A \iff A  $& (23) $(A\iff B) \iff (B \iff A)$\cr\smallskip
\centerline{(24) $((A\iff B) \land (B \iff C)) \to (A \iff C)$}
\line{\leaderfill} 
\+(25) $((A\land B)\land C) \iff (A\land(B\land C))$& (30) $((A\lor B)\lor 
C)\iff (A \lor (B\lor C))$\cr\smallskip 
\+(26) $(A \land B) \iff (B\land A)$& (31) $(A \lor B) \iff (B \lor 
A)$\cr\smallskip 
\+(27) $(A\land(B\lor C))\iff$& (32) $(A\lor (B\land 
C))\iff $\cr
\+$\ \ \ \ ((A\land B)\lor (A \land C))$&$\ \ \ \ ((A\lor B)\land (A \lor 
C))$\cr\smallskip 
\+(28) $(A\land A) \iff A$& (33) $(A\lor A) \iff A$\cr\smallskip 
\+(29) $(A\land (A\lor B))\iff A$& (34) $(A\lor (A \land B)) \iff 
A$\cr 
\line{\leaderfill}
\+(35) $(\neg (\neg A))\iff A$& (36) $\neg(A \land (\neg A))$\cr\smallskip 
\centerline{(37) $A\lor (\neg A)$}
\line{\leaderfill}
\+(38) $(\neg(A\lor B))\iff ((\neg A)\land (\neg B))$& (39) $(\neg(A \land B)) 
\iff ((\neg A)\lor (\neg B))$\cr\smallskip 
\centerline{(40) $(\neg (A \to B)) \iff (A \land ((\neg B))$}
\line{\leaderfill}
\+(41) $(A\lor B) \iff (\neg((\neg A) \land(\neg B)))$& (44) $(A \land B) \iff  
(\neg((\neg A) \lor(\neg B))) $\cr\smallskip 
\+(42) $(A \to B) \iff (\neg (A \land (\neg B)))$& (45) $(A \to B) \iff ((\neg 
A)\lor B)$\cr\smallskip 
\+(43) $(A \land B)\iff (\neg (A \to (\neg B)))$& (46) $(A\lor B) \iff ((\neg 
A)\to B)$\cr\smallskip
\centerline{ (47) $(A\iff B)\iff ((A\to B)\land (B \to A))$}\parm
\hrule
\smallskip
\hrule
\medskip 
\centerline{\bf EXERCISES 2.5}\medskip
\noindent (1) Use the truth-table method to establish that formula (1)$,$ (2)$,$ (21)$,$ 
(32)$,$ and (47) of theorem 2.4.3 are valid.\parm
\noindent (2) Determine by truth-table methods whether or not the following 
formula are contradictions.\pars
\+(a) $((\neg A)\lor (\neg B)) \iff$&(c) $(\neg (A \to B)) \iff ((\neg A)\lor 
B)$\cr
\+$\ \ \ \  (\neg((\neg A)\lor (\neg B)))$&\cr\smallskip
\+(b) $(\neg A) \to (A \lor B)$& (d) $(((A\lor (\neg B)) \land (\neg P))) \iff$ \cr  
\+& $\ \ \ \ (((\neg A)\lor B)\lor P)$\cr \s\ss
\noindent {\bf 2.6 Equivalent Formula.}\parm
As it will be seen in future sections$,$ one major objective is to investigate 
the classical human deductive processes and how these relate to the logic 
operator used in the solution to the General Grand Unification Problem. In order to 
accomplish this$,$ it has been discovered that humankind seems to believe that 
certain logical statements can be substituted completely for other logical 
statements without effecting the general ``logic'' behind an argument. Even 
though this fact will not be examined completely in this section$,$ we will 
begin its investigation.\pars
Throughout mathematics the most basic concept is the relation which is often 
called {{\it equality}}. In the foundations of mathematics$,$ there's a 
difference between equality$,$ which means that objects  are identically the 
same (i.e. they cannot be distinguished one from another by any property for
the collection  that contains them)$,$ and certain relations which {\it behave} 
like equality but that can be distinguished one from another and do not allow 
for substitution of one for another.  \parm
{\bf Example 2.6.1} When you first defined the rational numbers from the 
integers you were told the strange fact that 1/2$,$ 2/4$,$ 3/6 are  ``equal'' but 
they certainly appear to be composed of distinctly different symbols and would 
not be identical from our logical viewpoint. \parm
We are faced with two basic difficulties. In certain areas of mathematical 
logic$,$ it would be correct to consider the {\it symbols} 1/2$,$ 2/4$,$ 3/6 as 
names for a single unique object. On the other hand$,$ if we were studying the 
symbols themselves$,$ then 1/2$,$ 2/4$,$ 3/6 would not be consider as names for the 
same object but$,$ rather$,$ they are distinctly different symbols. These 
differences need not be defined specifically but can remain on the intuitive 
level for the moment. The are two types of  ``equality'' relations. One type 
simple behaves like equality but does not allow for substitution$,$ in general. 
But then we have another type that behaves like equality and does allow 
for substitution with {\it respect to certain properties.} This means that 
these two objects are identical as far as these properties are concerned. Or$,$ 
saying it another way$,$ a set of properties cannot distinguish between two such 
objects$,$ while another set of properties can distinguish them one from 
another. \parm
Recall  that a {{\it binary relation}}  ``on'' any set $X$ can be thought of as 
simply  a set of ordered pairs $(a,b)$ that$,$ from a symbol string viewpoint$,$ has 
an ordering. The first {\it coordinate} is the element you meet first in 
writing this symbol from left to right$,$ in this case the $a.$ The second 
coordinate is the next element you arrive at$,$ in this case the $b.$ Also recall 
that two ordered pairs are  identical (you can substitute one for another 
throughout your mathematical theory) if their first coordinates are 
identical and their second 
coordinates are identical (i.e. can not be distinguished one from another 
by the defining properties for the set in which they are contained.) The 
word  ``on'' means that the set of all first coordinates is $X$ and the set of 
all second coordinates is $X.$  
Now there are two ways of symbolizing such a binary relation$,$ either by 
writing it as a set of ordered pairs $R$ or by doing the following:\parm
\hrule
\smallskip
{\bf Definition 2.6.1} (Symbolizing ordered pairs.) Let $R$ be a nonempty set 
of ordered pairs. Then $(a,b) \in R$ if and only if $a \ R \ b.$ The 
expression $a \ R\ b$ is read  ``$a$ is $R$ related to $b$'' or similar types 
of expressions. \pars
\hrule
\medskip
In definition 2.6.1$,$ the reason the second form is used is that many times 
it's simply easier to write a binary relation's defining properties when the 
$a \ R\ b$ is used. It's this form we use to define a very significant binary 
relation that gives  the concept of behaving like ``equality.''\parm
\hrule
\smallskip 
{\bf Definition 2.6.2} (The equivalence relation.) A binary relation $\equiv$ 
\underbar{on} a set 
$X$ is called an {{\it equivalence relation}} if for each $a,\ b,\ c \in X$ 
it follows that\pars
(i)  $a \equiv a$ (Reflexive property).\pars
(ii) If $a \equiv b,$ then $b \equiv a.$ (Symmetric property). \pars
(iii) If $a \equiv b$ and $b \equiv c,$ then $a \equiv c.$ (Transitive 
property).\pars
\hrule\medskip
Now when we let $X=L,$ then the only identity or equality we use is the 
intuitive identity. Recall that this means that two symbol string are 
recognized as congruent geometric configurations or are intuitively 
similar strings of symbols. This would yield a trivial equivalence relation.  
As a set of ordered pairs$,$ an identity relation is $\{(a,a)\mid a \in X\},$ 
which is (i) in definition 2.6.2. Parts (ii)$,$ (iii) also hold for this 
identity relation. \pars
For the next theorem$,$ please recall that the validity of a formula $A$ does 
not depend upon an assignment $\ass$ that contains MORE members than the 
number of atoms contained in $A.$ Such an assignment is used by restricting 
the Ts and Fs to those atoms that are in $A.$ \parm
{\bf Theorem 2.6.1}\pars
{\sl Let $A,\ B  \in L$ and $\ass$ an arbitrary assignment to the atoms that 
are in  $A$ 
and $B.$\pars
(i)  Then $v(A\iff B,\ass) = T$ if and only if $v(A,\ass) = 
v(B,\ass).$\pars
(ii) $\mod A \iff B$ if and only if for any assignment $\ass$ to the atoms 
that are in $A$ and $B,$  $v(A,\ass) = v(B,\ass).$}\pars
Proof. Let $A,\ B \in L.$\pars
(i) Then let the ${\rm size}(A\iff B)=n\geq 1.$ Then 
$A,\ B \in L_{n-1}.$ This result now follows from the general truth-tables on 
page 22. \pars
(ii) Assume that $\mod A \iff B$ and let $\ass$ be an arbitrary 
assignment to the atoms that are contained in $A$ and $B.$ Then
$v(A\iff B,\ass) = T$ if and only if $v(A,\ass) = v(B,\ass)$ from part 
(i). Conversely$,$ assume that $\ass$ is an assignment for the atoms 
in $A$ and $B.$ Then $\ass$  also determines an assignment for $A$ and $B$ 
separately. Since$,$ $v(A,\ass) = v(B,\ass)$ then from (i)$,$ it follows that 
$v(A \iff B,\ass) = T.$ But $\ass$ is arbitrary$,$ hence$,$ $\mod A \iff B.$
\qed
\hrule
\smallskip
{\bf Definition 2.6.3} (The logical equivalence relation $\equiv$.) 
Let $A,\ B \in L.$ Then define $A\equiv B$ iff $\mod A\iff B.$\pars
\hrule
\medskip
Notice that definition 2.6.3 is easily remembered by simply dropping the 
$\mod$ 
and replacing $\iff$ with $\equiv.$ Before we proceed to the study of 
equivalent propositional formulas$,$ I'll 
anticipate a question that almost always arises after the next few theorems.
What is so important about equivalent formula? When we study the actual 
process of logical deduction$,$ you'll find out that within any classical 
propositional deduction  a formula $A$ can be substituted for an equivalent 
formula and this will in no way affect the deductive conclusions. What it may 
do is to present a more easily followed logical process. 
This is exactly what  happens if one truly
 wants to understand real world logical arguments. For example$,$ take a look 
at theorem 2.5.3 parts (29) and (34) and notice how logical arguments can be 
made more complex$,$ unnecessarily$,$ by adding some rather complex statements$,$ 
statements that include totally worthless additional statements such as any
additional statement $B$ that might be selected simply to confuse the 
reader. \parm 
{\bf Theorem 2.6.2} 
{\sl The relation $\equiv$ is an equivalence relation defined on $L.$}\pars 
Proof. Let $A,\ B,\ C \in L.$ From the list of valid formula that appear in theorem 2.5.3$,$ formula 
(22) yields that for each $A \in L,$ $\mod A\iff A.$ Hence$,$ $A \equiv A.$\pars
Now let $A, \ B \in L$ and assume that $A \equiv B.$ Then for any assignment 
$\ass,$ $\mod A\iff B$ implies that $v(A,\ass) = v(B,\ass)$ from theorem 
2.6.1. 
Since the 
equality means identically the same symbol (the only equality for our 
language)$,$ it follows that $v(B,\ass) = v(A,\ass).$ Consequently$,$ 
$B \equiv A.$ \pars
Now assume that $A\equiv B,\ B \equiv C.$ Hence$,$ $\mod A \iff B,\ \mod B 
\iff C.$ Again by application of theorem 2.6.1 and the definition of $\mod,$
it follows that $\mod A \iff C.$ Consequently$,$ $A \equiv C$ and $\equiv$ is an 
equivalence relation.\qed
We now come to the very important substitution theorem$,$ especially when 
deduction is concerned. It shows that substitution is allowed throughout the 
language $L$ and yields a powerful result.\parm
\hrule\smallskip
 {\bf Definition 2.6.4} (Substitution of formula). Let $C\in L$ be any formula
and $A$ a formula which is a composite element in $C.$ Then $A$ is 
called a {{\it subformula}}. Let $C_A$ denote the formula $C$ with the 
subformula $A$ specifically identified. Then the substitution process states 
that if you substitute $B$ for $A$ then you obtain the $C_B,$ where you have 
substituted for the specific formula $A$ in $C$ the formula $B.$\pars
\hrule\medskip
{\bf Example 2.6.2} Suppose that $C = (((\neg P) \lor Q)\to ((P \lor S) \iff 
S)).$ Let $A = (P \lor S)$ and consider $C_A.$  Now let $B = (S \land (\neg P)).$ Then
$C_B =  (((\neg P) \lor Q)\to (\underline{(S \land (\neg P))} \iff S)),$ 
where the substituted formula is identified by the underline.\parm
{\bf Theorem 2.6.3} {\sl If $A,\ B,\ C \in L$ and $A \equiv B,$  then
$C_A \equiv C_B.$}\pars
Proof. Let $A \equiv B.$ Then $\mod A \iff B$. Let $\ass$ be any assignment to 
the atoms in $C_A$ and $C_B.$ Then $\ass$ may be considered as an assignment 
for $C_A,$ $C_B,$ $A,\ B.$ Let ${\rm size}(C_A) = n.$ In the truth-table 
calculation 
process (or formal process) for $v(C_A,\ass)$ there is a step when we 
(first) calculate $v(A,\ass).$ Let ${\rm size}(A) = k\leq n.$ If ${\rm size} 
(A) = n,$ then $A = C$ and $C_B = B$ and we have nothing to prove. Assume that 
$k < n.$ Then the calculation of $v(C_A,\ass)$ at this specific level only 
involves the calculation of $v(A,\ass)$ and the other components and other 
connectives not in $A.$ The same argument for $C_B$ shows that calculation for 
$C_B$ at this level uses the value $v(B,\ass)$ and any other components and 
other connectives in $C$ {\it which are all the same} as in $C_A.$ However$,$ 
since $A  \equiv B,$ theorem 2.6.1 yields $v(A,\ass) = v(B,\ass).$ Of course$,$ 
the truth-values for the other components in $C_A$ that are the same as the 
other components in $C_B$ are equal since these components are the exact same 
formula. Consequently$,$ since the computation  of the truth-value for $C_A$ and 
$C_B$ now continue from this step and all the other connective are the same 
from this step on$,$ then $C_A$ and $C_B$ would have the same truth-value.
Hence $v(C_A,\ass) = v(C_B,\ass).$ But $\ass$ is arbitrary; hence,
 $\mod C_A \iff C_B.$ Thus $C_A \equiv C_B.$ \pars
{\bf Corollary 2.6.3.1} {\sl If $\mod A \iff B$ and $\mod C_A,$ then
$\mod C_B.$}\pars
Proof. From the above theorem $\mod C_A  \iff C_B,$ it follows that for any 
assignment $\ass$ to the atoms in $C_A$ and $C_B$$,$ $v(C_A,\ass) = 
v(C_B,\ass).$ However$,$ $v(C_A,\ass) = T.$ Moreover$,$ all of the assignments for 
the atoms in $C_A$ and $C_B$ will yield all of the assignments 
$\underline{b}$ for the atoms in $C_B$ as previously mentioned. Hence$,$ if 
$\underline{b}$ is any assignment for the atoms in $C_B,$ then 
$v(C_B,\underline{b})=T$ and the result follows.\qed
[Note: It follows easily that if $C,\ A,\ B$ are written in formula variables and$,$ hence$,$ represent a hidden atomic structure$,$ then Theorem 2.6.3 and its corollary will also hold in this case.]\pars
The next result seems to fit into this section. It's importance cannot be 
over-emphasized since it mirrors our major rule for deduction. For this 
reason$,$ it's sometimes called the  {{\it semantical modus ponens}} 
result.\parm
{\bf Theorem 2.6.4} {\sl If $\mod A$ and $\mod A \to B,$ then $\mod B.$}\pars
Proof. Suppose that $A,\ B\in L.$ Let $\mod A,\ \mod A \to B$ and $\ass$ be any 
assignment to the atoms in $A,\ B.$ Then $v(A,\ass) = T= v(A\to B, \ass).$
Thus $v(B,\ass) =T.$  Since $\ass$ is any assignment$,$ then$,$ as used 
previously$,$ using the set of all assignments for $A,\ B,$  we also obtain all of 
the assignments for $B.$ Hence $v(B,\underline{b})=T$ for any assignment for 
$B$s atoms and the result follows.\qed
\medskip
\hrule
\smallskip
\hrule
\medskip
\centerline{\bf EXERCISES 2.6} \parm
\baselineskip =12pt
\noindent 1. There is a very important property that shows how equivalence 
relations can carve up a set into important pieces$,$ where each piece contains 
only equivalent elements. Let $\equiv$ be {\it any} equivalence relation 
defined on  
the non-empty set $X.$ This equivalence relation can be used to 
define a subset of $X.$ For every $x \in X,$ this subset is denoted by $[x].$ 
Now to define this very special and important set. For each 
$x\in X,$ let $[x] = \{y \mid y \equiv x \}.$ Thus if you look at one of these 
sets$,$ say $[a],$ and you take any two members$,$ say $b,  \ c \in [a],$ it 
follows that $b \equiv c.$ Now see if you can establish by a simple logical 
argument using the properties (i)$,$ (ii)$,$ (iii) of definition 2.6.2 that:\parm
\baselineskip =14pt
(A) If there is some $z \in X$ such that $z \in [x]$ and $z \in [y],$ then 
$[x] =[y].$ (This equality is set equality$,$ which means that $[x]$ is a subset 
of $[y]$ and $[y]$ is a subset of $[x].$)\pars
(B) If $x \in X$$,$ then there exists some $y \in  X$ such that $x \in 
[y].$\pars\baselineskip=12pt
\noindent 2. (A) Of course$,$ there are usually many interesting binary relations defined 
\underbar{on} 
a non-empty set $X.$ Suppose that you take {\it any} binary relation $B$ 
defined on $X$ and you emulate the definition we have used for $[x].$ Suppose 
that you let $(x)= \{ y \mid y\  B\  x\}.$ Now what if properties (A) and (B) 
and the reflexive property (i) of definition 2.6.2 hold true for this 
relation. Try and give a simple argument that shows in this case that
 $B$ is$,$ indeed$,$ 
an equivalence relation. \pars
(B) In (A) of this problem$,$ we required that $B$ be reflexive. 
Maybe we can do without this 
additional requirement. Try and show that this requirement is necessary by 
looking at a set that contains two and only two elements $\{a,b\}$ and find 
a set of 
ordered pairs$,$ using one or both of its members$,$ that yields a binary 
relation on $\{a,b\}$ such that (A) and (B) of problem 1 
hold but (i) of definition 2.6.2$,$ the reflexive property$,$ does not hold.
If you can find one$,$ this is an absolute counter-example that establishes that 
the reflexive property is necessary.\pars\baselineskip=14pt
\noindent 3.  One of the more important properties of $\equiv$ is the 
transitive property (iii). For example$,$ if $C_A \equiv C_B$ and $C_B \equiv 
C_D,$ then $C_A \equiv C_D.$ Now this can be applied over and over again a 
finite number of times. Notice what can be done by application of theorem 
2.5.3 parts (26) and (31). Suppose that you have a formula $C$ containing a 
subformula $(A \lor B),$ where $(A \lor B) \equiv (B \lor A)$ or  $(A \land 
B)$$,$ where $(A\land B) \equiv (B \land 
A).$ Letting $H = (A \lor B)$ and $K = (B \lor A),$  then $C_H \equiv C_K.$
Now recall that$,$ in general$,$ on a set $X$ where an equality is defined$,$ 
an operation$,$ say $\Delta,$ is commutative if for each $x,\ y \in X,$ it 
follows that $x\  \Delta \ y = y\  \Delta \ x.$ Thus for the operation and the 
(not equality) equivalence relation $\equiv$ the same type of commutative law 
for $\lor$ holds. In the following$,$ using if necessary the transitive 
property$,$ establish that (A)$,$ (B)$,$ (C)$,$ (D) (E) hold by stating the particular
valid formula(s) from theorem 2.5.3 that need to be applied. \parm
(A) Given $C_D$ where $D = (A \lor (B \lor C)),$ then $C_D \equiv C_E,$ where 
$E = ((A \lor B) \lor C).$ This would be the  {\it associative law} for 
$\lor.$ Now establish the associative law for $\land.$\pars
(B) Given $C_H,$ where $H = (A\lor B).$  Show that $C_H \equiv C_K$ where 
$K$ only contains the $\neg$ and $\to$ connectives.\pars
(C) Given $C_H,$ where $H = (A\land B).$  Show that $C_H \equiv C_K$ where 
$K$ only contains the $\neg$ and $\to$ connectives.\pars     
(D) Given $C_H,$ where $H = (A\iff B).$  Show that $C_H \equiv C_K$ where 
$K$ only contains the $\neg$ and $\to$ connectives.\pars     
(E) Given $C_H,$ where $H = (\neg(\neg(\neg \cdots A\cdots ))).$  (i.e. the 
formula has  ``n'' $\neg$ to its left.) Show that $C_H \equiv C_K$ where 
$K$ only contains one and only one $\neg$ or no $\neg.$\pars 
\noindent 4. Using the results from problem (3)$,$ and using$,$ if necessary a 
finite number of transitive applications$,$ re-write each of the following 
formula in terms of an equivalent formula that contains only the $\neg$ 
and $\to$ connectives. (The formulas are written in the allowed slightly 
simplified form.) \parm
\+ (a) $(\neg(A \lor B)) \to (B \land C)$& (c) $((A \lor B)\lor C) \land D.$\cr
\+ (b) $A \iff (B \iff C)$& (d) $\neg ((A \lor (\neg B)) \lor (\neg(\neg 
D)))$\cr  \pars
\noindent 5. In mathematics$,$ it's the usual practice to try and weaken 
hypotheses as much as possible and still establish the same conclusion.
Consider corollary 2.6.3.1. I wonder if this can be weakened to the theorem
``if $\mod A  \to B$ and $\mod C_A,$ then $\mod C_B$? Try to find an 
explicit formula such that $\mod A  \to B$ and $\mod C_A,$ 
but $\not\mod C_B.$ If you can find such a formula in $L,$ then this would 
mean that the hypotheses cannot be weakened to $\mod A  \to B$ and $\mod C_A.$\pars
\ss
\noindent {\bf 2.7 The Denial$,$ Full Disjunctive Normal Form$,$ Logic 
Circuits.}\parm
{\it From this point on in this chapter$,$ when we use language symbols such as 
P$,$ Q$,$ R$,$ S$,$ A$,$ B$,$ C$,$ D$,$ E$,$ F and the like$,$ it will always be assumed that they 
are members of L. This will eliminate repeating this over and over 
again.}\pars
As was done in the exercises at the end of the last section$,$ checking theorem
2.5.3$,$ we see that $(A\to B) \equiv ((\neg A)\lor B).$ 
Further$,$ $(A \iff B) \equiv ((A \land B)\lor ((\neg A)\land (\neg B))).$ 
Consequently$,$ for any $C,$ we can take every subformula that uses connective 
$\to$ and $\iff,$ we can use the substitution process and$,$ hence$,$ express $C$ in a 
equivalent form $D$ where in $D$ only the connectives $\neg,\ \lor,$ and 
$\land$ appear. Obviously$,$ if $D$ is so expressed with at most these three 
connectives$,$ then $\neg C \equiv \neg D$ and $\neg D$ is also expressed with at 
most these three connectives. Further$,$ by use of the valid formula theorem,
any formula with more that one $\neg$ immediately to the left (e.g. 
$(\neg(\neg(\neg A)))$) is equivalent to either a formula for no $\neg$ 
immediately to the left$,$ or at the most just one $\neg $ immediately on the 
left. Since $\neg(A \lor B) \equiv ((\neg A) \land (\neg B))$ and 
$\neg(A \land B) \equiv ((\neg A)\lor (\neg B))$ then$,$ applying the above 
equivalences$,$ we can express any formula $C$ in an equivalent form $D$ with the 
following properties.\parm
(i) $D$ is expressed entirely in atoms.\pars
(ii) Every connective in $D$ is either $\neg,\ \lor, \land.$\pars
(iii) And$,$ when they appear$,$ only single $\neg$'s appear immediately to the 
left of atoms.\parm 
\hrule
\smallskip
{\bf Definition 2.7.1} (The denial.) Suppose that $A$ is in the form $D$ with 
properties (i)$,$ (ii)$,$ (iii) above. Then the {{\it denial}} $A_d$ of $A$ is 
the formula obtained by\parm 
(a) dropping the $\neg$ that appears before any atom.\pars
(b) Placing a $\neg$ before any atom that did not have such a connective 
immediately to the left.\pars
(c) Replacing each $\lor$ with $\land.$\pars
(d) Replacing each $\land$ with $\lor$.\pars
(e) Adjust the parentheses to make a correct language formula.  \pars
\hrule
\medskip
{\bf Example 2.7.1} Let $A = ((\neg P) \lor (\neg Q)) \land (R \land (\neg S))  
.$ Then $A_d = (P \land  Q) \lor ((\neg R) \lor  S).$ Notice were the 
parentheses have been removed and added in this example.\parm
{\bf Theorem 2.7.1} {\sl Let $A$ be a formula containing only atoms$,$ the 
connective $\neg$ appearing only to the immediate left of atoms$,$ if at all$,$ 
and any other connectives are $\land$ and/or $\lor.$ Then $ \neg A\equiv  
A_d.$}\pars
Proof. (This is the sort of thing where mathematicians seem to be proving the 
obvious since we have demonstrated a way to create the equivalent formula. 
The proof is a formalization  of this process for ANY formula based upon one of 
the most empirically consistent processes known to the mathematical 
community.  The process is called induction on the natural numbers$,$ 
in this case the 
unique natural number we call the size of a formula.)\par
 First$,$ we must show the 
theorem holds true for a formula of size 0. So$,$ let ${\rm size}(A) = 0.$ Then 
$A =P \in L_0$ and is a single atom. Then $A_d = \neg A.$ Further$,$ $\neg A = 
\neg P.$ We know that for any formula $D$$,$ $D \equiv D.$ Hence$,$ $\neg A \equiv 
A_d$ for this case. \pars
Now (strong) induction proofs are usually done by assuming 
that the theorem holds for all $A$ such that ${\rm size}(A) \leq n,$ where $n 
>0.$ Then using this last statement$,$ it is  shown that one method will yield 
the theorem's conclusions for ${\rm size}(A) = n+1.$ However$,$ this specific 
procedure may not work$,$ yet$,$ since there is not one simple method unless we start 
at $n >1.$ Thus let ${\rm size}(A) = 1.$ Then there are three possible forms.
(i) Let $A = \neg P.$ Then$,$ from theorem 2.5.3$,$ it follows that $\neg A =
\neg (\neg P) \equiv P.$ Now notice that if $A = B,$ then $\mod A \iff B.$ 
Further$,$ $A_d= P$ implies that $\neg A \equiv A_d.$   For the cases where  $A = 
P\land Q$ or $A = P \lor Q,$ the result follows from theorem 2.5.3$,$ parts (38)$,$ 
(39).\pars
Now assume that the theorem holds for a formula $A$ of size $r$ such that
$1 \leq r \leq n.$ Let ${\rm size}(A) = n+1.$  Also we make the following observation. Because of the 
structure of the formula $A,$ $A$ cannot be of the form $\neg B$ where 
${\rm size}(B) \geq 1.$ Indeed$,$ if the formula has a $\neg$ and a $\lor$ or an 
$\land,$ then ${\rm size}(A) >1.$  Consequently$,$ there will always be two 
and only two cases.\pars
Case (a). Let $A = B \lor C.$ Consider $\neg A = \neg (B \lor C).$ From the 
above discussion$,$ observe that ${\rm size}(\neg B) \leq n.$ Hence$,$ by 
the induction hypothesis$,$ $\neg 
B\equiv B_d$ and$,$ in like manner$,$ $\neg C\equiv C_d.$ Theorem 2.5.3$,$ shows 
that $\neg (B \lor C) \equiv (\neg B) \land (\neg C).$ Since equal formula are 
equivalent$,$ then substitution yields$,$ $\neg (A \lor C) \equiv B_d \land C_d.$ 
Again$,$ since $A_d = B_d \land C_d,$ and equal formula are equivalent$,$ 
substitution yields $\neg A \equiv A_d.$\pars
Case (b). Let $A = B \land C.$ This follows as in case (a) from theorem 
2.5.3. Thus the theorem holds for ${\rm size}(A) = n+1.$ From the induction 
principle$,$ the theorem holds for any (specially) constructed formula since 
every 
such formula has a unique size which is a natural number.\qed
On page 28$,$ I mentioned how you could take certain valid formula and find a 
correct set-theoretic formula. The same can be done with the denial of the 
special form $A.$ If you have any knowledge in this area$,$ the $\neg$ is 
interpreted as set-complementation with respect to the universe. We can get another one of D'Morgan's Laws$,$ for 
complementation using $\neg A \equiv A_d.$\pars
Since any formula is equivalent$,$ by theorem 2.5.3 part (29)$,$ to infinitely many 
different formula$,$ it might seen not to intelligent to ask whether or not a 
member of $L$ is equivalent to a formula that is unique in some special way? 
Even if this is true$,$ is this uniqueness useful?  So$,$ the basic problem is to  
define the concept of a {{ \it unique equivalent form}} for any give 
formula. \pars
Well$,$ suppose that $A$ is a contradiction and $P$ is any atom. Then 
$A \equiv (P\land (\neg P)).$ And if $B$ is any valid formula$,$ then
$B \equiv (P \lor (\neg P)).$ Hence$,$ maybe the concept of a unique equivalent 
form is not so easily answered. But$,$ we try anyway.\pars
Let $A$ be a formula that is in atomic form and contains only the atoms 
$P_1,\ldots,P_n.$ We show that there is a formula equivalent to $A$ that uses 
these are only these atoms and that does have an {\it almost} unique form.
The formula we construct is called the {{\bf full disjunctive normal form}} 
and rather than put this into a big definition$,$ I'll slowly described the 
process by the truth-table procedure. \pars
Let the distinct atoms $P_1,\ldots,P_n$ be at the top of a truth-table and in 
the first  ``n'' columns. Now observe that when we calculate the truth-values
for a formula $A \land (B \land C)$ we have also calculated the truth-values
for the formula $(A\land B)\land C)$ since not only are these formula 
equivalent$,$ but they use the same formula 
$A, \ B, \ C,$  the exact same number and type of connective$,$ in the exact 
same places. Indeed$,$ only the parentheses are in different places. For this 
reason$,$ we often {\it drop the parentheses} \underbar{in this case} when we 
are calculating the truth value for $A \land (B \land C).$ Now for the 
procedure. Consider any row $k$$,$ where $1\leq k \leq 2^n.$\parm
(a) For each $T$ that appears in that row under the atom $P_i,$ write done 
 the symbol $P_i.$ \parm     
(b) For each $F$ that appears in that row under the atom $P_j,$ write done 
the symbol $(\neg P_j).$\parm  
(c) Continue this process until you have used (once) each truth value that 
appears in row $k$ making sure you have written down all these symbols in a 
single row and have left spaces between them. \parm
(d) If there is more than one symbol$,$ then between each symbol put 
a $\land$ and insert the outer most 
parentheses. \parm
(e) The result obtained is called the {{\bf fundamental conjunction.}}\parm
{\bf Example 2.7.2} Suppose that the $k$th row of our truth table looks 
like\pars
\newbox\medstrutbox
\setbox\medstrutbox=\hbox{\vrule height14.5pt depth9.5pt width0pt}
\def\medstrut{\relax\ifmmode\copy\medstrutbox\else\unhcopy\medstrutbox\fi}

\bigskip

\newdimen\leftmargin
\leftmargin=0.0truein
\newdimen\widesize
\widesize=3.0truein
\advance\widesize by \leftmargin       
\hfil\vbox{\tabskip=0pt\offinterlineskip
\def\tablerule{\noalign{\hrule}}

\halign to \widesize{\medstrut\vrule#\tabskip=0pt plus2truein

&\hfil\quad#\quad\hfil&\vrule#
&\hfil\quad#\quad\hfil&\vrule#
&\hfil\quad#\quad\hfil&\vrule#

\tabskip=0pt\cr\tablerule


&$P_1$&&$P_2$&&$P_3$&\cr\tablerule
&$T$&&$F$&&$F$&\cr\tablerule
}}
\parm
\noindent Then$,$ applying (a)$,$ (b)$,$ (c)$,$ we write down\parm
\centerline{$P_1\ \ \ \ \ (\neg P_2)\ \ \ \ \ (\neg P_3)$}
\noindent Next we do as we are told in (d). This yields\parm
\centerline{$(P_1\land (\neg P_2)\land (\neg P_3)).$} \pars
Now the assignment for the $k$ row is $\ass = (T,F,F).$ Notice the important 
fact that $v((P_1\land (\neg P_2)\land (\neg P_3)),\ass) = T.$ What we have 
done to remove any possibility that the truth-value would be $F.$ But$,$ also 
it's significant$,$ that if we took any other {\it distinctly different} 
assignment $\underline{b},$ then $v ((P_1\land (\neg P_2)\land (\neg 
P_3)),\underline{b}) =F.$ These observed facts about this one example can be 
generally established.\parm

{\bf Theorem 2.7.2} {\sl Let $k$ be any row of a truth-table for the distinct  
set of atoms $P_1,\ldots,P_n.$ Let $\ass$ be the assignment that this row 
represents. For each $a_i=T,$ write down $P_i$. For each $a_j= F,$ write down 
$(\neg P_j).$ Let $A$ be the formula obtained by placing conjunctions between 
each pair of formula if there exists more than one such formula. Then 
$v(A,\ass) = T,$ and for any other distinct assignment $\underline{b},$ 
$v(A,\underline{b}) = F.$}\pars
Proof. See theorems 1.6 and 1.7 on pages 13$,$ 14 of the text  ``Boolean 
Algebras  and Switching Circuits,'' by Elliott Mendelson$,$ Schaum's Outline 
Series$,$ McGraw Hill$,$ 1970.\parm
Now for any formula $C$ composed of atoms $P_1,\ldots, P_m$ and which is 
{\it not a contradiction} there will be$,$ at the least$,$ one row assignment 
$\ass$ such that $v(C,\ass)=T.$ We now construct the formula that is 
equivalent to $C$ that we has an almost unique form.\parm
{\bf Definition 2.7.1} (Full disjunctive normal form.) Let $C$ not be a 
contradiction. \parm
(a) Take every row $k$ for which $v(C,\ass) =T.$ \parm
(b) Construct the fundamental conjunction for each such row.\parm
(c) Write down all such fundamental conjunctions and between each pair$,$ if 
any$,$ place a $\lor.$\parm
(d) The result of the construction (a)$,$ (b)$,$ (c) is called the 
{\bf full disjunctive normal form} for $C.$ This can be denoted by 
fdnf($C$). \parm

{\bf Example 2.7.2} Suppose that $C = P\iff (Q \lor R).$ Now 
consider the truth-table.\parm\vfil\eject

\newbox\medstrutbox
\setbox\medstrutbox=\hbox{\vrule height14.5pt depth9.5pt width0pt}
\def\medstrut{\relax\ifmmode\copy\medstrutbox\else\unhcopy\medstrutbox\fi}

\bigskip

\newdimen\leftmargin
\leftmargin=0.0truein
\newdimen\widesize
\widesize=2.5truein
\advance\widesize by \leftmargin       
\hfil\vbox{\tabskip=0pt\offinterlineskip
\def\tablerule{\noalign{\hrule}}

\halign to \widesize{\medstrut\vrule#\tabskip=0pt plus2truein

&\hfil\quad#\quad\hfil&\vrule#
&\hfil\quad#\quad\hfil&\vrule#
&\hfil\quad#\quad\hfil&\vrule#
&\hfil\quad#\quad\hfil&\vrule#

\tabskip=0pt\cr\tablerule


&$P$&&$Q$&&$R$&&$\ \ C\ \ $&\cr\tablerule
&$T$&&$T$&&$T$&&$\ \ T\ \ $&\cr\tablerule
&$T$&&$T$&&$F$&&$\ \ T\ \ $&\cr\tablerule
&$T$&&$F$&&$T$&&$\ \ T\ \ $&\cr\tablerule
&$T$&&$F$&&$F$&&$\ \ F\ \ $&\cr\tablerule
&$F$&&$T$&&$T$&&$\ \ F\ \ $&\cr\tablerule
&$F$&&$T$&&$F$&&$\ \ F\ \ $&\cr\tablerule
&$F$&&$F$&&$T$&&$\ \ F\ \ $&\cr\tablerule
&$F$&&$F$&&$F$&&$\ \ T\ \ $&\cr\tablerule}}
\parm
Now we construct fundamental conjunctions for rows 1$,$ 2$,$ 3$,$ 8. This yields
$(P\land Q \land R),\ (P\land Q \land (\neg R)),\ (P\land(\neg Q) \land R),
\  ((\neg P)\land(\neg Q)\land(\neg R)).$ Notice that I have not 
included for the $\land$ the internal parentheses since they can be placed 
about any two expressions. \pars
\noindent Then putting these together we have\pars
\centerline{fdnf$(C) = (P\land Q \land R)\lor (P\land Q \land (\neg R))\lor$}
\centerline{$(P\land(\neg Q) \land R)\lor ((\neg P)\land(\neg Q)\land(\neg 
R)).$}\parm
{\bf Theorem 2.7.3} {\sl Let non-contradiction $C\in L.$ Then fdnf$(C) \equiv C.$}\pars
Proof. See the reference given in the proof of theorem 2.7.2.\qed The is no example 2.7.3.
The fdnf is unique in the following sense.\parm
{\bf Theorem 2.7.4}. {\sl  If formula $B,\ C$ have the same 
number of atoms which we can always denote by the same symbols 
$P_1,\ldots,P_m.$ If fdnf$(B)$ and fdnf$(C)$ have the same set of fundamental 
conjunctions  except for a change in order of the individual conjuncts$,$ 
then  fdnf$(B)\equiv$ fdnf$(C).$}\pars
Proof. I'm sure you could prove this from our validity theorem$,$ substitution$,$ 
and the properties of $\equiv.$\qed\parm
Usually$,$ when elementary concepts in logic are investigated$,$ the subject of 
computers is often of interest. The reason for this is that$,$ technically$,$ 
computers perform only a very few basic underlying processes all related to 
propositional logic. So$,$ for a moment$,$ let's look at some of the basic {{\it 
logic circuits,}} many of which you can construct. Such circuits are extensions 
of the (switching) relay circuits where we simply suppress the actual device 
the functions in the fashion diagrammed on pages 17$,$ 18. Since any formula $A$ 
is equivalent to its fdnf$,$ its the fdnf that's used as a bases for these 
elementary logic circuits. An important procedure within complex logic 
circuits is simplification or minimizing techniques. Simplification does not 
necessarily mean fewer devices. The term simplification includes the concept of 
something being more easily constructed and/or less expensive. We will have not 
interest is such simplification processes.\pars
Looking back at pages 17$,$ we have three logic devices. The 
{{\it or-gate}} \orgate$,$ diagram (i); the {{\it and-gate}} \andgate$,$ 
diagram (ii)$,$ and an 
{{\it inverter}} {\inverter}$,$ which is the combination of the two diagrams above diagram 
(i). The inverter behaves as follows: when current goes in one end$,$ it opens 
and no current leaves the exit wire$,$ the output. 
But$,$ when input is no current$,$ then output is current. (We need not use 
current$,$ of course. Any two valued physical event can be used.) In the 
following diagrams the current direction and what is called the {{\it 
logical flow}} is indicated by the line with an arrow. If there is no flow$,$ 
then no arrows appear. Now each gate has at least two inputs and one output in 
our diagrams except the inverter.\pars
It's very$,$ very easy to understand the behavior of the or-gates and the 
and-gates. They have the same current flow properties as a corresponding 
truth-table$,$ where $T$ means current flows and $F$ means no current flows. The 
basic theorem used for all logic circuits is below.\parm
{\bf Theorem 2.7.5} {\sl If $A \equiv B,$ then any logic circuit that 
corresponds to $A$ can be substituted for any logical circuit that corresponds 
to $B$.}\pars
 Proof. Left to you. \parm
{\bf Example 2.7.4} Below are diagrams for two logic circuits. For the first 
circuit$,$ note that if no circuit flows into lines $A$ and $B,$ then there is a 
current flowing out the $C$ line. In current flows in the $A$ and $B$ line 
from left to right$,$ 
then again current flows out the $C$ line. But if current 
flows in the $A$ and 
not in the $B$$,$ or in the $B$ and not in the $A,$ then no current flows out 
the $C$ line. I'll let you do the ``flow'' analyze for the second diagram. \parm 
In the first diagram, the symbol $\searrow${\kern -.3em}$\rm _x$ means that the arrow has been removed from the pathway indicator. (This is done to minimize the storage space required when processing this monograph.) The diagrams only show what happens when both A and B have current.  \parm

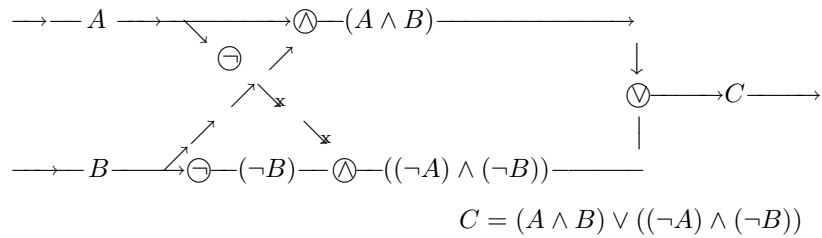
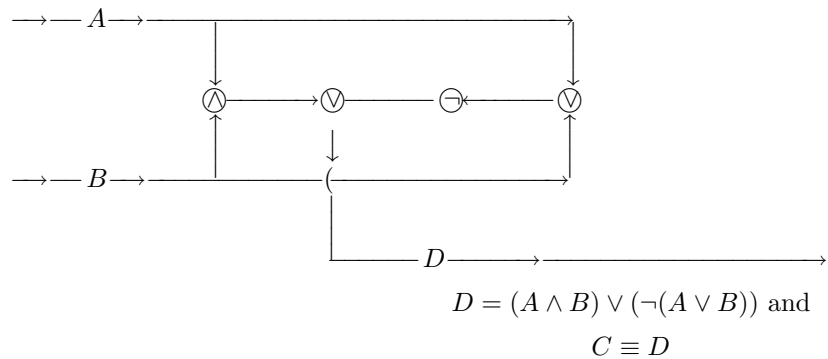
\line{{\hbox to 0.2in {\rightarrowfill}}{\hbox to 0.2in {\noarrowfill}}$A$    
{\hbox to 0.25in{\rightarrowfill}}{\kern -.5em} 
{\hbox to 0.7in{\rightarrowfill}}\andgate{\hbox to 
0.15in{\noarrowfill}}$(A\land B)${\hbox to 
1.05in{\rightarrowfill}}\hfil}
\vskip -0.14in
\line{\hskip 0.9in$\searrow$\hfil}
\vskip -0.05in
\line{\hskip 1.08in \inverter \hskip 0.15in$\nearrow$ \hskip 1.7in$\big\downarrow$\hfil}
\line{\hskip 1.15in $\nearrow$$\searrow${\kern -.3em}$\rm _x$\hskip 1.77in\orgate {\hbox to 0.4in {\rightarrowfill}}$C${\hbox to 0.4in {\rightarrowfill}}\hskip 0.3in}
\line{\hskip 0.94in$\nearrow$\hskip 0.45in$\searrow${\kern -.3em}$\rm _x$\hskip 1.58in$\big\vert$}
\vskip -0.06in
\line{\hskip 0.8in$\nearrow$}
\vskip -0.13in
\line{{\hbox to 0.26in {\rightarrowfill}}{\kern -.3em}{\hbox to 
0.19in{\noarrowfill}}$B${\hbox to 0.4in{\rightarrowfill}}\inverter{\hbox  to 
0.15in{\noarrowfill}}$(\neg B)${\hbox to 0.18in{\noarrowfill}}\andgate{\hbox to 
0.15in{\noarrowfill}}$((\neg A)\land (\neg B))${\hbox to 
0.52in{\noarrowfill}}\hfil}
\medskip
\centerline{$C= (A \land B) \lor ((\neg A) \land (\neg B))$}
\vskip 0.5in
\line{{\hbox to 0.2in 
{\rightarrowfill}}{\hbox to 0.2in {\noarrowfill}}$A${\hbox to 0.2in 
{\rightarrowfill}}{\hbox to 2.25in{\rightarrowfill}}
\hfil} 
\vskip -0.05in 
\line{\hskip 1.03in$\bigg\downarrow$\hskip 1.78in $\bigg\downarrow$}   
\vskip -0.0in
\line{\hskip 1.0in{\andgate}{\hbox to 0.5in{\rightarrowfill}}\orgate{\hbox to 
0.5in{\noarrowfill}}\inverter{\hbox to 0.5in{\leftarrowfill}}\orgate\hfil}
\line{\hskip 1.03in$\bigg\uparrow$\hskip 0.52in $\big\downarrow$\hskip 1.15in $\bigg\uparrow$}
\vskip -0.07in
\line{{\hbox to 0.2in 
{\rightarrowfill}}{\hbox to 0.2in {\noarrowfill}}$B${\hbox to 0.2in 
{\rightarrowfill}}{\hbox to 0.93in {\noarrowfill}}(\kern -0.2em {\hbox to 1.26in{\rightarrowfill}} 
\hfil}
\line{\hskip 1.66in $\bigg |$}    
\vskip -0.07in
\line{\hskip 1.66in {\hbox to 0.5in{\noarrowfill}}$D${\hbox to 
0.5in{\rightarrowfill}}{\hbox to 1.5in{\rightarrowfill}}  \hfil}\pars
\centerline{$D= (A\land B)\lor(\neg(A \lor B))$ and}\smallskip
\centerline{$C \equiv D$}\parm

An important aspect of logic circuits is that they can be so 
constructed so that they will do binary arithmetic. Here is an example of 
binary arithmetic. \par
$$\matrix{\ \ 0\ 1\cr
          +\underline{1\ 1}\cr
          1\ 0\ 0\cr}$$
The process goes like this. First $1 + 1 =1\ 0.$ Thus you get a 0 = $S$
with a carry over digit = 1 = $C.$ The carry over digit is then added to 
the next column digit 1 and you get $1\ 0,$ which is a 0 with a carry over of 
1. The following logic circuit does a part of this arithmetic. If current in 
$A,\ B$ indicates 1$,$ no current indicates 0. This represents the first step 
and yields the basic $S$ and the basic cover number $C.$ (Insert Figure 2 below.)\parm
\line{{\hbox to 0.2in 
{\rightarrowfill}}{\hbox to 0.2in {\noarrowfill}}$A${\hbox to 0.2in 
{\rightarrowfill}}{\hbox to 2.23in{\rightarrowfill}}
\hfil}
\vskip -0.05in
\line{\hskip 1.02in $\bigg\downarrow$\hskip 1.76in$\bigg\downarrow$}    
\line{\hskip 1.0in{\orgate}{\hbox to 0.5in{\rightarrowfill}}\andgate{\hbox to 
0.5in{\noarrowfill}}\inverter{\hbox to 0.5in{\leftarrowfill}}\andgate{\hbox to 
0.15in{\noarrowfill}}$C${\hbox to 0.5in{\rightarrowfill}}\hfil}
\line{\hskip 1.02in $\bigg\uparrow$\hskip 0.55in $\big \vert$\hskip 1.17in $\bigg\uparrow$}
\vskip -0.06in
\line{{\hbox to 0.2in 
{\rightarrowfill}}{\hbox to 0.2in {\noarrowfill}}$B${\hbox to 0.2in 
{\rightarrowfill}}{\hbox to 0.93in {\noarrowfill}}(\kern -.2em {\hbox to 1.27in{\rightarrowfill}} 
\hfil}
\line{\hskip 1.66in $\bigg |$}    
\vskip -0.06in
\line{\hskip 1.66in {\hbox to 0.5in {\noarrowfill}}$S${\hbox to 
1.95in{\noarrowfill}} \hfil}
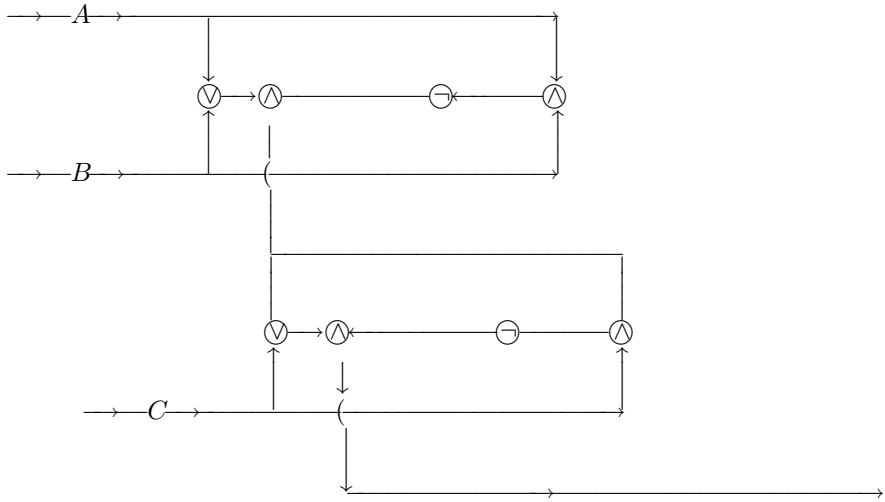
\medskip

\hrule\smallskip\hrule\medskip
\centerline{\bf EXERCISES 2.7}\parm
\indent 1. Express each of the following formula as an equivalent formula in 
terms of$,$ at the most$,$ one $\neg$ on the left of an atom$,$ and the connectives
$\lor$ or $\land.$\parm
\+(a) $P \iff(A \to (R\lor S))$&(d) $((\neg P)\iff Q) \to R$\cr
\+(b) $((\neg P) \to Q) \iff R$&(e) $(S \lor Q) \to R$\cr
\+(c) $(\neg((\neg P)\lor (\neg Q)))\to R$&(f) $(P\lor(Q\land S)) \to R$\cr
\parm
\noindent 2. Write the denial for each of the following formula.\parm
\+(a) $((\neg P)\lor Q)\land (((\neg Q)\lor P) \land R)$&(c) $((\neg R)\lor 
(\neg P)) \land (Q \land P)$\cr
\+(b) $((P\lor (\neg Q))\lor R) \land$&(d) 
$(((Q\land (\neg R))\lor Q) \lor (\neg P))\land (Q \lor R)$\cr 
\+\hfil $(((\neg P)\lor Q)\land R)$& \cr\parm
\noindent 3. Write each of the following formula in its fdnf$,$ if it has 
one.\parm
\+(a) $(P \land (\neg Q))\lor (P \land R)$& (c) $ (P\lor Q) \iff (\neg R)$\cr
\+(b) $ P \to (Q \lor (\neg R))$&(d) $(P \to Q) \to ((Q \to R) \to (P \to 
R))$\cr\parm
\noindent 4. Using only inverters$,$ or-gates and and-gates$,$ diagram the logic 
circuits for the following formula. Notice that you have three inputs.\parm
\+ (a) $(A\lor (\neg B))\lor (B \land (C \lor(\neg A)))$&(b) $(A \to B)\lor (\neg 
C)$\cr\parm
\noindent 5. The following two diagrams correspond to an output which is a 
composite formula in terms of $A,\ B,\ C.$ Write down this formula. \parm
\vskip 0.25in
\noindent (a) \medskip
\line{{\hbox to 0.2in 
{\rightarrowfill}}{\hbox to 0.2in {\noarrowfill}}\kern -0.3em$A$\kern -0.2em
{\hbox to 1.25in{\rightarrowfill}}\inverter{\hbox to 1.0in 
{\noarrowfill}}\hfil}
\vskip -0.06in 
\line{\hskip 1.03in $\bigg\downarrow$\hskip 1.65in $\bigg |$}   
\line{\hskip 1.0in{\orgate}{\hbox to 0.5in{\leftarrowfill}}\andgate{\hbox to 
0.4in{\noarrowfill}}\kern -0.2em{\hbox to 0.63in{\leftarrowfill}}\kern -0.2em\orgate\hfil}

\line{\hskip 1.04in $\big \downarrow$\hskip 0.52in $\bigg\uparrow$\hskip 1.01in $\bigg \uparrow$}
\vskip -0.07in
\line{{\hbox to 0.2in 
{\rightarrowfill}}\kern -0.2em{\hbox to 0.2in {\noarrowfill}}\kern -0.2em$B$\kern -0.2em{\hbox to 0.2in 
{\rightarrowfill}}\kern -0.2em{\hbox to 0.45in{\noarrowfill}$($\kern -0.2em{\hbox to 1.73in {\rightarrowfill}}} 
\hfil} 
\vskip -0.03in
\line{\hskip 1.03in $\big\downarrow$}   
\vskip -0.05in
\line{ {\hbox to 0.2in 
{\rightarrowfill}}\kern -0.2em{\hbox to 0.23in {\noarrowfill}}\kern -0.2em$C$\kern -0.1em{\hbox to 0.49in
{\rightarrowfill}}\andgate{\hbox to
0.5in{\rightarrowfill}}\kern -0.2em{\hbox to 1.5in{\rightarrowfill}}  \hfil}
\vfil\eject

\noindent (b)
\vskip 0.5in
\line{{\hbox to 0.2in 
{\rightarrowfill}}\kern -0.2em{\hbox to 0.2in {\noarrowfill}}\kern -0.2em$A$\kern -0.2em{\hbox to 0.2in 
{\rightarrowfill}}\kern -0.2em{\hbox to 2.31in{\rightarrowfill}}
\hfil} \vskip -0.05in
\line{\hskip 1.02in $\bigg\downarrow$\hskip 1.73in$\bigg\downarrow$}   
\line{\hskip 1.0in{\orgate}{\hbox to 0.2in{\rightarrowfill}}\andgate{\hbox to 
0.8in{\noarrowfill}}\kern -0.2em\inverter{\hbox to 0.5in{\leftarrowfill}}\kern -0.2em\andgate\hfil}
\line{\hskip 1.02in $\bigg \uparrow$\hskip 0.25in$\big |$\hskip 1.44in $\bigg\uparrow$}
\vskip -0.08in
\line{{\hbox to 0.2in 
{\rightarrowfill}}\kern -0.2em{\hbox to 0.2in {\noarrowfill}}\kern -0.2em$B$\kern -0.2em{\hbox to 0.2in 
{\rightarrowfill}}\kern -0.2em{\hbox to 0.77in{\noarrowfill}}\kern -0.2em(\kern -0.2em{\hbox to 1.53in{\rightarrowfill}} 
\hfil}  
\line{\hskip 1.34in \kern 0.2em$\bigg |$}  
\vskip -0.22in
\bigskip
\line{\hskip 1.38in{\hbox to 1.87in{\noarrowfill}}
\hfil} 
\vskip -0.05in
\line{\hskip 1.37in $\bigg |$\hskip 1.79in$\bigg |$}  
\vskip -0.01in
\line{\hskip 1.35in{\orgate}{\hbox to 0.2in{\rightarrowfill}}\andgate{\hbox to 
0.8in{\leftarrowfill}}\kern -0.2em\inverter{\hbox to 0.5in{\noarrowfill}}\kern -0.2em\andgate\hfil}
\line{\hskip 1.36in$\bigg\uparrow$\hskip 0.27in $\big \downarrow$\hskip 1.37in$\bigg\uparrow$}
\vskip -0.07in
\line{\hskip 0.4in{\hbox to 0.2in 
{\rightarrowfill}}\kern -0.2em{\hbox to 0.2in {\noarrowfill}}\kern -0.2em$C$\kern -0.2em{\hbox to 0.2in 
{\rightarrowfill}}\kern -0.2em{\hbox to 0.76in{\noarrowfill}}\kern -0.2em(\kern -0.2em{\hbox to 1.49in{\rightarrowfill}} 
\hfil}
\line{\hskip 1.74in $\bigg \downarrow$}    
\vskip -0.05in
\line{\hskip 1.78in{\hbox to 
1.1in{\rightarrowfill}}\kern -0.2em{\hbox to 1.75in{\rightarrowfill}}  \hfil}\pars
\ss
\noindent {\bf 2.8 The Princeton Project$,$ Valid Consequences (in General).}\parm
After the methods were discovered that use logical operators to generate a 
solution to the General Grand Unification Problem$,$ and that give an answer to 
the questions  ``How did our universe come into being?'' and  ``Of what is empty 
space composed?'' I discovered that the last two questions were attacked$,$ in 
February --- April 1974$,$ by John Wheeler and other members of the Physics and 
Mathematics Department at Princeton University. Wheeler and Patton write$,$ 
``It is 
difficult to imagine a simpler element with which the construction of physics 
might begin  than the choice  yes--no or true--false or open 
circuit--closed circuit. . . . which is isomorphic [same as] a proposition in 
the propositional calculus of mathematical logic.'' [ Patton and Wheeler$,$ 
{\it Is Physics Legislated by a Cosmogony}$,$ in Quantum Gravity$,$ ed. Isham$,$ 
Penrose$,$ Sciama$,$ Oxford University Press$,$ Oxford (1975)$,$ pp. 538--605.]
These basic concepts are exactly what we have just studied. \pars
These individuals$,$ one of the world's foremost group of scientists$,$ attempted 
to solve this problem by a statistical process$,$ but failed to do so. Because 
they failed$,$ they rejected any similar approach to the problem. They seemed to 
be saying that ``If we can't solve these problems$,$ then no one can.'' They 
were wrong in their rejection of the propositional calculus as a useful aspect  
for 
such a solution. But$,$ the solution does not lay with this two valued 
truth--falsity model for the 
propositional logic. The solution lies with the complementary aspect we'll 
study in section 2.11 called  {\it proof theory}.  
Certain proof theory concepts correspond to simple 
aspects of our truth--falsity model. 
In particular$,$ we are able to determine by assignment and truth-table 
procedures whether or not a logical argument is 
following basic human reasoning processes (i.e. 
classical propositional deduction).  
However$,$ what you are about to study will not 
specifically identify what the brain is doing$,$ but it will determine whether 
or not it has done its deduction in terms of the classical processes must 
easily comprehend by normal human beings.\pars
Let $\{A_1, \ldots, A_n\}$ be a finite (possibly empty) set of formula. 
These formula represent the {{\it 
hypotheses}} or {{\it premises}} for a logical argument. For convenience$,$ it 
has become common place to drop the set-theoretic notation $\{$ and $\}$ from
this notation. Since these are members of a set$,$ they are all distinct in 
form. Again from the concepts of set-theory$,$ these formula are not considered 
as   ``ordered'' by the ordering of the subscripts.  \parm\vfil\eject
\hs
{\bf Definition  2.8.1} (Valid consequence.) 
A formula $B$ is a {{\it valid consequence}} (or simply a consequence) of a 
set of premises $A_1,\ldots,A_n$ if for any assignment $\ass$ to the atoms in 
each $A_1,\ldots,A_n$ \underbar{AND} $B$ such that $v(A_1,\ass) = \cdots =
v(A_n,\ass)=T,$ then $v(B,\ass) = T.$ If $B$ is a valid consequence of
$A_1,\ldots,A_n,$  then this is denoted by  $A_1,\ldots,A_n \mod B.$\pars
\hm
As usual$,$ if  some part of definition  2.8.1 does not hold$,$ then $B$ is an
invalid (not a valid) consequence of $A_1,\ldots,A_n.$ This is denoted by
$A_1,\ldots,A_n \not\mod B.$  It is important to notice that definition 2.8.1 
is a conditional statement. This leads to a very interesting result.\parm
{\bf Theorem 2.8.1} {\sl Let $A_1,\ldots,A_n$ be a set of premises and there 
does not exist an assignment to the atoms in  $A_1,\ldots,A_n,$ such that the 
truth-values $v(A_i,\ass)=T$ for each $i,$ where $1\leq i\leq n.$ Then$,$ 
for ANY formula $B \in L,$ $A_1,\ldots,A_n\mod B.$}\parm
Proof. In a true conditional statement$,$ if the hypothesis is false$,$ then 
the conclusion holds. \qed
The conclusion of theorem 2.8.1 is so significant that I'll devote an entire 
section (2.10) to a more in-depth discussion. Let's continue with more facts 
about the valid consequence concept. We will need two terms$,$ however$,$ for the 
work in this section that are also significant for section 2.10.\parm
\hs
{\bf Definition 2.8.2} (Satisfaction.) If given a set of formula $A_1, \ldots,
A_n\ (n\geq 1),$  there exists an assignment $\ass$ to all the atoms in 
$A_1, \ldots, A_n$ such that $v(A_1,\ass) = \cdots = v(A_n,\ass)=T,$ then 
the set of premises are said to be {{\it satisfiable.}} The assignment itself 
is said to {{\it satisfy}} the premises. One the other hand$,$ if such an 
assignment does not exist$,$ then the premises are said to be {{\it not 
satisfied.}}\pars
\hm
{\bf Theorem 2.8.2} (Substitution of equivalence.) {\sl  If $A_n \equiv C$ and    
$A_1,\ldots,A_n \mod B,$ then  $A_1,\ldots,A_{n-1},C \mod B.$   If $B \equiv C$ and    
$A_1,\ldots,A_n \mod B,$ then  $A_1,\ldots,A_n \mod C.$}\parm
Proof. Left to you.\qed   
Are many logical arguments that seem very complex$,$ in reality$,$ a  
disguised simple deduction? Conversely$,$ can we take a simple deduction and 
make it look a little more complex? The following theorem is not the last word 
on this subject and is very closely connected with what we mean when we say 
that such and such is a set of premises.  When we write the premises
$A_1,\ldots,A_n,$ don't we sometime (all the time?) say  ``and'' when we write 
the  comma  ``$,$''? Is this correct?\parm
{\bf Theorem 2.8.3} (The Deduction Theorem) {\sl Let $\Gamma$ be any finite 
(possible empty) set of formula and $A,\ B \in L.$\parm
(i) $\Gamma,\ A \mod B$ if and only if $\Gamma \mod A  \to B.$\parm
(ii) Let $A_1, \ldots, A_i,\ldots A_n$ be a finite (nonempty) set  
of formula$,$ where $1<i \leq n.$ Then 
$A_1, \ldots, A_n\mod B$ if and only if $A_1, \ldots, A_i\mod (A_{i+1} \to 
(\cdots \to (A_n \to B)\cdots )).$\parm                                                                    
(iii) Let $A_1, \ldots, A_i,\ldots A_n$ be a finite (nonempty) set   
of formula$,$ where $1<i \leq n.$ Then 
$A_1, \ldots, A_n\mod B$ if and only if $(A_1\land \cdots \land A_i),\ldots, 
A_n \mod 
B.$\parm                                                                       
(iv) $A_1, \ldots, A_n\mod B$ if and only if $\mod (A_1\land \cdots \land 
A_n) \to B.$} \parm
Proof. (i) Assume that $\Gamma,\ A \mod B.$  Let $\ass$ be an assignment to 
the set of atoms in  $\Gamma, A, B.$  If $\ass$ satisfies $\Gamma,\ A,$ then
$v(A,\ass) = T.$  Also from the hypothesis$,$ $v(B,\ass) = T.$ Hence$,$ 
$v(A \to B, \ass) = T.$ Now if $\Gamma$ is not satisfied (whether or not $A$ 
is)$,$  then from theorem 
2.8.1$,$ $\Gamma \mod A\to B.$ If $\ass$ satisfies $\Gamma$ and does not satisfy 
$A,$ then $v(A,\ass) = F.$ Thus$,$ $v(A\to B, \ass) = T.$ All the cases have 
been covered; hence$,$ in general$,$ $\Gamma \mod A \to B.$ \pars
Conversely$,$ assume that $\Gamma \mod A\to B$  and let $\ass$ be as previous 
defined. If $\ass$ does not satisfy $\Gamma,$ then $\ass$ does not satisfy
$\Gamma,\ A.$ If $\ass$ satisfies $\Gamma$ and does not satisfy $A,$ then 
$\ass$ does not satisfy $\Gamma,\ A.$ Hence$,$ we only need to consider what
happens if $\ass$ satisfies $\Gamma,\ A.$ Then$,$ in this case$,$ $v(A,\ass) = T.$
Since $v(A\to B,\ass) = T,$ then $v(B,\ass) = T.$  Therefore$,$ $\Gamma,\ A \mod 
B.$\parm
(ii) (By induction on the number $m$ of connectives $\to$ placed between the 
formula on the right of $\mod$.)\pars
(a) The $m =1$ case is but part (i). Assume theorem holds for $m$ connectives. 
Then one more applications of (i) shows it holds for $m+1$ connectives $\to.$
Hence$,$ the result holds in general. Similarly the converse holds.\parm
(iii) (By induction on the $m$ number of connectives $\land$ placed between
$A_1,\ldots,\ A_i.)$  Suppose that  $A_1, A_2,\ldots, A_n\mod B.$ \pars
(a) Let $m = 1.$ Suppose that  $A_1, A_2,\ldots, A_n\mod B.$ Let $\ass$ be an 
assignment to all the atoms in $(A_1\land A_2),\ldots, A_n, B$ that satisfies 
$(A_1\land A_2),\ldots, A_n.$ Then $v(A_1 \land A_2,\ass) = T.$ Hence,
$v(A_1,\ass) = v(A_2,\ass) =T.$ 
Then 
$\ass$ is an assignment to all the atoms in $A_1, A_2,\ldots, A_n, B.$ 
Now $v(A_1,\ass) = v(A_2,\ass) = T = \cdots = v(A_n,\ass).$ This implies that 
$v(B,\ass) = T.$ Hence$,$ $(A_1\land A_2),\ldots, A_n\mod B.$ \pars
(b) Now assume theorem holds for $m$ or less connectives $\land$ and suppose 
that $i = m+2.$ We know that $A_1, A_2,\ldots, A_i,\cdots, A_n \mod B.$ Now 
from (ii)$,$ we know that $A_1, A_2,\ldots, A_{m+1} \mod (A_{m+2} \to (\cdots (A_n \to 
B)\cdots ))$ from (ii). The induction hypothesis yields
$(A_1\land A_2\land\ldots\land A_{m+1}) \mod (A_{m+2} \to (\cdots (A_n \to 
B)\cdots )).$ From this we have $(A_1\land A_2\land\ldots\land A_{m+1}), 
A_{m+2}\mod (A_{m+3} \to (\cdots (A_n \to B)\cdots )).$ Application 
of (i) and (ii) yields $(A_1\land \cdots \land A_{m+2}),\ldots, 
A_n \mod B.$ The general result follows by induction and the converse follows in a 
similar manner.\parm    
(iv) Obvious from the other results.\qed 
When you argue logically for a conclusion$,$ it's rather obvious that one of 
your hypotheses is a valid conclusion. Further$,$ if you start with a specific 
set of hypotheses and obtained a finite set of logical conclusions. You then often 
use these conclusions to argue for other consequences. Surely it should be 
possible to go back to your original hypotheses and argue to you final 
conclusions 
without going through the intermediate process. \parm
{\bf Theorem 2.8.4}\parm
 {\sl (i) $A_1,\ldots,A_n \mod A_i$ for each $i= 
1,\ldots,n.$\parm                       
(ii) If $A_1,\ldots, A_n \mod B_j,$ where $j = 1, \ldots, p,$ and $B_1, 
\ldots, B_p \mod C,$ then $A_1,\ldots, A_n \mod C.$}\parm
Proof. (i) Let $\ass$ be any assignment that satisfies  $A_1,\ldots,A_n.$ Hence$,$ 
$v(A_i,\ass) = T,$ for each $i= 
1,\ldots,n.$ Thus $A_1,\ldots,A_n \mod A_i$ for each $i= 
1,\ldots,n.$\parm  
(ii) Let $\ass$ be an assignment to all the atoms in $A_1,\ldots, A_n, B_i,\ldots, B_p, C.$ Of 
course$,$ this is also an assignment for each member of this set. Conversely$,$ 
any assignment to any of the formula in this set can be extended to an 
assignment to all the atoms in this set. Suppose that $v(A_i,\ass) = T,$ for 
each $i$ such that $1 \leq i \leq m.$ Since $A_1,\ldots, A_n \mod B_j,$ 
where $j = 1, \ldots, p,$ then $v(B_j,\ass) = T,$ where $j = 1, \ldots, p.$  
This $\ass$ satisfies $B_j$ for each $j = 1, \ldots, p.$ Hence$,$ from the remainder 
of the hypothesis$,$ $v(C,\ass) = T$ and the proof is complete.\qed
[For those who might be interested. This is not part of the course. First$,$ we 
know that if $A \equiv B,$ then we can substitute throughout the process 
$\mod$ anywhere $A$ for $B$ or $B$ for $A.$ Because of theorem 2.8.3$,$ all 
valid consequences can be written as $A \mod B,$ where $A$ is just one 
formula. (i) Now $A \mod A,$ and (ii) if $A\mod B,\ B \mod C,$ then $A \mod C.$ Thus
$\mod$ behaves almost like the partial ordering of the real numbers. 
If you substitute $\leq$ for $\mod$ you have $A\leq A,$ if $A \leq B,\ B \leq 
C,$ then $A \leq C.$ In its present form it does not have the requirement
that (iii) if $A\leq B, \ B \leq A,$ then $A = B.$ However$,$ notice that if
$A\mod B,\ B \mod A,$ then $\mod A \to B,\ \mod B \to A.$ This implies that 
$A\equiv B.$ Hence we could create a language by taking one and only one member 
of from each equivalence class [A] (see problem 1$,$ Exercise 2.6.) 
and only used these for logical deduction. (A very boring way to communicate.)
Then (iii) would hold. Hence$,$ in this case$,$ everything known about a partial 
ordering should hold for the $\mod.$] \parm
\hs\hm
\centerline{\bf EXERCISE 2.8}\parm
\noindent 1. First$,$ note the discussion at the top of page 57 as to how to use truth-tables 
to determine if a consequence is valid. Now use the truth-table method to 
determine whether the following consequences are valid consequences from the 
set of premises.\par
\noindent (a) $P\to Q, \ (\neg P)\to Q \mod Q.$\parm
\noindent (b) $P\to Q,\ Q\to R,\ P\mod R.$\parm
\noindent (c) $(P\to Q)\to P,\ \neg P \mod R.$ \parm
\noindent (d) $(\neg P)\to (\neg Q),\ P\mod Q.$\parm
\noindent (e) $(\neg P)\to (\neg Q),\ Q \mod P.$\pars
\ss
\noindent {\bf 2.9 Valid Consequences --- Model Theory and Beyond.} \parm
The most obvious way to show that $B$ is a valid consequence of 
$A_1,\ldots,A_n$ is by the truth-table method.\parm
(i)  Simply set up a
truth-table for the formulas  $A_1,\ldots,A_n,\ B.$\parm
(ii) Look at every row$,$ where under each $A_i$ there is a $T,$ and if in 
that row the $B$ is  a $T,$ then $B$ is a valid consequence of 
$A_1,\ldots,A_n.$\parm
Of course$,$ if $A_1,\ldots,A_n$ is not satisfiable$,$ then it is automatically 
the case that $B$ is a valid consequence. BUT$,$ in most cases you would have a 
very large truth-table. For example$,$ for the arguments of section 1.1$,$ you 
would need a truth-table with $2^6 +1=65$ rows and about 8 columns. 
You might get a research grant$,$ of some sort$,$ so that you could get the 
materials for such a truth-table construction. This is the strict {{\it model 
theory}} approach. Why? Because the truth-table is a  ``model'' (i. e. not the real 
thing) for classical human propositional deduction. When I deduce what I hope 
is a logical conclusion$,$ I don't believe I construct a truth-table in my mind.
Maybe some do$,$ but I don't. So$,$ is there another method that is 
model theory viewed in  a different way  that \underbar{may be}
a shorter method? The method comes from theorem 2.8.3 part  (iv). \parm
\hs {\bf NOTATION CHANGE.}  It seems pointless to keep saying ``Let 
$\ass$ be an assignment to the atoms in $A.$ Then $v(A,\ass) = T\ {\rm or} \ 
F.$'' Why not 
do the following: just write $v(A,\ass)= v(A) = T.$ When we see $v(A) = T\ 
{\rm or}\ F,$ we know that there is some assignment to the atoms that makes it 
so.\pars\hm
Now what theorem 2.8.3 tells us is that all we need to do is to show that 
$\mod (A_1 \land \cdots \land A_n) \to B.$ 
But$,$ if one selects an  assignment such that $v(B) = F$$,$ then  
$\mod (A_1 \land \cdots \land A_n) \to B$ if and only if $v((A_1 \land \cdots 
\land A_n)) = F.$ We need not look at the case when $v(B) = T.$ (Why not ?) 
The method is a natural language 
algorithm. This means that I use some of the terms that I've previously 
introduced and ordinary English to give a series of repeatable instructions.
Now whether or not this method is shorter than the truth-table method depends 
upon how clever you are in a certain selection process. Only experience 
indicates that it is often much shorter. The instructions themselves are not 
short in content. But remember these methods took over 2,000 years to 
develop.\parm

\hs
{\bf Special Method 2.9.1} (To show that $\mod (A_1 \land 
\cdots \land A_n) \to B.$)  In what follows$,$ the symbol $\Rightarrow$ 
represents that word {\it forces}. \pars
(1) First$,$ let $v(B) = F.$ In the simplest case$,$ this truth-value for $B 
\Rightarrow$ a fixed truth-value on each of the atoms of $B.$ Let these 
atoms have these forced truth-values.\pars
 {\bf Example 2.9.1.1} Let $B = P \to (Q \to R).$ Then $v((P \to (Q \to R)))= 
F \Rightarrow v(P) = T;\ v(Q) = T; \ v(R) = F$ and these are the only 
possibilities. \parm
(2) See if these forced and fixed atom truth-values forces any of the premises 
to have a fixed truth-value. \parm
(3) If (2) occurs and the value of the premise is $F,$ then the process stops 
and you have a {\bf valid} consequence. \pars
{\bf Example 2.9.1.2} Suppose that one of the premises is $A=P\land R.$ Then 
for example 2.9.1.1 atomic values you would have that $v(A) = F.$ You can stop 
the entire process. $B$ is an valid consequence of the premises.\parm
(4) If (2) occurs and the truth-value is $T$ for a premise$,$ then simply write it down as its 
value. \parm
(5) If (4) occurs and there are more premises that are NOT forced to take on a 
specific truth value by the forced atomic truth-values for $B,$ then you can 
select any of the remaining premises$,$ usually those with the fewest 
non-forced atoms (but not always)$,$ and set the truth value of the selected 
premise as $T.$ \pars
{\bf Example 2.9.1.3} Suppose that process did not stop at step (3). Let one 
of the premises be $A_1 = S \lor R.$ Then setting $v(A_1) = T \Rightarrow
v(S) = T.$\parm
(6) You now have some premises with forced or selected values of $T.$ Now use 
the forced atomic values and begin again with (2) for the remaining process. 
\parm
(7) If this the simplest part of the process stops$,$ then it will either force 
a premise to be $F$ and you may stop and declare the consequence valid or all 
the premises will be selected or forced to be $T$ and you have found one 
assignment that proves that the consequence is {\bf invalid}.\pars
\hm
{\bf Example 2.9.1.4} Suppose that you what to answer the question
$P_1 \to (P_2 \to P_3),\ (P_3 \land P_4) \to P_5,\ (\neg P_6) \to (P_4\land 
(\neg P_5)){\buildrel \rm ?\over \mod}P_1 \to (P_2 \to P_6).$\pars   
(i)  Well$,$ 
let $v(P_1 \to (P_2 \to P_6)) = F.$ Then $v(P_1) = v(P_2) = T,\ v(P_6) = F.$ 
These values do not force any of the premises to be any fixed value.\pars
(ii) Select the first premise and let $v(P_1 \to (P_2 \to P_3))= T.$ From (i) 
$\Rightarrow v(P_3) = T.$ The values that have been forced do not force the 
remaining premises to take any fixed value.\pars
(iii) Let $v((\neg P_6) \to (P_4 \land (\neg P_5)))= T\ \Rightarrow v(P_4 
\land (\neg P_5)) = T$ from (1). Hence $v(P_4) = T,\ v(P_5) = F.$ But these 
forced atoms $\Rightarrow v((P_3 \land P_4) \to P_5) = F.$\pars 
(iv) Hence$,$ $\mod$ 
holds. Notice that for this example a truth table requires 65 rows.\parm
{\bf Example 2.9.1.4} Suppose that you what to answer the question
$P \to R,\ Q \to S,\ (\neg R)\lor (\neg S){\buildrel \rm ?\over \mod}
P \lor (\neg Q).$\pars 
(i) Let $v(P \lor (\neg Q))= F,\ \Rightarrow v(P) = F,\ v(Q) = T.$\pars
(ii) Let $v(Q \to S) = T,\ {\rm from\ (i),} \Rightarrow v(S) = T.$ \pars
(iii) Let $v((\neg R)\lor (\neg S)) = T,\ {\rm from\ (ii),} \Rightarrow v(R) = 
F.$\pars
(iv) Now (i) and (iii) $\Rightarrow v(P\to R) = T.$   \pars
(v) Since all premises were either selected or forced to be $T,$ then $B$ is 
an invalid consequence from the premises.\parm
Is all of this important? Well$,$ suppose that you were given a set of orders by 
your commanding officer. You tried to follow these orders but could not do so.
Why can't they be carried out? You discover$,$ after a lot of work$,$ 
that the consequence your 
commanding officer claimed was a result of the set of premises he gave is 
invalid. Next you must prove this fact at a  court-martial. Yes$,$ it could be 
very important. But the above method need not be as straightforward as the 
examples indicate. \parm
{\bf Special Method 2.9.2 (Difficulties  with Method 2.9.1.)} (Case studies.) This special method$,$ can brake down 
and become very complex in character for one basic reason. Either the 
selection of the (might be) consequence as an $F$ or the selection of any of 
the premises as a $T$ need not produce fixed values for the atoms. Now what do 
you do? For the case study difficulties$,$ its easier to establish INVALID
consequences. \pars
(1) Suppose that your assumption that the (may be) consequence $B$ has 
truth-value 
$F$ does not yield unique atomic truth-values. Then you must brake up the 
problem into all the cases produced by all the different possible truth-values 
for the atoms in $B.$ \pars
{\bf Example 2.9.2.1} Let $B = P\iff Q.$ Then for $v(P \iff Q) = F$ 
there are the following two cases. (a) $v(P ) = F,\ v(Q) = T.$ (b)
$v(P) = T,\ v(Q) = F.$ \parm
(2) For premises that are not forced to have specific truth-values$,$ then your 
selection of a truth-value for a (possible) premise $A$ need not yield unique 
atomic truth-values. This will lead to more case studies. Indeed$,$ possible 
case studies within case studies. \parm
(3) During any of the specific case studies if the truth-values of 
\underbar{all} possible premises yields $T,$ then you may stop for you have an 
invalid consequence.\parm
(4) If during any case study you get one or more of the assume premises to be 
forced to be $F,$ then this does NOT indicate that you have a valid 
consequence. You must get an $F$ for some assumed premise \underbar{for all}
possible case studies before you can state that it is a valid argument.\parm
{\bf Example 2.9.2.2} Suppose that you what to answer the question
$P \to R,\ Q \to S,\ (\neg R)\lor (\neg S){\buildrel \rm ?\over \mod}
P\land Q.$\pars
(i) Let $v(P\land Q) = F.$ You have three cases. (a) $v(P) =F,\ v(Q) = F.$  
(b) $v(P) =T,\ v(Q) = F.$  (c)  $v(P) =F,\ v(Q) = T.$\pars
Case (a). No assumed premise is forced to be anything. So$,$ select $v(Q \to S) 
= T.$ This yields two subcases. (a$_1$) $v(S) = F,$ (a$_2$) $v(S) = T.$ \pars
Case (b). Again let $v(Q \to S)= T.$ Again we have two subcases.  
(b$_1$) $v(S) = F,$ (b$_2$) $v(S) = T.$\pars
Case (c).  Again let $v(Q \to S)= T.$ Now this $\Rightarrow v(S) = T.$ \pars
(ii) Now we would go back and select another assumed premises such as $P\to R$ 
and set its value to $T.$ Then assuming cases (a)$,$ (a$_1$) see what happens to 
the atoms in $P \to R.$ This might produce more cases such as (a$_{11}$) and 
an (a$_{12}$). We would have a lot to check if we believed that $B$ might be a 
valid consequence.\pars
(iii) Notice that we have only one more assumed premises remaining $(\neg R) 
\lor (\neg S).$  We can assume that $v(P \to R) = T$ and there are atomic 
values that will produce this truth-value. Now $S$ does not appear in this 
formula hence under condition  (b$_1$) $v(S) = F. \ v((\neg R) 
\lor (\neg S))=T.$ Hence we have found a special assignment that shows that
$B$ is an invalid consequence. \pars
(iv) Lets hope the method doesn't lead to many case studies since it might be 
better to use truth-tables. \pars
NOTE ON FORMULA VARIABLES. You do not need to use the atomic form 
of a formula when a valid consequence is being determined. What you actually 
can do is to substitute for \underbar{every} specific atom in every place it 
appears a formula variable symbol$,$ and distinct variables for distinct atoms.  
If for each variable formula symbol you substitute a fixed formula in 
atomic form$,$ then in the valid consequence truth-table the various levels at 
which the premises are $T$ is only dependent on the connectives that are in 
the original formula prior to substitution if the premise is not a single 
atom. If it is a single atom$,$ then the atom as a premise still only depends 
upon it 
being give a $T$ value. The same would be true for any formula substituted 
throughout the premises and assumed consequence for that atom. Thus$,$ in  
exercise 1 below formula variable symbols have been used. Simply consider 
them to behave like atoms. Each time you 
determine that the indicated formula is a valid consequence$,$ then 
you have actually determined the case for infinitely many formula. But if 
it is an invalid consequence$,$ then you cannot make such a variable 
substitution. Such an  ``invalid'' result only holds for atoms. 
\parm
\hs\hm
\centerline{\bf EXERCISE 2.9}\parm
\noindent 1. Using the special method 2.9.1 or 2.9.2 (in the formula variable 
form) 
to determine whether or 
not\parm
(a) $(\neg A) \lor B,\ C\to (\neg B) \mod A \to C.$\pars
(b) $A\to (B\to C),\ (C\land D)\to E,\ (\neg G)\to (D \land (\neg E))\mod A\to 
(B\to G).$ \pars
(c) $(A\lor B)\to (C\land D),\ (D\lor E) \to G \mod A \to G.$\pars
(d) $A \to (B\land C),\ (\neg B) \lor D,\ (E \to (\neg G)) \to (\neg D),\ B\to 
(A \lor (\neg E)) \mod B \to E.$\parm
\noindent 2. Translate the following natural language arguments into 
propositional formula using the indicated propositional symbols and determine 
by method 2.9.1 or 2.9.2 whether or not the argument is valid or invalid. 
(Remember that you do NOT need to know what the terms in a phrase mean to 
check the validity of an argument.)\parm
(a) Either I shall go home (H)$,$ or stay and study (S). I shall not go home. 
Therefore I shall stay and study.\pars
(b) If the set of real numbers is infinite (I)$,$ then it has cardinality $c$ (C). 
If the set of real numbers is not infinite$,$ then it forms a finite set (D). 
Therefore$,$ either the set of real numbers has cardinality c or it forms a 
finite set.\pars
(c) A Midshipman's wage may sometime increase (S) only if there is inflation 
(I). If there is inflation$,$ then the cost of living will increase (C). Now and 
then a Midshipman's wage has increased. Therefore$,$ the cost of living has 
increased. \pars
(d) If 2 is a prime number (P)$,$ then it is the least prime number (L). If 2 is 
the least prime number$,$ then 1 is not a prime number (N). The number 1 is not 
a prime number. Therefore$,$ 2 is a prime number.\pars
(e) Either the set of real numbers is well-ordered (W) or it contains a 
well-ordered subset (C). If the set of real numbers is well-ordered$,$ then every 
nonempty subset contains a first element (R). The natural numbers form a 
well-ordered subset of the real numbers (N). Therefore$,$ the real numbers are 
well-ordered.\pars
(f) If it is cold tomorrow (C)$,$ then I'll wear my heavy coat (I) if the sleeve 
is mended (M). It will be cold tomorrow and the sleeve will not be mended. 
Therefore$,$ I'll not wear my heavy coat.\pars
(g) If the lottery is fixed (L) or the Colts leave town again (C)$,$ then the 
tourist trade will decline (D) and the town will suffer (S). If the 
tourist trade decreases$,$ then the police force will be more content (P).
The police force is never content. Therefore$,$ the lottery is fixed.\pars
\ss
\noindent {\bf 2.10 Satisfaction and Consistency.}\parm
As mentioned previously$,$ the concept of when a set of formula is satisfied is 
of considerable importance. Suppose that we assume that a set of premises 
refer to  ``things'' that occur in reality. As you'll see$,$ in order for a set 
of premises to differentiate between different occurrences it must be 
satisfiable. We recall the definition. \parm
\hs
{\bf Definition 2.10.1} (Satisfaction.) A nonempty (finite) set of premises
$A_1,\ldots,A_n$ is {{\it satisfiable}} if there exists an assignment $\ass$ 
to all the atoms that appear in the premises such that $v(A_i) =  T$ for each 
$i$ such that $1  \leq i \leq n.$\pars \hm
One way to attack the problem of satisfaction  is to make a truth-table. If 
there is a row such that under every $A_i$ there is a $T,$ then the set of 
premises is satisfiable. As we did in the previous section$,$ there is a short 
way to do this without such specific truth-tables. But$,$ for the 
propositional calculus$,$ why is satisfaction so important?\parm\hs
{\bf Definition 2.10.2} A set of formula $A_1,\ldots,A_n$ is 
{{\it consistent}} if for each $B\in L,$ 
$A_1,\ldots,A_n\not\mod B \land (\neg B).$ A set of premises is {{\it 
inconsistent}} if there exist some $B \in L$ such that 
$A_1,\ldots,A_n \mod B \land (\neg B).$\pars\hm
Actually definition 2.10.2 is technical in character since if we only had 
this definition it might never (in time) be possible  to know whether a set of premises 
is consistent. Also$,$ as yet consistency may not seem as an important 
property. Notice that that formula $B \land(\neg B)$ is a contradiction. Thus 
sometimes a set of premises that is {{\it inconsistent}} are also said to be 
{{ \it contradictory}}.
Notice that definition 2.10.2 includes the possible empty set of 
premises. This yields the pure validity concept. The next result shows that 
our pure validity concept is consistent. \parm
{\bf Theorem 2.10.1} {\sl If $B \in L,$ then $\not\mod B \land (\neg 
B).$}\pars 
Proof. Let (for an appropriate assignment $\ass$) $v(B) = T.$ Then
$v(B\land (\neg B)) = F.$  One the other hand$,$ if $v(B) = F,$ then
$v(B\land (\neg B)) = F.$ Since every assignment to $B \land (\neg B)$ is an 
assignment to $B$ and conversely$,$ the result follows.\qed
As mentioned above$,$ it is assumed by many individuals$,$ although it cannot be 
established$,$ that human deduction corresponds to a humanly comprehensible  
``occurred in reality'' concept. I won't discuss the philosophical aspects of 
this somewhat dubious 
assumption$,$ but even if it's$,$ at the least$,$ partially true the concept of 
consistency is of paramount importance. Theorem 2.10.1 gives a slight 
indication of what is going on. Not every formula in our language is a valid 
formula. The concept of {{\it simply consistent}} is defined by the 
statement that a set of premises is simply consistent if not all formula are 
consequences of the premises. Obviously$,$ by theorem 2.10.1$,$ there is no 
difference between the two concepts for an empty set of premises. \pars
The worst thing that can happen for any nonempty set of premises 
$A_1,\ldots, A_n$ in the scientific or 
technical areas is that $A_1,\ldots, A_n\mod B,$  where $B$ is ANY member of
$L.$ Why? This would mean that all formula including contradictions are valid 
consequences. Now if we associate with $A_1,\ldots A_n \mod B,$  the notion 
that if each $A_i$ occurs in reality$,$ then $B$ will occur in reality$,$ then 
this worst case scenario says  ``all things $B$ will occur in reality.'' 
Intuitively$,$ this just doesn't imply that any theory based upon this set of 
premises cannot differentiate between occurrences$,$ but ``true'' could not be 
differentiated from  ``false.'' But how can we know when a set of premises has 
this worst case scenario property? \parm
{\bf Theorem 2.10.2} {\sl A nonempty set of premises $A_1,\ldots, A_n$ 
is inconsistent if and only if $A_1,\ldots, A_n\mod B$ for every $B \in L.$}
\pars
Proof. Let $A_1,\ldots, A_n$ be inconsistent and any $B\in L$. Then there is some $C \in L$ 
such that  $A_1,\ldots, A_n\mod C \land (\neg C).$ Considering \underbar{any} 
assignment $\ass$ to the atoms in $C$ and $B,$ then $v(C \land (\neg C))=F.$ 
Hence$,$ $C \land( \neg C) \mod B.$ Application of theorem 2.8.4 (ii) yields 
$A_1,\ldots, A_n\mod B.$\par
Conversely$,$ simply let the formula in $L$ be $C \land(\neg C).$ Then 
$A_1,\ldots, A_n\mod C \land(\neg C)$ satisfies the definition.\qed
{\bf Corollary 2.10.2.1} {\sl A nonempty set of premises $A_1,\ldots, A_n$ 
is consistent if and 
only if there exists some $B \in L$ such that $A_1,\ldots, A_n\not\mod 
B.$}\parm
Well$,$ the above definitions and theorems$,$ although they give us information 
about the concept of inconsistency$,$ DO NOT GIVE any actual way to determine 
whether or not a set of premises is consistent. The theorems simply 
say that we need to check valid consequences for infinity many formula. Not an 
easy thing to do. For over 2$,$000 years$,$ there was no way to determine 
whether or 
not a set of premises was consistent except to show that human propositional 
deduction leads to a \underbar{specific} contradiction. \parm
{\bf Theorem 2.10.3} {\sl A nonempty set of premises $A_1, \ldots,A_n$ is 
inconsistent if and only if it is not satisfiable.}\pars
Proof. Assume that $A_1, \ldots,A_n$ is inconsistent.  Thus there is some
$B \in L$ such that $A_1, \ldots,A_n\mod B \land (\neg B).$  Hence$,$
 $\mod (A_1\land  \cdots \land A_n)\to (B\land (\neg B))$ by the Deduction theorem. 
But for any assignment to atoms in $A_i$ and $B,$ $v(B \land (\neg B)) =F.$ 
Now let $A_1, \ldots,A_n$ be satisfiable. Hence $,$ there is an assigmment $v$ such that $v(A_i) = T,\ i = 1,\ldots,n.$ Extend this assignment to $v'$ so that it is an assignment for any different atoms that might appear in $B.$ Thus$,$ $v'(A_1\land  \cdots \land A_n)\to (B\land (\neg B)) = F.$ This contradicts the stated Deduction theorem. Hence$,$ inconsistency implies not satisfiable.  
  \pars
Conversely$,$ suppose that $A_1, \ldots,A_n$  is not satisfiable. Then 
for any $B\in L,$  $A_1, \ldots,A_n\mod B$ from the definition of satisfiable 
(theorem 2.8.1). 
By theorem 2.10.2$,$ $A_1, \ldots,A_n$ is inconsistent.\pars
{\bf Corollary 2.10.3.1} {\sl A nonempty set finite of premises is consistent if and
only if it is satisfiable.}\parm
Since the concept of satisfaction is dependent upon the concept of valid 
consequence$,$ the same {{\it variable substitution}} 
process (page 59) can be used for 
the atoms that appear in each premise.  Thus when you consider
 the formula variables 
as behaving like atoms$,$ when inconsistency is determined$,$ 
you have actually shown that infinitely many sets of premises are  
inconsistent. Consistency$,$ however$,$ only holds for atoms and not for the * 
type of variable substitution. As an example$,$ the set $P\to Q,\ Q$ is a 
consistent set$,$ but $P\to ((\neg R) \land R),\ (\neg R) \land R$ is 
an inconsistent set. Now$,$
it is theorem 2.10.3 and its corollary that gives a specific and FINITE method 
to determine consistency. A (large sometimes) truth-table will do the job. But 
we can also use a method similar to the forcing method of the previous 
section. \parm
\hs
{\bf Special Method 2.10.1} The idea is to try and to pick out a specific 
assignment that will satisfy a set of premises$,$ or to show that when you 
select a set of premises to be $T,$ then this $\Rightarrow$ a premise to be 
$F$ and$,$ hence$,$ the set would be inconsistent.\parm
(1) First$,$ if the premises are written in formula variables then either 
substitute atoms for the variables or$,$ at the least$,$ consider them to be 
atoms.\parm
(2) Now select a premise$,$ say $A_1$ and let $v(A_1) = T.$ If possible select a 
premise that forces a large number of atoms to have fixed values. If this is 
impossible$,$ then case studies may be necessary. \parm
(3) Now select another premise that uses the maximum number of the forced 
atoms and either show that this premise has a value $T$ or $F.$ If it has a 
value $F$ and there are no case studies then the set is inconsistent. If it is 
$T$ or you can select it to be $T,$ then the process continues.\parm
(4) If the process continues$,$ then start again with (3). Again if a premise if 
forced to be $F,$  then the set is inconsistent. This comes from the fact that 
the other premises that have thus far been used are FORCED to be $T.$ \parm
(5) If the process continues until all premises are forced to be $T,$ then 
what has occurred is that you have found an assignment that yields that the 
set is consistent.\parm
(6) If there are case studies$,$ the process is more difficult. A case study is 
produced when a $T$ value for a premise has non-fixed truth-values for the 
atoms. You must get a forced $F$ for each case study for the set to be 
inconsistent. If you get all premises to be $T$ for any case study$,$ then the 
set is consistent.\parm
(7) Better still if you'll remember what you're trying to establish$,$ then 
various short cuts can be used. For consistency$,$ we are trying to give a
metalogic argument that 
there is an assignment that gives a $T$ for all premises. Or$,$ for 
inconsistency$,$ show that under 
the assumption that some premises are $T,$ then this will  force$,$ in all 
cases$,$ some other premise to be an $F.$\pars
\hm 
{\bf Example 2.10.1} Determine whether or not the set of premises 
$(A \lor B)\to (C \land D),\ (D \lor F) \to G, \ A \lor (\neg G)$  
(written in formula variable form) is consistent. First$,$
re-express this set in terms of atoms. $(P\lor Q) \to (R \land S),\ (S \lor 
S_1) \to S_2,\ P \lor (\neg S_2).$ \parm     
(a) Let $v(P \lor (\neg S_2)) = T.$ Then there are three cases.\parm
\indent\indent (a$_1$)$,$ $v(P) = T,\ v(S_2) = F;$\par 
\indent\indent (a$_2$)$,$ $v(P) = T,\ v(S_2) = T;$\par 
\indent\indent (a$_3$)$,$ $v(P) = F,\ v(S_2) = F.$\parm 
(b) But$,$ consider the first premise. Then $v((P\lor Q) \to (R \land S)) = T.$
Now consider case (a$_1$). Then $v(P\lor Q) = T \Rightarrow v(R \land S) = T 
\Rightarrow v(R) = v(S) = T.$ Still under case (a$_1),\ \Rightarrow 
v((S\lor S_1 ) \to S_2) =F.$ But we must continue for the other cases for 
this (b) category. Now for (a$_2$)$,$ we have that $v((S\lor S_1 ) \to S_2) 
=T.$ Hence$,$ we have found assignments that yield $T$ for all premises and the 
set is consistent.\parm
{\bf Example 2.10.2}  Determine whether or not the set of premises 
$A \iff B,\ B \to C, \ (\neg C) \lor D, \ (\neg A) \to D,\ \neg D.$  
(written in formula variable form) is consistent. 
First$,$ re-express this set in terms of atoms. 
$P \iff Q,\ Q \to R, \ (\neg R) \lor S, \ (\neg P) \to S,\ \neg S.$ \parm     
(a) Let $v(\neg S) = T \Rightarrow S = F.$\parm
(b) Let $v((\neg R)\lor S) = T \Rightarrow v(R) = F.$\parm
(c) Let $v(Q \to R) = T \Rightarrow v(Q) = F.$ \parm
(d) Let $v(P\iff Q) = T \Rightarrow v(P) = F \Rightarrow v((\neg P)\to S) = 
F.$ Thus the set is inconsistent and$,$ hence$,$ we can substitute the original 
formula variables back into the set of premises. \parm
If $\mod$ is associated with ordinary propositional deduction$,$ then it should 
also mirror the propositional metalogic we are using. One of the major 
metalogical methods we are using is  called   ``proof by contradiction.'' This 
means that you assume as an additional premise the negation of the conclusion. 
Then if you can establish a contradiction of anything$,$ then the given hypotheses
can be said to ``logically'' establish the conclusion. The next theorem about 
valid consequence mirrors this notion. \parm
{\bf Theorem 2.10.4} {\sl For any set of premises$,$ $A_1,\ldots, A_n \mod B$ if 
and only if $A_1,\ldots, A_n,\neg B\mod C\land (\neg C)$ for some $C \in         
L.$}\pars   
Proof. Let $A_1,\ldots, A_n \mod B.$  If $A_1,\ldots, A_n$ is inconsistent$,$ 
then $A_1,\ldots, A_n\mod C \land (\neg C)$ for all $C \in L.$ Adding any 
other premise such as $\neg B$ does not alter this. So$,$ assume that 
$A_1,\ldots, A_n$ is consistent. Hence$,$ consider any assignment $\ass$ to 
all the atoms such that $v(A_i) = T,\ 1 \leq i \leq n$ and $v(B) = T$  and 
such an assignment exists. Thus for any such assignment $v(\neg B) = F.$ 
Consequently$,$ for any assignment $\underline{b}$ either $v(A_j) = F$ for some 
$j$ such that $1 \leq j \leq n,$ or all $v(A_i) = T$ and $v(\neg B) = F.$ 
Hence$,$ $A_1,\ldots, A_n,\neg B$ is not satisfied. Hence$,$ for some (indeed$,$ 
any) $C \in L,\ A_1,\ldots, A_n,\neg B\mod C\land (\neg C).$ \pars
Conversely$,$ let $A_1,\ldots, A_n,\neg B\mod C\land (\neg C)$ for some $C \in         
L.$ Then $A_1,\ldots, A_n,\neg B$ is inconsistent and thus given any 
assignment $\ass$ such that $v(A_i) =T$ for each $i$ such that $1 \leq i \leq 
n,$ then $v(\neg B) = F.$ Or$,$ in this case$,$ $v(B) = T.$ Consequently$,$ 
$A_1,\ldots, A_n \mod B.$ \qed
\hs\hm
\centerline{\bf EXERCISES 2.10}\parm
\noindent 1. By the special method 2.10$,$ determine if the indicated set of 
premises is consistent.\parm
(a) $A \to (\neg(B \land C)),\ (D \lor E) \to G,\ G \to (\neg (H \lor I)),\ 
(\neg C)\land E \land H.$\pars
(b) $(A\lor B) \to (C \land D),\ (D\lor E) \to G,\ A \lor (\neg G).$\pars
(c) $(A\to B)\land (C \to D),\ (B \to D)\land((\neg C)\to A),\ (E \to G)\land 
(G \to (\neg D)),\ (\neg E)\to E.$\pars
(d) $(A\to (B\land C))\land(D\to (B\land E)), ((G \to (\neg A))\land H)\to I,
 (H \to I)\to (G\land D),  \neg((\neg C)\to E).$\pars
\ss
\vfil
\eject
\noindent {\bf 2.11 Proof Theory -- General Concepts.}\parm
Actually$,$ we are investigating the propositional logic in a historically 
reversed order. The concept of {\it proof theory} properly began with Frege's 
{\it Begriffschrift} which appeared in 1879 and was at its extreme with {\it 
Principia Mathematica} written by Whitehead and Russell from 1910 --1913. 
Philosophically$,$ proof  ``theory$,$'' even though it can become very difficult$,$ 
seems to the nonmathematician to be somewhat more  ``pure'' in its foundations
since it seems to rely upon a very weak mathematical foundation if$,$ indeed$,$ it 
does have such a foundation. However$,$ to the mathematician this could indicate 
a kind of weakness in the basic tenants associated with the concept of the 
{{\it formal proof.}} Notwithstanding these philosophical differences$,$ basic 
proof theory can be investigated without utilizing any strong mathematical 
procedures. Since we are studying mathematical logic$,$ we will not restrict our 
attention to the procedures used by the formal logician but$,$ rather$,$ we are 
free to continue to employ our metamathematical concepts.\pars
The concept of  {\it formalizing}  a logical argument goes back to Aristotle.
But it again took 
over 2$,$000 years before we could model$,$ in another formal way$,$ the various 
Aristotle formal logical arguments. Proof theory first requires a natural 
language algorithm that gives simple rules for writing a formal proof. 
The rules that are used can be checked by anyone with enough training to apply 
the same rules. But$,$ in my opinion$,$ the most important part of proof theory
is the fact that we view the process externally while the philosophical 
approach views the process internally. (We study the forest$,$ many philosophers 
study the trees$,$ so to speak.) For example$,$ we know that every propositional 
formula can be replaced by an equivalent formula expressed only in the 
connectives $\neg$ and $\to.$ Thus we can simplify our language construction 
process considerably.\parm\hs
{\bf Definition 2.11.1} (The Formal Language $L^\prime$.)  \parm
(1) You start at step (2) of definition 2.2.3.\parm
(2) Now you go to steps (3) --- (6) of definition 2.2.3 but only use the 
$\neg$ and $\to$ connectives. \pars
\hm
Notationally$,$ the language levels constructed in definition 2.11.1 are denoted 
by $L^\prime_n$ and $B  \in L^\prime$ if and only if there is some $n  \in 
\nat$ such that $B \in L^\prime_n.$ Also don't forget that level $L^\prime_n$ 
contains all of the previous levels. Clearly $L^\prime$ is a proper subset of 
$L.$ Also$,$ all the definitions of size$,$ and methods to determine size$,$ and 
simplification etc. hold for $L^\prime.$  \pars
Now proof theory does not rely upon any of the concepts of truth---falsity$,$ 
what will or will not occur in reality$,$ and the like we  modeled previously. 
All the philosophical 
problems associated with such concepts are removed. It is pure$,$ so to speak. 
But as mention it relies upon a set of \underbar{rules} that tells us in 
a step-by-step manner$,$ hopefully understood by all$,$ exactly what formulas
we are allowed 
to write down and exactly what manipulations we are allowed to perform with 
these formulas. It's claimed by some that$,$ at the least$,$ it's mirroring some 
of the procedures the human beings actually employ to obtained a logical 
conclusion from a set of premises. Proof theory is concerned with the 
notion of logical argument. It$,$ of course$,$ is also related to the metalogical 
methods the human being uses. We start with the rules for a {\it formula
theorem} that is obtained from an empty set of premises. The word  
``theorem'' as used here does not mean the metalanguage things called
``Theorems'' in the previous sections. It will mean a specially located 
formula.  We can use formula variables and connectives to write formula 
schemata or schema. Each schemata would represent infinitely many specific 
selections from $L^\prime.$ We will simply call these {\it formulas}
 and use the 
proper letters to identify them as formula variables.\parm\hs
{\bf Definition 2.11.2} (A formal proof of a theorem.)\parm
(1) A {{\it formal proof}} contains two  FINITE COLUMNS$,$ one of formulas 
$B_i$$,$ the other column 
stating the reason $R_i$ you placed the formula 
in that step.\parm
\line{\hskip0.75in$B_1$\hskip 2.0in $R_1$\hfil}
\line{\hskip0.75in$B_2$\hskip 2.0in $R_2$\hfil}
\line{\hskip0.75in ...\hskip 2.1in ...\hfil}
\line{\hskip0.75in ...\hskip 2.1in ...\hfil}
\line{\hskip0.75in ...\hskip 2.1in ...\hfil}
\line{\hskip0.75in ...\hskip 2.1in ...\hfil}
\line{\hskip0.75in$B_n$\hskip 2.0in $R_n$\hfil}
\parm
\noindent The formula used represent any members of $L^\prime$ and each  must 
be obtained in the 
following manner. \parm
(2) A step $i$ in the formal proof corresponds to a specific formula $B_i.$
It can be a specific instance of one of the following axioms. For {\it any}
formula $A,\ B,\ C\,$ where as before $A,\ B,\ C$ are formula variables (we note that because of the way we have of writing subformula$,$ the $A,\ B,\ C$ can also be considered as expressed in formula variables)$,$ write\parm
\indent\indent $P_1:\ \ \ A \to (B \to A),$ \parm 
\indent\indent $P_2:\ \ \ (A  \to (B \to C)) \to ((A\to B)\to (A \to C)),$ \parm  
\indent\indent $P_3:\ \ \ ((\neg A) \to (\neg B)) \to (B \to A).$ \parm 
(3) Formula can only be obtained by one other procedure. It's called 
{{\it modus ponens}} and is abbreviated by the symbol MP. It's the only way 
that we 
can obtain a formula not of the type in (2) and is considered to be our one 
rule of logical inference. \pars
Step $B_j$ is obtained if there are two PREVIOUS steps $B_i$ and 
$B_k$ such that \parm
\indent\indent (i) $B_i$ is of the form $A,$\parm
\indent\indent (ii) $B_k$ is of the form $A\to B.$  Then \parm
\indent\indent (iii) $B_j$ is of the form $B$\parm
\noindent Note that the order of the two needed previous steps $B_i$ and $B_k$ 
is not specified. All that's needed is that they come previous to $B_j.$ \parm
(4) The  {\it last step} in the finite column is called a (formal) {\it 
theorem} and the total column of formulas is called the (formal) {\it proof of 
theorem.}\parm
(5) If a proof of a theorem $E$ exists$,$ then this is denote by $\vd E.$\pars
\hm
Since we are using variables$,$ any proof of a theorem is actually a proof for 
any formula you substitute for the respective variables (i.e. infinitely many 
proofs for specific formula.) In our examples$,$ problems and the like$,$ I won't 
repeat the statement  ``For any $A  \in L^\prime.$'' The fact that I'm using 
variable symbols will indicate this. \parm
{\bf Example 2.11.1} \parm
Show that $\vd A \to A.$\parm
Proof.\pars
\line{\indent (1) $(A \to ((A\to A)\to A)) \to ((A\to (A \to A)) \to (A \to 
A))$\hfil$P_2$}
\line{\indent (2) $A \to((A \to A) \to A)$\hfil $P_1$}  
\line{\indent (3) $(A \to (A\to A)) \to (A \to A)$\hfil $MP(1,2)$}
\line{\indent (4) $A \to (A \to A)$\hfil $P_1$}
\line{\indent (5) $A \to A$\hfil $MP(3,4)$}\parm
The MP step must include the previous step numbers used for the MP step. Now 
any individual who can follow the rules can check that this last formula was 
obtained correctly. Also$,$ there are many other proofs that lead to the same 
result. The idea is the same as used to demonstrate a Euclidean geometry 
proof$,$ but there is no geometric intuition available. Now at any point in some
other proof if you should need the statement $A \to A$ as one of you steps 
(of course$,$ you can use another symbol for $A$) then you could substitute the 
above finite proof just before the step is needed and re-number all steps. BUT$,$ 
rather than do this all you need to do is to write the symbol $\vd A \to A$ 
since this indicates that a proof exists. By the way$,$ there is stored at a university 
in Holland literally thousand upon thousands of formal proofs. However$,$ in 
this course almost all of the formal proofs exhibited will be needed to 
establish our major results. These are to show that the formal proof of a formal 
theorem is equivalent to the modeling concept we call validity.  \pars
One of the major results observed by Aristotle was the logical argument he 
called {{\it hypothetical syllogism}} or HS. We next use the above procedure 
to establish this. But it's a method of introducing a formal proof and 
our  result is a metatheorem proof. \parm
{\bf Theorem 2.11.1} {\sl Assume that you have two steps in a 
possible proof of the 
form (i) $A \to B$ and (j) $B \to C.$ Then you can write down at a larger step 
number (iii) $A \to C.$}\pars
Proof.\parm 
\line{\indent (1) \leaderfill}
\line{\indent ...\hfil}
\line{\indent ...\hfil}
\line{\indent (i) $A \to B$\hfil $R_i$}
\line{\indent ...\hfil}
\line{\indent ...\hfil}
\line{\indent (j) $B \to C$\hfil $R_j$}
\line{\indent (j + 1) $(B \to C) \to (A \to (B \to C))$\hfil $P_1$}
\line{\indent (j + 2) $(A \to (B \to C)$\hfil $MP(j, j+1)$}
\line{\indent (j + 3) $(A\to (B\to C)) \to ((A\to B)\to (A\to C))$\hfil $P_2$}
\line{\indent (j + 4) $(A\to B) \to (A \to C)$\hfil$MP(j+2,j+3)$}
\line{\indent (j + 5) $A \to C$\hfil$MP(i,j+4)$}\parm
\noindent Whenever we use HS$,$ the reason is indicated by the symbol $HS(i,j)$
in the same manner as is done for MP. 
IMPORTANT FACT. In the above metatheorem that yields HS$,$ please note that 
the step (j) can come before step (i). \parm
{\bf Example 2.11.2} \parm
\noindent $\vd(\neg (\neg A)) \to A$\parm
\line{\indent (1) $(\neg(\neg A))\to ((\neg(\neg(\neg(\neg A))))\to (\neg(\neg 
A)))$\hfil $P_1$}
\line{\indent (2) $( (\neg(\neg(\neg(\neg A))))\to (\neg (\neg A))) \to ((\neg 
A)\to 
(\neg (\neg (\neg A))))$ \hfil ....}
\line{\indent (3) $(\neg (\neg A)) \to ((\neg A) \to (\neg (\neg (\neg 
A))))$\hfil $HS(1,2)$}
\line{\indent (4) $((\neg A) \to (\neg (\neg (\neg A))))\to ((\neg (\neg 
A))\to A)$\hfil ....}
\line{\indent (5) $(\neg (\neg A)) \to ((\neg (\neg A)) \to A)$\hfil $HS(\ ,\ )$}
\line{\indent (6) $((\neg (\neg A))\to ((\neg (\neg A)) \to A))\to$\hfil} 
\line{\indent\indent $(((\neg(\neg A)) \to (\neg (\neg A))) \to ((\neg (\neg 
A) \to A)))$\hfil $P_2$}
\line{\indent (7) $((\neg (\neg A))\to (\neg (\neg A))) \to ((\neg (\neg A)) \to
A)$\hfil $MP(\ ,\ )$}  
\line{\indent (8) $((\neg (\neg A))\to (((\neg(\neg A))\to (\neg(\neg A))) \to 
(\neg(\neg A)))) \to (((\neg(\neg A))\to$\hfil}
\line{\indent\indent $((\neg(\neg A)) \to (\neg (\neg A))))\to ((\neg(\neg A)) 
\to (\neg (\neg A))))$\hfil $P_2$}
\line{\indent (9) $(\neg(\neg A)) \to (((\neg (\neg A)) \to (\neg(\neg A))) 
\to (\neg (\neg A)))$ \hfil ....}
\line{\indent (10) $((\neg(\neg A)) \to ((\neg(\neg A)) \to (\neg(\neg A)))) 
\to$\hfil} 
\line{\indent\indent $((\neg(\neg A)) \to (\neg(\neg A)))$\hfil $MP(\ ,\ )$}
\line{\indent (11) $(\neg(\neg A)) \to ((\neg(\neg A)) \to (\neg(\neg A)))$\hfil ....}
\line{\indent (12) $(\neg(\neg A)) \to (\neg(\neg A))$\hfil ....}
\line{\indent (13) $(\neg(\neg A)) \to A$\hfil $MP(\ ,\ )$} \parm
\hs\hm
\centerline{\bf EXERCISES 2.11}\parm
\noindent (1) Rewrite the proof of $\vd (\neg(\neg A)) \to A$ as it appears in 
example 2.11.2 filling in the missing reasons. Some may be examples. \parm
\noindent (2) Give formal proofs of the next two theorems. 
You may use any previously 
established $\vd,$ or method. You must state all your reasons.\parm
(a) $\vd A \to (\neg(\neg A))$.\parm
(b) $\vd (\neg B) \to (B \to A).$\pars
\ss
\noindent {\bf 2.12 Demonstrations$,$ Deductions From Premises.}\parm
Given a nonempty set premises $\Gamma$ (not necessarily finite)$,$ then what is 
the classical procedure employed to deduce a formula $A$ from $\Gamma$? 
Actually$,$ we are only allowed to do a little bit more than we are allowed to 
do in order to give a proof of a formal theorem. The next definition gives the 
one
additional rule that is assumed to be adjoined to definition 2.11.2 (formal 
proof of a  theorem) in order to demonstrate that a formula is deducible 
from a given set of premises.\parm
\hs
{\bf Definition 2.12.1} (Deduction from premises.)  A formula $B \in L^\prime$ is 
said to be {{\it deduced from a set of $\Gamma$}} or {{\it a consequence of 
$\Gamma$}}$,$ where $\Gamma$ is a not necessary finite (but possibly empty) set 
of premises$,$ if $B$ is the last step in a FINITE column of steps and reasons 
as described in definition 2.11.2 where you are allowed one additional 
rule.\parm
(1)  You may write down as any step a single instance of a formula as it 
appears in $\Gamma,$ where the premises are written in formula variable 
form. The reason given is {\it premise}.\parm
(2) The finite column with reasons is called a {{\it demonstration}}.\pars
\hm
Observe that whenever we have written a demonstration$,$ because of the use of 
variables$,$ we have actually written infinitely many demonstrations. Deduction 
from a set of premises is what one usually considers when  one uses the 
terminology  ``a logical argument.'' Now a proof of a theorem is a 
demonstration from an empty set of premises. When $B$ is deduced from the set 
$\Gamma,$ then this is symbolized by writing $\Gamma \vd B.$ The very 
straightforward definition 2.12.1 yields the following simple results.\parm
{\bf Theorem 2.12.1} {\sl Assume that $\Gamma$ is any set of formula.\parm
(a) If $A \in \Gamma,$ or $A$ is an instance of an axiom$,$ then $\Gamma \vd 
A.$\parm
(b) If $\Gamma \vd A$ and $\Gamma \vd A \to B,$ then  $\Gamma \vd B.$\parm
(c) If $\vd A,$ then $\Gamma \vd A.$\parm
(d) If $\Gamma$ is empty and $\Gamma \vd A,$ then $\vd A.$\parm
(e) If $\Gamma \vd  A$ and $D$ is any set of formula$,$ then $\Gamma \cup D \vd 
A.$\parm
(f) If $\Gamma \vd A,$ then there exists some finite subset $D$ of $\Gamma$ 
such that $D \vd A.$}\pars
Proof. (A) Let $A \in  \Gamma.$ Then one step consisting of the line
``(1) $A.\ .\ .\ .\ .\ .\ .\ $ Premises'' is a demonstration for $A.$ If $A$ is 
an axiom$,$ then one step from definition 2.11.2 yields a demonstration.\parm
(b) Let $B_1, \ldots, B_m$ be the steps in a demonstration that $A$ is 
deducible from $\Gamma.$ Then $B_m = A.$ Let $C_1,\ldots, C_k$ be a 
demonstration that $A \to B$ is deducible from $\Gamma.$ Now after writing the 
steps $B_{m+i} = C_1,\ i = 1,2,\ldots,k$ add the step $C_{k+1}: B 
.\ .\ .\ .\ .\ .\ .\ MP(m,m+k).$ Then this yields a demonstration that $B$ is  
deducible from $\Gamma.$\parm
(c) Obvious.\parm
(d) Since $\Gamma$ is empty$,$ rule (1) of definition 2.11.2 has no application. 
Hence$,$ only the rules for a proof of a formal theorem have been used and this result 
follows.\parm
(e) We have$,$ possibly$,$ used formulas from $\Gamma$ to deduce $A.$ We have$,$ 
possibly$,$ used formulas from $\Gamma \cup D$ to deduce $A$.\parm
(f) Assume that $\Gamma\vd A.$ Now let $D$ be the set of all premises that have 
been utilized as a specific step marked ``premise'' that appears in the 
demonstration $\Gamma \vd A.$ Obviously$,$ we may replace $\Gamma$ with $D$ and 
have not altered the demonstration. \qed
{\bf Example 2.12.1} To show that $(\neg(\neg A)) \vd A.$\parm
\line{(1) $(\neg(\neg A))$\hfil Premise}
\line{(2) $(\neg(\neg A)) \to ((\neg (\neg (\neg (\neg A)))) \to (\neg(\neg 
A)))$\hfil .....}
\line{(3) $ (\neg (\neg (\neg (\neg A)))) \to (\neg(\neg A))$\hfil $MP(\ ,\ )$}
\line{(4) $((\neg(\neg(\neg(\neg A))))\to (\neg(\neg A))) \to ((\neg A) \to 
(\neg(\neg(\neg A))))$\hfil .....}
\line{(5) $(\neg A) \to (\neg(\neg(\neg A)))$\hfil $MP(\ ,\ )$}
\line{(6) $((\neg A) \to (\neg(\neg(\neg A))))\to ((\neg(\neg A)) \to A)$\hfil 
 .....}
\line{(7) $(\neg(\neg A)) \to A$\hfil $MP(\ ,\ )$}
\line{(8) $A$\hfil $MP(\ ,\ )$}
\parm
\baselineskip 12pt
Notice that example 2.12.1 and example 2.11.3  show that $\vd (\neg(\neg 
A)) \to A$ involve similar formulas. Moreover$,$ we have not used HS in example 
2.12.1. Indeed$,$ if we had actually made the HS substitution into steps of 
example 2.11.3$,$ then we would have at the least 21 steps in example 2.11.3. 
Even though we have not considered as yet any possibilities that the 
semantical methods which have previously been employed might be equivalent in 
some sense to the pure proof methods of this and section 2.11$,$ it would 
certainly be of considerable significance if there was some kind of deduction 
theorem for our proof theory. We can wonder if it might be possible to 
substitute for $\mod$ in the metatheorem $\mod (\neg (\neg A)) \to A$ if and 
only if $ (\neg (\neg A))\mod  A$ the symbol $\vd$? Further$,$ it's possible to
stop at this 
point and$,$ after a lot of formal work$,$ introduce you to the actual 
type of deduction process used to create a universe. But to avoid the formal 
proofs of the needed special formal theorems$,$ I'm delaying this introduction 
until we develop the complete equivalence of valid formulas and formal 
theorems. I'll answer the previous question by establishing one of the most 
powerful procedures used to show this equivalence.   
\parm\baselineskip 14pt
\hs\hs
\centerline{EXERCISES 2.12}\parm
\noindent (1) Rewrite example 2.12.1 filling in the missing reasons.\parm
\noindent (2) Complete the following deductions from the indicated set of 
premises. Write down the missing steps and/or reasons.\parm
\noindent Show that $A \to B,\ B \to C \vd A\to C.$\parm
\line{(1) $(B \to C) \to (A \to (B \to C))$\hfil $P_1$}
\line{(2) $B \to C$\hfil .....}
\line{(3) \hfil $MP(1,2)$}
\line{(4) $(A\to (B \to C)) \to ((A \to B) \to (A \to C))$\hfil .....}
\line{(5) $(A \to B) \to (A \to C)$\hfil .....}
\line{(6) \hfil Premises}
\line{(7) $A \to C$\hfil .....}\parm
\noindent (b) Use 2(b) of exercise 2.11$,$ to show in THREE steps that
$(\neg A) \vd A \to B.$\parm
\noindent (c) Use the fact that $\vd (B\to A) \to ((\neg A) \to (\neg B))$ and any 
previous results to show that $\neg(A \to B) \vd B \to A.$\parm
\line{(1) $\neg(A \to B)$\hfil Premise}
\line{(2) $\vd (B\to (A \to B)) \to ((\neg (A \to B)) \to (\neg B))$\hfil Given}
\line{(3) \hfil $P_1$}
\line{(4) \hfil .....}
\line{(5) $(\neg B)$\hfil .....}
\line{(6) $(\neg B) \to (B \to A)$\hfil .....}
\line{(7) $B \to A$\hfil .....}\pars
\ss
\noindent {\bf 2.13 The Deduction Theorem.} \parm
The Deduction Theorem is so vital to logic that some mathematicians$,$ such as 
Tarski$,$  
use it as 
a basic axiom for different logical systems. The following metaproof was first 
presented in 1930 and it has been simplified by your present author. This 
simplification allows it to be extended easily to cover other types of logical 
systems. What this theorem does is to give you a step by step process to 
change a formal proof into another formal proof. \parm
{\bf Theorem 2.13.1} (The Deduction Theorem) {\sl Let $\Gamma$ a collection of 
formula from $L^\prime$ written in formula variable form.
Assume that formula variables $A ,B$ represent arbitrary members of 
$L^\prime.$ Let $\Gamma \cup \{A\}$ be the set of premises the contains and 
only contains the members of $\Gamma$ and $A.$ Then $\Gamma \cup \{A\} \vd B$ 
if and only if $\Gamma \vd A \to B.$}\pars
Proof. For the necessity$,$ assume that $\Gamma \cup \{A\} \vd B.$
Then there exists a finite $A_1,\ldots, A_n \in \Gamma$ such that
$A_1,\ldots, A_n,A \vd B.$ (Of course no actual member of $A_1,\ldots, A_n$ 
need to be used.) Assume that every step has been written with the 
reasons as either a premise$,$ or an axiom$,$ or MP. \pars
(a) Renumber$,$ say in red$,$ all of the steps that have reasons stated as 
premises 
or axioms. Assume that there are k such steps. Note that there must be at 
least three  steps for MP to be a reason$,$ and then it could only occur in 
step 3. We now construct a new demonstration. \parm
\indent\indent Case $(m = 1).$ There is only one case. Let the formula be 
$B_1.$ Hence$,$ $B_1 = B.$ \parm
\indent\indent\indent Subcase (1). Let $B_1$ be a premise $A_j.$ Assume that
$A_j \not=A.$ Now \underbar{keep}  $B_1 = A_j$ in the new demonstration and
insert the two 
indicated steps (2) (3). \pars
\line{\hskip 1.0in (1) $B_1 = A_j$\hfil premise}
\line{\hskip 1.0in (2) $A_j \to (A \to A_j)$\hfil $P_1$}
\line{\hskip 1.0in (3) $A \to A_j = A \to B$\hfil $MP$}\par
Now notice that the $A$$,$ in these three steps does not appear as a single 
formula in a step. \parm
\indent\indent\indent Subcase (2). Assume that $B_1 = A.$ Then 
 \underbar{insert} the five steps that yields 
$\vd A \to A$ (Example 2.11.1.) into the new demonstration. Now \underbar{do 
not place} the step $B_1 = A$ in the new demonstration.
Again we have as the last step $A\to A$ and $A$ does not appear as any step 
for in this new demonstration. \parm
\indent\indent\indent Subcase (3). Let $B_1 = C$ be any axiom. Simply 
repeat subcases (1) (2) and we have as a last step $A \to C$ and $A$ does not appear 
in any step. \pars
Now we do the remaining by induction. Suppose that you have 
used m of the original renumbered steps to thus far construct the new 
demonstration.  
Now consider step number $B_{m+1}.$ Well$,$ just apply the same subcase 
procedures that are used under case $m = 1.$ The actual induction hypothesis 
is vacuously employed (i.e. not employed). Thus by induction$,$ for any $k$ 
steps we have found a  way to use all non-MP steps of the original 
demonstration to construct a new demonstration in such a manner that no 
step in the new demonstration has$,$ thus far$,$  
only $A$ as its formula. Further$,$ all original non-MP steps $B_i$ are 
replaced in the new demonstration by 
$A \to B_i.$ Thus if the last formula $B_n$ in the original demonstration 
is a non-MP step$,$ 
then $B_n=B,$ and the last step in the new demonstration is $A \to B.$  We 
must now use the original MP steps to continue the new demonstration 
construction since $A$ might be in the original demonstration 
an MP step or the last step of 
the original demonstration might be an MP step. \pars
\underbar{After} using all of the original non-MP steps$,$ we now renumber$,$ for 
reference$,$ all the original MP steps. We now define by induction on the 
number $m$ of MP steps a 
procedure which will complete our proof.\parm
\indent Case $(m = 1)$ We have that there is only one application of 
MP. Let $B_j$ be this original MP step. Then two of the original 
\underbar{previous} steps in the 
original demonstration$,$ say $B_g$ and $B_h,$ were employed and $B_h = B_g \to 
B_j.$ We have$,$ however$,$ in the new demonstration used these to construct
the formula and we have not used the original $B_j$ step. Now use the two 
new steps 
(i) $A \to B_g$$,$ and (ii) $A \to (B_g \to B_j)$ and insert 
\underbar{immediately after}  
 new step (ii)  the following:\pars
\line{\hskip 1.0in (iii) $(A \to (B_g \to B_j))\to ((A \to B_g) \to (A \to 
B_j))$\hfil $P_2$}
\line{\hskip 1.0in (iv) $(A\to B_g) \to (A \to B_j)$\hfil $MP(ii,iii)$}
\line{\hskip 1.0in (v) $A\to B_j$ \hfil $MP(i,iv)$}\pars
Now \underbar{do not include} the  original $B_j$ in the new demonstration 
just as in case where $B_j = A.$\pars
Assume the induction hypothesis that we have $m$ of the 
original MP steps used in the original 
demonstration now 
altered so that the appear as $A \to \ldots$. \pars
\indent Case $(m+1).$ Consider MP original step $B_{m+1}.$ 
Hence$,$ prior steps in the 
old demonstration$,$ say $B_g,\ B_h = B_g \to B_{m+1}$ are utilized to obtain  
the $B_{m+1}$ formula. However$,$ all the original steps up to but not including 
the $B_{m+1}$ have been replaced by a new step in our new demonstration 
in such a manner they now look 
like (vi) $A \to B_g$ and (viii) $A \to (B_g \to B_{m+1})$ no matter how these 
step were originally obtained. (The induction hypothesis is necessary at this 
point.)
Now follow the exact same insertion process as in the case $m = 1.$ This 
yields a new step $A \to B_{m+1}$ constructed from the original $B_{m+1}.$
Consequently$,$ by induction$,$ we have defined a procedure by which all the 
original steps $B_k$ have been used to construct a new demonstration
and if $B_n$ was one of the original formula$,$ then it now appears in the new 
demonstration as $A \to B_n.$ Further$,$ no step in the new demonstration  
is  the single formula $A$. This yields a newly 
constructed demonstration that $A_1,\ldots,A_n \vd A \to B.$ \parm
For the sufficiency$,$ assume that $\Gamma \vd A \to B.$ The final step is 
$B_k = A \to B.$ Now add the following two steps. $B_{k+1} = A$ (premises)$,$ 
and $B_{k+1}=B,\ MP(k,k+1).$ This yields $\Gamma \cup \{A\}\vd B.$ \qed
{\bf Corollary 2.13.1.1} {\sl $A_1,\ldots,A_n \vd B$ if and only if
$\vd (A_1 \to (A_2 \to \cdots (A_n \to B) \cdots )).$}\parm
Now in the example on the next page$,$ I follow the rules laid out within the proof of the 
Deduction Theorem. The formula from the original demonstration used as a step 
in the new construction  are represented in roman type. Formula from the 
original demonstration used to obtain the new steps BUT not included in the 
new demonstration are in Roman Type$,$ BUT  are placed between square brackets 
[\quad ]. In the case $\vd A \to A,$ I will not include all the steps. \parm
{\bf Example 2.13.1} \parm
\centerline{(I) $A \to B, \ B \to C,\ A \vd C.$}\parm
\line{(1) $A \to B$\hfil Premise}
\line{(2) $B \to C$\hfil Premise}
\line{(3) $A$ \hfil Premise}
\line{(4) $B$ \hfil $MP(1,3)$}
\line{(5) $C$ \hfil $MP(2,4)$}\parm
We now construct from (I) a new demonstration that\par
\centerline{(II) $A \to B, \ B \to C, \vd A \to C.$}\parm
\line{(1) ${\rm A \to B}$\hfil Premise}
\line{(2) $(A \to B) \to (A\to (A \to B))$\hfil $P_1$}
\line{(3) $A \to (A \to B)$\hfil $MP(1,2)$}
\line{(4) ${\rm B\to C}$\hfil Premise}
\line{(5) $(B\to C) \to (A \to (B\to C))$\hfil $P_1$}
\line{(6) $A \to (B\to C)$ \hfil $MP(4,5)$}
\line{(7) $\vd A \to A$\hfil Example 2.11.1}
\line{${\rm [B_g = A,\ B_j = B]}$\hfil }
\line{(8) $(A \to (A \to B))\to ((A \to A) \to (A \to B))$\hfil $P_2$}
\line{(9) $(A \to A) \to (A \to B)$\hfil $MP(3,8)$}
\line{(10) $A \to B$ \hfil $MP(7,9)$}
\line{${\rm [B_g =B,\ B_j =C]}$\hfil }
\line{(11) $(A \to (B \to C)) \to((A \to B) \to (A \to C))$\hfil $P_2$}
\line{(12) $(A \to B) \to (A \to C)$\hfil $MP(6,11)$}
\line{(13) $A \to C$\hfil $MP(10,12)$}\parm
I hope this example is sufficient. But note that although this gives a 
demonstration$,$ it 
need not give the most efficient demonstration.\parm
\hs\hm
\centerline{EXERCISES 2.13}\parm
\noindent 1. Give a reason why we should NOT use the Deduction theorem  
as a reason that from $A \vd A,$ we have $\vd A \to A.$\parm
\noindent 2. Complete the following formal proofs of the indicated theorems
 by application of the Deduction Theorem in order to insert premises.
You will also need to insert$,$ in the usual manner$,$ possible $\vd$ statements
obtained previously. Please give reasons.\parm
\hs
\line{(A) Show that $\vd (B\to A) \to ((\neg A)\to (\neg B))$\hfil}\parm
\line{(1) $B\to A$ \hfil Premise and D. Thm.}\smallskip
\line{(2) $\vd (\cdots\cdots)\to B$\hfil Ex. 2.12.1 and D. Thm.}\smallskip
\line{(3) $(\neg (\neg B)) \to A$\hfil $HS(\ ,\ )$}\smallskip
\line{(4)\hfil}\smallskip
\line{(5) $(\neg(\neg B)) \to (\neg(\neg A))$\hfil .....}\smallskip
\line{(6) $((\neg(\neg B))\to (\cdots\cdots))\to ((\neg A) \to (\neg 
B))$\hfil .....}\smallskip
\line{(7) $\cdots\cdots\cdots$\hfil .....}\smallskip
\line{(8)  $\cdots\cdots\cdots$\hfil D. Thm.}\medskip
\hs
\line{(B) Show that $\vd ((A \to B)\to A) \to A$\hfil}\parm
\line{(1) $(A \to B) \to A$\hfil .....}\smallskip
\line{(2) $\cdots\cdots\cdots$\hfil D.Thm. and Ex 2.12.2b}\smallskip
\line{(3) $(\neg A) \to A$\hfil .....}\s
\line{(4) $(\neg A)\to (\neg(\neg((\neg A) \to A)) \to (\neg A))$\hfil ....}\s
\line{(5) $(\neg(\neg((\neg A) \to A)) \to (\neg A)) \to \cdots $\hfil $P_3$}\s
\line{(6) $(\neg A) \to (A \to \cdots ((\neg A)\to \cdots\cdots ))$\hfil .....}\s
\line{(7) $((\neg A)\to (A \to (\neg((\neg A) \to A))))\to(((\neg A)\to A)\to$\hfil}
\line{\hskip 1.5in $((\neg A) \to (\neg((\neg A) \to A))))$\hfil .....}\s
\line{(8) $\cdots\cdots\cdots$\hfil $MP(\ ,\ )$}\s
\line{(9) $\cdots\cdots\cdots$\hfil $MP(3 ,8)$}\s
\line{(10) $((\neg A) \to (\neg(\cdots\cdots)))\to (((\neg A) \to \cdots) \to 
\cdots\cdots )$\hfil $P_3$}\s
\line{(11) $\cdots\cdots\cdots$\hfil $MP(\ ,\ )$}\s
\line{(12) $A$ \hfil .....}\s
\line{(13) $\vd\cdots\cdots\cdots$ \hfil D. Thm.}\parm\hs
\noindent [3] Use the following demonstration 
that $A \to B,\ A \vd B$ and construct$,$ as in example 2.13.1$,$  
by use of the procedures within the metaproof of the Deduction Theorem a 
formal demonstration that $A \to B  \vd A \to B.$ [It's obvious that this 
will not give the most efficient demonstration.]\parm 
\hm
\centerline{$A \to B,\ A \vd B$}\s
\line{(1) $A$\hfil Premise}\s
\line{(2) $A \to B$ \hfil Premise}\s
\line{(3) $B$ \hfil $MP(1,2)$}\s
\ss
\noindent {\bf 2.14 Deducibility Relations or $\mod$ implies $\vd$ 
almost.}\parm
The Deduction Theorem for formal demonstrations seems to imply that the 
semantical concept of $\mod$ and the proof-theoretic concept of $\vd$ are 
closely related for they share many of the same propositional facts.
In this section$,$ we begin the study which will establish exactly how these two 
seeming distinct concepts are related. Keep in mind$,$ however$,$ that $\mod$ 
depends upon the mirroring of the classical truth-falsity$,$ will occur-won't 
occur notion while $\vd$ is dependent entirely upon the strict formalistic 
manipulation of formulas. We are in need of two more formal proofs.\parm
{\bf Example 2.14.1} $\vd A \to ((A \to B ) \to B).$ \pars      
\line{(1) $A$\hfil Premise}
\line{(2) $A \to B$ \hfil Premise}
\line{(3) $B$ \hfil $MP(1,2)$}
\line{(4) $\vd A \to ((A \to B ) \to B)$\hfil D. Thm.}\parm
{\bf Example 2.14.2} $\vd A \to ((\neg B) \to (\neg(A \to B))).$ \s
\line{(1) $\vd A \to ((A \to B ) \to B)$\hfil Ex. 2.14.1}
\line{(2) $((A \to B) \to B) \to ((\neg B)\to (\neg (A \to B)))$\hfil Exer.
2.13 (A)}
\line{(3) $A \to ((\neg B) \to (\neg(A \to B)))$\hfil $HS(1,2)$} \pars
The theorem we will need is that of Example 2.14.2. As done previously$,$ the 
formal proofs or demonstrations are done in formula variables. They hold for any 
formula consistently substituted for the variables. In the semantics sections,
most but not all deduction-type concepts such as validity and valid 
consequences also hold for formula variables. Non-validity was 
a notion that did not hold in formula variable form. In what follows$,$ we again 
assume that we are working in formula variables. Of course$,$ they also hold for 
atoms substituted for these variables. We use the notation for 
formula variables. Let $A$ be a formula written in the following 
manner. The formula variables $A_1, \ldots, A_n$      
and only these formula variables are used to construct $A$ with the $L^\prime$ 
propositional connectives $\neg ,\ \to.$ Thus $A$ is written in 
formula variables. Since $L^\prime \subset 
L,$ the truth-table concept can be applied to $L^\prime.$
 Now let 
$\ass$ be an assignment to the atoms that would appear in each $A_i,\ i = 
1,\ldots, n$ when specific formula are substituted for the formula 
variables.\parm
\hs
{\bf Definition 2.14.1}  For each $i,$ we define a formula $A_i^\prime$ as 
follows:\par
(i)  if $v(A_i) = T,$ then $A_i^\prime = A_i.$\pars
(ii) If $v(A_i) = F,$ then $A_i^\prime =(\neg A_i).$\pars
(iii) If $v(A) = T,$ then $A^\prime = A.$\pars
(iv) If $v(A) = F,$ then $A^\prime = (\neg A).$
\pars   \hm

In the truth-table for $A,$ reading from left to right$,$ you have all the atoms 
then the formula $A_i$ and finally the formula $A.$ Then we calculate the $A$ 
truth-value from the $A_i,$ with possibly additional columns if needed. 
\pars
\hs
{\bf Definition 2.14.2} (Deducibility relations.) For each row $j$ of the 
truth-table$,$ there are truth-values for each of the $A_i$ and $A.$ 
These generate
the formula $A_i^\prime$ and $A^\prime,$ where the $A_i$ are the formula that 
comprise $A.$\pars
(i) The $j$th {\it Deducibility relation} is $A^\prime_1,\ldots, A^\prime_i 
\vd A^\prime.$ \pars
\hm
{\bf Theorem 2.14.1.} {\sl Given $A_1,\ldots,A_n$ and $A$ as defined  
above. Then for any row of the corresponding truth-table for $A$ generated by 
the truth-values for $A_1,\ldots,A_n,$ we have that $A_1^\prime, 
\ldots,A_n^\prime \vd A^\prime.$}\pars
Proof. First, assume that each $A_i$ is an atom and $A$ is expressed in atomic form.  This allows for induction and for all the possible 
truth-values for non-atomic $A_i$. We now do an induction proof on the 
size of $A.$\pars
Let ${\rm size}(A) = 0.$ The $A = A_i$ for some $i$ (possible more than one). As 
we know $A_1, \ldots A_n \vd A_i$ for any $i,\ i = 1,\ldots,n.$ Thus the 
result follows by leaving the original $A_i$ or replacing it with $\neg A_i$ 
as the case may be. \pars
Assume the induction hypotheses (in strong form) that the result holds for any 
formula $A$ of size $\leq m$ where $m >0.$ Let ${\rm size}(A) = m+1.$ \pars
There are two cases. (i) The formula $A =  \neg B,$ or 
case (ii) $A = B \to C,$ where ${\rm size}(B),\ {\rm size}(C) \leq m.$ 
Note that by the induction hypothesis if $q_1,\ldots, q_k$ are the atoms in 
$B,$ then $q_1^\prime,\ldots,q_k^\prime \vd B^\prime.$ But adding any other 
finite set of atoms does not change this statement. Hence,
$A_1^\prime,\ldots, A_n^\prime \vd B^\prime.$ In like manner,
$A_1^\prime, \ldots, A_n^\prime \vd C^\prime.$ \pars
Case (i). Let $A =(\neg B)$ Then (a), suppose that $v(B) = T.$ Then
$B^\prime = B$ and $v(A) = F.$ Hence, $A^\prime = (\neg A) = (\neg(\neg B)).$
To the demonstration that $A_1^\prime,\ldots, A_n^\prime \vd B^\prime=B$ 
adjoin the proof that $\vd B\to (\neg (\neg B)).$ Then consider one MP step. 
This yields the formula $(\neg (\neg B)) = A^\prime.$ Consequently,
$A_1^\prime,\ldots, A_n^\prime \vd A^\prime.$ \pars
Subcase (b). Let $v(B) = F.$ Then $v(A) = T$ yields that $A = A^\prime.$
Hence, $A_1^\prime,\ldots, A_n^\prime \vd B^\prime = (\neg B) = A = A^\prime.$
\pars
For case (ii), let (a), $v(C) = T.$ Thus $v(A) = T$ and $C^\prime = C,\ 
A^\prime = B \to C.$ Hence, $A_1^\prime,\ldots, A_n^\prime \vd C.$ Now add the 
steps $\vd C \to (B \to C)$ and MP yields $B \to C.$ Consequently, 
$A_1^\prime,\ldots, A_n^\prime \vd B\to C = A^\prime.$ \pars
Now, let (b) $v(B) = F.$ Then $v(A) = T,\ B^\prime = (\neg B),\ A^\prime = A 
=B \to C.$ Then consider $A_1^\prime,\ldots, A_n^\prime \vd B^\prime = (\neg 
B)$ and adjoin to this proof the proof of $\vd (\neg B) \to (B \to C).$
(Exercise 2.12.2b) Then MP yields $A_1^\prime,\ldots, A_n^\prime \vd B\to C = 
A^\prime.$ \pars
Next part (c) requires $v(B) = T,\ v(C) = F.$ Then $v(A) = F,\ B^\prime = B,\  
C^\prime = (\neg C),\ A^\prime = (\neg A) = (\neg (B\to C)).$ Using both 
demonstrations for $A_1^\prime,\ldots, A_n^\prime \vd B^\prime$ and for 
$A_1^\prime,\ldots, A_n^\prime \vd C^\prime$ and the result of Example 2.14.2
that $\vd B \to ((\neg C) \to (\neg (B \to C)))$, two applications of MP yields
$A_1^\prime,\ldots, A_n^\prime \vd (\neg (B \to C)) = A^\prime.$ By induction, 
the proof is complete for atomic $Ai$. Note that the * substitution process holds if done throughout each step of the formal proof. Applying this process, the result holds in general. \qed
{\bf Examples 2.14.3} Let $A = B \to (\neg C).$ Then 

\newbox\medstrutbox
\setbox\medstrutbox=\hbox{\vrule height14.5pt depth9.5pt width0pt}
\def\medstrut{\relax\ifmmode\copy\medstrutbox\else\unhcopy\medstrutbox\fi}

\bigskip

\newdimen\leftmargin
\leftmargin=0.0truein
\newdimen\widesize
\widesize=2.5truein
\advance\widesize by \leftmargin       
\hfil\vbox{\tabskip=0pt\offinterlineskip
\def\tablerule{\noalign{\hrule}}

\halign to \widesize{\medstrut\vrule#\tabskip=0pt plus2truein

&\hfil\quad#\quad\hfil&\vrule#
&\hfil\quad#\quad\hfil&\vrule#
&\hfil\quad#\quad\hfil&\vrule#

\tabskip=0pt\cr\tablerule


&$B$&&$C$&&$A$&\cr\tablerule
&$T$&&$T$&&$F$&\cr\tablerule
&$T$&&$F$&&$T$&\cr\tablerule
&$F$&&$T$&&$T$&\cr\tablerule
&$F$&&$F$&&$T$&\cr\tablerule}}
\parm
\centerline{Deducibility relations}\pars
\line{\hskip 0.5in (a) $B,\ C \vd (\neg(B \to (\neg C))),$\qquad  (b) $B,\ (\neg C) \vd B \to (\neg C),$\hfil}
\line{\hskip 0.5in (c) $(\neg B),\ C \vd B \to (\neg C),$ \qquad (d) $(\neg B),(\neg C) \vd B \to (\neg C).$\hfil}\parm
\hs\hm
\centerline{\bf EXERCISES 2.14}\parm
\noindent 1. For the following formula$,$ write as in my example all of the 
possible deducibility relations.\parm
(a) $A = (\neg B) \to (\neg C).$ 
(b) $A = B \to C.$
(c) $A = B\to (C \to B).$
(d) $A = B \to (\neg(C\to D)).$\pars
\ss
\noindent {\bf 2.15 The Completeness Theorem.}\parm
One of our major goals is now at hand. We wish to show that $\mod$ and
that $\vd$ mean the same thing. As I mentioned$,$ I'm presenting only the 
necessary formal theorems that will establish this fact. We need just one more.
This will be example 2.15.4. We'll do this in a few small steps. \parm
{\bf Example 2.15.1}\parm
\centerline{$(\neg A) \to A,\ (\neg A) \vd B$ or $\vd ((\neg A)\to A) \to 
((\neg A)\to B).$}\parm 
\line{(1) $(\neg A)\to A$\hfil Premise}\s
\line{(2) $(\neg A)$\hfil Premise}\s
\line{(3) $A$ \hfil $MP(1,2)$}\s
\line{(4) $\vd (\neg A)\to (A \to B)$ \hfil Exer. 2.11 (2b)}\s
\line{(5) $A \to B$\hfil $MP(2,4)$}\s
\line{(6) $B$ \hfil $MP(3,6)$}\par
\bigskip
{\bf Example 2.15.2}\parm
\centerline{$B\to C \vd (\neg(\neg B)) \to C$ or $\vd (B \to C) \to 
((\neg(\neg B)) \to C)$}\parm
\line{(1) $B \to C$\hfil Premise}\s
\line{(2) $(\neg (\neg B)) \to B$\hfil Exam. 2.11.3}\s
\line{(3) $(\neg (\neg B)) \to C$\hfil $HS(1,2)$}\par
\bigskip
{\bf Example 2.15.3}\parm
\centerline{$(\neg A) \to A \vd A$ or $\vd ((\neg A)\to A) \to A$}\parm
\line{(1) $(\neg A) \to ((A \to A)\to (\neg A))$\hfil $P_1$}\s
\line{(2) $((\neg (\neg (A\to A)) \to (\neg A)) \to (\neg A))\to (A \to (\neg (A\to 
A)))$\hfil $P_3$}\s
\line{(3) $\vd ((A \to A)\to (\neg A))\to (\neg (\neg (A\to A)) \to (\neg 
A)))$\hfil Exam. 2.15.2}\s
\line{(4) $((A \to A)\to (\neg A)) \to (A\to (\neg(A \to A)))$\hfil $HS(3,2)$}\s
\line{((5) $(\neg A) \to (A \to (\neg (A\to A)))$ \hfil $HS(1,4)$}\s
\line{(6) $\vd ((\neg A)\to A) \to ((\neg A)\to (\neg(A \to A)))$\hfil Exam. 
2.15.1}\s
\line{(7) $(\neg A)\to A$ \hfil Premise}\s
\line{(8) $(\neg A)\to (\neg (A\to A))$\hfil $MP(6,7)$}\s
\line{(9) $((\neg A)\to (\neg (A\to A))) \to ((A\to A)\to A)$\hfil $P_3$}\s
\line{(10) $(A \to A) \to A$\hfil $MP(8,9)$}\s
\line{(11) $\vd A \to A$\hfil Exam. 2.11.1}\s
\line{(12) $A$\hfil $MP(10,11)$}\par
\bigskip
{\bf Example 2.15.4}\parm
\centerline{$A \to B, (\neg A) \to B \vd B$ or $\vd (A \to B)\to (((\neg A)\to 
B)\to B)$}\parm
\line{(1) $A \to B$ \hfil Premise}\s
\line{(2) $(\neg A) \to B$\hfil Premise}\s
\line{(3) $\vd B \to (\neg(\neg B))$ \hfil Exer. 2.11. 2(a)}\s
\line{(4) $(\neg A) \to (\neg (\neg B))$\hfil $HS(2,3)$}\s
\line{(5) $((\neg A) \to (\neg (\neg B)))\to ((\neg B) \to A)$\hfil $P_3$}\s
\line{(6) $(\neg B)\to A$\hfil $MP(4,5)$}\s
\line{(7) $(\neg B)\to B$\hfil $HS(6,1)$}\s
\line{(8) $((\neg B) \to B)\to B$\hfil Exam. 2.15.3}\s
\line{(9) $B$ \hfil $MP(7,8$}\par
\bigskip
The next two results completely relate $\mod$ and $\vd.$\parm
{\bf Theorem 2.15.1} (Completeness Theorem) {\sl If $A \in L^\prime$ and 
$\mod A$ (in L)$,$ then $\vd A.$}\parm
Proof. Note that $L^\prime \subset L.$ Now assume $A \in L^\prime$ and $\mod 
A.$ We use an illustration that shows exactly how an explicit proof for $\vdash A$ can be constructed. The process is a reduction process and this illustration can be easily extended to a formally stated reduction processes. Let $P_1,P_2,P_3 $ be the atoms in $A$ and ${\cal T} =\ \vd (P_i \to B)\to (((\neg P_i)\to 
B)\to B).$ The ordering of the construction process is {\bf 0}, followed by {\bf 1}, then {\b 2}, then \b 3 and finally followed by {\bf 4}. First, write done the formal steps that lead to each of the following deducibility relations and then construct, in each case using the Deduction Theorem, the forms {\bf 0}. \parm
\noindent (1) $P_1,P_2,P_3 \vdash A {{\rm D. Thm}\atop\Rightarrow} \ \b 0\   \vdash P_1 \to (P_2 \to (P_3 \to A))$\pars
\noindent (2) $P_1,P_2,\neg P_3 \vdash A {{\rm D. Thm}\atop\Rightarrow} \ {\bf 0}\ \vdash P_1 \to (P_2 \to ((\neg P_3) \to A))$\pars

\noindent (3) $P_1,\neg P_2,P_3 \vdash A {{\rm D. Thm}\atop\Rightarrow} \ {\bf 0}\ \vdash P_1 \to ((\neg P_2) \to (P_3 \to A))$\pars
\noindent (4) $P_1,\neg P_2,\neg P_3  \vdash A {{\rm D. Thm}\atop\Rightarrow} \ {\bf 0}\ \vdash P_1 \to ((\neg P_2) \to ((\neg P_3) \to A))$\pars
\noindent (5) $\neg P_1,P_2,P_3 \vdash A {{\rm D. Thm}\atop\Rightarrow} \ {\bf 0}\ \vdash (\neg P_1)\to (P_2 \to (P_3 \to A))$\pars
\noindent (6) $\neg P_1,P_2,\neg P_3 \vdash A {{\rm D. Thm}\atop\Rightarrow} \ {\bf 0}\ \vdash (\neg P_1)\to (P_2 \to ((\neg P_3) \to A))$\pars
\noindent (7) $\neg P_1,\neg P_2,P_3 \vdash A {{\rm D. Thm}\atop\Rightarrow} \ {\bf 0}\ \vdash (\neg P_1)\to ((\neg P_2)\to (P_3 \to A))$\pars
\noindent (8) $\neg P_1,\neg P_2,\neg P_3 \vdash A {{\rm D. Thm}\atop\Rightarrow} \ {\bf 0}\ \vdash (\neg P_1)\to ((\neg P_2)\to ((\neg P_3) \to A))$\pars
\line{\hskip .5in {\bf 1} Using (1) (5), insert  ${\cal T}$, MP, MP (a) $\Rightarrow \vdash P_2 \to (P_3 \to A)$\hfil}\smallskip
\line{\hskip .5in {\bf 1} Using (2) (6), insert ${\cal T}$, MP, MP (b) $\Rightarrow \vdash P_2 \to ((\neg P_3) \to A))$\hfil}\smallskip
\line{\hskip .5in {\bf 1} Using (3) (7), insert ${\cal T}$, MP, MP (c) $\Rightarrow \vdash (\neg P_2) \to (P_3 \to A)$\hfil}\smallskip
\line{\hskip .5in {\bf 1} Using (4) (8), insert ${\cal T}$, MP, MP (d) $\Rightarrow \vdash (\neg P_2)\to ((\neg P_3) \to A))$\hfil}\smallskip
\line{\hskip 1.0in {\bf 2} Using (a) (c), insert ${\cal T}$, MP, MP (e) $\Rightarrow \vdash P_3 \to A$\hfil}\smallskip
\line{\hskip 1.0in {\bf 2} Using (b) (d), insert ${\cal T}$, MP, MP (f) $\Rightarrow \vdash (\neg P_3) \to A$\hfil}\smallskip

\line{\hskip 1.5in {\bf 3} Using (e) (f), insert ${\cal T}$, MP, MP $\Rightarrow \vdash A$\hfil}\parm
Using this illustration, it follows that if $P_1,P_2,\ldots,P_k$ are the atoms in $A$ and we let $P_0 = A$, then, using the $2^k$ deducibility relations to obtain the combined steps that yield the statements \b 0,  each of the steps of the form $\b j \leq k,$ will reduce the problem, by including additional steps, to one where only the forms $\vd P'_{k- \b j} \to (\cdots (P'_k \to P_0)\cdots )$ occur. In this case, each $P'_m$, $\b j\leq m \leq k,$ will be either a $P_m$ or a $\neg P_m$ and reading from left-to-right will be, for each row of a standard truth table, the same forms $P'_{k-\b j},\ldots, P'_k$ as constructed by Definition 2.14.2. After $\b k$ applications of this process we have a formal proof, without any need for the insertion of premises, that has the last step $P_0=A.$
\qed    
{\bf Theorem 2.15.2} (Soundness Theorem) {\sl If $A \in L^\prime$ and $\vd 
A,$ then $\mod A.$}\parm
Proof. Note that each instance of the axioms $P_1,\ P_2,\ P_3$ is a valid 
formula. Also$,$ we have that $\mod A$ and $\mod A \to B,$ then $\mod B.$ 
Thus at each  step in the  proof of for $\vd A,$ we can insert correctly to the left of the 
formula the symbol $\mod.$ Since the last step in the proof is $A,$ then
we can correctly write $\mod A.$ \qed
{\bf Corollary 2.15.2.1} {\sl Let $\{A_1,\ldots,A_n,A\} \subset L^\prime.$ Then 
$A_1,\ldots,A_n \vd A$ if and only if
$A_1,\ldots, A_n \mod A.$}\parm
Proof. By repeated application  of  deduction theorem and theorem 2.8.1 
(a).\qed
Due to these last two theorems$,$ we can identify the connectives which we have 
used in $L$ but not in $L^\prime$ with equivalent formula from $L^\prime.$ 
Hence define $A \lor B$ by $(\neg A)\to B,$ $A\land B$ by $(\neg (A\to (\neg 
B)))$ and $A \iff B$ by $(A\to B) \land (B \to A).$ \pars
There is a slight difference between the concept for $L^\prime$ we denote by
$\Gamma \vd A$ and the concept $\Gamma = \{A_1,\ldots,A_n\} \mod A.$ The concept 
$\Gamma \vd A$ includes the possibility that $\Gamma$ is infinite not just 
finite. Shortly$,$ we'll be able to extend $\Gamma \mod A$ to the possibility 
that $\Gamma$ is infinite. Before we do this however$,$ I can$,$ at last$,$ introduce 
you to the very first basic steps in the generation of a logical operator 
that mirrors a physical-like process that will create universes.\pars
\ss
\noindent{\bf 2.16 Consequence Operators}\parm
In set theory$,$ if you are given  a set of anything $A,$ like a set of formula 
from $L^\prime,$ then another set is very easily generated. The set is 
denoted by $\power A.$  This set is the set of all subsets of $A.$ In the 
finite case$,$ suppose that $A =\{a,b,c\}.$ Letting $\emptyset$ denote the 
empty set$,$ then the set of all subsets of $A$ is $\power A = \{\emptyset, A, 
\{a\},\{b\},\{c\}, \{a,b\}, \{a,c\}, \{b,c\}\}.$ Notice that this set has 8 
members. Indeed$,$ if any set has $n$ members$,$ then $\power A$ has $2^n$ 
members. Obviously$,$ if $A$ is not finite$,$ then $\power A$ is not finite. \pars
Most logical processes$,$ like $\vd,$ satisfy a very basic set of process 
axioms. Notice that you can consider $\Gamma \vd B$ as a type of function.
First$,$ the entire process that yields a demonstration is done upon subsets
of $A.$ When the subset is $\emptyset,$  you get a formal theorem. The set of 
ALL theorems or deductions from a give $\Gamma$ is a subset of $A.$ Let 
$\Delta$ denote the set of all the deductions that can be obtained from 
$\Gamma$ by our entire deductive process. Then can we get 
anything new $B \not\in \Delta$ by considering the premises $\Gamma\cup 
\Delta$? Well$,$ suppose we can and$,$ of course$,$ $B$ is not a theorem. 
Then there is a finite set of formula from 
either $\Gamma$ or $\Delta$ that is used in the demonstration. Suppose that 
from $\Delta$ you have $D_1,\ldots, D_m,$ where we may assume that these are 
not theorems and  they are used in the proof. But 
each of these comes from a proof using members from $\Gamma.$ So$,$ 
just substitute for each occurrence of $D_i$ its proof. Then you have a 
proof of $B$ using only members of $\Gamma.$ Thus you cannot get any more 
deductions by adjoining the set of all deductions to the original set of 
premises.\pars
So$,$ what have we determined about the $\vd$ process? Well$,$  to express what we 
have learned mathematically$,$ consider a function $C$ with the domain the set 
of all premises. But$,$ this is just $\power {L^\prime}.$ Then the codomain$,$ the set 
of all deductions $\Delta$ is also a subset of $\power {L^\prime}.$\par\bigskip
(1) $C$ has as its domain $\power {L^\prime}$ and its range is contained in $\power 
{L^\prime}.$\parm
(2) Since a one step demonstration yields a premise$,$ then 
for each ${\cal B} \in \power {L^\prime}, \ {\cal B} \subset C({\cal B}).$\parm
(3) From the above discussion$,$ for each ${\cal B}\in \power {L^\prime},\ C(C({\cal B})) 
= C({\cal B}).$\parm
(4) Theorem 2.12.1 (f) or our above discussion  states that if
$B \in C({\cal A})$ then there exists a finite subset ${\cal D}$ of ${\cal A}$ 
such that 
$B \in C({\cal D}).$ \parm
Any function that satisfies$,$ (1)$,$ (2)$,$ (3)$,$ (4) is called a {{\it consequence 
operator}}. The important thing to know is that $\vd$ can be replace by such a 
consequence operator with additional axioms. For example$,$ 
(5) for each $A,\ B,\ C \in 
{L^\prime},$ the set of all $P_1,$ the set of all $P_2,$ and the set of all $P_3$ form 
the set $C(\emptyset).$ Then (6) for each ${\cal A} \in \power {L^\prime},$ if 
$A,\ A \to B \in C({\cal A}),$ then $B \in 
C({\cal A}).$ \parm  
What mirrors the physical-like behavior that creates a 
universe are  very special type of consequence operators$,$ one of which is 
denoted by  
$\Hyper S.$ 
Operator  $\Hyper S$ is basically  
determined by a very simple logical process $S$. 
It's basis uses our description for $\vd$ including (6) (our MP). But a 
different set of axioms. These axioms are actually four very simple  
theorems from the language $L^\prime$ with the definition for $\land.$
The completeness theorem tells us that they are theorems in $\power 
{L^\prime}.$ Specifically$,$ they are\parm
\indent\indent (1) $(A \land (B \land C)) \to ((A \land B)\land C).$\pars
\indent\indent (2) $((A \land B) \land C) \to (A \land (B\land C)).$\pars 
\indent\indent (3) $(A \land B) \to A.$\pars
\indent\indent (4) $(A \land B) \to B.$\parm
It turns out that every know propositional deduction used throughout all the 
physical sciences$,$ if they are different from the one we are studying$,$ have 
(1) -- (4) as theorems. \parm
Consequence operator  $S,$ can be generated by  the consequence operators
$S_n,$ where the only difference between $S$ and $S_n$ is that for $S_n$ the 
MP step is restricted to level $L^\prime_n.$  Although for ${\cal A} \in 
\power {L^\prime}$ if $A\in S({\cal A})$ there exists some $n$ such that
$A \in S_n({\cal A}),$ it is also true that for each $n\geq 3$ $S_n({\cal 
A})\subset S({\cal A})$ and $S_n({\cal A})\not= S({\cal A}).$  When 
a consequence operator like $S_n$ has this property then $S$ said to 
be {{\it stronger than}}  $S_n.$\pars
Now I can't go any further in discussing the very special  consequence 
operator that generates a universe. Why? Since the language $L^\prime$ and 
the deductive process $\vd$ must be greatly expanded so that it is {{\it more 
expressive}}. Indeed$,$ so that we can express almost everything within 
mathematics with our language. But$,$ after we have done this$,$ then in the very 
last section of this book$,$ I'll be able to show the mathematical existence of$,$ 
at the least$,$ one of these universe generating consequence operators.  \parm
\hs\hs
\centerline{\bf EXERCISE 2.16}\parm
\noindent In the following$,$ let $C$ be a consequence operator defined on
$\power L^\prime.$ See is you can give an argument that establishes the 
following additional consequence properties based$,$ originally$,$ upon  the 
axioms. \parm
\noindent 1. Let ${\cal A},\ {\cal B}$ be two sets of premises taken from $L^\prime.$ 
Suppose that ${\cal A} \subset {\cal  B}.$ Show that $C({\cal A}) \subset 
C({\cal B}).$\parm
\noindent 2. Recall that ${\cal A} \cup {\cal B},$ the ``union'' set$,$ 
is the set of all formula a formula $A\in {\cal A} \cup {\cal B}$  if and only 
if $A \in {\cal A}$ or $A \in {\cal B}.$ Suppose that ${\cal A} \cup {\cal 
B}\subset L^\prime$. Show that ${\cal A} \subset C({\cal B})$ if and only if
$C({\cal A}) \subset C({\cal B}).$\parm
\noindent 3. Suppose that ${\cal A} \cup {\cal B}\subset L^\prime$. Show that
$C({\cal A} \cup {\cal B}) = C({\cal A} \cup C({\cal B})) = 
C(C({\cal A} \cup C({\cal B})).$ \parm
\hs\hs
\centerline{\bf Some other properties of idempotent operators.}\parm
We use the consequence operator as our prototype. Recall that an operator $C$ (function, map, etc.) is {\it idempotent} if for each $X\in \power {L},\ C(C(X)) = C(X).$\pars
Let $S_1 = \{C(X)\mid X\in \power {L}\},\ S_2 = \{Y\mid Y = C(Y) \in \power {L}\}.$\parm 
{\bf Theorem 1.} {\it The sets $S_1 = S_2.$}\pars
Proof. Let $Y = C(X) \in S_1$. Then $Y \in \power {L}$ and $Y = C(C(X))= C(Y).$ Thus $Y \in S_2.$ \pars
Conversely, let $Y=C(Y) \in \power {L}.$ Then $Y \in S_1.$ Hence, $S_1 = S_2.$\qed 
We can ask if $C(X_1) =C(X_2)$, does this matter? The answer is no since if $Y = C(X_1)= C(X_2) = C(C(X_1)) = C(C(X_2))= C(Y) \in S_2.$ There is a significant unification $\cal U$ for any collection of physical theories The definition of $\cal U$ required that we consider the set $\{Y\mid X \subset Y=C_1(Y) = C_2(Y)\}.$ \pars
Note the identity map $I(X) = X$ for each $X \in \power {L}$ is idempotent. (Indeed, a consequence operator if we are considering only these
objects.) So, one can inquire as to when a given $C$ has an inverse $C^{\leftarrow} = C^{-1}.$ As usual $C^{-1}$ is an inverse if $C^{-1}(C(X))= I(X).$ \parm
{\bf Theorem 2.} {\it The idempotent operator $C$ has an inverse if and only if $C= I.$}\parm
Proof. Suppose that idempotent $C$ has an inverse $C^{-1}.$ Then for each $X \in \power {L},\ X = I(X) =(C^{-1}C)(X) = C^{-1}(C(X))=C^{-1}(C(C(X)))=(C^{-1}C)(C(X)) = I(C(X)) = C(X).$ Thus, from the definition of the identity operator, $C =I$. \qed
{\bf Corollary 2.1}. {\it The only idempotent operator $C$ that is one-to-one is the identity.}\parm
This corollary is an interesting result for consequence operators since all science-community theories are generated by  logic-systems which generate corresponding consequence operators. (You can find the definition of a logic-system in my published paper ``Hyperfinite and standard unifications for physics theories,'' Internat. J. Math. Math. Sci, 28(2)(2001), 17-36 with an archived version at http://www.arxiv.org/abs/physics/0105012). So, for a the logic system used by any science-community, there are always two distinct sets of hypotheses $X_1,\ X_2$ such that $C(X_1)= C(X_2).$\parm
Another example of idempotent operators are some matrices. 
Can we apply these notions to functions that take real or complex numbers and yield real or complex numbers? There is a result that states that the only non-constant idempotent linear real valued function defined on, at least, $[b,d],\ b ,0,\ d>0$ and continuous at $c \in (b,d)$ is $f(x) = x.$ Can you ``prove'' this? It also turns out in this case that if you want this identity form (i.e. $f(x) = x$), then continuity at some $x = c$ is necessary. In fact under the axiom of set theory called The Axiom of Choice, there is a function defined on all the reals that has either a rational value or is equal to zero for each $x$, is, clearly, not constant, is linear, idempotent and not continuous for any real $x.$ This makes it somewhat difficult to ``graph.''\pars
\ss
\noindent {\bf 2.17 The Compactness Theorem.}\parm
With respect to $\vd$$,$ consistency is defined in the same manner as it was 
done for $\mod$. We 
will use the defined connectives $\land,\ \lor,\ \iff.$\parm
\hs
{\bf Definition 2.17.1} (Formal consistency.) A nonempty set of premises 
$\Gamma$ is {{\it formally consistent}} if there does not exist a formula $B 
\in L^\prime$ such that $\Gamma \vd B \land (\neg B).$\pars
\hm
For finite sets of premises $\Gamma$ the Completeness and Soundness  Theorems
show that definition 2.17.1 is equivalent to consistency for $\mod.$ In the 
case of a set of finitely many premises$,$ then all of the consistency results 
relative to $\mod$ can be transferred. But$,$ what do we do when $\Gamma$ is an 
infinite set of premises? \pars
Well$,$ we have used the assignment concept for finitely many premises. We then  
used the symbolism $v(A,\ass)$ for the truth-value for $A$ and the assignment 
to a finite set of atoms that includes the atoms in $A.$ But if you check the
proof in the appendix that such assignments exist in general$,$ you'll find out 
that we have actually define a truth-value function on all the formula in $L.$ 
It was done in such a manner$,$ that it preserved all of the truth-value 
properties required for our connectives for $L.$ Obviously$,$ we could reduce 
the number of connectives and we would still be able to construct a function 
$v$ that has a truth-value for each of the infinitely many atoms and preserves 
the truth-value requirements for the connectives $\neg$ and $\to.$ What has 
been proved in the appendix is summarized in the following rule.\parm
\hs
{\bf Truth-value Rule.} There exist truth-value functions$,$ $v,$ defined for each 
$A \in L^\prime$ such that for each atom $A \in L^\prime$:\pars
 (a) $v(A) = T$ or $v(A) =  F$ not both. \pars
(b) For any $A \in L^\prime,$  $v(A) = T$ if and only if $v(\neg A) = F.$\pars
(c) For any $A,\ B \in L^\prime,$ $v(A \to B) = F$ if and only if $v(A) = T$ 
and $v(B) = F.$\pars
(d) A truth-value function $v$ will be called a {{\it valuation}} function 
and any such function is unique in the following sense. Suppose $f$ and $v$ 
are two functions that satisfy (a) and for each $P \in L_0$ $f(P) = v(P).$   
Then if (b) and (c) hold for both $f$ and $v,$ then $f = v.$ \parm
\hm
Because of part (d) of the above rule$,$ there are  many different valuation 
functions. Just consider a different truth-value for some of the atoms in 
$L_0$ and you have a different valuation function. In what follows$,$ we let 
${\cal E}$ be the set of all valuation functions. \parm
\hs
{\bf Definition 2.17.2} (Satisfaction) \parm
(a) Let $\Gamma \subset L^\prime.$ If there 
exists a $v  \in {\cal E}$ such that for each $A \in \Gamma,\ v(A) = T,$ then
$\Gamma$ is {{\it satisfiable}}.\pars
(b) A formula $B \in L^\prime$ is a {{\it valid consequence of $\Gamma$}} if 
for every $v  \in {\cal E}$ such that for each $A \in \Gamma$ $v(A) = T,$ then 
$v(B) = T.$\parm
\hm
It's obvious that for finite sets of premises the definition 2.17.2 is the same 
as our previous definition (except for a simpler language). Thus we use the 
same symbol $\mod$ when definition 2.17.2 holds. The proofs of metatheorems for 
this extended concept of truth-values are slightly different than those for the 
finite case of assignments. Indeed$,$ it would probably have been better to have 
started with this valuation process and not to have considered the finite 
assignment case. But$,$ the reason I did not do this was to give you a lot a 
practice with the basic concepts within elementary mathematical logic so as to build up a 
certain amount of intuition. In all that follows$,$ all of our formula 
variables are considered to be formula in $L^\prime.$ I remind you$,$ that all of our previous 
results that used assignments hold for this extended concept of $\mod$ if the 
set of premises is a finite or empty set.\parm
\vfil\eject
{\bf Theorem 2.17.1}\parm
{\sl (a)  $\mod A$ if and only if $\{(\neg A)\}$ is not satisfiable.\parm
(b) A single formula premise $\{A\}$  is consistent if and only if $\not\vd 
(\neg A).$ \parm
(c) The Completeness Theorem is equivalent to the statement that 
``Every consistent formula is satisfiable.'' \parm
(d) The set of premises $\Gamma \mod A$ if and only if $\Gamma \cup \{(\neg 
A)\}$  is not satisfiable. \parm
(e) If the set $\Gamma$ is formally consistent and $C \in \Gamma,$ then
$\Gamma \not\vd (\neg C).$}\parm
Proof. (a)$,$ (b)$,$ (c) are left as an exercise.\pars
(d) First$,$ let $\Gamma \not\mod A.$ Then there is a $v \in {\cal E}$ such that 
$v(C) = T$ for each $C \in \Gamma$ but $v(A) = F.$ This $v((\neg A)) = T.$ 
Hence$,$ $\Gamma \cup \{(\neg A)\}$ is satisfiable.\pars
For the converse$,$ assume that $\Gamma \mod A.$  Now assume that  there exists 
some $v \in {\cal E}$ such that  $v(C) = T$ for each $C \in \Gamma.$ Then for 
each such $v$ $v(A) = T.$ In this case$,$ $v((\neg A)) = F.$ Now if no such $v$ 
exists such that $v(C) = T$ for each $C \in \Gamma,$ then $\Gamma$ is not 
satisfiable. Since these are the only two possible cases for 
$\Gamma \cup \{(\neg A)\},$ it follows that  $\Gamma \cup \{(\neg A)\}$ is not 
satisfiable. \parm
(e) Assume that $\Gamma$ is formally consistent$,$ $C \in \Gamma$ and $\Gamma 
\vd (\neg C).$ This yields that $\Gamma \vd C.$ Now to the demonstration 
add the step $\vd (\neg C) \to (C \to A).$ Then two MP steps$,$ yields that 
$\Gamma \vd A.$ Since $A$ is any formula$,$ simply let $A = D \land (\neg D).$ 
Hence $\Gamma \vd D \land (\neg D).$ This contradicts the consistency of 
$\Gamma.$ \parm
The important thing to realize is that to say that $\Gamma$ is consistent says 
that no \underbar{finite} subset of $\Gamma$ can yield a contradiction. But if $\Gamma$ 
is not itself finite$,$ then how can we know that no finite subset of premises 
will not yield a contradiction? Are there not just too many finite 
subsets to check out? We saw that if $\Gamma$ is a finite set$,$ then all it 
needs in order  to be consistent is for it to be satisfiable due to the Corollary 
2.15.2.1. The next theorem states that for infinite $\Gamma$ the 
converse of what we really need holds.   
But$,$ just wait$,$ we will be able to show that certain infinite sets of 
premises are or are not consistent.\parm
{\bf Theorem 2.17.2} {\sl The set $\Gamma$ is formally consistent if and only if $\Gamma$ is 
satisfiable.}\parm
Proof.  We show that if $\Gamma$ is consistent$,$ then it is satisfiable. 
First$,$ note that we can number every member of $L^\prime.$ We can 
number them with the set of natural numbers $\nat.$ Let 
$L^\prime = \{A_i\mid i \in \nat \}$ Now let $\Gamma$ be given. We extended 
this set of premises by the method of induction. (An acceptable method within 
this subject.) \parm
\indent (1) Let $\Gamma = \Gamma_0.$\pars
\indent (2) If $\Gamma_0 \cup \{A_0\}$ is consistent$,$ then let $\Gamma_1 = 
\Gamma_0 \cup \{A_0\}.$  If not$,$ let $\Gamma_1 = \Gamma_0.$ \pars
\indent (3) Assume that $\Gamma_n$ has been defined for all $n \geq 0.$\pars
\indent (4) We now give the inductive definition. For $n+1$$,$ let 
$\Gamma_{n+1} = \Gamma_n \cup \{A_{n}\}$ if $\Gamma_n \cup \{A_{n}\}$ is 
consistent. Otherwise$,$ let $\Gamma_{n+1} = \Gamma_n.$\parm
It follows by the method of definition by induction$,$ that the entire set of
$\Gamma_n$s has been defined for each $n \in \nat$ and each of these sets 
contains the original set $\Gamma.$ We now define that set $\overline{\Gamma}$ 
as follows: $K \in  \overline{\Gamma}$ if and only if there is some $m\in\nat$ 
such that 
$K \in \Gamma_m.$\parm
 We show that $\overline \Gamma$ is a consistent set. Suppose that 
$\overline \Gamma$ is not consistent. Then for some $C \in L^\prime$ and 
finite a  subset $A_1, \ldots, A_n$ of $\overline\Gamma$$,$ it follows that 
$A_1, \ldots, A_n\vd C\land (\neg C).$ However$,$ since $\Gamma \subset \Gamma_1 
\subset \Gamma_2 \subset \cdots \subset  \Gamma$ etc. and since   $A_1, \ldots, 
A_n$ is a finite set$,$ there is some $\Gamma_j$ such that 
$A_1, \ldots, A_n$ is a subset of $\Gamma_j.$ But this produces a 
contradiction that $\Gamma_j$ is inconsistent. Thus $\overline\Gamma$ is 
consistent. \pars
We now show that it is the  ``largest'' consistent set containing $\Gamma.$
Let $A \in L^\prime$ and assume that $\overline{\Gamma} \cup \{A\}$ is 
consistent. If $A \in \Gamma$$,$ then $A \in \overline\Gamma.$ If $A \notin 
\Gamma,$ then we know that there is some $k$ such that $A = A_k.$ 
But $\Gamma_{k} \cup \{A_k\} \subset \overline\Gamma \cup \{A_k\}$ 
implies that $\Gamma_{k} \cup \{A_k\}$ is consistent. But then $A_k \in 
\Gamma_{k+1} \subset \overline\Gamma.$ Hence$,$ $A \in \overline\Gamma.$ We need a 
few additional facts about $\overline\Gamma.$\parm
(i) $A \in \overline\Gamma$ if and only if $\overline\Gamma \vd A.$
(Such a set of formulas is called a {{\it deductive system}}.) First$,$ the 
process $\vd$ yields immediately that if $A \in \overline\Gamma,$ then 
$\overline\Gamma \vd A.$ Conversely$,$ assume that $\overline\Gamma \vd A.$ Then 
$F_1 \vd A$ for a finite subset of $\overline{\Gamma}.$ We show that 
$\overline\Gamma \cup \{A\}$ is consistent. Assume not. Then there is some 
finite subset $F_2$ of $\overline\Gamma$ such that $F_2 \cup \{A\}\vd 
C\land(\neg C)$ for some $A \in L^\prime.$ Thus $F_1 \cup A_2 \vd C\land(\neg 
C).$ But this means that $\overline\Gamma \vd C\land(\neg 
C).$ This contradiction implies that $\overline\Gamma \cup \{A\}$ is 
consistent. From our previous result$,$ we have that $A \in 
\overline\Gamma.$\pars
(ii) If $B\in L^\prime,$ then either $B \in \overline\Gamma$ or $(\neg B) \in 
\overline\Gamma.$ (When a set of premises has this property they are said to 
be a {{\it (negation) complete set.}}) From consistency$,$ not both $B$ and $\neg B$ can 
be members of $\overline\Gamma$\pars
(iii) If $B\in \overline\Gamma$$,$ then $A \to B \in \overline\Gamma$ for each $
A \in L^\prime.$\pars
(iv) If $A \not\in\overline\Gamma,$ then $A \to B \in \overline\Gamma$ for 
each $B\in L^\prime.$\pars
(v) If $A \in \overline\Gamma$ and $B \not\in\overline\Gamma,$ then  
$A \to B \not\in\overline\Gamma.$\pars
(Proofs of (ii) -- (v) are left as an exercise.)\parm
We now need to define an valuation on all of $\overline\Gamma.$ Simple to do. 
Let $v(A) = T$ if $A \in \overline\Gamma$ and $v(B) =F$ if $B \not\in 
\overline\Gamma.$ Well$,$ does this satisfy the requirements of a valuation 
function? First$,$ it is defined on all of $L^\prime$? \pars
(a) By (ii)$,$ if $v(A) = T,$ then $A\in \overline\Gamma$ implies that 
$(\neg A) \not\in \overline\Gamma.$ Hence$,$ $v(\neg A) = F.$\pars
(b) From (iii) --(v)$,$ it follows that $v(A \to B) =F$ if and only if $v(A) = 
T$ and $v(B) = F.$ Thus $v$ is a valuation function and $\overline\Gamma$ is 
satisfiable. But $\Gamma \subset  \overline\Gamma$ implies that $\Gamma$ is 
satisfiable and the proof is complete.\parm
Now to show that if $\Gamma$ is satisfiable$,$ then it is consistent. Suppose 
that  $\Gamma$ is satisfiable but not consistent. Hence there exists a finite 
$F \subset \Gamma$ such that $F \vd C \land (\neg C)$ for some $C \in 
L^\prime.$ But then $F\mod C \land (\neg C),$ by corollary 2.15.2.1 and 
is not satisfiable. Hence since $F \subset \Gamma,$ 
then $\Gamma$ is not satisfiable. This contradiction yields this result.\qed
{\bf Theorem 2.17.3} {\sl The following statement are equivalent.\parm
(i) If $\Gamma \mod B,$ then $\Gamma \vd B.$ (Completeness)\pars
(ii) If $\Gamma$ is consistent$,$ then $\Gamma$ is satisfiable.}\parm
Proof. Assume that (i) holds and that $\Gamma$ is consistent. If $C \in 
\Gamma,$ then by theorem 2.17.1 part (e) $\Gamma\not\vd (\neg C).$ By the 
contrapositive if (i)$,$ then $\Gamma \mod (\neg C).$ But by (d) of theorem 2.17.1$,$ 
$\Gamma \cup \{(\neg(\neg C))\}$ is satisfiable. Let $v$ be the valuation.
Then $v(A) = T$ for each $A \in \Gamma$ AND $v((\neg(\neg C))) = T$. Hence$,$
$v(C) = T.$ Thus $\Gamma$ is satisfiable.\parm
Now assume that (ii) holds and let $\Gamma \mod B.$  Then by theorem 2.17.1 
part (d) $\Gamma \cup \{(\neg B)\}$ is not satisfiable. Hence $\Gamma \cup 
\{(\neg B)\}$ is inconsistent. Consequently$,$ there is some $C$ such that 
$\Gamma \cup \{(\neg B)\}\vd C \land(\neg C).$ By Corollary 2.15.2.1$,$ 
$C \land (\neg C) \vd A$  for any $A \in L^\prime.$  So$,$ let $A = (\neg B) \to 
B.$ So$,$ in the demonstration that $\Gamma \cup \{(\neg B)\}\vd C \land(\neg 
C)\vd (\neg B) \to B,$  one MP step yields $B$. Thus $\Gamma \vd B.$ \qed
From theorem 2.17.3$,$ since (ii) holds$,$ then the extended completeness theorem 
(i) holds. But it really seems impossible to show that an infinite
$\Gamma$ is inconsistent unless we by chance give a demonstration that $\Gamma 
\vd C \land (\neg C)$ or to show that it is consistent by showing that $\Gamma$ is 
satisfiable. For this reason$,$ the next theorem and others of a similar 
character are of considerable importance.\parm
\vfil\eject
{\bf Theorem 2.17.3} (Compactness) {\sl A set of formulas $\Gamma$ is 
satisfiable if and only if every finite subset of $\Gamma$ is 
satisfiable.}\parm
Proof. Assume that $\Gamma$ is satisfiable. Then there is some $v \in \cal E$ 
such that $v(A) = T$ for each $A \in \Gamma.$ Thus for any subset $F$ of $\Gamma$ 
finite or otherwise $v(B) = T$ for each $B\in F.$ \parm
Conversely$,$ assume that $\Gamma$ is not satisfiable. Then from theorem 2.17.2$,$
$\Gamma$ is not consistent. Hence there is some $C\in L^\prime$ and finite 
$F\subset \Gamma$ such that $F \vd C\land (\neg C).$ By the soundness theorem$,$
$F \mod C\land (\neg C).$ Thus $F$ is not satisfiable and the proof is 
complete.\qed
{\bf Corollary 2.17.4.1} {\sl A set of premises $\Gamma$ is consistent if and 
only if every finite subset of $\Gamma$ is consistent.}\parm
{\bf Example 2.17.1} Generate a set of premises by the following rule. 
Let $A_1 = (\neg A),\ A_2 = (\neg A) \land A, A_3 = (\neg A)\land A \land A,$ 
etc. Then the set $\Gamma =\{A_i \mid i \in \nat \}$ is inconsistent since
$A_2$ is not satisfiable.\parm
{\bf Example 2.17.2} Generate a set of premises by the following rule.
Let $A_1 = P,\ A_2 = P\lor P_1,\ A_3 = P \lor P_1 \lor P_2, \ldots, 
A_n = P \lor \cdots \lor P_{n-1}.$ Then the set $\Gamma = \{A_i \mid i \geq
1,\ i  \in \nat \}$ is consistent. For consider a nonempty finite subset $\c F$ 
of $\Gamma.$ Then there exists a formula $A_k \in \c F$ (with maximal subscript) 
such that if $A_i \in \c F,$  then $1 \leq i \leq k.$ Let $v \in \c E$ be a 
valuation such that $v(P) = T.$  Obviously such a valuation exists. The 
function $v$ also gives truth-values for all other members of $\c F.$ 
But all other formula in $\c F$ contain
$P$ and are composed of a formula $B$ such that $A_i= P \lor B,\  2\leq i \leq 
k.$ But $v(P \lor B) = T$ independent of the values $v(B).$ Hence $\c F$ is 
satisfied and the compactness theorem states that $\Gamma$ is consistent.\parm
\hs\hs
\centerline{\bf EXERCISES 2.17}\parm
\noindent 1. Prove properties (ii)$,$ (iii)$,$ (iv)$,$ (v) for $\overline{\Gamma}$ 
as they are stated in the proof of theorem 2.17.2.\parm
\noindent 2. Prove statements (a)$,$ (b)$,$ (c) as found in theorem 2.17.1.\parm
\noindent 3. Use the compactness theorem and determine whether or not the 
following sets of premises are consistent.\pars
(a) Let $A \in L^\prime.$ Now $\Gamma$ contains $A_1 = A \to A,\ A_2= (A \to 
(A \to A),\  A_3 = A \to (A\to (A \to A)),\ $ etc. \pars
(b) Let $A \in L^\prime.$ Now $\Gamma$ contains $A_1 = A \to A, \ A_2 =
(\neg (A \to A)),\ A_3 = A \to (\neg(A\to A))$ etc.\pars
(c) Let $A_1 = P_1 \iff P_2,\  A_2 = P_1 \iff (\neg P_2),\  A_3 = P_2 \iff 
P_3,\ A_4 = P_2 \iff (\neg P_3),$ etc. 
\vskip 0.25in
\centerline{The consistency of the process (1)(2)(3)(4) used in Theorem 2.17.2}\par\bigskip
We do not need Theorem 2.17.2 if $\Gamma$ is a non-empty finite set of premises since then Corollary 2.15.2.1 applies. However, what follows holds for any set of hypotheses. Clearly, if the process (1)(2)(3)(4) is consistent it is not an effective process since there are no rules given to determine, in a step-by-step process, when a set such as $\Gamma \cup {A}$, where $A \notin \Gamma,$ is consistent. On the other hand, there are various ways to show convincingly that the process itself is consistent. This is done by the intuitive method of re-interpretation and modeling. Two examples are described below. Later it will be shown more formally that such statements that characterize such processes are consistent if and only if they have a set-theoretic model. \pars
(A) Call any $A \in L'$ a ``positive formula'' if the 
$\neg$ symbol does not appear in the formula.  Let $\Gamma$ be a set of positive formula. Apply the process (1)(2)(3)(4) where we substitute for ``consistent'' the phrase ``a set of positive formula.'' Note that since we started with a set of positive formula, then we can determine whether the set $\Gamma_n \cup \{A_n\}$ is a set of positive formula just by checking the one formula $A_n.$ What has been done in this example is that the original process (1)(2)(3)(4) has been re-interpreted using a determining requirement that can actually be done or a requirement that most exist. The fact that an actual determination can be made in finite time is not material to the consistency of the process. All that is required is that each formula be either positive or negative, and not both. This re-interpretation is called a ``model'' for this process and implies that the original process is consistent relative to our intuitive metalogic. This follows since if the process is not consistent, then a simple metalogical argument would yield that there is an actual formula in $L'$ that is positive and not positive.\footnote*{Formally, (1) $ B\land (\neg B)\equiv \neg(B \to B).$ (2) $\vd (\neg(B \to B) \to ((B\to B)\to C_i)).$ (3) $(B\to B)\to C_i.$ (4) $\vd B \to B.$ (5) $C_i$. Just take $C_1 = $ ``$A$ is positive'' and $C_2 =$ ``$A$ is not positive.''} All of this is, of course, based upon the acceptance that the processes that generated $L'$ are also consistent.\pars
(B) We are using certain simple properties of the natural numbers to study mathematically languages and logical processes. It is assume that these natural number processes are  consistent. Using this assumption, we have the following model that is relative to the a few natural number properties.  Interpret $\Gamma$ as a set of even numbers and interpret each member of $L'$ as a natural number. Substitute for ``consistent'' the phrase ``is a set of even natural numbers.'' Thus, all we need to do is to determine for the basic induction step is whether $\Gamma_n \cup \{A_n\}$ is a set of even numbers. Indeed, all that is needed is to show that $A_n$ can be divided by 2 without remainder, theoretically a finite process. This gives a model for this induction process relative to the natural numbers since a natural number is either odd or even, and not both. Hence, we conclude that the original process is consistent relative to our metalogic. 
\vfil\eject\centerline{NOTES}\vfil\eject 
\centerline{\bf Chapter 3 - PREDICATE CALCULUS\hfil}\parm

\noindent {\bf 3.1 A First-Order Language.}\parm
The propositional language studied in Chapter 2 is a fairly expressive formal 
language. Unfortunately$,$ we can't express even the most obvious scientific or 
mathematical statements by  using only a propositional language. \parm
{\bf Example 3.1.1} Consider the following argument. Every natural number is a 
real number. Three is a natural number. Therefore$,$ three is a real number. 
\parm
The propositional language does not contain the concept of ``every.''    Nor 
does it have a method to go from the idea of ``every'' to the specific 
natural number 3. Without a formal language that mirrors the informal idea of 
``every'' and the use of the symbol ``3'' to represent a specific natural 
number$,$ we can't follow the procedures such as those used in chapter 2 
to analyze the logic of this argument. In this first section$,$ we construct 
such a formal language. Accept for some additional symbols and rules as to 
how you construct new formula with these new symbols$,$ the method of construction 
of language levels is exactly the same as used for $L$ and $L^\prime.$ So that 
we can be as expressive as possible$,$ we'll use the connectives $\lor,\ 
\land,$ and $\iff$ as they are modeled by the connectives $\neg$,$ \to$. We 
won't really need the propositions in our construction. But$,$ by means of a 
special technique$,$ propositions can be considered as what we will call
{{\it predicates}} with {{\it constants inserted}}. This will be seen 
 from our 
basic construction. This construction follows the exact same pattern as 
definition 2.2.3. After the construction$,$ I will add to the English language 
interpretations of definition 2.2.4 the additional interpretations for the new 
symbols. \parm
\hs
{\bf Definition 3.1} (The First-Order Language $Pd$)\parm
(1) A nonempty set of symbols written$,$ at the first$,$ with missing pieces (i.e. holes 
in them). They look like the following$,$ where the underlined portion means a 
place were something$,$ yet to be described$,$ will be inserted.\pars
\indent\indent (a) The {{\it 1-place predicates}}\hfil\break $P(\underline{\quad}),\  
Q(\underline{\quad}),\ R(\underline{\quad}),\ S(\underline{\quad}),\ 
P_1(\underline{\quad}),\ P_2(\underline{\quad}),\ldots.$\pars                   
\indent\indent (b) The {{\it 2-place predicates}}\hfil\break 
$P(\underline{\quad},\underline{\quad}),\  
Q(\underline{\quad},\underline{\quad}),\ 
R(\underline{\quad},\underline{\quad}) S(\underline{\quad},\underline{\quad}),\ 
P_1(\underline{\quad},\underline{\quad}),\ P_2(\underline{\quad},\underline{\quad}),\ldots.$\pars  
\indent\indent (c) The {{\it 3-place predicates}}\hfil \break
$P(\underline{\quad},\underline{\quad},\underline{\quad})$,$  
Q(\underline{\quad},\underline{\quad}$,$ \underline{\quad}),
R(\underline{\quad},\underline{\quad},\underline{\quad}),\ldots.$\pars
\indent\indent (d) And so-forth$,$ continuing$,$ if necessary$,$ through any
{{n-place predicate}}$,$ a symbol with a  ``('' followed by $n$  ``holes$,$'' 
where the 
underlines appear$,$ followed by a  ``)''. \parm
(2) An infinite set of {{\it variables}}$,$ $V= 
\{x_1,x_2,x_3,\ldots \},$ where we let the symbols $x$,$ \ y,\ z$,$ w$ represent 
{\it any} distinct members of $V.$\parm
(3) A nonempty set of constants ${\cal C}  = 
\{c_1,c_2,c_3,c_4,\ldots\}$ and other similar notation$,$ where we let $a,\ b,\ c,\ d$ represent 
distinct members of ${\cal C}.$ \parm
(4) We now use the stuff in (1)$,$ (2)$,$ (3) to construct a set of {{\it 
predicates}}$,$ ${\c P},$ which are the atoms of $Pd.$ In each of the underlined 
places in (1)$,$ insert  either a single variable from $V$ 
or a single constant from ${\cal C}.$ Each time you construct one of these new 
symbols$,$ it is called a $1,2,\ldots,$ n-place predicate. 
{\bf The language $Pd$ will use one or more of these predicates.} \parm              
(5) The non-atomic formula are constructed from an {\it infinite} set of 
connectives. The u-nary connective $\neg$ and the previous binary
connectives $\lor,\land,\ \to,\ \iff.$  You place these in the same 
positions$,$ level for level$,$ as in the construction of definition 2.2.3.\pars
(6) Before every member of $V$$,$ the variables$,$ will be placed immediately
to them left the symbol $\forall.$ Then before every member of $V$ will be 
placed immediately to the left the symbol $\exists.$ This gives infinitely 
many symbols of the form ${\c U}=\{ \forall x_1,\ \forall x_2,\ \forall x_3,\ 
\forall x_4$,$ 
\ldots\}$ and of the form ${\c EX}= \{ \exists x_1,\ \exists x_2,\ 
\exists x_3,\ \exists x_4$,$ 
\ldots\}$ where we use the symbols $\forall x,\ \forall y,\ \forall x$,$ \ 
\exists x,\ \exists y,\ \exists z$ and the like to represent {\it any} of 
these symbols. Each  $\forall x$ is called an {{\it universal quantifier}}$,$ 
each  $\exists x$ is  called an {{\it existential quantifier}}.\pars
(7) Now the stuff constructed in (6) is an infinite set of u-nary operators 
that behavior in the 
same manner as does $\neg$ relative to the method we used to construct our 
language. \pars
(8) Although in an actual first-order language$,$ we usually need only a 
few predicates 
and a few universal or existential quantifiers or a few constants$,$ it is 
easier to simply construct the levels$,$ after the first level that contains all 
the predicates$,$ in the following manner using parentheses about the outside 
of every newly constructed formula each time a formula 
from level 2 and upward is constructed. Take every formula from the {\it 
preceding  
level} and 
place immediately to the left a 
universal quantifier$,$ take the same 
formula from the previous level and put immediately to the left a 
existential 
quantifier$,$ repeat these constructions for all of the universal and 
existential quantifiers. {\bf Important:  The formula to the right of  
the quantifier used for this  construction is called the} {{\it scope}} 
{\bf of the quantifier.} Now 
to the same formula from the  preceding level  
put immediately to left 
the $\neg$ in the same manner as was done in definition2.2.3.\pars
(9) Now following  definition 2.2.3$,$ consider every formula from the previous 
level and construct the  
formula using the binary connectives $\lor,\ \land,\  \to,\ \iff$ making sure 
you put parentheses about the newly formed formula. AND after a level is 
constructed with these connectives we adjoin the previous level to the one 
just constructed. \par
(10) The set of all formula you obtain in this manner is our language 
$Pd.$\pars
(11) Any predicate that contains only constants in the various places behaves 
just like a proposition in our languages $L$ and $L^\prime$ and can be 
considered as forming a propositional language. What we have constructed 
is an extension of these propositional languages. \pars
(12) I point out that there are other both less formal and more formal ways 
to construct the language $Pd.$\pars
\hm 
\centerline{\bf This is not really a difficult construction.} \parm
{\bf Example 3.1.1} (i) The following are predicates.
 $P(x),$ $ P(c),$ $ P(x,c),$ $ Q(x,x),$ $ Q(c,c_1,z).$\pars
(ii) The following are not predicates. $P(),$ $ P(\underline{\quad}),$ $ 
P(P,Q),$ $ P)(,$ $ Px$,$ $ $Pc.$\pars
(iii) The following are members of $Pd.$  $(P(c) \iff Q(x)),$ $ 
(\forall x P(c)),$ $(\exists yP(x)),$ $ (\forall xP(x,y)),$  
$(\neg(\forall xP(z))),$ $ 
(\exists x(\forall x(\exists y(\forall yP(c,d,w))))),$ $((\exists x(P(x) \iff 
Q(d)))).$\pars
(iv) The following are not formula. $P(\exists x),$ $(\forall 
P(x)),$ $(\forall 
c P(x)),$ $(\exists x(\forall y)P(x,y)),$ $((P(d) \iff Q(x) \exists x).$\parm
As was done with the language $L$$,$ we employ certain terminology to discuss 
various features of $Pd.$ We employ$,$ whenever possible$,$ the {\bf same 
definitions.} Each element of $Pd$ is called a {{\it formula}}. The first 
level $Pd_0$ is the set of all the {{\it atoms}}. Any member of $Pd$ which is 
not an atom is called a {{\it composite formula}}. A formula expressed 
entirely in terms of atoms and connectives is called an {{\it atomic 
formula.}} We let$,$ as usually$,$ the symbols $A,\ B,\ C$,$ \ldots$ be {{\it formula 
variables}} and represent arbitrary members of $Pd.$ \pars
If $A \in Pd,$ then there is a smallest $n$  such that $A \in Pd_n$ and   
$A \notin Pd_k$ for any $k <n.$ The number $n$ is the {{\it size}} of $A.$ 
All the parenthesis rules$,$ the common pair rules hold$,$ where $\forall x$ and 
$\exists x$ behave like the $\neg$. However$,$ the parentheses that form the 
predicates are never included in these rule processes. \pars
The English language interpretation rules for the predicates and the 
quantifiers take on a well-know form$,$ while the interpretation is the same for 
all the other connectives as they appear in 2.1. For the following$,$ we utilize 
the {{\it variable predicate forms}} $A(x)$ and the like. These may be 
thought of as 1-place predicates for the moment. Since this interpretation is 
an inductive interpretation$,$ where we interpret only the immediate symbol$,$ we 
need not know the exact definition for $A(x)$ at this time.\parm
\hm
\+ (i) $\lint P(\underline{\quad})\rint:\ \rint\ 
P(\underline{\quad},\underline{\quad})\rint:\ \lint 
P(\underline{\quad},\underline{\quad},\underline{\quad})\rint:$&\hfil\cr
\noindent These are usually considered to be simple declarative sentences that 
relate various noun forms. In order to demonstrate this$,$ we use the 
``blank'' word notation where the blanks are to be understood by the 
``--'' symbol and the ``blanks'' must be filled in with a variable or a 
constant. Here is a list of examples that may be interpreted as 1-place$,$ 
2-place$,$ 3-place and 4-place predicates. ``-- is lazy.'' ``-- is a man.''
``-- plays --.'' ``-- is less than --.'' ``-- $>$ --.'' ``-- plays -- with --.'' 
``-- plus -- equal --.'' ``-- + -- = --.'' ``-- plays -- with -- at --.'' 
``-- + -- + -- = --.''\parm
\hm
\noindent (ii) The constants $C$ are interpreted as identifying ``names''
for an element of a domain$,$ a proper name and the like$,$ where if there is a common or required name such as ``0$,$'' ``1$,$'' or ``sine$,$'' then this common name is used as a constant.\parm

\hm
\+(iii) $\lint \forall x A(x)\rint:$\hfil\cr\s
\+For all  $x,\ \lint A(x) \rint:$&For arbitrary $x,\ \lint A(x) \rint:$\cr\s
\+For every $x,\ \lint A(x) \rint.$&For each $x,\ \lint A(x) \rint.$\cr\s
\+Whatever $x$ is $\lint A(x)\rint.$&(Common language)\cr\s
\+&Everyone is $\lint A(x) \rint.$\cr\s
\+&$\lint A(x) \rint$ always holds.\cr\s
\+&Each one is $\lint A(x) \rint.$\cr\s
\+&Everything is $\lint A(x)\rint.$\cr\medskip
\hm 
\+(iv) $\lint \exists x A(x)\rint:$&\cr\s
\+There exists an $x$ such that $\lint A(x)\rint.$&There is an $x$ such that $\lint A(x) 
\rint.$\cr\s
\+For suitable $x,\ \lint A(x) \rint.$&There is some $x$ such that $\lint A(x) \rint.$\cr\s
\+&For at least one $x,\ \lint A(x) 
\rint.$\cr\s
\+&(Common language)\cr\s
\+&At least one $x$ $\lint A(x) \rint.$\cr\s
\+&Someone is $\lint A(x) \rint.$\cr\s
\+&Something is $\lint A(x) \rint.$\cr\medskip
\hm
Great care must be taken when considering a negation and quantification as is 
now demonstrated since a  ``not'' before a quantifier or after a quantifier 
often yields different meanings.\medskip
\hm
\+(v) $\lint \neg(\forall x A(x))\rint:$&\cr\s
\+Not for all $x$$,$ $\lint A(x)\rint.$&$\lint A(x) \rint$ does not hold for all 
$x.$\cr\s
\+$\lint A(x)\rint$ does not always hold.&Not everything is $\lint A(x) 
\rint.$\cr\medskip
\hm
\+(vi) $\lint \forall x(\neg A(x))\rint:$&\cr\s
\+For all $x$$,$ not $\lint A(x) \rint.$&$\lint A(x) \rint$ always fails.\cr\s
\+&Everything is not $\lint A(x) \rint.$\cr\medskip
\hm
\+(vii) $\lint \neg(\exists x A(x)) \rint:$&\cr\s
\line{There does not exists an $x$ such that $\lint A(x) \rint.$\hfil}\s
\line{There does not exists any $x$ such that $\lint A(x) \rint.$\hfil}\s
\+There exists no $x$ such that $\lint A(x)\rint.$&Nothing is $\lint 
A(x)\rint.$\cr\s
\+There is no $x$ such that $\lint A(x) \rint.$&No one is $\lint A(x)\rint.$\cr\s
\+There isn't any $x$ such that $\lint A(x)\rint.$&\cr\s
\hm
\+(viii) $\lint \exists x(\neg A(x))\rint:$&\cr\s
\+For some $x$$,$ not $\lint A(x) \rint.$&Something is not $\lint A(x) 
\rint.$\cr\s
\centerline{(All other $\exists x$ with a not just prior to $\lint A(x) 
\rint.$}\s
\rm
There are many other possible English language interpretations for the symbols 
utilized within our first-order language but the ones listed above will 
suffice.\parm
\hs\hm
\centerline{\bf EXERCISES 3.1}\parm
\centerline{NOTE: Outer parenthesis simplification 
may have been applied.}
\noindent 1. Let $A$ represent each of the following strings of symbols. 
Determine if $A \in Pd$ or $A\notin Pd.$ \pars
\+(a) $A= (\exists x(\neg(\forall xP(x,y))))$&(d) $A = (P(c) \lor \exists 
x\forall y)Q(x,y))$\cr\s
\+(b) $A =\forall x(P(x) \to (\exists cP(c,x)))$&(e) $A= \exists x(\forall y(P(x) 
\to Q(x,y)$\cr\s
\+(c) $\forall x(\exists y(\forall z (P(x)))).$&\cr\medskip
\noindent 2. Find the size of $A.$\pars
\+(a) $A = \exists x(\forall y(\exists zP(x,y,z))).$&\cr\s
\+(b) $A = (P(c) \to (\exists y(\neg(\forall xP(x,y))))).$&\cr\s
\+(c) $A = \forall x(P(x) \to ((\exists yP(x,y))\to (P(c) \lor Q(c))))$&\cr\s
\+(d) $A = (P(c) \lor (\exists y(Q(x) \to P(y))))\to R(x,y,z).$&\cr\parm
\noindent 3. Use the indicated predicate symbols$,$ which may appear in previous 
lettered sections$,$  
and any arithmetic symbols 
such as $\underline{\quad} +\underline{\quad}$ for  ``sum of -- and --'' 
or $\underline{\quad}<\underline{\quad}$ for  ``-- is less than --$,$'' or 
$\underline{\quad} = 
\underline{\quad}$ for ``-- is equal to --'' and the like and 
translate the following into symbols from the language $Pd.$ You may need to 
slightly re-write the English language statement into one with the same 
intuitive meaning prior to translation \parm
(a) If the product of finitely many factors is equal to zero$,$ then at least 
one of the factors is zero. [Let $P(\underline{\quad})$ be the predicate
``-- is the product of finitely many factors$,$'' $Q(\underline{\quad 
},\underline{\quad })$ be the predicate   ``-- is a factor of --.'']\parm
(b)  For each real number$,$ there is a larger real number. [Let 
$R(\underline{\quad})$ be  ``-- is a real number.'']\parm
(c) For every real number $x$ there exists a real number $y$ such that for 
every real number $z,$ if the sum of $z$ and 1 is less than $y,$ then the sum 
of $x$ and 2 is less than 4. [Note that  ``$- + 1 < -$'' is a three place 
predicate with a constant ``1'' in the second place.]\parm
(d) All women who are lawyers admire some judge. [Let $W(\underline{\quad})$ 
be 
``-- is a women$,$'' let $L(\underline{\quad})$ be  ``-- is a lawyer$,$'' let
$J(\underline{\quad})$ be ``-- is a judge$,$'' and let  $A(\underline{\quad},
\underline{\quad})$ be ``-- admires --.'']\parm
(e) There are both lawyers and shysters who admire Judge Jones. [Let 
$S(\underline{\quad})$ be ``-- is a shyster$,$'' and jet $j$ symbolize the name
``Jones.'']\parm
(f) If each of two persons is related to a third person$,$ then the first person 
is related to the second person. [Let $R(\underline{\quad},\underline{\quad})$ 
be ``-- is related to --$,$'' and let $P(\underline{\quad})$ be  ``-- is a 
person.'']\parm
(g) Every bacterium which is alive in this experiment is a mutation.
[Let   $B(\underline{\quad})$ be ``-- is a bacterium$,$ $A(\underline{\quad})$ be  
``-- is alive in this experiment$,$'' and $M(\underline{\quad})$ be  ``-- is a 
mutation.'']\parm
(h) The person responsible for this rumor must be both clever and 
unprincipled. [Let $P(\underline{\quad})$  be ``-- is a person responsible for 
this rumor$,$'' $C(\underline{\quad})$ be  ``-- is clever$,$'' and $U(\underline{\quad})$
be  ``-- is unprincipled.'']\parm
(i) If the sum of three equal positive numbers is greater than 3 and the sum 
of the same equal positive numbers is less than 9$,$ then the number is greater 
than 1 and less than 3.\parm
(j) For any persons $x$ and $y$$,$ $x$ is a brother of $y$ if and only if
$x$ and $y$ are male$,$ and $x$ is a different person than $y$$,$ and $x$ and $y$ 
have the same two parents. [Let (again) $P(\underline{\quad})$ be ``-- is a 
person$,$'' $B(\underline{\quad},\underline{\quad})$ be  ``-- is a brother of --
$,$'' $M(\underline{\quad})$ be ``-- is a male$,$'' $ -- \not= --$ be 
``-- is a different person than --$,$'' $Q(\underline{\quad},\underline{\quad})$
be ``-- and -- have the same two parents.'']\parm
\noindent 4. Let $P(\underline{\quad})$ be  ``-- is a prime number$,$'' 
$E(\underline{\quad})$ be  ``-- is an even number$,$''
$\r O(
\underline{\quad})$ be  ``-- is an odd number$,$'' $D(\underline{\quad},\underline{\quad})$
be ``-- divides --.'' Translate each of the following formal sentences into 
English language sentences.\parm
(a) $P(7) \land \r O(7).$\pars
(b) $\forall x(D(2,x) \to E(x)).$\pars
(c) $(\exists x(E(x) ) \land P(x))) \land (\neg(\exists x(E(x) \land P(x))))
\land (\exists y (x \not= y))\land E(y)\land P(y).$\pars
(d) $\forall x(E(x) \to (\forall y(D(x,y)\to E(y)))).$\pars
(e) $\forall x (\r O(x) \to (\exists y(P(y) \to D(y,x)))).$\pars
\ss
\noindent {\bf 3.2. Free and Bound Variable Occurrences.}\parm
Prior to the model theory (the semantics) for the language $Pd,$
a further important concept must be introduced. It's interesting to note that 
the next concept dealing with variables took a considerable length of time to 
formulate in terms of an easily followed rule. As usual$,$ we are fortunate that 
much of the difficult work in mathematical logic has been simplified so that 
we may more easily investigate these important notions. For simplicity$,$ the 
n-place predicates where all places are filled with constants will act as if 
they are propositions. As will be seen$,$ the presence of any quantifiers prior 
to such predicates will have NO effect upon their semantical meaning nor will 
such a case involve any additional nonequivalent formula in the formal 
deduction (proof theory) associated with $Pd.$\parm \hs
{\bf Definition 3.2.1} (Scope of a quantifier.)\pars
Assume that a quantifier is symbolized by $Q$. Suppose that $Q$ appears in a 
formula $A \in Pd.$ Then in that quantifier appears one and only one variable$,$
say $v \in\c V.$  Immediately to the right of the variable $v,$ two and only 
two mutually distinct symbols appear. \parm
(1) A left parenthesis. If this is the case$,$ then the subformula from that 
left parenthesis to its common pair parenthesis us called the {{\it scope of 
that quantifier.}}\pars
(2) No parenthesis appears immediately on the right. In this case$,$ the 
predicate that appears immediately on the  right is called the {{\it scope of 
that quantifier.}}\pars
\hm
In general$,$ to identify quantifiers and scopes$,$ the quantifiers are counted by 
the natural numbers from left to right. In what follows$,$ we place the 
``quantifier'' count as a subscript to the symbols $\forall$ and $\exists$ rather 
than as subscripts to the variables. Note: if the same quantifier occurs more 
than once$,$ then it is counted more than once. Why? Because it is the location 
of the quantifier within a formula that is of first importance.\parm
{\bf Example 3.2.1}\pars
(a) Let $A = (\forall_1 xP(x)) \land  Q(x).$ Then the scope of quantifier 
(1) is the predicate $P(x).$\parm
(b) Let $A = \exists_1 y (\forall_2 x(P(x,y) \to (\forall_3 zQ(x)))).$ The scope 
of the quantifier (1) is $(\forall x(P(x,y) \to (\forall z Q(x)))).$ 
The scope for 
the quantifier (2) is $(P(x,y) \to (\forall z Q(z))).$ The scope for 
quantifier (3) is $Q(x).$\parm
(c) Let $A = \forall_1 x(\forall_2 y(P(c) \to Q(x,y)))\iff (\exists_3 
xP(c,d)).$ Then the scope for quantifier (1) is $(\forall y(P(c) \to 
Q(x,y))).$ The scope for quantifier (2) is $(P(c) \to Q(x,y)).$ The scope for 
quantifier (3) is $P(c,d).$\parm
(d) Let $A = (\forall_1 x ((R(x) \land (\exists_2 xQ(x,y)))\to (\exists_3 y 
P(x,y)))) \to Q(x,z).$ The scope of quantifier (1) is $((R(x) \land (\exists 
xQ(x,y)))\to (\exists y P(x,y))).$ The scope of quantifier (2) is $Q(x,y).$ 
The scope of 
quantifier (3) is $P(x,y).$ Notice that the second occurrence of $Q(x,y)$ 
is not in the scope of any  quantifier.\parm
Once we have the idea of the scope firmly in our minds$,$ then we can define the 
very important notion of the  ``free'' or  ``bound'' occurrence of variable. 
We must have 
the formula in atomic form to know exactly what variables behave in these two 
ways. This behavior is {\it relative to their positions} within a formula. We call 
the position where a variable or constant explicitly appears as an {\it occurrence} of 
that \underbar{specific} variable or constant with in a \underbar{specific}
formula. \parm
{\bf Definition 3.2.2} (Bound or free variables) \parm
(1) Every place a constant appears in a formula is called a {{\it bound 
occurrence of the constant.}}\pars
(2) You determine the free or bound occurrences after the scope for each 
quantifier has been determined. You start the determination at level $Pd_0$ 
and work your way up each level until you reach the level at which the formula 
first appears in $Pd.$ This is the  ``size'' level.\pars
(3) For level $Pd_0$ every occurrence of a variable is a free occurrence.\pars
(4) For each succeeding level$,$ an occurrence of a variable is a {\it bound}
 occurrence if that variable occurs in a quantifier or is not previously marked as a bounded occurrence in the scope of the 
quantifier. \pars
(5) After all bound occurrences have been determined$,$ than any occurrence
of a variable that is not bound is called a {\it free}
occurrence. \pars
(6) A bound occurrence of a specific variable is determine by a quantifier in 
which the variable appears. This is called the {\it bounding quantifier.} 
\pars
(7) Any variable that has a free occurrence within a formula is said to be 
{{\it free in the formula}}.\pars
\hm
There are various ways to diagram the location of a specific variable 
and its bound occurrences. One method is a line diagram and the other is to use the 
quantifier number and the {\it same} number for each bound occurrence of the 
variable. I illustrate both methods. It turns out that for one concept the 
line diagram is the better of the two. \parm
{\bf Example 3.2.2} This is the example of the line diagram method of showing 
the location of each bound variable. \parm
\baselineskip=9pt
\indent \hskip .4in{\vrule height0pt width1pt depth6pt}\kern-.1em\vbox{\hrule width 0.33in}{\vrule height0pt width1pt depth6pt}\kern -.1em\vbox{\hrule width .65in}{\vrule height0pt width1pt depth6pt}\kern -.1em\vbox{\hrule width 1.0in}{\vrule height0pt width1pt depth6pt}\par
\indent (a) $(\forall x((P(x)\land(\exists zQ(x,z))) \to (\exists yM(x,y))))\land Q(x).$\par
\hskip 1.15in{\vrule height6pt width1pt depth0pt}\kern-.1em\vbox{\hrule width .40in}{\vrule height6pt width1pt depth0pt}\hskip .52in {\vrule height6pt width1pt depth0pt}\kern-.1em\vbox{\hrule width .42in}{\vrule height6pt width1pt depth0pt}\par
\vskip -4pt\indent \hskip .4in{\vrule height0pt width1pt depth6pt}\kern-.1em\vbox{\hrule width 0.33in}{\vrule height0pt width1pt depth6pt}\kern -.1em\vbox{\hrule width .65in}\kern -.1em\vbox{\hrule width 1.02in}{\vrule height0pt width1pt depth6pt}\par
\indent (b) $(\forall y((P(y)\land(\exists xQ(x,z))) \to (\exists zM(y,z))))\land Q(z).$\par
\hskip 1.15in{\vrule height6pt width1pt depth0pt}\kern-.1em\vbox{\hrule width .25in}{\vrule height6pt width1pt depth0pt}\hskip .73in {\vrule height6pt width1pt depth0pt}\kern-.1em\vbox{\hrule width .38in}{\vrule height6pt width1pt depth0pt}\par
\vskip -4pt\indent \hskip .4in{\vrule height0pt width1pt depth6pt}\kern-.1em\vbox{\hrule width 0.33in}{\vrule height0pt width1pt depth6pt}\kern -.1em\vbox{\hrule width .65in}\kern -.1em\vbox{\hrule width 1.02in}{\vrule height0pt width1pt depth6pt}\par
\indent (c) $(\forall x((P(x)\land(\exists xQ(x,z))) \to (\exists yM(x,y))))\land Q(z).$\par
\hskip 1.15in{\vrule height6pt width1pt depth0pt}\kern-.1em\vbox{\hrule width .25in}{\vrule height6pt width1pt depth0pt}\hskip .73in {\vrule height6pt width1pt depth0pt}\kern-.1em\vbox{\hrule width .38in}{\vrule height6pt width1pt depth0pt}\par
\indent \hskip .4in{\vrule height0pt width1pt depth6pt}\kern-.1em\vbox{\hrule width 0.35in}{\vrule height0pt width1pt depth6pt}\kern -.1em\vbox{\hrule width .78in}{\vrule height0pt width1pt depth6pt}\kern -.1em\vbox{\hrule width 0.87in}{\vrule height0pt width1pt depth6pt}\par
\indent (d) $(\forall z((P(z)\land(\exists xQ(x,z))) \to (\exists yM(z,y))))\land Q(x).$\par
\hskip 1.15in{\vrule height6pt width1pt depth0pt}\kern-.1em\vbox{\hrule width .25in}{\vrule height6pt width1pt depth0pt}\hskip .73in {\vrule height6pt width1pt depth0pt}\kern-.1em\vbox{\hrule width .38in}{\vrule height6pt width1pt depth0pt}\par
\baselineskip=14pt
\medskip
In these line diagrams$,$ each vertical line segment attached to a line segment
identifies the bound occurrences for a specific quantifier. Notice that in (a) 
the $x$ has a bound  occurrence and is free  in the formula. Thus the concepts 
of free in a formula and bound occurrences are not mutually exclusive. 
In (b)$,$ the 
$z$ has a bound occurrence and free in the formula. 
In (c)$,$ $z$ is the only variable free in the 
formula. In (d)$,$ $x$ has both a bound occurrence and is free in the 
formula.\parm
{\bf Example 3.2.3} The use of subscripts to indicate the bound occurrences
of a variable AND its bounding quantifier.\parm 
(a) $\forall x_1((P(x_1) \land (\exists z_2Q(x_1,z_2)))\to (\exists y_3 
M(x_1,y_3)))\land 
Q(x).$\parm
(b) $\forall y_1((P(y_1) \land (\exists x_2Q(x_2,z)))\to (\exists z_3  
M(y_1,z_3)))\land 
Q(z).$\parm
 The line diagram or number patterns themselves turn out to be of 
significance. \parm
\hs
{\bf Definition 3.2.3} (Congruent Formula) Let $A,\ B \in Pd$ and $A,\ B$ are in 
atomic form. Assume that all the bound occurrences of the variable in both $A$ 
and $B$ have been determined either by the line segment method or numbering 
method. If in both $A$ and $B$ all the bound variables and only the bound 
variables are erased and the remaining geometric form for A and for B are 
exactly 
the same (i.e. congruent)$,$ then the formula are said to be {{\it congruent.}}\parm
\hm
The basic reason the concept of congruence is introduced is due to the 
following theorem. After the necessary machinery is introduce$,$ it can be 
established. The $\iff$ that appears in this theorem behaves exactly as it 
does in the language $L.$\parm
{\bf Theorem 3.2.1} {\sl If $A,\ B\in Pd$ and $A,\ B$ are congruent$,$ then
$\vdash A\iff B.$}\parm
{\bf Example 3.2.3} In this example$,$ we look at previous formula and determine 
whether they are congruent. \parm
(a) In example 3.2.2$,$ formula (a) is not congruent to (b) since the variable 
that occurs free in $Q(\underline{\quad})$ in (a) is $x$ while the variable 
that occurs in $Q(\underline{\quad})$ in (b) is $z.$ \pars
(b) Formulas (b) and (c) in example 3.2 are congruent.\parm
\hs
{\bf Definition 3.2.4} (Sentences) Let $A \in Pd.$ Then $A$ is a {{\it 
sentence}} or {{\it is a closed formula}} if there are NO variables that are 
free in the formula $A.$ 
\pars
\hm
For this text$,$ sentences will be the most important formula in $Pd.$  Although 
one need not restrict investigations to sentences only$,$ these other 
investigations are$,$ usually$,$ only of interest to logicians. Indeed$,$ most 
elementary text books concentrate upon the sentence concept due to 
their significant applications. Further$,$ you can always assign the concept of 
``truth'' or  ``falsity'' (occur or won't occur) to sentences. The  ``always 
assign'' in the last sentence does not hold for formula in general.\parm
\hs
\hm 
\centerline{\bf EXERCISES 3.2}\parm
\noindent 1. List the scope for each of the numbered quantifiers. \parm
(a) $(\forall_1x (\exists_2x Q(x,z))) \to (\exists_3xQ(y,z)).$\pars
(b) $(\exists_1x(\forall_2y(P(c) \land Q(y)))) \to (\forall_3xR(x)).$\pars
(c) $(P(c) \land Q(x)) \to (\exists_1y(Q(y,z) \to (\forall_2xR(x)))).$\pars
(d) $\forall_1z((P(z) \land(\exists_2xQ(x,z)))\to(\forall_3z(Q(c) \lor 
P(z)))).$\parm
\noindent 2. Use the subscript or line segment method to display the variables 
that are bound by a specific symbol $\forall$ or $\exists.$\parm
(a) $(\forall z(\exists y(P(z,y)\land (\forall zQ(z,x))) \to M(z))).$\pars
(b) $(\forall x(\exists y(P(x,y)\land (\forall yQ(y,x))) \to M(x))).$\pars
(c) $(\forall z(\exists x(P(z,x)\land (\forall zQ(z,y))) \to M(z))).$\pars
(d) $(\forall y(\exists z(P(y,z)\land (\forall zQ(z,x))) \to M(y))).$\pars
(e) $(\forall y(\exists z(P(z,y)\land (\forall zQ(z,x))) \to M(y))).$\pars
(f) $\exists x(\forall z(P(x,z) \lor (\forall u M(u,y,x)))).$\pars
(g) $\exists y(\forall x(P(z,x) \lor (\forall x M(x,u,y)))).$\pars
(h) $\exists y(\forall x(P(y,x) \lor (\forall x M(x,y,z)))).$\pars
(i) $\exists z(\forall x(P(z,x) \lor (\forall x M(x,y,z)))).$\pars
(j) $\exists x(\forall x(P(z,x) \lor (\forall z M(x,y,z)))).$\parm
\noindent 3. For the formula in [2] list as ordered pairs all formula that 
are congruent.\parm
\noindent 4. In the following formula$,$ make a list as follows: write down the 
formula identifier$,$ followed by the word ``free$,$'' followed by a list of the 
variables that are free in the formula$,$ followed by the word ``bound$,$'' 
followed by the variables that have bound occurrences. (e.g.
(g) Free$,$ $z,x,z$; Bound $x,z.$)\parm
(a) $(P(x,y,z) \lor (\exists x(P(y) \to Q(x)))) \land (\forall zR(z)).$\pars
(b) $(\forall xP(x,x,x))\to (Q(z) \land (\forall zR(z))).$\pars 
(c) $(\forall x(P(c) \land Q(x))\to M(x))\iff P(y).$\pars
(d) $(\forall x(\forall y(\forall z(P(x) \to Q(z)))))\lor Q(z).$\pars
(e) $(P(c) \land Q(x)) \to ((\exists x Q(x))\land (\exists y(Q(y,x) \to 
M(d)))).$\pars
(f) $\forall x(\exists y((P(x) \land Q(y)) \to (P(c) \land Q(c)))).$\parm
\noindent 5.\quad (a) Which formula in problem 1 are sentences?\pars
(b) Which formula in problem 2 are sentences?\pars
(c) Which formula in problem 4 are sentences?\pars 
\ss
\noindent {\bf 3.3 Structures}\parm
In the discipline of mathematical logic$,$ we use the simplest and most 
empirically consistent processes known to the human mind to study the human 
experience of communicating by strings of symbols and logical deduction. 
Various human abilities are necessary in order to investigate logical 
communication. One must be able to recognize that the symbols within the 
quotation marks ``a'' and ``b'' are distinct. Moreover$,$ one must use various 
techniques of ``ordering'' in order to communicate in English and most other 
languages. An individual must intuitively know which string of symbols starts 
a communication and the direction to follow in order to obtain the next string of 
symbols. In the English language$,$ this is indicated by describing the 
direction as from ``left to right'' and from ``top to bottom.'' Without a 
complete intuitive knowledge of these  ``direction concepts'' no 
non-ambiguous English 
language communication can occur. These exact same intuitive processes are 
necessary in the mathematical discourse. Mathematics$,$ in its must fundamental 
form$,$ is based upon such human experiences. \pars
In 1936$,$ Tarski introduced a method to produce conceptually  the semantical 
``truth'' or  ``falsity'' notion for a predicate language $Pd$ that mirrors the 
concepts used for the propositional language $L\p.$ The Tarski ideas 
are informal in character$,$ are based upon a few aspects of set-theory$,$ but 
these aspects 
cannot be considered as inconsistent$,$ in any manner$,$ since by the above 
mentioned human abilities and experiences concrete symbol string examples can
be given for these specific Tarski notions. The only possible difficulty with 
the Tarski concepts lies in the fact that in order for certain sentences 
in $Pd$ to 
exhibit ``truth'' or ``falsity,'' the concept would need to be extended to 
non-finite sets. When this happens$,$ some mathematicians take the view that
the method is weaker than the method of the formalistic demonstration. In 
this case$,$ however$,$ most mathematicians believe that the Tarski method 
is not weaker than formal demonstrations  
since no inconsistency has occurred in using the most simplistic of the 
non-finite 
concepts. This simplistic non-finite notion has been used for over 3,000 
years. This non-finite notion is the one associated with the set 
called the natural numbers $\nat$.\pars
 In this text$,$ the natural numbers$,$ beginning with 
zero$,$ are considered as the most basic of intuitive mathematical 
objects and will not be formally discussed. This Tarski approach can be 
considerably simplified when restricted to the 
{{\it set of all sentences}} $\c S.$ It's this simplification that we  
present next. \pars
We first need to recall some of the most basic concepts$,$ concepts that appear 
even in many high school mathematics courses. Consider the following set of two 
alphabet letters $A=\{a,a,b\}.$ First$,$ the letters are not assumed  
to be in anyway different simply due to the left to right order in which they 
have been written or read. There is an ``equality'' defined for this set. Any 
two letters in this set are equal if they are intuitively the same letter. One 
sees two ``a'' symbols in this set. In sets that contain identifiable objects
such as the two ``a'' symbols$,$ it's an important aspect of set-theory that 
only one such object should be in the set. Hence$,$ in general$,$ 
specific objects 
in a set are considered as unique or distinct. On the other hand$,$ variables can 
be used to  ``represent'' objects. In this case$,$ it makes sense to write such 
statements such as  ``x = y.'' Meaning that both of the variables ``x'' and ``y'' 
represent a unique member of the set under discussion. \pars
In mathematical logic$,$ another refinement is made. The set $A$ need not be 
composed of just alphabet symbols. It might be composed of midshipmen. In this 
case$,$ alphabet symbols taken from our list of constants are used to name the 
distinct objects in the set. Take an $a \in C.$ In order to differentiate 
between this name for an object and the object it represents in a set$,$ we use the new 
symbol $a\p$ to denote the actual object in the set being named by $a$. \parm
Only one other elementary set-theory notion is necessary before we can define 
the Tarski method. This is the idea of the  {\it ordered pair.} This is where$,$ 
for the intuitive notion$,$ the concept of left to right motion is used. Take 
any one member of $A$. Suppose you take $a$. Then the symbol $(a,a)$ is an 
ordered pair. Using the natural number counting process$,$ starting with 1$,$ the 
set members between the  left parenthesis  ``('' and the right parenthesis
``)'' are numbered. The first position$,$ moving left to right$,$  is called the   {\it first coordinate}$,$ 
the second position is called the {\it second coordinate.} In this example$,$ 
both coordinates are equal under our definition of equality.\pars
 We can construct 
 the set of ALL ordered pairs from the members in $A.$ This set is denoted by 
$A^2$ or $A \times A$ and contains the distinct objects 
$\{(a,a),(b,b),(a,b)(b,a)\}.$  What one needs to do is to use the equality 
defined for the set $A$ to {\it define} an equality for ordered pairs. This 
can be done must easily by using a word description  for this ``ordered pair 
equality.'' Two ordered pairs are ``equal'' (which we denote by the symbol =) 
if and only if the first coordinates are equal AND the second coordinates are 
equal.\pars
Well$,$ what has been done with two coordinates can also be done with 3$,$ 4$,$ 5$,$ 
25$,$ 100$,$ 999 coordinate positions. Just write down a left  ``('' a member of 
$A,$ followed by a comma$,$ a member of $A,$ followed by a comma$,$ continue the 
described process$,$ moving left to right$,$ until you have$,$ say 999 coordinates but follow the last one 
in your symbol not by a comma but by a right ``)''. What this is called is a
999-tuple. The set of ALL 999-tuples so constructed from the set $A$ is denoted 
by $A^{999}$ rather than  writing the symbol $A$ 999 times and putting 998
$\times$ symbols between them. For any natural number $a$ greater than 1,
$A^n$ denotes the set of ALL n-tuples that can be constructed from the members 
of $A.$ \pars
Notice that we have defined equality of order 
pairs by using the coordinate numbers. To define equality for n-tuples$,$ where 
$n$ is any natural number greater than 1$,$ simply extend the above 
2-tuple (ordered pair) definition in terms of the position numbers. Thus two
n-tuples are ``equal'' if and only if the first coordinates are equal AND
the second coordinates are equal AND . . . . AND the nth coordinates are 
equal.\pars
The last set theory concept is the simple {\it subset} concept.  A set $B$ is 
a subset of the set $A$ if and only if every object in $B$ is in $A,$ where 
this statement is symbolized by $B \subset A.$ The symbols ``in $A$'' mean 
that you can recognize that the members of $B$ are explicitly members of $A$ 
by their properties. There is 
a special set called the {\it empty set} (denoted by $\emptyset$)  
that is conceived of as have NO members. From our definition for ``subset,''
the empty set is a subset of any set since the $\emptyset$ contains no 
members and$,$ thus$,$ it certainly follows that all of its members are 
members of any set. Lastly$,$ a subset $R$ of $A^n,$ for $n >1,$ will be called 
an n-place relation or an n-ary relation. \parm
\hs  
{\bf Definition 3.3.1} (Structures) Let $Pd$ be a first-order language 
constructed from a nonempty finite or infinite set $\{P_1,P_2,\ldots\}$ 
of predicates and an appropriate set $C$ of constants that satisfies (iv). For our basic application$,$ a {{\it structure with an interpretation}} 
is an object ${\c M}= \langle D, \{P_1\p, P_2\p, 
\ldots\}\rangle$ where\pars
(i)  $D$ is a nonempty set 
called the {{\it domain}}. 
\pars
(ii) each n-place predicate $P_i,$ where $n > 1$$,$ corresponds to an n-place 
relation $P_i\p \subset D^n,$ and \pars
(iii) the 1-place predicates correspond to specific subsets 
of $D$ and \pars
(iv) there is a function $I$ that denotes the correspondence in (ii) and (iii) 
and that corresponds constants $c_i\in N \subset C$ 
to 
members $c_i\p$ in $D,$ (this function $I$ is a naming function).  Constants in $N$ correspond to one and only one member in $D$ and each member of $D$ corresponds to one and only one member of $N.$ Due to how axioms are stated$,$ the structure symbol ${\c M}= \langle D, \{P_1\p, P_2\p, 
\ldots\}\rangle$ may include a set of distinguished elements that correspond to distinguished constants in $\cal C$. In this case$,$ these constants are considered as contained in $N$ and correspond under $I$ to fixed members of $D.$ [Only members of $\cal C$ appear in any member of $Pd.$]\pars 
The rule you use to obtain the correspondences between the language 
symbols and the set theory objects is called an {{\it interpretation.}} In general,  ${\c M}= \langle D, \{P_1\p, P_2\p, 
\ldots\}\rangle$ is a {\it structure} for various interpretations and includes distinguished elements.\pars
\hm
   
{\bf Example 3.3.1} Let our language be constructed from  a 1-place 
predicate $P,$ a 2-place predicate $Q$ and two constants $a,b.$ For the domain 
of our structure$,$ let $D =\{a\p,b\p,c\p\}.$ Now correspond the one place 
predicate $P$ to the subset $\{a\p,b\p\}.$ Let the 2-place predicate $Q$ 
correspond to $Q\p = \{(a\p,a\p),(a\p,b\p)\}.$ Finally$,$ let $a$ correspond to 
$a\p$ and $b$ correspond to $b\p.$ The interpretation $I$ can be symbolized as 
follows: $I$ behaves like a simple function. $I$ takes the atomic portions of 
our language and corresponds then to sets. $I(P) = P\p,\ I(Q) = Q\p,\ 
I(a)=a\p,\ I(b) =b\p,$  $a,\ b \in C.$\pars
Since the only members of $Pd$ that will be considered as having a truth-value 
will be sentences$,$ the usual definition for the truth-value for members in 
$Pd$ can be simplified. A level by level inductive definition$,$ similar to that 
used for the language $L$ is used to obtain the truth-values for the 
respective sentences. This will yield a valuation function $v$ that is 
dependent upon the structure. \pars
There is one small process that is needed to define properly this valuation 
process. It is called the {{\it free substitution operator.}} We must know the 
scope of each quantifier and the free variables in the formula that is the 
scope.\pars 
\hs
{\bf Definition 3.3.2} (Free substitution operator) Let the symbol 
$S^x_\lambda$ have a language variable as the superscript and a language 
variable or language constant as a subscript. For any formula $A\in Pd,$
$S^x_\lambda A]$ yields the formula where $\lambda$ has been substituted for 
every free occurrence of $x$ in $A.$\pars
\hm
{\bf Example 3.3.2} Let $A = \forall y(P(x,y) \to Q(y,x)).$  \pars
(i) $S^x_y A] = \forall y(P(y,y) \to Q(y,y))$\pars
(ii) $S^x_x A] = \forall y(P(x,y) \to Q(y,x)).$\pars
(iii) $S^x_c A] = \forall y(P(c,y) \to Q(y,c)).$\parm  
\hs
{\bf Definition 3.3.3} (Structure Valuation for sentences). Given a structure 
$\c M,$ with domain $D,$ for a 
language $Pd.$  \pars
(i) Suppose that $P(\underline{\quad},\ldots,\underline{\quad}) \in Pd_0$ is an n-place 
predicate $n \geq 1,$ that contains only constants $c_i \in N$ in each of the places. Then  ${\c M}\mod P(c_1,\ldots,c_n)$ if and only if
$(c_1\p,\ldots,c_n\p) \in P\p,$ where for any $c_i \in N,$  
$c_i^\prime = I(c_i),$ and constants that denote {\it special} required objects such as ``0'' or ``1'' and the like denote fixed members of $D$ and are fixed members of $D$ throughout the entire valuation
process. Note that $(c_1') \in P'$ means that $c_1' \in P'.$  [Often in this case$,$ such a set constants is said to ``satisfy'' (with respect to $\c M$,) 
the predicate(s) or formula.]  Since we are only 
interested in  ``modeling'' sentences$,$ we will actually let the ``naming'' subset $N$ of $C$ vary in such a manner such that $I(N)$ varies over the entire set $D.$  
 \pars
(ii) If (i) does not hold$,$ then ${\c M}\not\mod P(c_1,\ldots,c_n).$   \pars
(iii) Suppose that for a level m$,$ the valuation $\mod$ or $\not\mod$ has been 
determined for formula $A$ and $B,$ and specific members $c_i^\prime$ of 
$D.$\pars
\indent\indent (a) ${\c M}\not\mod A\to B$ if and only if ${\c M}\mod A$ and 
${\c M}\not\mod B.$ In all other cases$,$ ${\c M}\mod A \to B.$\pars
\indent\indent (b) ${\c M}\mod A \iff B$ if and only if ${\c M}\mod A$ and
${\c M}\mod B,$ or ${\c M}\not\mod A$ and  ${\c M}\not\mod B.$ \pars
\indent\indent (c) ${\c M}\mod A\lor B$ if and only if ${\c M}\mod A$ or ${\c 
M}\mod B.$\pars
\indent\indent (d) ${\c M}\mod A\land B$ if and only if ${\c M}\mod A$ and ${\c 
M}\mod B.$\pars
\indent\indent (e) ${\c M} \mod (\neg A)$ if and only if ${\c M}\not\mod 
A.$\parm
Take note again that any constants that appear in the original 
predicates have been assigned FIXED members of $D$ and never change throughout this valuation for a given structure. In what follows$,$ $d$ and $d\p$ are used to denote arbitrary members of $N$ and the corresponding members of $D.$ Steps (f) and (g) below show how the various constants are obtained that are used for the previous valuations.\parm
\indent\indent (f) For each formula $C=\forall x A,$ ${\c M}\mod \forall x A$ 
if and only 
if for every $d\p \in D$ it follows that ${\c M}\mod S^x_d A].$ Otherwise,
${\c M}\not\mod \forall x A.$ \pars
\indent\indent (g) For each formula $C=\exists xA,$ ${\c M} \mod \exists x A$    
if and only 
if there is some $d\p \in D$ such that ${\c M}\mod S^x_d A].$ Otherwise$,$ 
${\c M} \not\mod \exists x A.$ \pars
(iv) Note that for 1-place predicate $P(x)$ (f) and (g) say$,$ that 
${\c M} \mod \forall x P(x)$ if and only if  for each $d\p \in D,\ d\p \in P\p$ and 
${\c M} \mod \exists x P(x)$ if and only if there exists some $d\p \in D,$ such that 
$d\p \in P\p.$ \parm                                                     
\hm
Please note that$,$ except for the quantifiers$,$ the valuation process follows 
that same pattern as the $T$s and $F$s follow for the truth-value valuation of 
the language $L.$ Simply associate the symbol ${\c M}\mod A$ with the $T$ and 
the ${\c M}\not\mod A$ with the $F.$ Since every formula that is valuated is a 
sentence$,$ we use the following language when ${\c M}\mod A.$ The structure $\c 
M$ is called a {{\it model}} for $A$ in this case. If ${\c M}\not\mod A,$ 
then we say that ${\c M}$ is {{\it not a model}} for $A.$ This yields the 
same logical pattern as the  ``occurs'' or  ``does not occur'' concept that 
can be restated as  ``as $\c M$ models'' or  ``as $\c M$ does not 
model.'' {\bf Definition 3.3.3 is informally applied in that metalogical arguments are used relative to the informal notions of ``there exists'' and ``for each'' as they apply to members of sets.} For involved sentences$,$ 
this process can be somewhat difficult. \pars
In the appendix$,$ it is shown that the process described in definition 3.3.3 is unique for 
every structure in the following manner. Once an interpretation $I$ is 
defined$,$ then all valuations that proceed as described in definition 3.3.3 
yield the exact same results. Definition 3.3.3 is very concise. To illustrate what it means$,$ let $D= \{a\p, 0\p\},\ A = P(x,y) \to Q(x)$ and $P\p, Q\p  \in {\c M}.$ It is determine whether ${\c M}\mod P(a,0), \  {\c M}\mod P(a,a), \ {\c M}\mod P(0,a), \ {\c M}\mod P(0,0), \ {\c M}\mod Q(a),\ \ {\c M}\mod Q(0).$ Now we use these results to determine whether ${\c M}\mod \forall x(P(x,0) \to Q(x)).$ This occurs only if ${\c M}\mod P(a,0) \to Q(a),\ {\c M}\mod P(0,0) \to Q(0),$ where we use the previous determinations. \parm
{\bf Example 3.3.3} (a) Let $A= (\forall x(\exists y P(x,y))) \to (\exists y(\forall 
xP(x,y))).$ Let the domain $D= \{a\p\}$ be a one element set. We 
define a structure for $A= (\forall x(\exists y P(x,y))) \to (\exists y(\forall 
xP(x,y))).$ Consider the 2-place relation $P\p = \{(a\p,a\p)\}$ and the 
interpretation $I(P) = P\p.$ As you will see later$,$ one might try to show 
this by simply assuming that  ${\c M} \mod (\forall x(\exists y P(x,y))).$ But 
this method will only be used where $\c M$ is not a specific structure. For 
this example$,$ we determine whether ${\c M} \mod (\forall x(\exists y 
P(x,y)))$ and whether ${\c M} \mod (\exists y(\forall 
xP(x,y)))$  and use part (a) of definition 3.3.3.
The valuation process proceeds as follows: first we have only the one statement that ${\c M} \mod P(a,a).$ Next$,$ given each $c^\prime \in D$$,$ is 
there some $d^\prime \in D$ such that $(c^\prime,d^\prime) \in P^\prime$?
Since there is only one  member $a^\prime$ in $D$ and $(a^\prime,a^\prime) \in 
P^\prime,$ then ${\c M} \mod (\forall x(\exists y P(x,y))).$ 
 Now we test the statement $\exists y(\forall x P(x,y)).$ 
Does there exist a $c\p\in D$ such that ${\c M} \mod S^y_c \forall x P(x,y)] = 
\forall x P(x,c)$; which continuing through the substitution process$,$ does 
there exists a $c\p\in D$ such that for all $d\p\in D$ in $(d\p,c\p) \in P\p.$ 
The answer is yes$,$ since again there is only one element in $D.$ Thus ${\c M} \mod \exists y(\forall x P(x,y)).$ Hence$,$ part 
(v) (a) of definition 3.3.3 implies that ${\c M}\mod A.$ Clearly$,$ the more 
quantifiers in a formula the more difficult it may be to establish that a 
structure is a model for a sentence. \parm
(b) A slighter weaker approach is often used. The informal theory of natural numbers is used as the basic mathematical theory for 
the study of logical procedures. It's considered to be a consistent theory 
since no contradiction has been produced after thousands of years of theorem 
proving. Thus a structure can be constructed from the theory of natural 
numbers. Let $D =\nat$ and there is two constants $1,2$ and one 3-placed
relation $P\p$ defined as follows: $(x,y,z) \in P\p$ if and only if 
$ x \in \nat,\ y \in \nat$ and $z \in \nat$ and $x + y = z.$ Using the theory 
of natural numbers$,$ one can very quickly determine whether ${\c M} \mod A$ for 
various sentences.\pars
Consider $A = \forall x P(c_1,c_2,x).$ Then let $c_1 = 1,\ c_2 = 2$. [Note that the subscripts for the constants are not really natural numbers but are tick marks.] Now consider the requirement that 
for all $a\p \in \nat,$ ${\c M} \mod P(1,2,a).$ Of course$,$ since 
$(1\p,2\p,5\p)\notin P\p,$ then ${\c M}\not\mod A.$ You will see shortly that this would imply that $A$ is not valid.\parm
(c) Using the idea from (b)$,$ the sentence $\exists x \forall y P(x,y)$ has ${\c M} =<\nat, \leq >$ as a model. For there exists a natural number $\r x = 0$ such that $0\leq \r y$ for each natural number $\r y$.\parm{\bf Example 3.3.4} The following shows how the instructions for structure valuation can be more formally applied. \pars
For ${\c M} = \langle D, P \rangle,$ let $D = \{a',b'\},$ and $P'= \{(a',a'),(b',b')\}.$ we want to determine whether ${\c M}\mod (\exists y(\forall xP(x,y))) \to (\forall x(\exists y P(x,y))),$ and whether ${\c M}\mod (\forall x(\exists y P(x,y)))\to (\exists y(\forall xP(x,y)))$ and also if ${\c M}\mod \forall x(\exists x P(x,x).$ \pars
First we look at whether ${\c M}\mod (\exists y(\forall xP(x,y)).$ This means that first we establish that ${\c M}\mod P(a,a),\ {\c M}\not\mod P(a,b),\ {\c M}\mod P(b,b), \ {\c M}\not\mod P(b,a).$\pars
For the next stage we must determine whether $({\rm a})\ {\c M}\mod \forall xP(x,a),\ 
(\r b)\ {\c M}\mod \forall xP(x,b).$ Under the substitution requirement, for (a) that (i) ${\c M}\mod P(a,a)$ and (ii) ${\c M}\mod P(b,a).$ However, ${\c M}\not\mod P(b,a).$ Thus for (a) ${\c M}\not\mod \forall xP(x,a).$ In the same manner, it follows that for (b) ${\c M}\not\mod \forall xP(x,b).$ Under the substitution requirement $S^y_d]\forall xP(x,y)$ produces the two formula $\forall xP(x,a),\ \forall xP(x,b)$ and if one or the other or both satisfy $\mod$ then we know the ${\c M}\mod$ holds. But, from (a) and (b) we  
know this is not the case so ${\c M}\not\mod \exists y(\forall x P(x,y)).$ Now from our understanding of the $\to$ connective this implies that ${\c M}\mod (\exists y(\forall xP(x,y))) \to (\forall x(\exists y P(x,y))).$\pars
Now consider $\forall x(\exists y P(x,y)).$ The first step looks at $\exists y P(a,y)$ and $\exists y P(b,y).$ In both cases, we have that ${\c M}\mod \exists y P(a,y)$ and ${\c M}\mod \exists y P(b,y).$ Hence, we have that ${\c M}\mod \forall x(\exists x P(x,y)).$ But this now implies that ${\c M}\not\mod (\forall x(\exists y P(x,y)))\to (\exists y(\forall xP(x,y))).$ \parm 

This is how one ``thinks'' when making and arguing for these valuations and you can actually make the substitution if $D$ is a finite set and get a formula in $Pd.$ The $\forall$ can be replaced with a set of $\land$ symbols, one for each member of $D$ and the $\exists$ can be replaced with a set of $\lor$ symbols one for each member of $D.$ It follows easily that ${\c M} \mod \forall x(\exists y P(x,y))$ if and only if ${\c M} \mod (P(a,a)\lor P(a,b))\land (P(b,a)\lor P(b,b)).$]\pars
What about  $\forall x(\exists x P(x,x))$? The substitution process
says that under the valuation process this $\forall x(\exists x P(x,x))$ is the same as $\exists x P(x,x).$ The valuation process is ordered, first we do the $S^x_d P(x,x).$ This gives under the first step the formulas $P(a,a),\ P(b,b).$ The second valuation process does not apply since there are no free variables in $\exists x P(x,x).$ The ``or'' idea for $\exists$ yields that 
${\c M}\mod \forall x(\exists x P(x,x).$ \pars
We now change the structure to $D = \{a',b'\},$ and $P'= \{(a',b'),(b',a')\}$ and interpret this as same before. One arrives at the same conclusions that ${\c M}\not\mod (\exists y(\forall xP(x,y))$ and that ${\c M}\mod \forall x(\exists x P(x,y))$. But, for this structure 
 ${\c M}\not\mod \forall x(\exists x P(x,x))$. Hence, both structure are models for the sentences ${\c M}\mod (\exists y(\forall xP(x,y))) \to (\forall x(\exists y P(x,y))),\ \forall x(\exists x P(x,x)).$ But for these three sentences, the second structure is only a model for $\exists y(\forall xP(x,y))) \to (\forall x(\exists y P(x,y))).$ \pars
Note: A single domain and collection of n-place relations, in general, may have many interpretations. This follows from considering interpretations that use different members of $D$ for the distinguished constants, if any, and the fact that any number of n-placed predicates can be interpreted as the same n-place relation. Thus for the above two structures if you have formula with five two-placed predicates, then you can interprete them all to be $P'.$ However, if you are more interested in the logical behavior relative to a fixed structure, then in this case, such a formula with the five n-place predicates holds in this structure if and only if the sentence you get by replacing the five predicates with one predicate $P(x,y)$ holds 
in the structure. [If you have different distinguished constants in these predicates they still are different in the place where you changed the predicate symbol name to $P$.] The obvious reason for this is that the actual valuation does not depend upon the symbol used to name  the n-placed predicate (i.e. the $P$, $Q$, $P_1$ etc.) but only on the function and the relation that interpretes the predicate symbol. For this reason, sometimes you will see the definition for an interpretation state that the correspondence between predicate symbols and relations be one-to-one.\parm   
  
\hs\hm
\centerline{\bf EXERCISES 3.3}\parm
\noindent 1. In each of the following cases$,$ write the formula that is the 
result of the substitution process.\parm
\+(a) $S^x_a\ (\exists x P(x)) \to R(x,y)]$&(d) $S^x_aS^x_b (\exists x P(x)) \to 
R(x,y)]]$\cr
\+(b) $S^y_x\ (\exists y R(x,y))\iff (\forall x R(x,y))]$&(e) $S^x_a S^y_x 
(\exists y R(x,y))\iff (\forall x R(x,y))]]$   \cr
\+(c) $S^y_a\ (\forall x P(y,x)) \land (\exists y R(x,y))]$&(f) $S^x_aS^y_b 
(\forall z P(y,x)) \land (\exists y R(x,y))]]$\cr \parm
\noindent 2. Let $A \in Pd.$ Determine whether 
the following is ALWAYS true or not. If the statement does not hold for all 
formula in $Pd,$ then give an example to justify your claim.\parm
\+(a) $S^x_aS^y_b A]] = S^y_bS^x_a A]]$&(c) $S^x_yS^z_w A]]=S^z_wS^x_y A]]$\cr
\+(b) $S^x_y A] = S^y_x A]$&(d) $S^x_x S^y_y A]] = A$\cr
\noindent 3. Let $D = \{a\p,b\p\},$ the 1-place relation $P\p = \{a\p\}$ and the 2-place relation $Q\p = \{(a\p,a\p),(a\p,b\p)\}.$ 
For each of the following sentences and interpretation of the constant $c$$,$ 
$I(c) = a\p,$
determine whether ${\c M} \mod A,$ where ${\c M} = \langle D,\ a\p,\ P\p,\ Q\p 
\rangle.$\parm
\noindent (a) $A = (\forall x(P(c) \lor Q(x,x))) \to  (P(c) \lor \forall xQ(x,x)).$\pars
\noindent (b) $A = (\forall x(P(c) \lor Q(x,x))) \to  (P(c) \land \forall xQ(x,x)).$\pars
\noindent (c) $A = (\forall x(P(c) \lor Q(x,x))) \to  (P(c) \land \exists xQ(x,x)).$\pars 
\noindent (d) $A=  (\forall x(P(c) \land Q(x,x)))\iff (P(c) \land \forall xQ(x,x)).$\pars 
\noindent (e) $A = (\forall x(P(c) \land Q(c,x)))\iff (P(c) \land \forall xQ(x,x)).$\pars 
\ss
\vfil\eject
\noindent {\bf 3.4 Valid Formula in $Pd.$.}\parm
Our basic goal is to replicate the results for $Pd$$,$ whenever possible$,$ that 
were obtained for $L$ or $L\p.$ What is needed to obtain$,$ at the least$,$ one 
ultralogic is a compactness type theorem for $Pd.$\pars
Since we are restricting the  model concept to the set of sentences ${\c S}$ 
contained in  
$Pd,$ 
we are in need of a method to generate the simplest sentence for formula 
that is not sentence. {\bf Always keep in mind that fact that a structure is 
defined for all  predicates and constants in a specific language} $Pd.$ In 
certain cases$,$ it will be necessary to consider special structures with 
special properties. \parm
\hs
{\bf Definition 3.4.1} (Universal closure) For any $A \in Pd,$  with free 
variables $x_1,\ldots, x_n,$ (written in subscript order the {{\it 
universal closure}} is denoted by $\forall A$ and $\forall A = \forall 
x_1(\cdots (\forall x_n A)\cdots ).$ \pars
\hm
Obviously$,$ if there are no free variables in $A$$,$ then $\forall A = A.$ In any 
case$,$ $\forall A \in \c S.$ This fact will not be mentioned when we only 
consider members of $\c S.$ We will use the same symbol $\mod,$ as previously 
used$,$ 
to represent 
the concept of a  ``valid'' formula in $Pd.$ It will be seen that it's the 
same concept as the $T$ and $F$ concept for $L.$ \parm
\hs
{\bf Definition 3.4.2} (Valid formula in $Pd.$) A formula $A \in Pd$ is a 
{{\it valid}} formula (denoted by $\mod A$) if for every domain $D$ and 
every structure $\c M$ with an 
interpretation $I$ for each constant in $A$ and each n-placed predicate in 
$A$$,$ ${\c M}\mod \forall A.$\pars
A formula $B$ is a {{\it contradiction}} 
if for every domain $D$ and every structure $\c M$ with an 
interpretation $I$ for each constant in $A$ and each n-placed predicates in 
$A$$,$ ${\c M}\not\mod \forall A.$\hm
We show that this is an extension of the valid formula concept as defined 
for $L.$ Further$,$ since every structure is associated with an interpretaion,
this will be denoted  by $\str.$\parm
{\bf Theorem 3.4.1} {\sl Let $A,\ B,\ C \in \c S.$ Suppose that $\mod A;\ 
\mod A \to B.$ Then $\mod B.$}\pars
Proof. Suppose that $\mod A$ and $ \mod A \to B.$ Let $\str$ be any 
structure for $Pd.$  From the hypothesis$,$ $\str \mod A$ and $\str 
\mod A\to B$ imply that $\str\mod B$ from definition 3.3.3 (consider 
the size $(A\to B) = m$) part (a).\qed
In the metaproofs to be presented below$,$ I will not continually mention the 
size of a formula for application of definition 3.3.3.\parm
{\bf Theorem 3.4.2} {\sl Let $A,\ B,\ C \in \c S.$ Then 
(i) $\mod A \to (B \to A)$$,$ (ii) $\mod (A \to (B\to C)) \to ((A\to B)\to (A 
\to C)),$ (iii) $\mod ((\neg A)\to (\neg B)) \to (B \to A).$}\pars
Proof. (i)$,$ (ii)$,$ (iii) follow in the same manner as does theorem 3.4.1 by 
restricting structures to specific formulas. Notice that every formula in the 
conclusions are members of $\c S$. We establish (iii).\pars
(iii) Let $\str$ be a structure for $Pd.$ 
Assume that $\str\mod (\neg A)\to (\neg B)$ for otherwise the result would 
follow from  definition 3.3.3. Suppose that $\str \mod \neg A$ and 
$\str\mod \neg B.$  Hence from definition 3.3.3$,$
$\str \not\mod A$ and $\str\not\mod B.$ Consequently$,$ $\str \mod B \to 
A \Rightarrow \str\mod ((\neg A)\to (\neg B)) \to (B \to A)$ from 
definition 3.3.3. Suppose that $\str \not\mod (\neg A).$ Then $\str\mod A
\Rightarrow \str \mod B \to A$ and in this final case the result also 
holds. \qed
{\bf Theorem 3.4.3} {\sl Let $A$ be any formal theorem in $L.$ Let $\Hyper A$ 
be obtained  by substituting for each  proposition $P_i$ a member $A_i \in \c 
S.$ Then $\mod \Hyper A.$}\pars
Proof. This follows in the exact same manner as the soundness theorem 2.15.2 
for $L\p$ extended to $L$ along with theorem 3.4.2. \qed
{\bf Theorem 3.4.4} {\sl Let $A\in L$ and $\mod A$ as defined in $L\p.$ Then 
$\mod \Hyper A$ as defined for $Pd.$}\pars
Proof. This follows from the completeness theorem 
2.15.1 for $L\p$ extended to $L$ and theorem 3.4.3.\qed
What theorems 3.4.2 and 3.4.4 show is that we have a great many valid formula in $Pd.$
However$,$ is this where all the valid formula come from or are there many 
valid formula in $Pd$ that do not come from this simple substitution process? 
\parm
{\bf Example 3.4.1} Let $A = P(x) \to (\exists xP(x)).$ Then $\forall A = 
\forall x(P(x) \to (\exists xP(x))).$ Let $\str$ denote a structure for 
$Pd.$ Suppose that $\str \not\mod\forall A.$ Then 
there is some $d\p \in D$ such that $\str\not\mod P(d) \to \exists x 
P(x).$ Thus it most be that $\str \mod P(d)$ and $\str\not\mod \exists 
xP(x).$ But this contradicts definition 3.3.3 (v) part (g). Thus $\mod 
P(x) \to (\exists xP(x)).$ It's relatively clear$,$ due to the location of the 
$\forall$ in the universal closure in the formula$,$ that there is no formula 
$B \in L\p$ such that $\Hyper B = \forall A.$\parm
There are many very important valid formula in $Pd$ that are not obtained from 
theorems 3.4.2 and 3.4.4. To investigate the most important$,$ we use the following {{\it 
variable predicate}} notation. Let $A$ denote a formula from $Pd.$ Then 
there are always three possibilities for an $x \in \c V.$ Either $x$ does not 
appear in $A,$ $x$ appears in $A$ but is not free in $A,$ or $x$ is free in 
$A.$ There are certain important formula$,$ at least for the proof theory portion of 
this chapter$,$ that are valid and that can be expressed in this general 
variable predicate language. Of course$,$ when such metatheorems are 
established$,$ you need to consider these three possibilities.\parm

With respect to our substitution operator $S^x_\lambda$$,$ if $A$ either does 
not contain the variable $x$ or it has no free occurrences of $x$$,$ then  
$S^x_\lambda A] = A.$ Further it is important to note that the constant $d$ 
that appears is a general constant that is relative to a type of extended 
interpretation where it corresponds to some $d\p.$ But it is not part of the 
original interpretation. This difference must be strictly understood. \parm
\hs 
{\bf Definition 3.4.3} (Free for) Let $A \in Pd.$ Then a variable $v$ is
{{\it free for x in A}} if the formula $S^x_vA],$ at the least$,$ has free 
occurrences of $v$ in the same positions as the free occurrences of $x.$\pars
\hm
{\bf Example 3.4.2} Notice (i) that $x$ is free for $x$ in any formula $A \in 
Pd$ and $S^x_xA]=A.$ \pars
(ii) Further$,$ if $x$ does not occur in $A$ or is not free in $A$$,$ then any 
variable $y$ is free for $x$ in $A$ and $S^x_y A] = A$ for any
$ y\in \c V.$  The only time one gets a different (looking)  
formula that {\it preserves free occurrences} through the use of the substitution operator 
$S^x_\lambda A]$ is when $x$ is free in $A$ and  $\lambda\in 
\c V$ is free for $x$ and $\lambda \not= x.$ \pars
(iii) Let $A = \exists yP(y,x).$ Then $y$ is NOT free for $x$ since
$S^x_y A] = \exists yP(y,y).$ You get a non-congruent formula by this 
application of the substitution operator. \pars
(iv) But $z$ is free for $x$ since $S^x_z A] = \exists yP(y,z).$ \pars
(v) If $A= (\exists xP(x,y)) \to (\exists yQ(x)),$ then $y$ is NOT free for $x$ 
since $S^x_yA] = (\exists xP(x,y) \to (\exists yQ(y)).$ Again $z$ is free for 
$x.$\parm
{\bf Theorem 3.4.5} {\sl For any formula $A$ with variables $x_1,x_2, \ldots, 
x_n$ and only these free variables (where the subscripts only indicate the 
number of distinct variables and not their subscripts in the set $\c V$)$,$ then 
for any structure $\str,$ $\str\mod \forall A$ if and only if 
for each $c_1\p\in D$ 
and 
for each $c_2\p\in D$ and $\cdots$ and for each $c_n\p \in D,$  
$\str\mod S^{x_1}_{c_1}S^{x_2}_{c_2}\ldots S^{x_n}_{c_n} A]\ldots 
]].$}\pars
Proof. From the definition of universal closure and definition 3.3.3.\qed
{\bf Corollary 3.4.5.1} {\sl Under the same hypotheses as theorem 3.4.5,
$\str \mod S^{x_i}_{c_i}S^{x_j}_{c_j}\ldots S^{x_k}_{c_k} A]\ldots 
]]$ for any permutation $(i,j,\ldots, k)$ of the subscripts.}\parm
{\bf Theorem 3.4.6} {\sl Let $y$ be free for $x$ in $A.$ Then 
$S^y_d S^x_y A]] = S^x_dS^y_d A]].$}\pars
Proof. The major argument to establish must of our validity results is 
dependent upon a rewording of the substitution process. If there are any free 
occurrences of $y,$ then $x$ does not occur free at those places that $y$ 
occurs free in $A.$  Substituting $d$ for these specific free occurrences 
can be done first. Then$,$ each $y$ obtained by substituting for a free $x,$ 
due to the fact that $y$ is free for $x,$ can be changed to a $d$ by simply 
substituting the $d$ for the free occurrences $x$. This yields the left hand 
side of the equation where all free occurrence of $x$s are changed to $y$s$,$ 
any other free  occurrence of $y$ remains as it is$,$ and then all the free 
occurrences of $y$ are changed to $d.$\qed
{\bf Theorem 3.4.7} {\sl Suppose that $A \in Pd$ and $y$ is free for $x$ in 
$A$$,$ then \pars
(i) if $x$ does not occur free in $A,$ then for any structure $\str$ 
for $Pd,$ $\str\mod \forall (\forall xA)$ if and only if $\str\mod \forall A$ if and only if $\str\mod \forall x(\forall A).$\pars 
(ii) If $x$ does not occur free in $A,$ then for any structure $\str$ 
for $Pd,$ $\str\mod \forall (\exists xA)$ if and only if $\str\mod \forall A$ if and only if $\str\mod \exists x(\forall A)$\pars 
(iii) $\mod (\forall xA) \to S^x_yA]$\pars 
(iv) $\mod S^x_yA] \to (\exists x A).$}\pars
Proof. (i)  Let $\str$ be  a structure for $Pd.$ Suppose that $\str \mod \forall (\forall x A).$ 
Since $x$ is not free in $A,$ then from corollary 3.4.5.1$,$ 
$\str \mod (\forall x\forall A)\Rightarrow \str\mod \forall A,$ since 
under the substitution 
process there is no free $x$ for the substitution.\pars
Conversely$,$ suppose that $\str\mod \forall A.$ Since $x$ is not free in $A,$ 
then again from corollary 3.4.5.1$,$  
$\str\mod \forall A \Rightarrow \str\mod \forall (\forall xA)$ for the same 
reason. \pars
(ii) This follows in the same manner as (i).\pars
(iii) First suppose that $x$ is not free in $A$. Then $S^x_yA] = A.$  
Let $\str$ be a structure $Pd$ and consider the sentence 
$\forall ((\forall xA) \to A).$ 
Suppose that $\{x_1, \ldots, x_n\}$ are 
free variables in $A.$ From 3.4.5$,$ we must show that 
$\str\mod S^{x_1}_{c_1} \cdots  S^{x_n}_{c_n}((\forall xA) \to A),$ where $c_j \in D.$
However$,$ making the actual substitutions yields that 
$S^{x_1}_{c_1} \cdots  S^{x_n}_{c_n}((\forall xA) \to A)=
S^{x_1}_{c_1} \cdots  S^{x_n}_{c_n}(\forall xA)]\cdots] \to  
S^{x_1}_{c_1} \cdots  S^{x_n}_{c_n} A]\cdots ]= \forall(\forall xA) \to 
\forall A.$ We need only suppose that $\str\mod \forall(\forall xA)$. 
Then from (i)$,$   
$\str\mod \forall A.$ But then $\str\mod \forall A \to \forall A$ from our 
definition 3.3.3. 
\pars
Now assume that $x$ is free in $A.$ 
Note that $x$ is not free in $\forall x A,$ or in 
$S^x_yA].$ There are free occurrences of $y$ in $S^x_yA]$ and there may be 
free occurrences of $y$ in $\forall xA.$ Any other variables that occur free 
in  $\forall x A,$ or  $S^x_yA]$ are the same variables. Considering the 
actual substitution process for $\forall((\forall xA)\to S^x_yA])$ and using 
corollary 3.4.5.1$,$ we can permute all the other substitution processes for the 
other possible free variables \underbar{not $y$} to be done ``first.'' 
When this is done 
the positions that the \underbar{arbitrary} $d$s take yield $\forall xC$ 
and $S^x_yC,$ where $C$ contains the various symbols $d$ in the place of the other possible free variables and that correspond to 
the $d'\in D.$ Consider $\forall y((\forall x A) \to S^x_yC]).$ 
 The valuation process
for each $d\p\in D$ yields $S^y_d ((\forall x C)\to S^x_yC])] = 
(\forall x S^y_d C]) \to S^y_dS^x_y C]] = (\forall x (S^y_d C)) \to S^x_dS^y_d 
C]] = (\forall x B) \to S^x_d B,$ where $B = S^y_d C].$ Now simply assume that
$\str\mod \forall xB.$ Then for all $d\p \in D,\ \str\mod S^x_dB].$  
From definition 3.3.3$,$ $d\p \in D,$ $\str \mod (\forall x B)\to S^x_dB].$ Since $d\p \in D$ is 
arbitrary$,$ we have as this point in the valuation $\str \mod \forall y((\forall x A) \to S^x_yC]).$ Multiple applications of metalogic generalization as the d's associated with the other free variables vary over $D$ 
completes the proof of part (iii) since $\str$ is also arbitrary.\pars
(iv) The same proof from (iii) for the case that $x$ is not free in $A$ holds for this 
case. Again in the same manner as in the proof of part (iii)$,$ we need only 
assume that $\str$ is a structure for $Pd$ and for the formula  
$S^x_y C] \to (\exists x C),$  
 $x$ or $y$ are the only free variables in $C$ and $y$ is free for $x.$
Consider $\forall y ((S^x_y C] \to (\exists x C)).$
Then for arbitrary $c\p \in D,$    
$S^y_c(S^x_y C] \to (\exists x C)) = S^y_cS^x_y C]]\to S^y_c 
\exists xC]= 
S^x_cS^y_c C]] \to \exists x (S^y_cC])= S^x_c B \to \exists xB,$ where 
$B = S^y_c C].$ 
Now if $\str\mod \exists xB.$  Then there exists some $d\p \in D$ such that 
$\str\mod S^x_d B.$ Hence$,$ letting $c\p =d\p,$ we have 
$\str \mod S^x_c B \to \exists x B$. On the other hand$,$
if $\str \not\mod \exists xB,$ then for all $d\p \in D,$ $\str \not\mod
S^x_d B.$ This implies that $\str \not\mod S^x_c B.$ From definition 3.3.3$,$ 
 $\str\mod  S^x_c B] \to (\exists x 
B).$ Since $c$ is arbitrary$,$ then $\str \mod \forall y ((S^x_y C] \to (\exists x C)).$ Again by multiple applications of generalization and since 
$\str $ is an arbitrary structure$,$ the result follows. \qed
There are many formulas in $Pd$ that are not instances of valid propositional 
formula that may be of interest to the pure logician. It's not the purpose of 
this text to determine the validity of a member of $Pd$ that will not be of 
significance in replicating within $Pd$ significant propositional 
metatheorems. However$,$ certain formula in $Pd$ are useful in simplifying 
ordinary everyday logical arguments. The next two metatheorems relate to both 
of these concerns. \parm
{\bf Theorem 3.4.8} {\sl If $x$ is any variable and $B$ does not contain a 
free occurrence of $x$$,$ then\pars
(i) Special process and notation $\str\mod A.$ \pars
(ii) For any structure $\str$ for $Pd,$$,$ let $A$ have free 
variables $x_1, \ldots, x_n$. Then $\str \mod \neg (\forall x_1,\ldots, 
(\forall x_nA)\ldots )$ 
if and only if $\str\mod (\exists x_1,\ldots,(\exists x_n(\neg A))\ldots )$.
 \pars  
(iii) If $\str\mod B \to A,$ then $\str\mod B \to (\forall x A).$\pars
(iv) If $\str\mod A \to B,$ then $\str\mod (\exists x A) \to B.$}\pars
Proof. (i) To determine whether or not $\str \mod A,$ where $A$ is not a 
sentence$,$ we consider whether or not $\str \mod \forall A.$ First$,$ corollary 
3.4.5.1 indicates that the order$,$ from left to right$,$ in which we make the 
required substitution has no significance upon the whether or not $\str \mod 
\forall A.$ Note that after we write the, possibly empty, sequence of statements ``for each $d_1',$ for each $d_2', \cdots,\in D$'' {\bf the universal closure substitution operators $S_{d_1}^{x_1}, S_{d_2}^{x_2}, \cdots$ distribute over all of the fundamental connections $\lor,\ \land,\ \to,\ \iff,\ \neg.$} [Note that one must carefully consider the statements ``for each $d_1',$ for each $d_2', \cdots,\in D$.''] What happens is that when we have a variable that is not free in a subformula then the substitution process simply does not apply. Further if the statements ``for each $d_1',$ for each $d_2', \cdots,\in D$'' still apply to the entire formula and substitution operators have not been eliminated, then we can go from basic subformula that contain a universal closure substitution operators to the left of each subformula back to a universal closure for the entire composite formula. 
After making these 
substitutions$,$ we would have$,$ depending upon the domain $D,$ a large set of 
objects that now carry the d's (or c's) in various places and that act like sentences. If these sentences satisfy the requirements of $\str\mod$, 
then by the metalogical axiom of generalization $\str\mod$ holds for the universally closed formula.  \par
For example$,$ consider the hypotheses of this theorem and the formulas $B \to A,$  and
$B \to (\forall xA).$  To establish the result in (iii)$,$
$\str\mod B\to A$ means $\str\mod \forall(B\to A).$ Hence, we have the,  possibly empty, sequence of statements for each $d_1\cdots,d \in D,\ \str\mod S^{x_1}_{d_1} \cdots S^x_d (B \to A)$ holds. The valuation can be rewritten as $S^{x_1}_{d_1} \cdots S^x_d B \to S^{x_1}_{d_1} \cdots S^x_dA \Rightarrow  S^{x_1}_{d_1} \cdots B \to S^x_d(S^{x_1}_{d_1} \cdots A)= B' \to S^x_d A'.$ For the valuation process $B'$ acts like a sentence and $A'$ acts like a formula with only one free variable, $x$. Thus, from our observation about the metalogical process for quantification over the members of $D$$,$ if by simply consider $B$ to be a sentence$,$ and $A$ to have$,$ at the most$,$ one free variable $x,$ it can be shown that $\str\mod B \to \forall x A$, then we have established, in general, that  $\str\mod B \to \forall x A.$  \par 
Thus, this type of argument shows that in many of our following arguments$,$
relative to structures$,$ we can reduce the valuation process  
to a minimum number of 
free variables that are present within a specific formula.\pars
{\bf From this point on$,$ unless otherwise stated$,$ we will ALWAYS assume 
that $\str$ means a structure for $Pd$ and $\str\mod A$ means $\str\mod \forall A$}\pars
(ii) If we have no free variables in $A$$,$ then we have nothing to prove.  
Assume that $A$ has only one free variable $x$. Suppose that 
$\str \mod (\neg (\forall x A)).$ This implies that 
$\str\not\mod \forall x A.$ This means that there exists some $d\p \in D$ such 
that $\str\not\mod S^x_d A].$ This implies that there exists some $d\p \in D$ 
such that $\str\mod (\neg S^x_d A]) = S^x_d (\neg A)].$ Hence $\str\mod (\exists 
x(\neg A)).$ Now apply induction. \pars
Note: As well be seen in the following proofs, other conclusions hold that are not expressed in (iii) or (iv). These restricted conclusions are presented since these are basic results even where the universal closure is not used. \pars
(iii) By the special process$,$ assume that $B$ is a sentence and that $A$ has at the most one free 
variable $x$. Let $\str\mod B \to A.$ If $A$ has no free variables$,$  
then $A$ and $B$ are sentences. The statement that for 
each $d\p \in D,$ $\str \mod B \to A$ and the statement for each 
$d\p \in D,\  \str \mod B \to S^x_dA] = B\to A$ are identical and the result 
holds in this case. \pars
Suppose that $x$ is the only free variable in $A$ and $B$ is a sentence. Then $\str \mod B 
\to A$ means that, for each $d' \in D,$ $\str\mod S^x_d(B \to A)] = B \to S^x_dA].$  This is but the 
valuation process for the formula $B \to (\forall x A).$ Hence $\str \mod
B \to (\forall xA).$ (Assume $\str\mod B.$ Then $\str\mod S^x_dA]$ for each $d'\in D$ implies $\str\mod \forall x A.$)\pars
(iv) Assume that $B$ is a sentence and that $A$ has$,$ at the most$,$ one free 
variable $x$. As was done in (iii)$,$ if $x$ is not free in $A,$ the result 
follows. Assume that $x$ is free in $A$ and that $\str\mod  A \to B.$ Then 
this means that, for each $d'\in D,$ $\str\mod S^x_d (A \to B)] = (S^x_d A]) 
\to B.$  Hence$,$ since $D \not= \emptyset,$ that
 there exists some $d\p \in D$$,$ such that $\str \mod
(S^x_d A]) \to B.$ Consequently$,$ by the special process$,$ 
$\str\mod (\exists xA) \to B.$ This complete the proof.\parm
[Note: Parts (iii) and (iv) above do not hold if $B$ contains $x$ as a free variable. For an example$,$ let $P(x) =A =B,\ D =\{a,b\},\ P' = \{a\}.$ Assume that $\str\mod P(x) \to P(x).$ Thus$,$  for $a,$ $\str\mod S^x_a P(x)\to P(x)] = 
P(a) \to P(a)$ and$,$ in like manner$,$ $\str\mod P(b) \to P(b).$ Since 
$\str\not\mod P(b)$$,$ then $\str\not\mod \forall xP(x).$ Hence$,$ for the case $d = a,$ we have that $\str\not\mod P(a) \to \forall x P(x) \Rightarrow
\str\not\mod \forall (P(x) \to \forall x P(x)) \Rightarrow \str\not\mod P(x) \to \forall x P(x).$] \pars   
An important aspect of logical communication lies in the ability to re-write 
expressions that contain quantifiers into logically equivalent forms. The next 
theorem yields most of the principles for quantifier manipulation that are 
found in ordinary communication.\parm
{\bf Theorem 3.4.9} {\sl For formulas $A,\ B, C$ the following 
are all valid formulas$,$ where $C$ does not contain $x$ as a free 
variable.\pars
(i) $(\exists x(\exists y A(x,y)) \iff (\exists y(\exists x A(x,y)).$\pars 
(ii) $(\forall x(\forall y A(x,y)) \iff (\forall y(\forall x A(x,y)).$\pars 
(iii) $(\neg(\exists x A)) \iff (\forall x(\neg A))$.\pars
(iv) $(\neg(\forall x A)) \iff (\exists x(\neg A)).$\pars
(v) $(\exists x(A \lor B)) \iff ((\exists xA)\lor ((\exists xB)).$ \pars
(vi) $(\forall x(A\land B)) \iff ((\forall xA)\land (\forall xB)).$ \pars
(vii) $(\forall x(C \lor B))\iff (C \lor (\forall xB)).$\pars
(viii) $(\exists x(C \land B))\iff (C \land (\exists xB)).$}\pars
Proof. (i) This follows from Theorem 3.4.5 and the fact that the expression ``there exists some $d\p \in D$ and there exists some $c\p \in D$'' is metalogically equivalent to ``there exists some $c\p \in D$ and 
there exists some $d\p \in D$.''  
Consequently$,$ for arbitrary $\str\mod$, if $\str\mod (\exists x(\exists y A(x,y)))$, then $\str\mod(\exists y(\exists x 
A(x,y)))$ and conversely. \pars
(ii) This follows immediately by Corollary 3.4.5.1.\pars
(iii) (Special process.) This follows from the assumption that $x$ is the only possible free 
variable in $A$$,$ the propositional equivalent $\mod \neg(\neg A) \iff A,$ and 
theorem 3.4.8 part (ii). \pars
(iv) Same as in (iii).\pars
(v) (Special process.) We may assume that the only possible free variables in $A\lor B$ is the 
variable $x.$ Assume that for $\str,$ an arbitrary structure and 
$\str\mod (\exists x(A\lor B)).$ 
Then there exists some $d\p \in D$ such that 
$\str\mod S^x_d(A \lor B)] = (S^x_dA] \lor S^x_dB]).$ Note that it does not 
matter whether the variable $x$ is free or not in $A$ or $B$ since 
this still holds whether or not a substitution is made. Hence$,$
there exists 
some $d\p \in D$ such that $\str\mod S^x_dA]$ or there exists at the least 
the same $d\p \in D$ such that $\str\mod S^x_dB].$ Hence$,$
$\str\mod ((\exists xA) \lor (\exists xB)).$ Now assume that $\str\not\mod 
\exists x(A\lor B).$ Thus$,$ there does not exist any $d\p \in D$ such that 
$\str\mod S^x_d (A \lor B)].$ This means there does not exists any $d\p \in D$ 
such that $\str\mod (S^x_dA] \lor S^x_dB]).$ Therefore$,$ $\str\not\mod ((\exists 
xA) \lor (\exists xB)).$ This result now follows. \pars
(vi) Taking the proof of (v) and change the appropriate words and $\lor$ to 
$\land$ this proof follows.\pars
(vii) (Special process.) Assume that $C$ does not have $x$ as a free variable and that $B$ may 
contain $x$ as a free variable. Further$,$ it's assumed that there are no other 
possible free variables. $\forall x(C \lor B)$ is a sentence. Let $\str$ be an 
arbitrary structure. Assume that $\str\mod \forall x(C \lor B).$ 
Then for each $d\p \in D,\ \str\mod S^x_d(C \lor B)] = C\lor S^x_dB].$ Hence$,$
$\str\mod C$ or for each $d\p \in D,\ \str\mod S^x_d B].$  Hence$,$  
$\str\mod (C \lor (\forall x B)).$ Then in like manner$,$ since $x$ is not 
free in 
$C$$,$ $\str\mod (C \lor (\forall xB))\Rightarrow 
\str\mod \forall x(C \lor B).$ [Note that this argument 
fails if $C$ and $B$ both have $x$ as a free variable. Since if 
$\str\mod \forall x(C \lor B),$ then considering any $d\p \in D$ we have that 
$\str\mod S^x_d(C \lor B)]=S^x_dA]\lor S^x_dB] \Rightarrow \str\mod S^x_dC]$ or 
$\str\mod S^x_dB].$ But both or these statement need not hold for a specific 
$d\p \in D.$ Thus we cannot conclude that  $\str\mod \forall x C$ or $\str\mod \forall x 
B.$ ]\pars
(viii) Left as an exercised.\qed
Obviously$,$ theorems such as 3.4.9 would be very useful if the same type of 
substitution for valid formula with the $\iff$ in the middle holds for $Pd$ as 
it holds in $L$ or $L^\prime.$ You could simply substitute one for the other 
in various places. Well$,$ this is the case$,$ just by simple symbolic 
modifications of the proofs of the metatheorems 2.6.1$,$ 2.6.2$,$ 2.6.3$,$ and 
corollary 2.6.3.1. We list those results not already present 
as the following set of metatheorems for $Pd$ and for 
reference purposes. \parm
\hs
{\bf Definition 3.4.4} ($\equiv$ for $Pd.$) Let $A,\ B \in Pd.$ Then define
$A \equiv B$ if and only if $\mod A \iff B.$ [See note on page 138.]\pars
\hm
{\bf Theorem 3.4.10} {\sl The relation $\equiv$ is an equivalence 
relation.}\pars
Proof. See note on page 38.\parm
{\bf Theorem 3.4.11} {\sl If $A,\ B, C \in {Pd}$ and $A\equiv B,$ then 
$C_A \equiv C_B.$}\pars
Proof. Let $\mod \forall (A\iff B).$ It can be shown that this implies $\mod \forall (C_A \iff C_B)$. \pars
{\bf Corollary 3.4.11.1} {\sl If $A,\ B, C \in {Pd}$ and $A\equiv B,$ and $\mod C_{A},$ then
$\mod C_{B}.$}\parm
With respect to $A$ being congruent to $B$ recall that all the formula have 
the exact same form$,$ the exact same free variables in the exact same places$,$ 
and the bounded variables can take any variable name as long as the formula retain the same bound occurrence patterns. Hence, structure valuation would yield the same statement that either $\mod$ or $\not\mod,$ holds for $A$ and $B$. Also see note on page 38  \parm
{\bf Theorem 3.4.12} {\sl Suppose that $A,\ B \in Pd$ are congruent.
Then $A\equiv B.$}\parm
{\bf Theorem 3.4.13} {\sl Let $A \equiv B.$ 
Then $\str\mod A$ if and only if  $\str\mod B.$} \pars
Proof. From Definition 3.4.4, if $A \equiv B$, then $\mod \forall A \iff \forall B$ if and only if $\str\mod A$ implies $\str\mod B$ and $\str\mod B$ implies $\str\mod A$. \parm
\hs\hm
\centerline{\bf EXERCISES 3.4}\parm
\noindent 1. A formula $A$ is a said to be $n$-valid$,$ where $n$ 
is a natural number
greater than 0$,$ if for {\it any} structure $\str$$,$ with the 
domain containing $n$ and only $n$ elements$,$ $\str \mod A.$ \parm
\noindent (a) Prove that $A = (\forall x(\exists y P(x,y))) \to (\exists y(\forall x 
P(x,y)))$ is 1-valid.\pars
\noindent (b) A {{\it countermodel}} $\str$ must be used to show that a formula $A$ is not 
valid. You must define a structure $\str$ such that $\str \not\mod A.$ Show by 
countermodel  that the formula $A$ in (a) is not 2-valid.\parm
\noindent 2. For each of the following$,$ determine whether the indicated 
variable $\lambda$ is free for $x$ in the given formula $A.$\parm
\noindent (a) $A = \forall w(P(x) \lor (\forall xP(x,y))\lor P(w,x));\ \lambda = 
y.$\pars
\noindent (b) $A$ is the same as in (a) but $\lambda = w.$\pars
\noindent (c) $A = (\forall x(P(x) \lor (\forall yP(x,y)))) \lor P(y,x);\ \lambda  = x.$\pars
\noindent (d) $A$ is the same as in (c) and $\lambda = y.$\pars
\noindent (e) $A = (\forall x(\exists y P(x,y))\to (\exists yP(y,y,));\ \lambda = y.$\pars
\noindent (f) $A= (\exists zP(x,z))\to (\exists zP(y,z));\ \lambda = z.$\parm
\noindent 3. For the formula in question 2$,$ write $S^x_\lambda A]$ whenever$,$ 
as given in 
each problem$,$ $\lambda$ is free for $x$. \parm
\noindent 4. Give a metaproof for part (viii) of theorem 3.4.9.\pars
\noindent 5. Determine whether the formula are valid.\parm
\+(a) $Q(x) \to (\forall x P(x)).$ &(d) $(\exists x(\exists y P(x,y)))\to 
(\exists xP(x,x)).$\cr\s
\+(b) $(\exists xP(x))\to P(x).$&(e) $(\exists xQ(x)) \to (\forall 
xQ(x)).$\cr\s
\+(c) $(\forall x(P(x) \land Q(x)))\to ((\forall xP(x))\land (\forall 
xQ(x))).$&\cr\parm
\noindent 6. A formula $A$ is said to be in {{\it prenex normal form}} if 
$A = Q_1x_1(Q_2x_2(\cdots (Q_nx_n(A))\cdots)),$ where $Q_i,\ 1 \leq i \leq n$ is 
one of the symbols $\forall$ or $\exists.$ The following is a very important 
procedure. Use the theorems in this or previous sections relative to the
language $L$ or $Pd$ to re-write each of the following formula in a prenex 
normal form that is equivalent to the original formula. [Hint. You may need 
to use the congruency concept and change variable names. For example$,$ 
$\forall xP(x) \to \forall xQ(x) \equiv \forall xP(x) \to \forall y P(y)\equiv
\forall y(\forall xP(x) \to Q(y))\equiv \forall y(Q(y) \lor (\neg (\forall 
xP(x))))\equiv \forall y(Q(y) \lor (\exists x(\neg P(x))))\equiv \forall 
y(\exists x(Q(y) \lor (\neg P(x)))).$]\parm
\noindent (a) $(\neg (\exists x P(x))) \lor (\forall xQ(x)).$\pars
\noindent (b) $((\neg (\exists x P(x)))\lor (\forall xQ(x)))\land (S(c) \to 
(\forall xR(x))).$\pars
\noindent (c) $\neg (((\neg (\exists xP(x)))\lor(\forall xQ(x))) \land 
(\forall xR(x))).$\parm
\centerline{NOTE}\pars
It is not the case that $\mod \forall (A \iff B) \iff ((\forall A) \iff (\forall B)).$ Let $D = \{a',b'\},\ P'= \{a'\},\ Q'=\{b'\}.$ Then $\str\not\mod \forall x\forall y(P(x) \iff Q(y))$ since $\str\not\mod P(a)\iff Q(a)$. However, since $\str\not\mod P(b)$ implies that $\str\not\mod \forall x P(x)$ and $\str\not\mod Q(a)$ implies that $\str\not\mod \forall y Q(y),$ then $\str\mod (\forall x P(x)) \iff (\forall y Q(y)).$ \pars
Using material yet to come, a way to establish Theorems 3.4.10 and 3.4.11 is to use the soundedness and completeness theorem, the usual reduced language $Pd'$, and 14.7, 14.9 that appear in J.W. Robbin, Mathematical Logic a first course,W. A. Benjamin, Inc NY (1969) p. 48. I note that these results are established in Robbin by use of formal axioms and methods that are identical with the ones presented in the next section. [The fact that the definition of $\equiv$ via the universal closure is an equivalence relation follows from the fact that for any formula $C,\ \vd (\forall C) \to C$ and $\vd \forall(A \iff B)$ if and only if $\vd A \iff B.$ Also $\vd \forall (A \iff B) \to (\forall A) \iff (\forall B)$ implies that if $\mod \forall (A \iff B),$ then $\mod (\forall A) \iff (\forall B).$ From $\vd \forall (A \iff B) \to (C_A \iff C_B),$ if $\mod \forall (A \iff B),$ then $\vd  (C_A \iff C_B)$ implies   
$\vd \forall (C_A \iff C_B)$ implies $\mod \forall (C_A \iff C_B)$ implies $C_A\equiv C_B.$ Using $\vd\forall (C_A \iff C_B) \to ((\forall C_A)\iff (\forall C_B)),$ we have that if $A \equiv B,$ and $\str\mod C_A,$ then $\str\mod C_B.$ The fact that for congruent $A$ and $B$ that $A \equiv B$ also follows from $\vd A \iff B.$ ]\pars
\ss
\noindent{\bf 3.5 Valid Consequences and Models}\parm
As in the case of validity$,$ the (semantical) definition for the concepts of a 
valid consequence and satisfaction are almost identical to those used for the 
language $L.$ \parm
\hs
{\bf Definition 3.5.1} (Valid consequence for $Pd.$) A sentence $B$ in $Pd$ (
i.e. $B \in \c S$) is a {{\it valid consequence}} of a set of premises $\Gamma \subset\c S$$,$ which 
may be an empty set$,$
 if for any $\str$ for  $\Gamma \cup \{B\}$ 
whenever $\str\mod A$ for each $A \in \Gamma$$,$ 
then $\str \mod B.$ This can be most easily remembered by 
using the following notation. Let $\Gamma \subset \c S.$ Let $\str \mod 
\Gamma$ mean that $\str \mod A$ for each $A  \in \Gamma.$ Then $B$ is a valid 
consequence from $\Gamma$ if whenever $\str\mod \Gamma,$ then $\str 
\mod B,$ it being understood that $\str$ is a structure for the set 
$\Gamma \cup \{B\}.$ The notation used for valid consequence is $\Gamma \mod 
B.$ \pars
\hm
The concept of  ``satisfaction'' also involves models and we could actually do 
without this additional term.\parm
\hs
{\bf Definition 3.5.2} (Satisfaction for $Pd$) A $\Gamma \subset \c S$ 
is {{\it satisfiable}} if there exists a $\str$ such that $\str\mod 
\Gamma.$\pars
$\Gamma$ is not satisfiable if no such structure exists.\pars
\hm
As in the previous section$,$ all of the results of section 2.8$,$ 2.9 and 2.10 
that hold for $L$ also hold for this concept extended to $Pd.$ The following 
metatheorem contains processes that are not used in the metaproof of its 
corresponding $L$ language metatheorem. The remaining metatheorems follow in 
a manner very similar to their counterparts in sections 2.8 --- 2.10.
\parm
{\bf Theorem 3.5.1} {\sl (Substitution into valid consequences.) Let 
$A_n,\ C,\ B \in \c S.$\pars
(i) If $A_n \equiv C,$ and $A_1,\cdots,A_n, \cdots \mod B,$ 
then $A_1,\ldots,A_{n-1},C,\cdots \mod B.$\pars
(ii) If $B \equiv C,$ and $A_1,\cdots,A_n,\cdots \mod B,$ then  $A_1,\cdots,A_n,\cdots\mod 
C.$}\pars
Proof. (i) Suppose that $A_n \equiv C.$ Let $\str$ be a structure. The $\str$ 
is defined for $C, B, A_1,\ldots,A_{n-1},A_n,\cdots.$  
Suppose that $\str \mod A_i$ for $1 
\leq i$. Then $\str \mod B.$ But $\str\mod A_n$ if and only if  
$\str \mod C.$  Consequently$,$ if
$\str \mod \{A_1,\ldots,A_{n-1},C,\cdots\},$ then $\str \mod B.$ [Notice that these 
statements are conditional. If for any structure $\str$ you get for any 
sentence $A$ in the set of premises that $\str\not\mod A$ you can simply disregard 
the structure.] \pars
(ii) Suppose that $\str$ is a structure and that $\str$ is defined for
 $B, C, A_1,\ldots,A_n,\cdots.$   
If $\str \mod A_i$ for $1\leq i,$ then $\str\mod B.$ But $\str 
\mod C.$ From this the 
result follows.\qed
{\bf Theorem 3.5.2} {\sl (Deduction theorem) Let $A,\ B,\ A_i\in \c S$ for $1 
\leq i \leq n.$ \pars
(i) $A\mod B$ if and only if $\mod A \to B.$\pars
(ii) $A,\ldots, A_n \mod B$ if and only if $A_1\land \cdots\land A_n \mod 
B.$\pars                                                                
(iii) $A,\ldots, A_n \mod B$ if and only if $\mod (A_1\land \cdots\land A_n) \to 
B.$\pars 
(iv) $A,\ldots, A_n \mod B$ if and only if $\mod (A_1\to \cdots \to (A_n \to 
B) \cdots ).$}\parm                                                               
As with the language $L$$,$ consistency is of major importance. It is defined by 
the model concept. The metatheorems that follow the next definition are 
established  
the same manner as their counterparts for $L.$ \parm
\hs 
{\bf Definition 3.5.3} (Consistency) A nonempty  set of premises $\Gamma 
\subset \c S$ is {{\it consistent}} if there does not exist a $B\in \c S$ 
such that $\Gamma \mod B\land (\neg B).$ \pars
\hm
{\bf Theorem 3.5.4} {\sl If $B \in \c S,$ then $\not\mod B \land (\neg 
B).$}\parm
{\bf Theorem 3.5.5} {\sl A nonempty  finite set of premises $\Gamma \subset \c S$ is 
inconsistent if and only if $\Gamma \mod B$ for every $B\in  \c S.$} (Based on the proof method for 2.10.3.) \parm
{\bf Corollary 3.5.5.1} {\sl A nonempty finite set of premises $\Gamma \subset\c S$ 
is consistent if and only if there exists some $B \in \c S$ such that 
$\Gamma \not\mod B.$}\parm
{\bf Corollary 3.5.6.1} {\sl A nonempty finite set of premises $\Gamma \subset 
\c S$ is consistent if and only if it is satisfiable.}\parm

Even though we seem to have  strong results that can be used to determine 
whether a sentence is a valid consequence or that such a set is consistent$,$ it 
turns out that it's often very difficult to make such judgments for the 
language $Pd.$ The reasons for this vary in complexity. One basic reason is 
that some sets of premises tend to imply that the domain $D$ is not finite. 
Further$,$ when we made such determinations for $L,$ we have a specifically 
definable process that can be followed$,$ the truth-table method. It can 
actually be shown$,$ that there is no known method to describe one fixed process
that will enable us to determine whether a sentence is valid or whether
a set of premises is consistent. Thus we must rely upon ingenuity to establish 
by models these concepts. Even then some mathematicians do not accept such
metaproofs as informally correct since some claim that the model chosen is  
defined in an unacceptable manner. The next examples show how we must rely 
upon the previous definitions and metatheorems to achieve our goals of 
establishing valid consequences or consistency. In some cases$,$ however$,$ no 
amount of informal argument will establish the case one way or the other with 
complete assurance. I point out that these formal sentences are translations 
from what ordinary English language sentences.\parm
\hrule\smallskip
If you are given a finite set of premises $A_1,\cdots, A_n$ then there are 
various ways to show that $A_1,\cdots, A_n\mod B.$ One method is to use that 
Deduction Theorem and a general argument that $\mod A_1 \land \cdots\land A_n \to 
B$; the propositional method of assuming that for a structure $\str$ in 
general if $\not\mod B,$ then $\not\mod A_1$ or $\not\mod A_2.$ For the 
invalid consequence concept$,$ you have two choices$,$ usually. Consider 
ANY structure $\str$ such that $\str \not\mod B$ and show that this leads to 
each $\str \mod A_i,\ i= 1,\ldots,n.$ However$,$  
to show that 
$A_1,\cdots, A_n\not\mod B,$ it is often easier to define a structure 
$\str$ and show that $\str \mod A_i,\ i =1,\ldots, n,$ but that $\str \not\mod 
B$ (i.e. the definition for valid consequence does not hold.) Note that if there is no $\str$ such that $\str \not\mod B$, then $\mod B.$\pars
\hrule\medskip
{\bf Example 3.5.1.} Consider the premises $A_1 = \forall x(P(x) \to (\neg 
Q(x))),\ A_2 = \forall x(W(x) \to P(x)),\ B = \forall x(W(x) \to (\neg 
Q(x))).$ We now (attempt) to determine whether $\mod A_1\land A_2 \to B.$\pars
Let $\str$ be any structure. 
 First$,$ let $D$ be any nonempty domain and the $P',\ Q',\ W'$ be any 
subsets of $D$ (not including the possible empty ones).  
Suppose that $\str {\not\mod}  \forall x(W(x) \to (\neg 
Q(x))).$ Then there exists some $c\p \in D$ such that $\str {\not\mod} S^x_c(W(x) 
\to (\neg Q(x)) = W(c) \to (\neg Q(c)).$ From our definition of what 
$\str \not\mod$ means$,$ (*) $\str \mod W(c)$ (*) and that  
$\str \not\mod \neg Q(c)$; which means that  $\str\mod Q(c).$ 
This does not force any of the premises to have a specific $\mod$ or 
$\not\mod$. Thus as done for the propositional calculus$,$ let  
$\str \mod A_2 = \forall x(W(x) \to (\neg Q(x)).$ Hence$,$ for all $d' \in D$$,$ we 
have that $\str \mod W(d) \to P(d).$ Now $c$ is one of the $d$s. Hence$,$ we 
have that $\str \mod W(c) \to P(c).$ What does this do to $A_1$? We show that 
$\str \not\mod A_1.$ Assume that $\str\mod A_1.$ Hence$,$ for all $d'\in D,$
$\str \mod P(d) \to (\neg Q(d)).$ Again this would give $\str\mod P(c) \to (\neg 
Q(c)).$  We know that $\str \not\mod \neg Q(c).$  Hence $\str\not\mod P(c).$
Thus $\str {\not\mod} W(c).$ This contradicts the (*) expression. Hence$,$
$\str {\not\mod} A_1.$ Now consider an empty $W'$. Then $\str \mod \forall x(W(x) \to (\neg Q(x))).$ Consider $Q'$ empty and $W'$ not empty. Then $\str \mod \forall x(W(x) \to (\neg Q(x))).$ Hence, $B$ is a valid consequence
of $A_1,A_2.$\parm
{\bf Example 3.5.2.} Let $A_1 = \exists x(P(x) \to Q(x)),\ A_2 = \forall 
x(W(x) \to P(x)),\ B= \exists x(W(x) \to (\neg Q(x))).$ Is $B$ a valid 
consequence of the premises? Well$,$ we guess$,$ that this might be an invalid 
consequence. We try to find a structure $\str$ such that
$\str  \mod A_1,\ \str \mod A_2$ and $\str\not\mod B.$ Let $D=\{a\p\}, 
\ P\p=Q\p=W\p=D.$ There are no constants in the sentences $A_1,\ A_2$ or $B.$ $N$ contains only one constant $a$ which is interpreted as $a'$.
Clearly$,$ $\str \mod A_1$ and $\str \mod A_2.$ But$,$
$\str\mod W(a)$ and $\str\not\mod (\neg Q(a)).$ Consequently$,$ $\str\not\mod 
B.$ [Note the difficulty would be to first have a ``feeling'' that the 
argument from which these formula are taken is not logically correct 
and then construct an acceptable structure that 
establishes this ``feeling.'' Not an easy thing to do$,$ if it can be done at 
all.] \parm
{\bf Example 3.5.3.} Let $A = \exists y(\forall x P(x,y))$ and $B = \forall 
x(\exists y P(x,y).$ [In this case$,$ $P\p$ cannot be the empty relation$,$ for if 
this were the case$,$ then $\str \not\mod A$ for an appropriate structure.] We 
try to show that $A\mod B$ and get nowhere. So$,$ maybe the 
$A\not\mod B.$ Let's see if we can obtain a countermodel.  
Let $D =\{a\p,b\p\}.$ Let $P\p 
=\{(a\p,b\p),(b\p,a\p)\}.$ Now extend this construction to a structure $\str$ 
for $Pd.$ 
Then for each $d\p \in D$ there exists some $c\p 
\in D$ such that $(c\p,d\p) \in P\p,$ which implies that $\str \mod A.$ 
On the other hand$,$ there does not exists 
some $d\p \in D$ such that for each $c\p \in D,\ (c\p,d\p) \in P\p.$  Hence$,$ 
$\str \not\mod B$ for this structure. Thus$,$ $B$  is 
not a valid consequence of $A.$ However$,$ it can easily be shown that $A$ is 
a valid consequence of $B.$  \parm
{\bf Example 3.5.4.} Determine whether $A_1 = \forall x(P(x) \to (\neg Q(x))),\ 
A_2 = \forall x(W(x) \to P(x))$ is a consistent set of premises. This is not 
very difficult since we notice that if we construct a structure such that for 
each $d\p \in D$ $\str\not\mod P(d)$ and $\str\not\mod W(d),$ then this 
structure would yield that the set is consistent. Simply let $P\p$ and $W\p$ 
be the empty set (and extended the construction to all of $Pd$) and 
the conditions are met. Thus the set of premises is consistent.\parm
{\bf Example 3.5.5.} Add the premise $A_3= \forall x(\neg P(x))$ to the premises 
$A_1$ and $A_2$ in example 3.5.3. The same structure yields $\str\mod A_3.$ 
Thus this new set of premises is consistent.\parm
{\bf Example 3.5.6.} Determine whether 
$A_1 = \forall x(P(x) \land (\neg Q(x)),\ 
A_2 = \forall x(W(x) \land P(x))$ form a consistent set of premises. Clearly$,$ 
using only empty sets for the relations will not do the job. Let $D= \{d\p\}.$  
Now let $P\p = W\p = D$ and $Q\p$ 
is the empty set. 
Then for all $d\p \in  D$$,$ $\str \mod P(d),\ \str\mod W(d),\ 
\str \mod (\neg Q(d)).$ \parm
To show inconsistency$,$ constructing models will not establish anything since we 
must show things don't work for {\it any} structure. The argument must be a 
general argument. \parm
{\bf Example 3.5.7.} Let $A_1 = \exists x(\neg P(x)),\ A_2 = \forall x(Q(x) 
\to P(x))$ and $A_3 = \forall x Q(x).$ Let $\str$ be any structure such that
$\str \mod \forall xQ(x).$ Then whatever set $Q\p \subset D$ you select$,$ 
it must have the property that for all $d\p \in D,\ d\p \in Q\p.$ Hence$,$
$Q\p = D.$ Now if $\str \mod A_2,$ then for all $d\p \in D,\ \str \mod Q(p) \to 
P(d).$ Thus for all $d\p \in D,\ \str\mod P(d)$ implies that $d\p \in P\p.$ 
Thus $P\p = D.$ It follows that there does not exist a member $d\p$ of the set 
$P\p$ such that $d\p$ is not a member of $P\p.$ Thus $\str\not\mod A_1.$ 
Hence$,$ the set of premises is inconsistent. \parm
It might appear that it's rather easy to show that sets of premises are or 
are not consistent$,$ or that sentences are valid consequences of sets of 
premises. But$,$ there are sets of premises such as the 14 premises discovered 
by Raphael Robinson for which it can be shown that there is no known  way
to construct a model for these premises without assuming that the model is 
constructed by means that are either equivalent to the premises themselves or
by means that assume a set of premises from which Robinson's premises can be 
deduced. Such an obvious circular approach would not be accepted as an 
argument for the consistency of the Robinson premises. The 
difficulty in determining valid consequences or consistency occurs when 
nonempty $n$-place$,$ $n>1,$  and nonfinite  domains are needed to satisfy some 
of the premises. The only time you are sure that your model will be accepted by 
the mathematical community for consistency argument is when your  model uses a 
finite domain. Robinson's system$,$ if it has a model$,$ must have a domain that 
is nonfinite. Further$,$ the set of natural numbers satisfies the Robinson 
axioms. \pars 
Due to the above mentioned difficulties$,$ it's often necessary to consider a 
weaker form of the consistency notion -- a concept we term {{\it relative 
consistency}}. If you have empirical evidence that a set of premises $\Gamma$ 
is consistent$,$ then they can be used to obtain additional premises that we 
know are also consistent. After we show the equivalence of formal proof theory 
to model theory$,$ then these premises can be used to construct models.\parm
{\bf Theorem 3.5.7} {\sl If a nonempty  set of premises $\Gamma$ is 
consistent and $\Gamma \mod B_i,$ where each $B_i$
is  
sentence$,$ then the combined set of premises $\Gamma \cup \{B_i\}$ is 
consistent.}\pars
Proof. From consistency$,$ there is a structure $\str$ such that 
for each $A \in \Gamma,$  $\str \mod A$. But for $\str$ we also have that 
$\str \mod B_i$ for each $i.$ 
Consequently$,$ the structure $\str$ yields the requirements for 
the 
set of 
premises $\Gamma \cup \{B_i\}$ to be consistent. \par\bigskip

\hs\hm
\centerline{EXERCISES 3.5}
\noindent 1. The following are translations from what we are told are valid
English language argument. Use the model theory approach and determine if$,$ 
indeed$,$ the sentence $B$ is a valid consequence of the premises.\parm
\noindent (a) $A_1 = \forall x(Q(x) \to R(x)),\ A_2 = \exists x Q(x) \mod
B = \exists x R(x).$\pars
\noindent (b) $A_1 = \forall x(Q(x) \to R(x)),\ A_2 = \exists x (Q(x) \land 
Z(x)) \mod B = \exists x(R(x) \land Z(x)).$\pars
\noindent (c) $A_1 =\forall x(P(x) \to (\neg Q(x))),\  A_2 =\exists x (Q(x) 
\land R(x)) \mod B = \exists x ( R(x) \land (\neg Q(x))).$\pars
\noindent (d) $A_1 = \forall x (P(x) \to Q(x)),\ A_2 = \exists x(Q(x) \land 
 R(x)) \mod B = \exists x (R(x) \land (\neg Q(x))).$\pars
\ss
\noindent {\bf 3.6 Formal Proof Theory}\parm
I hope to make the formal proof theory for this language as simple as possible. 
Again we will only do formal proofs or demonstrations that will lead directly 
to our major conclusions. {\bf The language will be $Pd\p\subset Pd,$ where no
existential quantifiers appear.} Later we will simply use the abbreviation
$\exists x A$ for the expression $\neg (\forall x (\neg A)).$ This will allow 
us to extend all the proof theoretic notions for $Pd\p$ to $Pd.$ 
Also recall that 
the language $Pd$ can actually be considered as constructed with the 
connectives $\land,\lor,\ \iff$ since these were abbreviations for 
the equivalent propositional equivalent formula. (i) $A \lor B$ is an 
abbreviation for $(\neg A) \to B$. (ii) $A \land B$ is an abbreviation for
$(\neg (A\to (\neg B)))$. (iii) $A \iff B$ is an abbreviation for $(A\to B) 
\land (B\to A).$ From the semantics$,$ these formula have the same pattern for 
$\mod$ as they have for truth-values if it is assumed that $A,\ B \in L$ (say 
$A,\ B \in \c S$.)\parm
\hs
{\bf Definition 3.6.1} The language $Pd\p$ is that portion of $Pd$ in which no 
formula contains an existential quantifier.\pars
\hm
{\bf Definition 3.6.2} (A formal proof of a theorem) \parm
(1) Use the entire process as described in definition 2.11.2 for the language
$L\p$ with the following \underbar{additions} to the axioms $P_1,\ P_2,\ P_3$ 
and to the rule modus ponens.\parm
(2) For predicate variables $A,\ B$ (i.e. metasymbols for any formula) you may write down as any step the 
formula\parm
\indent\indent $P_4:\ \ \forall x(A \to B)\to (A \to (\forall xB))$ whenever
$x \in \c V$ AND $x$ is NOT free in $A$.\parm
\indent\indent $P_5:\ \ (\forall xA) \to S^x_\lambda A]$ whenever $\lambda \in 
\c V$ AND  $\lambda$ is free for $x$ 
OR $\lambda \in \c C$ (the set of constants$,$ if any.)\parm

(3) You add one more rule of inference$,$ the rule called {{\it generalization.}}
\indent\indent $G(i).$ Taking any previous step $B_i,$ you may write down as a 
new step with a larger step number the step $\forall xB_i,$ for any $x \in \c  
V.$
\pars
(4) If you follow the above directions and the last step in your formal proof 
is $A$ then we call $A$ a {{\it theorem}} for $Pd\p$ and denote this
(when no confusing will result) by $\vd A.$
\hm 
The additional steps that can be inserted when we construct a formal 
demonstration from a set of premises is exactly the same as in definition 
2.12.1.\parm
\hs 
{\bf Definition 3.6.3} Let $\Gamma$ be any subset of $Pd\p.$ Then a {{\it 
formal demonstration}} from $\Gamma$ follows the exact same rules as in 
definition 3.6.2 with the additional rule that we are allowed to insert as a 
step at any point in the demonstration a member of $\Gamma.$ If the last 
formula in the demonstration is $A,$ then  $A$ is a deduction from $\Gamma$ 
and we denote this$,$ when there is no confusion$,$ by the symbol $\Gamma \vd A.$ 
\pars
\hm
Of course$,$ a formal proof of a theorem or a demonstration for a formula 
contains only a finite number of steps. Do we already have a large 
list of formal proofs for formal theorems that can be used 
within the proof theory for $Pd\p$? \parm
{\bf Theorem 3.6.1} {\sl Consider any formal proof in $L\p$ for the $L\p$ 
theorem $A$. For each specific propositional variable within $A,$ substitute a 
fixed 
predicate variable at each occurrence of the propositional variable in $A$ 
to obtain the predict formula $\hat A.$ 
Then for $Pd\p$$,$ 
we have that $\vd \hat A.$}\pars
Proof. Simply take each step in the formal proof in $L\p$ that $\vd A$ and 
make the same corresponding consistent variable predicate substitutions. Since each 
axiom $P_1,\ P_2,\ P_3$ and the MP rule are the same in $Pd\p,$ this will give 
a formal proof in $Pd\p$ for $\vd \hat A.$\qed
{\bf Theorem 3.6.2} {\sl Let $\Gamma \subset Pd\p$ and $A, \ B\in Pd\p.$\pars
(a) If $A \in \Gamma$ or $A$ is an instance of an axiom $P_i,$ $1\leq i\leq 
5,$  then $\Gamma\vd A.$\pars
(b) If $\Gamma \vd A$ and $\vd A \to B,$ then $\Gamma \vd B.$\pars
(c) If $\vd A,$ then $\Gamma \vd A.$\pars
(d) If $\Gamma$ is empty and $\Gamma \vd A,$ then $\vd A.$\pars
(e) If $\Gamma \vd A$ and $D$ is any set of formula$,$ then $\Gamma \cup D \vd 
A.$\pars
(f) If $\Gamma \vd A,$ then there exists a finite $D \subset \Gamma$ such that 
$D \vd A.$}\pars
Proof. The same proof for theorem 2.12.1. \parm
We now give an example of a formula that is a member of $Pd\p$ and is not 
obtained by the method that yields $\hat A$ and is a formal theorem.\parm
{\bf Example 3.6.1} For any 1-place predicate $P(x) \in Pd_0\p,$ $x,\ y \in \c 
V,$ 
$\vd (\forall x P(x)) \to (\forall y P(y)).$\parm
\line{(1) $(\forall xP(x))\to P(y)$\hfil $P_5$}\s
\line{(2) $\forall y((\forall x P(x))\to P(y))$\hfil $G(1)$}\s
\line{(3) $\forall y((\forall xP(x))\to P(y)) \to ((\forall xP(x))\to 
(\forall yP(y))$\hfil $P_4$}\s
\line{(4) $(\forall x P(x)) \to (\forall y P(y))$ \hfil $M(1,3)$}\parm
 It should be noted$,$ that if $x$ is not 
free in $A$$,$ then $P_5$ yields $(\forall x A) \to S^x_yA]=(\forall x A) \to A,$
for any $\lambda \in {\c V}\cup {\c C}$$,$  
since the second $A$ has 
had no changes made. Further$,$ if $x$ is free in $A,$ then $P_5$ yields $(\forall xA) 
\to A.$ \parm    
{\bf Example 3.6.2} (Demonstration)\par 
\centerline{$A \to (\forall xB)\vd \forall x(A\to B)$}\s
\line{(1) $A \to (\forall xB)$\hfil Premise}\s
\line{(2) $(\forall x B) \to S^x_xB]=(\forall x B) \to B$\hfil $P_5$}\s
\line{(3) $A \to B$\hfil HS(1$,$2)}\s
\line{(4) $\forall x(A \to B)$\hfil G(3)}\parm
The rule that we have given for Generalization is not the only rule of this 
type that is used throughout the literature. In particular$,$ some authors 
define deduction from a set of premises in a slightly different manner than we 
have defined it. The reason for this depends upon how nearly we want the 
deduction theorem for $Pd\p$ to mimic the deduction theorem for $L\p,$ among 
other things.  Also$,$ note that the statement that y is 
free for x$,$ holds when $x$ does not occur in $A$ or does not occur free in $A$.
The following two results are useful 
in the sequel.   
\parm 
{\bf Theorem 3.6.3} {\sl If $y$ is free for $x$ in $A\in Pd\p$$,$ then 
$\forall xA \vd\forall yS^x_yA].$}\pars
\line{(1) $\forall x A$\hfil Premise}\s
\line{(2) $(\forall x A) \to S^x_yA]$\hfil $P_5$}\s
\line{(3) $ S^x_yA]$\hfil MP(1$,$2)}\s
\line{(4) $ \forall yS^x_yA]$\hfil G(3)}\parm
{\bf Corollary 3.6.3.1} {\sl If $y$ does not appear in $A$$,$ then $
\forall yS^x_yA]\vd \forall xA.$}\pars
Proof. This follows since$,$ in this case$,$ $y$ is free for $x$ for $y$ cannot be 
bound by any quantifier and $x$ is free for $y$ for $S^x_y A$.
 Further$,$ $A = S^y_xS^x_yA]].$ We now apply theorem 3.6.3. Hence$,$  $\forall y S^x_y A 
\vd \forall x S^y_xS^x_yA]]= \forall xA.$\parm \vfil\eject 
{\bf Theorem 3.6.4} {\sl Let $\{x_1,x_2,\ldots, x_n\}$ be any $n$ 
free variables in 
$A\in Pd\p.$ Then $\Gamma \vd A$ if and only if $\Gamma \vd \forall x_1(\forall 
x_2(\cdots(\forall x_n A)\cdots)).$}\parm
Proof. Let $n = 0.$ The result is clear.\pars
Let $n = 1.$ If $\Gamma \vd A$ has $x_1$ as free$,$ then by application of 
generalization on $x_1$ we have that $\Gamma \vd \forall x_1 A.$\pars
Suppose that result holds for $n$ free variables and suppose that 
$x_1,\ldots, x_{n+1}$ are free variables in $A.$ Then $x_1,\ldots,x_n$ are 
free variables in $\forall x_{n+1}A.$ Now if $\Gamma \vd A,$ then $\Gamma 
\forall x_{n+1}A$ 
by generalization. Hence$,$ from  the induction hypothesis$,$ $\Gamma \vd\forall x_1(\forall 
x_2(\cdots(\forall x_{n+1} A)\cdots)).$ Hence$,$ by induction$,$ for any $n$ free 
variables in $A$ the result holds.\pars
The converse follows by a like induction proof through application of 
$P_5$ and one MP. \qed
{\bf Corollary 3.6.3.1} {\sl For any $A \in Pd\p,$ $\Gamma \vd A$ if and only if  
$\Gamma \vd \forall A.$}\parm
{\bf Theorem 3.6.5} {\sl For any $A \in Pd\p,$ if $x\in \c V,$ 
then $\Gamma \vd A$ if and only if $\Gamma \vd \forall xA.$}\parm
Proof. Let $\Gamma \vd A$. Then one application of generalization yields 
$\Gamma \vd \forall x A.$\pars
Conversely$,$ suppose that $\Gamma\vd \forall xA.$ Then one application of 
$\forall xA \to S^x_xA = \forall x A \to A$ followed by one MP yields 
$A$.\parm
\hs\hm
\centerline{\bf EXERCISES 3.6}\parm
\noindent 1. Complete the following formal demonstrations. Note that whenever 
$P_5$ is applied$,$ we write $A$ for $S^x_xA$. 
\parm
\centerline{(A) $\forall x(A \to B),\ \forall x(\neg B)\vd \forall x(\neg A).$}\s
\line{(1) $\forall x(A \to B)$\hfil $\ldots\ldots$}\s
\line{(2) $\forall x(\neg B)$\hfil $\ldots\ldots$}\s 
\line{(3) $\ldots\ldots\ldots\ldots$\hfil $P_5$}\s
\line{(4) $A \to B$\hfil MP(\quad$,$\quad)}\s
\line{(5) $(A \to B) \to ((\neg B)\to (\neg A))$\hfil Exer. 2.13$,$ 2A.}\s
\line{(6) $\ldots\ldots\ldots\ldots$\hfil MP(\quad$,$\quad)}\s
\line{(7) $(\forall x(\neg B)) \to (\neg B)$\hfil $\ldots\ldots$}\s
\line{(8) $\neg B$\hfil $\ldots\ldots$}\s
\line{(9) $\neg A$\hfil $\ldots\ldots$}\s
\line{(10) $\ldots\ldots\ldots\ldots$\hfil $\ldots\ldots$}\medskip
\centerline{(B) $\forall x(\forall y A)\vd \forall y(\forall x A)$}\s
\line{(1) $\ldots\ldots\ldots\ldots$\hfil $\ldots\ldots$}\s
\line{(2) $(\forall x(\forall y A)) \to \forall y A$\hfil $P_5$}\s
\line{(3) $\ldots\ldots\ldots\ldots\ldots$\hfil MP(\quad$,$\quad)}\s
\line{(4) $\ldots\ldots\ldots\ldots$\hfil $P_5$}\s
\line{(5) $A$\hfil MP(3$,$4)}\s
\line{(6) $\forall x A$\hfil $\ldots\ldots$}\s
\line{(7) $\forall y(\forall xA)$\hfil $\ldots\ldots$}\medskip
\centerline{(C) $A,\ (\forall xA) \to C \vd \forall x C$}\s
\line{(1) $\ldots\ldots\ldots\ldots$\hfil $\ldots\ldots$}\s
\line{(2) $\forall xA$\hfil $\ldots\ldots$}\s
\line{(3) $\ldots\ldots\ldots\ldots$\hfil $\ldots\ldots$}\s
\line{(4) $C$ \hfil MP(2$,$3)}\s
\line{(5) $\forall xC$\hfil $\ldots\ldots$}\medskip
\centerline{(D) $\forall x(A \to B),\ \forall xA \vd \forall xB$}\s
\line{(1) $\ldots\ldots\ldots\ldots$\hfil Premise}\s
\line{(2) $\ldots\ldots\ldots\ldots$\hfil Premise}\s
\line{(3) $\ldots\ldots\ldots\ldots$\hfil $P_5$}\s
\line{(4) $A \to B$\hfil MP(\quad$,$\quad)}\s
\line{(5) $(\forall xA) \to A$\hfil $\ldots\ldots$}\s
\line{(6) $A$\hfil MP(2$,$5)}\s
\line{(7) $B$\hfil $\ldots\ldots$}\s
\line{(8) $\forall xB$\hfil $\ldots\ldots$}\pars
\ss
\noindent {\bf 3.7 Soundness and Deduction Theorem for $Pd\p.$}\parm
Although the model theory portions of this text have been restricted to sentences$,$ this 
restriction is technically not necessary. Many authors give valuation 
processes for any member of $Pd\p.$ This yields certain complications that 
appear only to be of interest to the logician. In the sciences$,$ 
a theory is determined by sentences. These are language elements that carry a 
definite  ``will occur'' or  ``won't occur'' content. As far as relations 
between proof theory concepts and the model theory concepts$,$ we consider 
 $\Gamma \vd A,$ where the elements need not be sentences. But for our model 
theory$,$ we$,$ at the least$,$ require that each member of $\Gamma$ is a sentence
and can always let $A$ be a sentence when $\Gamma \mod A$ is considered. When 
the notation $\str$ is employed$,$ it will be assume that $\str$ is a structure 
for $Pd\p.$ Further$,$ unless it is necessary$,$ we will not 
mention that our language variable are member of $Pd\p.$ Next is an important relation between a model and a formal demonstration. \parm
{\bf Theorem 3.7.1} ({\bf Soundness Theorem}) {\sl Suppose that $\Gamma \subset \c S,$  
$\str \mod \Gamma,$ and $\Gamma \vd A.$ Then $\str \mod A$ (i.e. $\str\mod \forall A.$)}\parm
Proof. In this proof$,$ we use induction on the steps $\{B_i\mid i=1,\ldots,n\}$ in a 
demonstration that $\Gamma \vd A$, and show that $\str\mod \forall B_i.$ \pars
Case $(n =1,\ 2,\ 3).$ Let $n =1.$ Suppose that $B_1$ is an 
instance of axioms $P_1,\ P_2,\ P_3.$ First$,$ recall the special process (i) of 
theorem 3.4.8. 
Thus$,$ to establish that $\str \mod P_i$$,$ we may 
consider the language variables as being sentences. 
Now simply replicate the proof for theorem 
3.4.2.\pars  
 Suppose that $C=\forall x B \to 
S^x_\lambda B]$ 
is an instance of axiom 
$P_5$. Again we need only assume that $x$ is the only free variable 
in $B$. Suppose that $\lambda= y$ is 
free for $x$ in $B.$ Then theorem 3.4.7 (iii) yields the result.  
Let the substitution in $P_5$ be a constant $\lambda= c.$ Suppose that $\str 
\mod \forall xB.$  Then for all $d\p \in D,$ $\str \mod S^x_dB.$ In 
particular$,$ 
$\str \mod S^x_c B.$ 
We leave it as an exercise$,$ to show that if $C$ is an instance of axiom 
$P_4$$,$ then $\str \mod C.$
Finally$,$ suppose that $C\in \Gamma.$
Then $\str \mod \Gamma$ means that $\str \mod C.$ \pars
Suppose $n = 2$ and $B_2= \forall x B.$ Then for the case $n = 1,$
$B_1 = B$ and $\str\mod B$ means that $\str\mod \forall B.$ Hence$,$ $\str\mod
\forall(\forall xB).$\pars
Suppose that $n = 3$ and $B_3=C$ is an MP step. Then for the two previous steps$,$ we have that
$\str\mod B$ and $\str\mod B\to C.$ From definition 3.3.3 and the special process$,$ $\str\mod C.$\pars
Case $(n+1>3).$ Assume the induction hypothesis that the result holds for any 
proof with $n$ or less steps. Consider a proof with $n+1$ steps. Now let $C = B_{n+1}.$ If $C$ is an instance of any of the axioms$,$ then the result follows as in 
the case where $n = 1.$ Suppose that $C$ comes from two previous steps$,$ 
$ B_i=B $ and $B_j=B \to C$ by MP.   By 
induction and the case $n =3$ process$,$ we have that $\str \mod C$. Now suppose that 
$C = \forall x B$ is an instance of generalization. 
 By induction and the case $n =2$ process$,$ we have that $\str \mod C.$ 
Finally$,$ for the case that $C \in \Gamma$$,$ the result follows as in case $n = 1.$ The complete result follows 
by induction. \qed
{\bf Corollary 3.7.1.1} {\sl If $\vd A,$ then $\mod A.$}\parm  
Of course$,$ the importance of theorem 3.7.1 lies in desire that whatever is 
logically produced by the mind using scientific logic will  also hold ``true''
for any model for which each of the hypotheses hold  ``true.'' {\bf This is the 
very basis for the scientific method assumption that human deductive 
processes correspond to natural system behavior.} Before we can establish a 
converse to theorem 3.7.1 and with it the very important ``compactness'' 
theorem$,$ a deduction-type theorem is required. This theorem$,$ however$,$ will not
completely mimic the deduction theorem for the language $L\p.$\parm
{\bf Theorem 3.7.2} {\sl (Deduction theorem for $Pd\p.$) Let $\Gamma \subset 
Pd\p,\ A,\ B \in Pd\p.$ If $\Gamma, \  A \vd B$ and the demonstration for $B$ 
contains no application of Generalization for a variable free in $A,$ then 
$\Gamma \vd A \to B.$}\pars                              
Proof. We use all the methods described in the metaproof of theorem 2.13.1$,$ 
the deduction theorem for $L\p.$ Proceed by induction exactly as done in 
theorem 2.13.1 and replace the axioms $P_1,\ P_2,\ P_3,\ P_4,\ R_5$ with  
steps that lead to $A \to C,$ no new step contains the single formula
$A.$ Now if $A$ appears as a step then this formula as been replaced by a $A 
\to A.$ This leaves the MP and Generalization steps to replace. Now re-number 
(i.e. count) these steps with their original ordering as $G_1,\ldots, G_n$. 
We now consider the induction process on these steps. \pars
Case $n = 1.$ Suppose that $G_1$ is an MP step. Then alter this step to 
produce other steps as done in the metaproof of theorem 2.13.1. Now suppose 
that $G_1 =\forall y B,$ where $B$ is from one of the previous original steps 
and was not obtained by generalization or MP. However$,$ the original steps have 
all be replaced by $A \to B.$ Between this new step and the step $G_1$ insert 
the following three steps.\pars
\line{\indent\indent (1) $\forall y(A \to B)$\hfil Generalization}\s
\line{\indent\indent (2) $\forall y(A \to B) \to (A\to (\forall yB))$
\hfil  $P_4$}\s
\line{\indent\indent (3) $A \to \forall yB$\hfil M(1$,$2)}\parm
\noindent [Notice that the insertion of step (2) requires that $y$ not be free 
in $A$.] Now \underbar{remove} the original step $G_1 = \forall y B.$ \pars
Case (n+1) Assume that all the alterations have been made for $G_1,\ldots,G_n$ 
steps. If $G_{n+1}$ is MP$,$ proceed in the same manner as in theorem 2.13.1.
Suppose that $G_{n+1} = \forall yB,$ where $B$ comes from an original step.  
However$,$ by induction all of the previous steps $B_i,\ i \leq n$ 
have been replaced by $A \to B_i.$ Consequently$,$ using this altered step$,$ 
proceed as in case $n = 1$ to obtain the formula $A \to \forall y B$ that 
replaces step $G_{n+1}.$ By induction$,$ all the original steps have been changed
to the form $A \to B$ and $A$ does not appear as a step in our new 
demonstration. The last step in our old demonstration was $B$ and it now has 
been changed to $A \to B.$ Hence the theorem has been proved.\qed
{\bf Corollary 3.7.2.1} {\sl Let $\{A_1\ldots, A_n\} \subset \c S$ and 
$A_1\ldots, A_n \vd B.$ Then $\vd (A_1 \to (A_2 \to \cdots (A_n \to B)\cdots 
)).$}\parm                                                                      
{\bf Theorem 3.7.3} {\sl Let $\Gamma \subset 
Pd\p,\ A,\ B \in Pd\p.$ If $\Gamma, \vd A\to B,$  then 
$\Gamma,\  A \vd B.$}\pars    
Proof. The proof is the same as in theorem 2.13.1.\qed
{\bf Corollary 3.7.3.1} {\sl If $\vd (A_1 \to (A_2 \to \cdots (A_n \to B)\cdots 
)),$ then $A_1,\ldots,A_n \vd B.$}\parm
{\bf Theorem 3.7.4} {\sl Let $\{A_1,\ldots,A_n\} \subset \c S.$ Then 
$A_1,\ldots,A_n \vd B$ if and only if $\vd (A_1 \to (A_2 \to \cdots (A_n \to B)\cdots 
)).$}\parm
We now show why the deduction theorem must be stated with the additional 
restriction. \parm
{\bf Example 3.7.1} First$,$ it is obvious that $P(x) \vd \forall x P(x).$
Just do the two step demonstration with the premise as the first step$,$ and 
generalization in $x$ as the second step. Now to show that $\not\vd P(x) \to 
\forall xP(x).$ Suppose that we assume that $\vd P(x) \to \forall x P(x).$ 
Consider the follows formal proof.\parm
\line{(1) $\vd P(x) \to \forall x P(x)$\hfil Given}\s
\line{(2) $\vd (B \to A) \to ((\neg A) \to (\neg B))$\hfil (A) page 83.}\s
\line{(3) $\neg (\forall x P(x)) \to (\neg P(x))$\hfil MP(1$,$2)}\s
\line{(4) $\forall x(\neg (\forall x P(x)) \to (\neg P(x)))$\hfil G(3)}\s
\line{(5)  $\forall x(\neg (\forall x P(x)) \to (\neg P(x)))\to(\neg (\forall x P(x)) \to (\forall x(\neg P(x))))$\hfil $P_4$}
\line{(6) $ \neg (\forall x P(x)) \to (\forall x(\neg P(x)))$\hfil MP(4$,$5)}
\parm
Now by the soundness theorem for $Pd\p$ (6) yields $\mod \neg (\forall x P(x)) \to (\forall x(\neg P(x))).$  Consider the structure $\str$ where $D 
=\{a\p,b\p\},$ and $P\p = \{a\p\}.$ Suppose that $\str\mod \neg (\forall x 
P(x)).$ This means that it is not the case that $ a\p\in P\p,$ or that 
$b\p\in P\p.$ Since $b\p \notin P\p,$ then $\str \mod \neg (\forall x 
P(x)).$ But $\str\not\mod \forall x(\neg P(x))$ for this means that 
for $a\p$ and $b\p$ we must have that $a\p \notin P\p$ and $b\p \notin P\p$ 
which is not the case. Hence$,$ $\not\mod \neg (\forall x P(x)) \to 
(\forall x(\neg P(x))).$  This contradiction shows that 
$\not\vd \neg (\forall x P(x)) \to (\forall x(\neg P(x))).$ Thus the use of 
the unrestricted deduction theorem to go from $P(x) \vd P(x),$ 
where we used 
generalization on $x,$ to obtain $P(x) \vd \forall x P(x),$ then to $\vd
P(x) \to \forall x P(x)$ is in error. \parm
\hs
\hm
\centerline{\bf EXERCISES 3.7}\parm
\noindent 1. Show that if $A \in Pd\p,$ and $A$ is an instance of axiom $P_4.$ 
Then $\mod A.$\pars\ss
\noindent {\bf 3.8 Consistency$,$ Negation completeness and Compactness.}\parm
Almost all of the sciences$,$ engineering and any discipline which must 
determine whether or not something ``will occur'' or  ``won't occur$,$'' use 
first-order languages to describe behavior.  The material in this section 
covers  
the most important of all of the first-order language concepts. 
The ramifications of these investigations$,$ 
 some completed only 2 years ago$,$ cannot be over stated. Many 
individuals who first obtained the results presented next$,$ became$,$ over night$,$
world famous figures within the scientific  community due to the significance of 
their findings. The conclusions 
in the next section on applications$,$ can only be rigorously obtained because of
the results we next present.\parm
\hs
{\bf Definition 3.8.1} A set of sentences $\Gamma \subset \c S$ is {{\it 
(formally) consistent}} if there does NOT exist a sentence $B \in \cal S$ 
such that 
$\Gamma \vd B$ and $\Gamma \vd \neg B.$\pars                                 
A set of sentences $\Gamma \subset \c S$ is {{\it 
inconsistent}} if it is not consistent.\pars                                 
A set of sentences $\Gamma \subset \c S$  is {{\it (negation) complete}}
if for every  $B\in \c S$ either $\Gamma \vd B$ or $\Gamma \vd \neg B.$\pars   
\hm
{\bf Theorem 3.8.1} {\sl A set of sentences $\Gamma$ is inconsistent if and only if
for every $B\in Pd\p,\ \Gamma \vd B.$}\parm
Proof. First$,$ suppose that $\Gamma$ is inconsistent. Then there is some
$B \in \cal S$ such that $\Gamma \vd B$ and $\Gamma \vd \neg B.$ Hence there are 
two demonstrations that use finitely many members of $\Gamma$ and yield the final 
steps $B$ and $\neg B$. Put the two demonstrations together$,$ and let $B$ appear 
as step $i$ and $\neg B$ appear as step $B_j.$ Now$,$ as the next step$,$ put the 
formal theorem (from the propositional calculus) $\vd \neg B \to (B \to A)$
from Exercise 2.11 problem 2b$,$ where $A$ is ANY member of $Pd\p.$ Now two MP 
steps yields the formula $A.$\pars
Suppose that for any formula $C \in Pd\p,$ it follows that 
$\Gamma \vd C$. Thus take any sentence $B,$ where $\neg B\in \cal S$ also.
Then $\Gamma \vd B$ and $\Gamma \vd \neg B.$ \qed
What theorem 3.8.1 shows is that if $\Gamma$ is inconsistent$,$ then there is an 
actual logical argument that leads to ANY pre-selected formula. The logical 
argument is correct. But the argument cannot differentiate between the 
concepts of  ``will occur'' or  ``won't occur.'' Indeed$,$ any scientific 
theory  using such an argument  
would be worthless as a predictor of behavior. Well$,$ is it possible that our 
definition for logical deduction without premises is the problem or must 
it be the premises themselves that lead to the serious theorem 3.8.1 
consequences?\parm
{\bf Theorem 3.8.2} {\sl There does not exist a formula $B \in Pd\p$ such that 
$\vd B$ and $\vd \neg B.$}\parm
Proof. Suppose that there exists a formula $B$ such that $\vd B$ and $\vd\neg 
B.$ Using the same process as used in theorem 3.8.1$,$ but only assuming that $B \in Pd\p$$,$ it follows that for 
any sentence $A \in Pd\p$$,$ $\vd A.$ Hence $\vd \forall A.$ Thus by corollary 
3.7.1.1$,$ for any structure $\str$ we have that $\str \mod \forall A.$ But the sentence
in step (5) of example 3.7.1$,$ does not have this property. The result follows 
from this contradiction.\qed
Thus the difficulties discussed after theorem 3.8.1 are totally caused by the 
premises $\Gamma$ and not caused by the basic logical processes we use. 
However$,$ if it's believed that a set of sentences in consistent$,$ then we 
might be able to obtain a larger consistent set. \parm
{\bf Theorem 3.8.3} {\sl Let $\Gamma \subset \c S$ be a consistent. Suppose 
that $\Gamma \not\vd B \in Pd\p.$ Then $\Gamma,\ (\neg \forall B)$ is 
consistent.}\parm
Proof. Assume that there exists some sentence $C$ 
such that (a) $\Gamma,\ (\neg    
\forall B)\vd C$ and (b) $ \Gamma,\ (\neg \forall B)\vd \neg C.$ 
Now from the propositional 
calculus (c) $\vd (\neg C \to (C \to \forall B)$ and$,$ hence$,$ this holds for 
$Pd\p.$ 
Using both demonstrations  (a) and (b)$,$ and inserting them as steps in a 
demonstration$,$ adjoin the step $(\neg C)\to (C \to \forall B).$ 
Two MP steps$,$ yield 
that $\Gamma,\ \neg \forall B\vd \forall B.$ Since $\neg \forall B$ is a 
sentence$,$ then the deduction
theorem yields $\Gamma \vd (\neg \forall B)\to \forall B.$ Now adjoin the 
steps for the 
$\vd ((\neg \forall B)\to \forall B) \to \forall B.$ One more MP step yields 
the contradiction
that $\Gamma \vd \forall B \vd B.$ \qed
{\bf Corollary 3.8.3.1} {\sl Suppose that $\Gamma$ is any set of 
sentences and $B$ is any sentence. If $\Gamma,\neg B$ is inconsistent$,$ then 
$\Gamma \vd B.$  If $\Gamma, B$ is inconsistent$,$ then $\Gamma \vd \neg B.$}\parm
Proof. If $\Gamma$ is inconsistent$,$ then the result follows. If $\Gamma$ is 
consistent$,$ then the result follows from the contrapositive of theorem 
3.8.3.  The second part follows from $\neg(\neg B) \vd B$ and $B\vd 
\neg(\neg B).$ \parm
{\bf Theorem 3.8.4} {\sl If consistent $\Gamma \subset \c S,$ then there 
is a language $Pd^{\prime\prime}$ that contains all of the symbols of 
$Pd\p,$ but with an 
additional set of constants and only constants adjoined$,$ and a set of 
sentences $\Gamma^{\prime\prime} \subset Pd^{\prime\prime}$ such that 
$\Gamma \subset \Gamma^{\prime\prime}$ 
and $\Gamma^{\prime\prime}$ is consistent and negation complete.}\pars
Proof. See appendix.\parm
We now come to the major concern of this section$,$ an attempt to mimic$,$ as 
close  as possible$,$ the completeness and consistency results for $L\p.$ In 
1931$,$ G\"odel$,$ in his doctoral dissertation$,$ convincingly established his 
famous  ``completeness'' theorem for $Pd\p$ by showing for the natural number 
domain that $\vd B$ if and only if $\mod B.$ Since that time G\"odel's methods 
have been highly refined and simplified. Indeed$,$ totally different methods 
have achieved his results as simple corollaries and these new methods have 
allowed for a greater comprehension of the inner workings of structures and 
models. \parm
\hs\hm
{\bf Theorem 3.8.5} {\sl If consistent $\Gamma \subset \c S$$,$ then there 
exists a structure $\str$ for $Pd\p$ such that $\Gamma \vd A$ if and only if 
$\str \mod \forall A.$ Further$,$ the domain of $\str$ is in one-to-one correspondence with the natural numbers.}\parm
Proof. See appendix. \pars
\hs\hm     
{\bf Corollary 3.8.5.1} (G\"odel) {\sl Let $A \in Pd\p.$ $\vd A$ if and only 
if $\mod A.$}\parm
Proof. Let $\Gamma = \emptyset.$ Assume that $\vd A.$  Corollary 3.7.1.1 to 
the soundness theorem yields $\mod A.$\pars
Conversely$,$ assume that $\not\vd A.$ Then $\not\vd \forall A.$ Thus the set
$\{\neg(\forall A)\}$ is consistent by theorem 3.8.3. Hence$,$ from 3.8.5$,$
there exists a structure $\str$ such that $\str\mod \{\neg(\forall A)\}.$ 
Consequently$,$ $\str \not\mod \forall A.$ Thus $\not\mod A.$ \qed
Due to corollary 3.8.5.1$,$ the abbreviations we have used that yield  
quantifiers $\exists x$$,$ and abbreviations for the connectives $\lor,\land,\ 
\iff$ all have the correct properties with respect to $\mod$ as they would 
under the truth-value definitions for $L.$ This allows us to use the more 
expressive language $Pd$ rather than $Pd\p.$ \parm

{\bf Theorem 3.8.6***} {\sl A nonempty set of sentences $\Gamma \subset Pd$ is consistent 
if and only if $\Gamma$ has a model.}\pars
Proof. Let $\Gamma$ be consistent. For each $A \in\Gamma$$,$ $\Gamma \vd A.$ 
Theorem 3.8.5 yields that there is a structure $\str$ such that 
$\str \mod A$ for each $A \in \Gamma.$ Hence $\str$ is a model for $\Gamma.$ \pars
Conversely$,$ assume that $\str$ is a model for $\Gamma$ and that $\Gamma$ is 
inconsistent. Then let $B$ be any sentence. Then $\Gamma \vd B$ and $\Gamma 
\vd \neg B.$ Hence the soundness theorem would imply that $\str \mod B$ 
and $\str \mod \neg B.$ This contradicts definition 3.3.3 part (e).\qed
{\bf Theorem 3.8.7} (Extended Completeness) {\sl If $\Gamma \subset Pd$ is a 
set of sentences$,$ $B$ is a sentence and $\Gamma \mod B,$ then
$\Gamma \vd B.$}\parm
Proof. Let $\Gamma \mod B$ and $\Gamma$ is not consistent. Then $\Gamma \vd 
B.$ So assume $\Gamma$ is consistent and that $\Gamma \not\vd B.$ Then from 
theorem 3.8.3$,$ the set of sentences $\{\Gamma, \ \neg B\}$ is consistent. 
Thus $\{\Gamma, \ \neg B\}$ has a 
model $\str.$ Therefore$,$ $\str \not\mod B.$ But $\str\mod \Gamma.$ Hence
$\str \mod B.$ This contradiction yields the result.\qed
We now come to one of the most important theorems in model theory. This 
theorem would have remained the only why to obtain useful models if other 
methods had not been recently devised. It's this theorem we'll use in the next 
section on applications.\parm
{\bf Theorem 3.8.8} (Compactness) {\sl A nonempty set of sentences $\Gamma \subset Pd$ has a model if 
and only if every finite subset of $\Gamma$ has a model.}\parm
Proof. Assume that $\Gamma$ has a model $\str.$ Then this is a model for any 
finite subset of of $\Gamma.$ [Note: it's a model for an empty set of 
sentences since it models every member of such a set. (There are none.)]\pars
Conversely$,$ suppose that $\Gamma$ does not have a model. Then $\Gamma$ is 
inconsistent. Hence$,$ taken any sentence $B \in Pd.$ Then there are two 
finite subsets $F_1,\  F_2$ of $\Gamma$ such that 
$F_1 \vd B$ and $F_2 \vd \neg B.$  Thus the finite 
subset $F=F_1 \cup F_2$ is inconsistent.
Hence $F$ has no model. This completes the proof.\qed
\centerline{\bf Important Piece of History}
It was convincingly demonstrated in 1931 by G\"odel that there probably 
is no formal 
way  to demonstrate that significant mathematical theories are consistent.
For example$,$ with respect to set-theory$,$ 
the only known method would be to establish a contradiction that is 
not forced upon set-theory by an actual demonstrable error on the part of an 
individual.  
All mathematical and scientific theories can be stated in the first-order
language of set-theory. There are approximately 50$,$000 research papers 
published each year in the mathematical sciences. No inconsistencies in modern 
set-theory itself has ever been demonstrated. The number theory we have used
to  establish all of the conclusions that appear in this book is thousands of 
years old and again no contradiction has ever been demonstrated. The theory 
of real numbers is hundreds of years old and no contradictions in the theory 
have ever been demonstrated. The specific axiomatic systems and 
logical procedures used to 
establish all of our results have never been shown to produce a contradiction.    
What all this tends to mean is that empirical evidence demonstrates
that our basic mathematical 
theories and logical processes are consistent. Indeed$,$ this evidence is more 
convincing than in any other scientific discipline. This is one reason why 
science tends to utilize mathematical models. \pars

{\bf Example 3.8.1} 
Let $\c N$ be the set of axioms for the natural numbers expressed in a set-theoretic 
language. Let $Pd$ be the first-order language that corresponds to $\c N$. 
Then the theory of natural numbers is the set of sentences $\Gamma =\{B\mid
B\in {\c S},\ {\c N}\vd B\}.$ This set $\Gamma$ is assumed$,$ from evidence$,$ to be consistent. 
Hence from theorem 3.8.6 there exists a structure $\str$ such that
$\str \mod \Gamma$. This means from the definition of $\str \mod \Gamma$ that 
there is a set of constants that {{\it name}} each member of the domain of the 
structure. This domain is denoted by the symbol $\nat.$ There are also n-ary 
predicate symbols for the basic relations needed for the axiom system and for 
many defined relations. \pars
We interpret $P(x,y)$ as the  ``='' of the natural numbers in $\nat.$ 
Depending upon the structure$,$ this could be the simple identity binary relation
$P(x,x).$ Now we come to the interesting part. Let our language constants 
$C$ be the constants 
naming all the members $c\p$ of $\nat$ and adjoin to $C$ a new constant $b$ 
not a member of $C.$ Consider the sentences 
$\Phi =\{\neg P(c,b)\mid c \in C\}$ that are members of a first-order language 
$Pd_b$ that is exactly the same as $Pd$ except for adjoining one additional 
constant $b.$ 
Now consider the entire set of 
sentences $\Gamma \cup 
\Phi.$ What we will do is to show that $\Gamma \cup \Phi$ has a model$,$ $\hyper 
{\str},$ by application of the compactness theorem.\pars
Take any finite set $A$ of members of  $\Gamma \cup \Phi$. If each member of 
$A$ is a member of $\Gamma,$ then $\str$ is a model for $A$. Indeed$,$ if 
$\{a_1,\ldots, a_n\}$ are all the members of $A$ that are members of $\Gamma,$
then again $\str \mod \{a_1,\ldots, a_n\}.$ So$,$ assume that $\{a_{n+1},\ldots, 
a_m\}$ are the remaining members of $A$ that may not be recognized as members 
of $\Gamma.$ Well$,$ there are only finitely many different constant $\{c_{n+1},
\ldots, c_m\}$ that are contained in these remaining  sentences in $A$. 
The theory 
of natural numbers states that given any nonempty finite set of natural 
numbers$,$ there is always another natural number $c\p$ that is not equal to any 
member of this finite set. These constants correspond to a finite set of 
natural numbers $\{c\p_{n+1},\ldots,c\p_m\}$ and interpret $b$ to be one of 
the natural numbers $c\p.$ {\bf This process of interpreting$,$ for a domain of 
the original structure$,$ the one additional constant as an appropriate domain 
member is a general procedure that is usually needed when the compactness 
theorem is to be used.} 
Thus each of the remaining sentences is modeled by 
$\str$. The compactness theorem states that there is a structure 
$\hyper{\str}$ such that $\hyper{\str}\mod \Gamma \cup \Phi$. 
Thus there is a domain $D$ and various n-place relations that behave exactly 
like the natural numbers since $\Gamma$ is the theory of the natural numbers
and $\Hyper \str \mod \Gamma.$   Indeed$,$ the interpretation $I$ restricted 
to the 
original set of constants yields an exact duplicate of the natural numbers and 
we denote this duplicate by same symbol $\nat.$ 
But there 
exists a member of the domain $D$ of this structure$,$ say $b\p,$ such that for 
each $c\p \in \nat,$  $b\p \not= c\p.$ But for every original $c\p$ we have that $c\p = c\p.$  
 Therefore$,$ we have a structure that {\it behaves}$,$ as described by 
$Pd_b,$ like the natural 
numbers$,$ but contains a new member that does not correspond to one of the 
original natural numbers. [End of example.]\parm
I give but two exercise problems for this section$,$ each relative to showing 
that
there exists mathematical structures that behave like well-known 
mathematical structures but that each contains a significant new member with a 
significant new property. \parm
\hs\hm
\centerline{\bf EXERCISES 3.8}\parm
\noindent 1. Modify the argument given in example 3.8.1 as follows: let 
$L(x,y)$ correspond to the natural number binary relation of ``less than''
(i.e.  $<$). Give an argument that shows that there is a structure 
$\hyper{\str}$ that behaves like the natural numbers but in which there 
exists a member $b\p$ that is ``greater than'' any of the original 
natural numbers. [Note: 
this solves the  ``infinite'' natural number problem by showing that there 
exists 
a mathematical object that behaves like a natural number but is  ``greater 
than'' every original natural number.\parm
\noindent 2. Let $\real$ denote the set of all real numbers. Let $C$ be a set 
of constants naming each member of $\real.$\footnote\ddag{The assumption is that the 
theory of real numbers is a  consistent theory and$,$ hence$,$ has a model. The proofs of Theorems 3.8.4 and 8.8.5 in the appendix$,$ show that such a theory has a model with a ``denumerable'' domain (i.e. can be put into a one-to-one correspondence with the natural numbers) and$,$ hence$,$ for this structure the 
set of constants used in our language for the real numbers and in definition 3.1 can be extended so as to name all the members of $\real.$}
 Suppose that $b$ is a constant 
not a member of $C$. Let $\Gamma$ be the theory of real 
numbers. Let $Q(0,y,x)$ be the 3-place predicate that corresponds to the 
definable real number 3-place relation  $0\p <c\p <d\p$$,$ where $0\p,\ c\p,\ 
d\p \in 
\real.$  Now in the real numbers there is a set of elements $G\p$ such that 
each member $c\p$ of $G\p$ has the property that $0\p< c\p$. Let $G$ be the set 
of constants that correspond to the members of the set $G\p.$ Consider the 
set of sentences $\Phi =\{Q(0,b,g)\mid g \in G\}$ in the language $Pd_b.$ 
Give an argument 
that  shows that there exists a structure $\hyper{\str}$ such that 
$\hyper{\str} \mod \Gamma \cup \Phi.$ That is there exists a mathematical 
domain $D$ that behaves like the real numbers$,$ but $D$ contains a member
$b\p$ such that $b\p$ is ``greater than zero'' but $b\p$ is  ``less than''
every one of the original positive real numbers. [Note: this is an example of an infinitesimal and 
solves the three hundred year old problem of Leibniz$,$ the problem to show the 
mathematical existence of the infinitesimals. These are objects 
that behave like real numbers but are not real 
numbers since they are  greater than zero but less than any positive real 
number.] \pars

The method used in this section and the next to obtain these ``nonstandard'' objects has been greatly improved over the years. The method presented here relies upon the assumption that the theory of real numbers and the theory of propositional deduction are consistent theories and$,$ hence$,$ have models. Also other processes have been employed that $,$ although correct$,$ may need further justification. Since the middle and late 1960s$,$ algebraic methods have been used to first obtain structures associated with first-order statements. One structure$,$ the standard structure ${\c M},$ would be based upon predicates that are present within statements about the real numbers or$,$ as explained in the next section$,$ the first-order theory of propositional deduction. Another structure that can be considered as containing $\c M$ and termed an ``enlargement$,$'' $\Hyper {\c M},$ is then constructed. These constructions do require the acceptance of an additional set-theoretical axiom. To obtain the various nonstandard objects discussed in this and the next section$,$ all that is needed is to have enough real number or propositional deduction statements hold in the standard structure. This will immediately show that the nonstandard objects$,$ the infinitesimals$,$ ultrawords and ultralogics$,$ mathematically exist within the $\Hyper {\c M}.$ 
\parm

Due to the significance of exercise problem 2 above, the following slightly more general result is established here. \pars
The following convention is used. The symbols used to denote constants and predicates will also be used to denote the members of the domain and  relations in a structure $\str.$ Let $\Gamma$ be the theory of real numbers $\real$. This theory includes all of the defined relations and all possible deduced sentences, etc. For this example, there is a 1-place predicate  $P(x)$ that is modeled by the set of all positive real numbers $P = \{x\mid (0< x)\land (x\in \real)\},$ where $0$ corresponds to the zero in $\real$. The basic predicate to be used here is $Q(0,x,y),$ where $Q$ will correspond to the relation $Q = \{(0,x,y)\mid (0<\vert x\vert <y)\land (x\in \real)\land (y\in P) \}.$\pars
The most basic assumption is that $\Gamma$ is consistent and from Theorem 1 it has a model $\str$, where $\real$ is the domain. The facts are that from a viewpoint external to this structure, it can be assumed that all we need to name every member of $\real$ is a set of constants that can be put into one-to-one relation with the even natural numbers, our most basic assumed consistent theory. So, let there be such a set of constants $\cal C$
and consider but one more constant $b$. Consider the following set of sentences
$${\cal A}= \{Q(0,b,a)\mid a \in P\subset {\cal C}\}.$$
Note that this does not add any new predicates. 
What we do is to show that the set of sentence $\Gamma \cup {\cal A}$ has a model. \parm

{\bf Theorem.} {\it There is a structure $\hyper \str$ that models $\Gamma \cup {\cal A},$ where the  domain $\hyperreal$  contains new objects not in $\real$ and the entire theory $\Gamma$ holds for $\hyperreal.$}\pars
Proof. We have that $\Gamma$ has a model $\str.$ Let 
finite $\{A_1,\ldots A_n\} \subset {\cal A}.$ If $\{A_1,\ldots A_{i}\} \subset \Gamma,$ then $\str\mod  \{A_1,\ldots A_{i}\}.$ Suppose that 
$\{A_{i+1},\ldots,A_n\}\subset {\cal A},$ where $A_k= Q(0,b,a_k),\ i+1\leq k \leq n.$ Then $\{a_{j+1},\ldots,a_n\}$ contains a minimal member in $\real$, $a_k$. Now let $b$ correspond to $a_k/2.$ Then for this member of $\real$ we have that $\str\mod \{A_{i+1},\ldots,A_n\}.$ Thus, by the Compactness Theorem, there is a structure $\hyper \str$ such that $\hyper \str\mod {\cal A}.$  \qed
Now there are no new predicates in $\cal A.$ The only thing that can be different is the domain $\hyperreal.$ Of course, all relations that model these predicates are now relative to $\hyperreal.$ What we are doing is metalogically investigate the two structures $\str$ and $\hyper \str$ from the external world rather than only internal to either of the structures. By convention, we denote some of the predicates as they are interpreted in $\hyper \str$ for $\hyperreal$ by the same letter when there is no confusion as to which domain they apply; otherwise, they will be preceded by a * (translated by the word ``hyper''). As pointed out, each member of the original $\real$ is named by a constant from $\cal C$ and these same constants name members of $\hyperreal.$ We call each such interpreted constant a {\bf standard} member of $\hyperreal.$ What happens is that there is a object in $\hyperreal$ that now is being named by $b$, and that we denote by $\eps.$ This $\eps \not= a$ for any standard member of $\hyperreal$. Why? Well, take any standard $r \in \hyperreal$. Using the theory $\Gamma$, if $\eps = r$, then we have contradicted that if $\eps < r,$ then $\eps \not= r.$ Thus, this is a new object that does not correspond to any standard element of $\hyperreal.$ Moreover, there are two such objects since $\Gamma$ holds for $\hyperreal$, and, hence, $-\eps \in \hyperreal$ and $\vert -\eps \vert = \vert \eps \vert.$ [As indicated the ``$\leq$,'' ``$\vert \cdot \vert,$'' and ``='' predicates must be interpreted in $\hyper \str$. This is the last time I'll mention this fact.]\pars
The objects $-\eps,\ \eps$ are called {\bf infinitesimals.} But do they behave like the {\it infinitely small or little} ideal numbers $o$ of Newton? Newton, as well as Leibniz, required the $o,$ on one hand, to behavior like a non-zero real number but then to behave like a zero. This is why  some people rejected the entire idea of an ``ideal'' number like $o$ since it contradicts real number behavior. For example, Newton divides the equation $3pox^2 +3p^2o^2x +p^3o^3 -2dqoy -dq^2o^3-abpo=0$ by $o$ where $o$ cannot act like a zero, and gets  $3px^2 +3p^2ox +p^3o^2 -2dqy -dq^2o-abp=0.$ But now, he treats $o$ as if it is zero and writes ``Also those terms are infinitely little where $o$ is. Therefore, omitting them there results $3px^2 -abp -2dqy = 0.$'' This is a direct contradiction as to the behavior of the $o$ and, indeed, one should be able to apply the ``omitting'' process to the first equation and this would only yield the identity $0 = 0.$\pars
First, let $\monad {0}$ be the set of all infinitesimals in $\hyperreal$ and include standard $0$ in this set. To determine that members of $\monad{0}$ have properties different from those of the original $\real$, we must investigate these from the ``meta'' viewpoint. The non-Greek lower case letters will always be the constants used to name the standard members of $\hyperreal$.\pars

(1) Let $0\not= \lambda \in \monad {0}.$ Let non-zero standard $r \in \hyperreal$ and arbitrary standard $x \in \hyper {P}.$ Then $0 < \vert \lambda \vert < x/\vert r\vert$ since $x/\vert r\vert$ is as standards member of $\hyper {P}.$  From the theory $\Gamma$, we have that $0 < \vert r\lambda \vert < x$. Thus $r\lambda \in \monad{0}$ since $x$ is an arbitrary standard member of $\hyper {P}$. [Note: Some would write this as $\hyper r\lambda \in \monad {0}.$] \pars 

(2) Using stuff from (1), let $\lambda,\ \gamma \in \monad{0}.$ We have that $0< \vert \lambda\vert  < x$ and $0 < \vert \gamma\vert < y$ where both $x$ and $y$ are arbitrary standard members of $\hyper P$. Hence, from $\Gamma$, we have that $0 < \vert \lambda + \gamma \vert \leq \vert \lambda \vert+\vert \gamma \vert < x + y.$ But $x + y =z$ is also an arbitrary standard member $\hyper P$. Thus, $\lambda + \gamma \in \monad{0}.$ In like manner $\lambda\gamma\in \monad{0}.$ \pars

(3) You can do all the ordinary real number algebra for the members of $\monad{0}$ since $\Gamma$ holds for them. However, if $0\not= \lambda,$ then $1/\lambda \not\in \monad {0},$ since for arbitrary standard $x \in \hyper P$, $x < 1/\lambda.$ \pars

(4) For $r$, let $\monad{r} = \{r + \lambda \mid \lambda \in \monad {0}\}.$ This is can a {\it monad about r}. Let $r_1 \not= r_2.$ Then $\monad{r_1}$ is completely disjoint from $\monad {r_2}.$ Why? Well, suppose not. Then there are two infinitesimals $\lambda,\ \gamma$ such that $\lambda + r_1 = \ 
 \gamma +r_2 \Rightarrow r_1 - r_2 = \gamma -\lambda \in \monad {0}.$ But $r_1 -r_2$ is a standard number and the only standard number in $\monad {0}$ is $0$. Thus, $r_1 = r_2;$ a contradiction.\pars

(5) Hence, every object in $\hyperreal$ that behaves like the original real numbers is ``surrounded,'' so to speak, by its monad. Because of this uniqueness, we can for every member of $\monad {r_1} \cup \monad {r_2}\cup \cdots \cup \monad {r_n} \cup \cdots= M_{fin},$ where the $r$ vary over all standard $\hyperreal,$ define an operator $\St$ with domain $M_{fin}$ that yields the unique standard $r.$ Although the notion of the ``limit'' need never be considered, the $\St$ operator mirrors  ``limit'' algebra and can be applied to the Newton material above. \pars
Let Newton's $o = \lambda \not= 0.$ Consider $3p\lambda x^2 +3p^2\lambda^2x +p^3\lambda^3 -2dq\lambda y -dq^2\lambda^3-abp\lambda=0.$  Now divide by $\lambda$ and get $3px^2 +3p^2\lambda x +p^3\lambda^2 -2dqy -dq^2\lambda-abp=0(1/\lambda)=0.$ Using the properties of the $\St$ operator applied to both sides of this equation, $3px^2 + 3p^2 \st {\lambda} x + p^3 \st {\lambda^2} -2dqy -dq^2 \st {\lambda} - abp = 3px^2 + 3p^2(0)x + p^3(0) - 2dqy -dq^2(0) - abp = 3px^2 -2dqy -abp =0.$ Which is the same result as Newton's and eliminates any contradiction. \pars\ss
\noindent {\bf 3.9 Ultralogics and Natural Systems.}\parm
In 1949 in the Journal of Symbolic Logic$,$ Leon Henkin gave a new proof for our 
theorem 3.8.5 from which the compactness theorem follows. Indeed$,$ the proof 
that appears in the appendix is a modification of Henkin's proof. I note 
that the appendix proof most also be modified if we wish to apply a 
compactness theorem to sets like the real numbers$,$ as needed for exercise 3.8 
problem 2. These modifications are rather simple in character and it's assumed 
that 
whenever the compactness theorem is used that it has been established for the 
language being discussed. \pars
What is particularly significant about the Henkin method is that it's hardly 
possible to reject the method. Why? Well$,$ as shown in the proof in the 
appendix the model is constructed by using the language constants and 
predicates themselves to construct the structure. What seems to be a very 
obvious approach was used 18 years after the first G\"odel proof. There is a 
certain conceptual correspondence between the use of the language itself and 
the material discussed in this last section.\pars
All of the metatheorems established throughout this text use a first-order 
metalanguage. These metatheorems describe various aspects of any formal
first-order language. But these metatheorems also apply to informal languages 
that can be represented or encoded by a formal language. That is such concepts 
as the compactness theorem can be applied to obtained models for various 
``formalizable'' 
natural languages and the logical processes 
used within science$,$ engineering and many other disciplines.
In 1963$,$ Abraham Robinson did just this with the first published paper 
applying a similar device as the compactness theorem to obtain some 
new models for the valuation process within formal languages. Your author  
has extended Robinson's work and has applied his  model theoretic methods 
to all natural languages such as  English$,$ French$,$ etc. The method used is the Tarski 
concept of the consequence operator. Before we apply the compactness theorem 
to obtain mathematically an ultralogic and an ultraword$,$ a few very simply 
communication concepts need to be discussed.\parm

\hs
{\bf Definition 3.9.1}  
A {{\it Natural system}} is a set or arrangement of physical entities 
that are so related or connected as to form an identifiable 
whole. Science specifically defines such collections of entities and gives them 
identifying names. The universe in which we dwell$,$ our solar system$,$
the Earth$,$ or a virus are Natural systems.\pars
\hm
The appearance of most Natural systems changes with  ``time.'' I'll not define 
the concept of time$,$ there are various definitions$,$ and I simply mention 
that the time concept  
can be replaced by something else called a universal event number if the time 
concept becomes a philosophical problem.  
One of the 
most important aspects of any science that studies the behavior of a
 Natural system is the communication of the predicted or observed behavior to 
other individuals. This communication can come in the form of word-pictures
or other techniques I'll discuss below. Even if science cannot predict the 
past or future behavior of an evolving or developing (i.e. changing in time) 
Natural system$,$ it's always assumed that at any instant of time the appearance 
can be described. \parm\vfil\eject
\hs
{\bf Definition 3.9.2} A {{\it Natural event}} is the actual and real objective 
appearance of a Natural system at a given instant whether or not it's 
described in any language or by other techniques.\pars
\hm  
From a scientific communication point of view$,$ 
a description is all that 
can be scientifically known about such a Natural event and the description is 
substituted for it. Relative to the behavior of a Natural system$,$
a  general scientific 
approach is taken and it's assumed that scientists are interested in various 
types of 
descriptions for Natural system behavior. It's not difficult to show that all 
forms of scientific description can be reduced to finitely long 
strings of symbols. 
Modern computer technology is used to produce an exact string of 
symbols that will reproduce$,$ with great clarity$,$ any photograph$,$ TV tape$,$ 
or sound.  
Today$,$ information is ``digitized.'' Relative to a visual instant-by-instant 
description$,$ television is used.  
Each small fluorescent region on a TV 
screen is given a coded location within a computer program. What electron 
beam turns 
on or off$,$ the intensity of the beam and the like is encoded in a series of 
binary digits. The computer software then decodes this information and the 
beam sweeps out a glowing picture on a TV tube. At the 
next sweep$,$ a different decoded series of digits yields a slightly different 
picture. And$,$ after many hundreds of these sweeps$,$ the human brain coupled 
with the eye's persistence of vision yields a faithful mental motion picture. 
Record companies digitize music in order to improve upon the reproduction 
quality. 
Schematics for the construction of equipment can be faithfully described in 
words and phrases if a fine enough map type grid is used. Thus$,$
the complete computer software expressed in a 
computer language$,$ the digitized 
inputs along with schematics of how to build the equipment to encode and 
decode digitized information$,$ taken together$,$ can be considered as an 
enormous symbol string the {\bf exact} content of which will be what you 
perceive
on ``the tube$,$'' hear from a CD player$,$ or other such devices.\par
Much of what 
science considers to be {\bf perception} may be replaced by a long 
exact string of symbols. All of the methods used by science to communicate 
descriptions for Natural system behavior will be called the {{\it 
descriptions.}} 
\par
What all this means is that the actual objective evolution of a Natural system 
can be replaced by descriptions for how such a system appears at 
specific instances during its evolution. The actual time differences between 
successive ``snap shorts'' will depend upon the Natural system being studied$,$ 
but they could be minuscule if need be. Just think of a developing Natural 
system as an enormous sequence of Natural events$,$ that have been replaced by an 
enormous sequence of descriptions. As discussed above 
these descriptions can be replaced by finitely long 
strings of symbols of one sort of another. This communication fact is the 
common feature of all scientific disciplines. \pars 
The basic object used to study the behavior of Natural systems by means of 
descriptions of such behavior is the consequence operator described in section 
2.16. Actually$,$ consequence operators are used more for mathematical convenience 
than any other reason. For the purposes of this elementary  exposition$,$ 
ALL of the consequence operator results can be re-expressed in the 
$\vd$  notation.    
However$,$ in the next definition we indicate how {{\it deductive processes}} 
such as the propositional process denoted by $\vd$ determine a consequence 
operator.  
\parm

\hs
{\bf Definition 3.9.3} For a propositional language $\c L$ and a propositional 
type of deduction from premises $\Gamma \vd_\ell B,$  a consequence operator $C$ is 
defined as follows: For every $\Gamma \subset L,$ $C(\Gamma) = \{A\mid A\in 
{\c L}\  {\rm and}\  \Gamma \vd_\ell A\}.$ Thus $C(\Gamma)$ is the set of all formula 
``deduced'' from $\Gamma$ by the rules represented by $\vd_\ell.$\pars
\hm
 Now recall the propositional language introduced in section 2.16. It's 
constructed from a set of atoms $\{P_0,P_1,\ldots \}$ in the usual manner$,$ but 
only using the connectives $\land$ and $\to.$ Of course$,$ this language $ 
L_S$ is a sublanguage 
of our basic language $L.$ For formal deductive process$,$ there are four axioms 
written in language variables.\parm
\indent\indent (1) $(A \land (B \land C)) \to ((A \land B)\land C).$\pars
\indent\indent (2) $(A \land B) \land C \to (A \land (B\land C)).$\pars
\indent\indent (3) $(A \land B) \to A.$\pars
\indent\indent (4) $(A \land B) \to B.$\pars
Notice that these axioms are theorems in $L$ (i.e. valid formula) and$,$ 
thus$,$ 
also in $L_S.$ The process 
of inserting finitely many premises in a demonstration is retained. The one 
rule of inference is {\it modus ponens} as before. We denote the deductive 
process these instructions yield by the usual symbol $\vd$. The 
consequence operator defined by the process $\vd$ is denoted by $S.$\pars 
There is actually an infinite collection of deductive processes that can be 
defined for $L$ and $L_S.$ This is done by restricting the modus ponens rule. \parm
\hs {\bf Definition 3.9.4} Consider all of the same rules as described above 
for $\vd$ except we use only the MP$_n$ modus ponens rule.  
Suppose two previous 
steps of a demonstration are of the form $A \to B$ and $A$ {\it and}  the 
${\rm size}(A \to B) \leq n.$  Then you can write down at a larger step 
number the formula $A.$ [This is not the only way to restrict such an 
MP process.] This is 
the only MP type process allowed for the $\vd_n$ deductive process.\pars
\hm 
The major reason the process $\vd_n$ is introduced is to simply indicate that 
there exist many different deductive processes that one can investigate. 
All of the 
metamathematical methods used to obtain the results in the previous sections 
of this text are considered as the most simplistic possible. They are the same 
ones used in the theory of natural numbers. Hence$,$ they are considered as 
consistent. Thus$,$ everything that has been done has a model. Let
$L = L'$ be the domain. 
 The language contains a set of atoms
$\{P_0,P_1 ,\ldots \}.$ Then we have the various 2-place relations between 
subsets of $L\p$ we have denoted by $\vd_n$ and $\vd.$ The first coordinate 
of each of these relations is a subset of $L\p,$  the premises$,$ and the 
second coordinate is a single formula deduced from the premises. If the set of 
premises is the empty set$,$ then the second coordinates are called theorems. 
To be consistent with our previous notation$,$ we would denote the 2-place 
relations for a deductive processes by $\vd_n\p$ and $\vd\p.$ This would not be 
the case for the members of $L\p.$\parm 
We now come to an important idea relative to the relation $\vd$ introduced, in 1978,  by the author of this book. The metatheory 
we have constructed in the past sections of this text$,$ 
the mathematical theory that tells us about the behavior of first-order
languages and various deductive processes$,$ is constructed from a first-order
language using the vary deductive processes we have been studying. Indeed$,$ all 
the proofs can be written in the exact same way as in the sections on formal 
deduction$,$ using a different list of symbols or even the same list but$,$ say$,$ 
in a different color. Another method$,$ the one that is actual used for the more 
refined and complex discusses of ultralogical and ultraword behavior$,$ is to 
use another formal theory$,$ first-order set theory$,$ to re-express all of these 
previous metatheorems. In either case$,$ the compactness theorem for first-order
languages and deduction would hold.\parm
There are two methods used to obtain models. If we assume that a mathematical 
theory is consistent$,$ then the theory can be used to define a structure$,$ which 
from the definition$,$ would be a model for the theory. Theory consistency would 
yield all of the proper requirements for definition 3.3.3. On the other hand$,$ 
you can assume that a structure is given. We assume that the n-place 
relations and constant named objects are related by a set of informal axioms 
that are consistent. Using this structure$,$ we develop new information about 
the structure. This information is obtained by first-order deduction and
is expressed in an informal first-order 
language. This informal first-order 
language can then be  ``formalized'' by substituting for the informal 
constants and relations$,$ formal symbols. This leads to a formal theory for 
the structure. Both of these methods yield what is called a {{\it standard}} 
structure and the formal theory is the  {{\it standard theory.}}\parm
For our propositional language $L,$ we have the entire collection 
$\c T$ of 
sentences that are established informally in chapter 2 about this language relative to 
various predicates and we 
translate these informal statements into a formal first-order theory $\cal E.$ One part of our mathematical analysis has been  associated with a standard 
structure. We have used first-order logic and a first-order language to 
investigate the propositional language and logic.
Thousands of  years of working with this structure has not produced a 
contradiction. The structure is composed of a domain $D,$ where  
 $D=L = C$ (the set of constants). There are various
n-place relations used. For example$,$ three simplest are 
the 2-place relation on this 
domain  $\vd\p,$ the 1-place relation ${\cal P}'$ that corresponds to the 
non-empty set of atoms $L_0 \subset L,$ and $x \in L$ that corresponds to $L(x).$ 
The relation $\vd'$ is defined as follows: $A \vd\p B,$ 
where $A,\ B \in L\p,$ if and only if $A \vd B.$ Then ${\cal P}'$ is defined 
as $P \in {\cal P}'$ if and only if ${\cal P}(P)$ (i.e. ``$P$ is an atom.'')
Although it is not necessary$,$ in order to have some ``interesting'' results$,$ we assume that there are as many atoms as there are natural numbers. [In this case$,$ it can be shown that there are also as many members of $L$ as there are natural numbers.]
Since $\c E$ is assumed consistent$,$ then there is a standard 
structure $\str =\langle D, L\p,\vd\p,{\cal P}',\ldots \rangle$ that models all of the n-placed predicates and constants that appear in $\c E,$ where we again note 
that for each $A\in L,$ $I(A) = A.$ The set $\c E$ is a subset of a first-order language $\c L,$ where $\c L$ is constructed from $C$ and all of n-placed predicates use in $\c E.$\parm  
 {\bf Theorem 3.9.1} {\sl Let $b$ be one new constant added to the 
constants ${C}.$  Construct a new first-order language ${\c L}\p$ with  
the set of constants $C \cup \{b\}$ and all of the n-placed predicates used to construct $\c L.$ Then ${\c L}\subset {\c L}\p.$  
Let $L_0$ be the set of atoms in $L\subset C,$ where  we consider  
them also as constants that name the atoms$,$ and let $ x \in L_0$ be the interpretation of the 1-place predicate ${\c P}(x).$ 
Consider the set of sentences $\Phi =\{{\cal P}(P) \land b \vd P \mid  
P \in L_0\}\subset {\c L}\p.$ Then there exists a model $\Hyper \str$ for 
${\c E} \cup \Phi.$}\par
Proof. This is established as in the example and exercises of section 3.8. Consider any non-empty finite $F \subset {\c E}\cup \Phi.$ Then $F = F_1 \cup F_2,$ where $F_1 \subset \c E$ and $F_2 \subset \Phi.$  Suppose that $F_1 \not= \emptyset.$ Then $\str\models F_1.$ 
 Suppose that $\{A_1,\ldots,A_n\}= F_2,\ n \geq 1.$ (Note: $A_i,\ 1 \leq i \leq n$ are all distinct.)  
  Consider the finite set 
$\{P_1,\ldots,P_n\}$  
of atoms that appear in $A_1,\ldots,A_n.$ (1) If $n = 1,$ then there is only one 
atom $P$ in this set$,$ and  $P \in D = L.$ Otherwise$,$ (2) consider the formula 
formed 
by putting $\land$ between each of the $P_i$ as follows: $(P_1 \land(\cdots(P_{n-1} \land P_n)$ and $(P_1 \land(\cdots(P_{n-1} \land P_n)\in D = L.$  
In case (1)$,$ let $b = P;$ in case (2)$,$ let 
$b = (P_1 \land(\cdots(P_{n-1} \land P_n).$
Then from the theory $\c E,$ we know that (1)
${\cal P}(P) \land b \vd P$ or (2) ${\cal P}(P_i)\land b \vd P_i$$,$ 
$i = 1,\ldots,n.$ Hence$,$ interpreting the $b$ in case (1) as $P$$,$ and $b$ in case (2) as $(P_1 \land(\cdots(P_{n-1} \land P_n),$ then   
$\str\models F_1.$  Consequently$,$ $\str$ models any finite subset of ${\c E} \cup \Phi.$ Hence by the compactness theorem there 
is a  model $\Hyper \str$ for ${\c E} \cup \Phi.$ \qed
Let $ \hyper {b}$ now denote that object in the domain of $\Hyper {\str}$ that satisfies 
each of the sentences in $\Phi,$ where $\Hyper {\vd}$ denotes the corresponding binary relation that corresponds to $\vd.$  
Since the formal theory $\c E$ corresponds to the informal theory $\c T$ and   $\Hyper {\str}\models {\c E},$ then the 
structure$,$ at the very least$,$  
behaves$,$ as described by the theory of chapter 2$,$ as a propositional logic.   
But does this new structure have additional properties? 
Note that $\Hyper {\vd}$  ``behaves'' 
like propositional deduction and $\hyper {b}$ ``behaves$,$'' thus far in this analysis$,$ simply like a formula in $L.$   
The object $\hyper {b}$ is 
called an {{\it ultraword}} and $\Hyper {\vd}$ is an example of a
(very weak) 
{{\it ultralogic.}} The reason why it is weak is that we have only related $\hyper {b}$ to the one relation $\Hyper {\vd}.$ Hence$,$ not much 
can be 
said about the behavior of $\hyper {b}$. However$,$ we do know that
$\hyper {b}$ does not behave like an atom in $\Hyper {\str}$
 for the following statement (3)
holds in 
$\str$; (3) $\neg \exists x({\cal P}(P_1)\land {\cal P}(P_2)\land {\c P}(x)\land x 
\vd P_1\land x\vd P_2).$ Hence$,$ (3) holds in $\Hyper {\str}.$ But the  $\neg \exists$ varies over the elements that$,$ at the least$,$ correspond to $C \cup \{b\}$ in its domain and$,$ hence$,$ this statement applies to the each member of $\{\hyper {P_1},\ \hyper {P_2}\},$ the $\Hyper {\str}$ interpreted constants $P_1,\ P_2,$ that are members of the $\Hyper {\cal P}.$ Hence$,$ whatever objects ``behavior'' like the atoms in $\Hyper {\str},$ $\hyper {b}$ is not 
one of them. If we consider other predicates$,$ then more information and properties will hold in another structure.  Suppose that 
$\lint C(\underline{\quad})\rint$: ``\underbar{\quad} is a consistent member of
 $L.$'' In 
$\Phi$$,$ replace ${\c P}(P)$ with ${\c P}(P) \land C(b)$ to obtain $\Phi\p.$ The same method as above shows that there is a structure $\Hyper {\str}_1$ such that $\Hyper {\str}_1 \models {\c E} \cup \Phi\p.$ The previous properties for $\hyper {b}$ still hold.  Using the sentence
$\neg \exists x(\forall y(({\cal P}(y) \land C(x)) \to x \vd y)),$ we are led to the conclusion that 
$\hyper{b} \Hyper{\vdash} \Hyper {P}$ in $\Hyper {\str}_1$ for all of the original atoms $P \in L_0.$ But$,$ when a comparison is made$,$ 
there is at least one other 
object that we name $d'$ that behaves like an ``atom'' and $\hyper{b} \Hyper{\not\vdash} d'$ in   
 $\Hyper {\str}_1.$\parm
{\bf Example 3.9.1}  Correspond each atom in $L$ to a 
description for 
the behavior of a Natural system at a specific moment. Let the ordering of the natural numbers to which each atom corresponds correspond to an ordering of an event sequence. Using the concepts of 
Quantum Logic$,$ one can interpret the 
ultralogic $\Hyper {\vd}$ as a physical-like process that when applied to
a single 
object $\hyper {b}$ yields each original interpreted description $\Hyper {P}.$ Using certain types of special constructions that yield $\Hyper {\str},$ one has that $\Hyper {P} = P.$ Hence$,$ under these conditions this $\hyper {b}$ and $\Hyper {\vd}$ yield the  
moment-by-moment event sequence that is the objectively 
real evolution of a Natural system. This idea applies to 
ANY Natural system including the universe in which we dwell and$,$ thus$,$ gives 
a describable process that can produce a universe. \parm
\baselineskip 12pt

Obviously$,$ the approach used to  obtain ultrawords and ultralogics as 
discussed 
in this section is very crude in character. There arise numerous questions 
that one would like to answer. By refining the above processes$,$ using 
set-theory$,$ consequence operators and other more complex procedures these 
questions have all been answered. I list certain interesting ones with  
(very) brief answer.\pars
(1) Can this process be refined so that the actual complex behavior expressed 
by each description is retained within the model while it's also being 
modeled by a proposition atom? (Yes)\pars
(2) Can the ultraword $\hyper {b}$ be analyzed? (Yes$,$ and they can have very interesting internal structures.)\pars
(3) Can you assign a size to $\hyper {b}$ and$,$ if so$,$ how big is it? (Yes$,$ and it is 
very$,$ very big. It is stuffed with a great deal of information.)\pars 
(4) If you take all the ultrawords  that generate all of the 
our Natural systems$,$ does there exist one ultraword that when $\Hyper {\vd}$
 is 
applied 
to this one ultraword then all the other ultrawords are produced and$,$ hence$,$ 
all of the 
event sequences for all of the Natural systems that comprise our universe are 
produced? (Yes. There is a ultimate ultraword $w$ such that when $\Hyper {\vd}$ 
is applied to $w$ all of the other ultrawords are produced as well as all the 
consequences produced by these other ultrawords. What this shows is 
that all natural event 
sequences are  
related by the physical-like process $\Hyper {\vd},$ among others$,$ applied to 
$w.$ This gives a solution to the General Grand Unification problem.) \pars
(5) Will there always be these ultranatural events no matter how we might 
alter the natural events? (Yes) \pars
(6) And many others. \pars
\vfil\eject
{\quad}

\centerline{\bf APPENDIX}
\bigskip\bigskip
\centerline{Chapter 2}\par\medskip
The major proof and definition method employed throughout this text is called 
{{\it induction}}. In the first part of this appendix$,$ we'll explore this 
concept which is thousands of years old.\parm
There are two equivalent principles$,$ the weak and the strong. In most cases$,$ 
the strong method is used. We use the natural numbers 
$\nat=\{0,1,2,3,4,\ldots \}$ either in their 
entirety or starting at some fixed natural number $m.$ There are two actual properties 
used. The first property is \parm
(1) Any nonempty finite set of natural numbers contains a maximal member. This 
is a natural number that is greater than or equal to every member of the set 
and is also contained in the finite set.\parm
(2) The ordering $<$ is a {{\it well-ordering.}} This means that any nonempty set of 
natural numbers$,$ finite or otherwise$,$ contains a first element. This means 
that the set contains a member that is less than or equal to every other 
member of the nonempty set. 
The actual induction property holds for many different sets of natural numbers 
and 
you have a choice of any one of these sets. Let $m\in \nat.$ Then we have 
the set $N_m$ of all natural numbers greater than or equal to $m.$ That is 
$N_m = \{n \mid n \in \nat,\ m \leq n\}.$  \parm

\hs
{\bf Principle 1.} (Strong Induction). Suppose that you take any nonempty $W 
\subset N_m.$ Then you show that a statement holds for \pars
(i) $m \in W.$\pars
(ii) Now if upon assuming that for a specific $n \in W,$ the statement we wish to 
establish holds for each $p \in \nat,$ where $m\leq p\leq n,$ 
you can show that $n+1 \in W$$,$ (this is called the strong induction hypothesis.) \pars
(iii) then you can declare that $W = N_m.$ (We say that the result follows by 
induction.) \pars                      
\hm
The major method in applying principle 1 is in how we define the set $W$. It 
is defined in terms of the statement we wish to establish. The 
set $W$ is defined by some acceptable description$,$ a set of rules$,$ 
that gives a method for 
{\it counting} objects. I have used the term ``acceptable.'' This means a 
method that is so clearly stated that almost all individuals having 
knowledge of the terms used in the description would be able to count the 
objects in question and arrive at the same count. Now there is another 
principle that may seem to be different from principle 1$,$ but it is actually 
equivalent to it.\parm
\hs
{\bf Principle 2.} (Ordinary (weak) induction.) 
Suppose that you take any nonempty $W 
\subset N_m.$ Then you show that a statement holds for \pars
(i) $m \in W.$\pars
(ii) Now if upon assuming that for a specific $n \in W,$ the statement we wish to 
establish holds $n,$ 
you can show that $n+1 \in W$$,$ (this is called the weak induction 
hypothesis) \pars
(iii) then you can declare that $W = N_m.$ (We say that the result follows by 
induction.) \pars                      
\hm
The difference in these two principles is located in part (ii). What you 
assume holds seems to be different. In principle 1$,$ we seem to require a 
stronger assumption than in principle 2. The next result shows that the two 
principles are equivalent. \parm
{\bf Theorem on the equivalence of the two principles 1 and 2.} {\sl Relative 
to the natural numbers and the subset $N_m$$,$ principle 1 holds if and only if principle 2 
holds.}\pars
Proof. We first show that principle 2 is equivalent to the fact that the 
simple ordering $<$ on $N_m$ is a well-ordering.  Assume principle 2 for 
$N_m$. Let nonempty $W \subset N_m.$ Suppose that $W$ does not have a first 
element with respect to the ordering $<.$  Then $W\not= N_m$ since $m$ is the 
first element in $N_m.$ Thus $m \notin W.$ We now define in terms of the ordering $<$ a relation 
$R$ where the second coordinate is the set $W$. We let $x\ R\ W$ if and 
only if for $x \in \nat,$ $x <y$ for each $y \in W.$ Let $W_1 = \{x\mid x \in 
N_m,\ x\ R\ W\}.$  From above$,$  we have since $W \subset N_m,$  that $m \in W_1.$ 
Also $W_1 \cap W = \emptyset.$ For if $a \in W_1 \cap W,$ then  $a < a;$ a 
contradiction. Assume that $p \in W_1$ and $q < p.$ Since $p <n$ for each $n 
\in W,$ then $q <n$ for each $n \in W.$ Hence $q \in W_1.$ Hence$,$ the 
nonempty natural number interval $[m,q] = \{x \mid m\leq x\leq q,\ x\in \nat 
\}\subset W_1.$ We now show that $q + 1 \in W_1.$ Assume that $q + 1 \notin 
W_1.$ Then there is some $y\in W$ such that $y\leq q+1$ from the definition of 
$W_1$. If $y \not= q+1,$ then $y \in [m,q]$ yields that $y \notin W.$ From 
this contradiction$,$ we have that $y =q+1$ and no $x \leq q < q+1$ is a member 
of $W$. Consequently$,$ $q + 1$ is a first element of $W$; a contradiction. Application of 
principle 2 implies that $W_1 = N_m.$ Thus yields the contradiction that 
$W = \emptyset.$ \parm
Now assume that $<$ is a well-ordering. Let nonempty $W \subset N_m.$ Assume 
that $m \in W$ and if arbitrary $n\in N_m$$,$ then $n+1 \in W,$ BUT $W \not= 
N_m.$ Consider $W_1 = N_m -W = \{x\mid x \in N_m,\  x\notin W\}.$ From our 
assumption$,$ $W_1 \not=\emptyset.$ By the well-ordering of $<,$  there exists 
in $W_1$ a first element $w_1.$ Since $m \notin W_1,$ it follows that 
$m \leq x_1 -1 \notin W_1.$ Thus $x_1 - 1 \in W.$ From our assumption$,$ 
$x_1 - 1 + 1 =x_1\in W.$ But by definition $W_1 \cap W = \emptyset.$ This 
contradiction yields the result. The fact that $<,$  with respect to the 
natural numbers$,$ is a well-ordering is equivalent to principle 2.\parm
Now we show that the well-ordering $<,$ restricted to any $N_m,$ implies  
principle 1 for $N_m.$  Let nonempty $W\subset N_m.$ Assume that $W \not= 
N_m.$ Assume that $m \in W$ and that if $x \in W$ and  $m\leq x \leq n,\ 
x\in N_m,$  then $n+1 \in W,$ but $W \not= N_m.$ Consider $W_1 = N_m - W_1.$ 
By well-ordering$,$ $W_1$ contains a first member $w_1.$ From the above 
assumption$,$ $w_1 \not= m.$ Hence$,$ we can express $w_1$ as $w_1 = n+1,$ for some
$n \in \nat.$ Since $w_1$ is a first element$,$ then each $x \in N_m$ such 
that $m\leq x\leq n,$ has the property that $x \in W.$ From our principle 1 
assumption$,$ this implies that $w_1 = n+1 \in W.$ This contradicts the 
definition of $W_1$ since $W_1  \cap W = \emptyset.$ The result follows.\parm
We now complete this proof by  showing that principle 1 and principle 
2 are equivalent. Given $N_m$ and assume principle 2. Then $N_m$ is 
well-ordered by $<.$ From above principle 1 holds.\pars
Now assume principle 1 holds. Let nonempty $W \subset N_m$ and suppose that 
$m \in W$ and if $n \in W,$ 
then $n+1 \in W.$ Principle 1 states that if $m \in W,$ and assuming that for 
each $r\in \nat$ such that $m\leq r\leq n,$ it can be established that $n+1 
\in W,$ then $W =N_m.$ However$,$ we are given that if $n \in W,$ then $n+1 \in 
W$ and $m\leq n\leq n.$ So$,$ trivially$,$ there is such an $r=n$ and
this yields 
that $n+1 \in W.$ Thus principle 1 implies that $W =N_m.$ \qed
There are two places that we use these equivalent induction processes. 
The first is called {{\it definition}} by induction.\parm
\hs
{\bf Principle 3} (Definition by induction) Consider a construction based upon 
the natural numbers $N_m.$ \pars
(i) Suppose that we describe a process for the case 
where $n = m$ (i.e. for step $m.$) \pars
(ii) Suppose that we assume that we have described a process for each $n$, 
where $n\geq m.$  (i.e. each step $n$.) \pars
(iii) Now use the $n$ notation and describe a fixed set of rules for the 
construction of the entity for the $n+1$ step. (The induction step.)\pars
Then you have described a process that obtains each step $n \in N_m.$\pars
\hm
{\bf Theorem on principle 3.} {\sl Principle 3 holds.}\pars
Proof. Suppose that you followed the rules in (i)$,$ (ii)$,$ (iii). Let $W \subset 
N_m$ be the set of all natural numbers great than or equal to $m$ for which 
the construction has been described. From (i)$,$ $m \in W.$ From (ii)$,$ we may 
assume that you have constructed step $n \in W.$ From (iii)$,$ you have 
described step $n+1$ from step $n$. Thus $n+1 \in W.$ Hence by principle 2$,$ 
$W = N_m.$\qed
Notice that the concept of what is an acceptable description for a 
construction by induction depends upon whether the description is so clear
that all individuals will obtain the same constructed object.\parm
{\bf Example 1.} (Definition by induction)   Let $m = 1$ and let $b$ be a 
positive real number. \pars
(1) Define $b^0 = 1$ and $b^1 = b$ each being a real number.\pars
(2) Assume that for arbitrary $n \in N_1,$ $b^n$ has been defined and is a 
real number.\pars
(3) Define $b^{n+1} = b\cdot (b^n),$ where $\cdot$ means the multiplication of 
real numbers. (Note since $b$ is a real number and by assumption $b^n$ is a 
real number$,$ then $b\cdot (b^n)$ is a real number.)\pars
(4) Hence$,$ by induction$,$ $b^n$ has been defined and is a real number 
for all $n\in N_1$ (and separately for $n = 0.$) \parm
In definition 2.2.3 that appears below for the language $L,$ we use the concept 
of definition by induction to obtain a definition for each language level 
$L_n.$ Although there are numerous examples of proof by induction within the 
main part of this text$,$ here is one more example.\parm
{\bf Example 2.} (Proof by induction) We show that if a natural number $n$ 
is greater than or equal to 2$,$ then there exists a prime number $p$ that is a 
factor of $n.$\pars
Proof. Let $m = 2$ and let $W \subset N_m$ such that each member $n \in W$ has 
a prime factor. Since $2$ is a prime factor of itself$,$ then $2 \in W.$ Assume 
that for arbitrary $n \in W$ and each $r\in N_2$ such that $2 \leq r\leq n,$ we 
know that $r \in W.$ Now consider $n+1.$ If $n+1$ has no nontrivial factor $b$
(i.e. not equal to 1 or $b$) such that $2\leq b< n+1,$ then $n+1$ is a 
prime number and$,$ hence$,$ contains a prime factor.
(Note: any such nontrivial factor would be less than $n+1.$) If $n+1$ has a 
nontrivial factor $b$$,$ then $2 \leq b < n+1.$ Thus for this case$,$ 
$2 \leq b \leq n.$ By the induction hypotheses$,$ $b$ has a prime factor $p$.
Hence $p$ being a prime factor of $b$ is also a prime factor of $n+1.$ Thus 
$n+1$ has a prime factor. Consequently$,$ $W = N_n$ by induction. \parm 
\hrule\smallskip
{\bf Definition 2.2.3.} The construction by induction of the propositional 
language $L.$ \pars 
(1) Let ${\cal A} =\{P,\ Q,\ R,\ S\}\cup (\cup \{P_i, \ Q_i, \ R_i,
\ S_i\mid i \in \nat -\{0\} \}.$ The set $\cal A$ is called a set of atoms. \pars
(2) Let $\emptyset \not=L_0 \subset {\cal A}.$\pars 
(3) Suppose that $L_n$ has been defined. Let $L_{n+1} =\{(\neg A)\mid A \in 
L_n\} \cup \{(A \land B)\mid A, b \in L_n\}\cup \{(A \lor B)\mid A, B \in 
L_n\} \cup \{(A \to B)\mid A, B \in L_n\} \cup\{(A \iff B)\mid A,B
 \in L_n\}\cup L_n.$ \par 
(4) Now let $L = \bigcup \{L_n \mid 0 \not= n \in \nat \}.$ \pars
(5) Note: We are using set theoretic notation. If one wants to formalize the 
above intuitive ideas$,$ the easiest way is to use a set theory with 
atoms where each member of $\cal A$ is an atom. 
Then consider various $L_{n+1}$ level  n-ary relations$,$ such as 
the $\land_{n+1}$ 
relation 
defined on the $L_{n}$ and for the other logical 
connectives. Then all nary $\land_n$ relations have the same properties and the 
same interpretation. Because$,$ they have 
the same properties and interpretation$,$ there is nothing gained by 
formalizing this construction. \pars
(6) These formulas can also be defined in terms of the class concept$,$ sequences 
of atoms$,$ trees$,$ closure concepts and a lot more stuff. Or$,$ just keep them 
intuitive in character. \pars
\hrule
\medskip
The next theorem holds obviously due to the inductive definition but I present 
it anyway.\parm
{\bf Theorem on uniqueness of size.} {\sl For any $A \in L$ there exists a 
natural number $n \in \nat$ such that $A\in L_n$ and if $m\in \nat$ and $m<n,$
then $A \not\in L_m.$}\pars
Proof. Suppose the $A \in L.$ Then from the definition of $L$ there exists 
some $n$
such that $A \in L_n.$ Let $K = \{k \mid A \in L_k\}.$ Then $\emptyset \not= K 
\subset \nat.$
Hence$,$ $K$ has a smallest member which by definition would be the size.\qed
\medskip
{\bf Theorem of the existence of a finite set of atoms  
that are contained in any formula $A$.} {\sl Let $A  \in L.$ Then there 
exists a finite set $A_1$ of atoms that contains all the atoms in $A.$}\pars
Proof.  We use strong induction.\pars
Let $A \in L.$ Then there exists a unique $n$ such that ${\rm size}(A) = n.$ 
\pars (1) Let ${\rm size}(A) = 0.$ Then $A \in L_0$ and$,$ hence$,$ $A$ is a single 
atom. The set that contains this single atom is a finite set of  
atoms in $A$. \pars
(2) Suppose that there exists a finite set that contains all the atoms for a 
formula of size  
$r \leq n$. Let 
${\rm size}(A) = n+1.$  Then from the definition of the levels either
(i) $A = \neg B,\ A = B\lor C,\ A = B\land C,\ A = B \to C,$ or $A = B \iff 
C,$ where size of 
$B\leq n$ and $C \leq n.$ Hence$,$ from the induction hypothesis$,$ 
there is a finite set of atoms $A_1$  that 
are contained in $B,$ and a finite set of  atoms $A_2$ that are contained 
in $C.$ Hence$,$ there  is a finite set that contains the atoms in $A.$ Thus by 
induction$,$ given any $A\in L,$ then there exists a finite set of atoms that 
contains all the atoms in $A.$.\qed
\medskip
{\bf Theorem on existence of a unique assignment dependent valuation 
function.} {\sl There exists a 
 function $v\colon L \to \{T,F\}$ such that \pars
(a) if $A\in L_0,$ then  $v(A) =F$ or $T$ not both.\pars
(b) If $ A = \neg B$$,$ then $v(A) = F$ if $v(B) = T,$ or $v(A) = T$ if $v(B) = 
F.$ \pars
(c) If $A = B\lor C,$ then  $v(A) = F$ if and only if $v(B)=v(C) = F.$ Otherwise $v(A) = 
T.$\pars
(d) If $A = B\land C,$ then $v(A) = T$ if and only if $v(B) = v(C) = T$. Otherwise $v(A) 
=F.$\pars
(e) If $A = B\to C,$ then $v(A) = F$ if and only if $v(B) = T,\ v(C) = F.$ Otherwise
$v(A) = T.$\pars
(f) If $A = B \iff C,$ then $v(A) = T \iff v(B) = v(C).$ Otherwise $v(A) = 
F.$\pars
The function $v\colon L\to \{T,F \}$ is unique in the sense that if any other 
function $f\colon L \to \{T,F\}$ has these properties and $f(A) = v(A)$ for each 
$A \in L_0$ then $f(A) = v(A)$ for each $A \in L.$}\pars
Proof. (By induction). We show that for every $n \in\nat$ there exists a 
function $v_n\colon L_n \to \{T,F\}$ that satisfies (a) --- (f) above. 
Let $n = 0$ and $A \in L_0.$ Then define $v_0$ by 
letting $v_0(A) = T$ or $v_0(A) =F$ not both. Then $v_0$ is a function. 
Now the other properties 
(b) --- (f) hold vacuously. 
Suppose there is a function $v_n \ (n \not= 0)$ defined that satisfies the 
properties. Define 
$v_{n+1}$ as follows: $v_{n+1} | L(n) = v_n.$ Now for $A\in L_{n+1}-L_n,$
then $A = \neg B,\ B\lor C,\ B\land C, b \to C, \ B \iff C.$   Thus define $v_{n+1}(A)$ in accordance with the 
requirements of (b)$,$ (c)$,$ (d)$,$ (e)$,$ (f) of the above theorem.  
This gives a function since $v_n,$ by the induction hypothesis$,$ is a function 
on $L_n$ and  $B,\ C\in L_n$ are unique members of  $L_n.$ \pars
Now to show that for all $n \in \nat,$ if $f_n \colon L \to \{T,F\}$ has the 
property that if $f_0 = v_0,$ and (a) --- (f) hold for $f$$,$ then $f_n =v_n.$
Suppose $n = 0$. Since $f_0=v_0$ property holds for $n = 0.$ Now suppose that
$v_n\ (n \not= 0)$ is unique (hence$,$ (a) --- (f) hold for $v_n$) Consider$,$ 
$f_{n+1} \colon L \to
\{T,F\}$ and (a) --- (f) hold for $f_{n+1}.$  Now $f_{n+1}| L_n$ satisfies 
(a) --- (f); hence $f_{n+1}| L_n=v_n.$ Now looking at $L_{n+1}-L_n$ and the 
fact that $f$ satisfies (b) --- (f)$,$ it follows that for each $A \in L_{n+1}-
L_n,\ f_{n+1}(A) = v_{n+1}(A).$ Consequently$,$ this part holds by 
induction.\pars
Now let $v = \bigcup \{v_n\mid n \in \nat \}.$ $v$ is a function since for 
every $n \leq m,\ v_n \subset v_m.$ We show the $v$ satisfies (a) --- (f).
Obviously $v$ satisfies (a) since $v|L_0 = v_0.$ Now if $A = \neg B,$  then 
there  exists $n \in \nat$ such that $B \in L_n$ and $A \in L_{n+1}$ from the 
existence of size of $A.$ Then $v(A) = v_{n+1}(A) = T,$ if $v_n(B) = F =v(B)$ or
$v_{n+1}(A)= F,$ if  $v_n(B) =T =v(B).$ Hence (b) holds. In like manner$,$ it 
follows that (c) -- (f) hold and the proof is complete.\qed   
It's obvious how to construct the assignment and truth-table concepts from
the above theorem. If $A$ contains a certain set of atoms$,$ then restricting
$v$ to this set of atoms gives an assignment $\ass.$ Conversely$,$ all 
assignments are generated by such a restriction. The remaining part of this 
theorem is but the truth-table valuation process restricted to all the formula 
that can be constructed from this set of atoms starting with the set 
$L_0.$\parm
\hs
{\bf Definition 2.11.1} The inductive definition for the sublanguage $L\p$ is 
the exact same as in the case for 2.2.3.\pars
\hm
\bigskip\bigskip
\centerline{Chapter 3}\parm
\noindent {\bf Important Note: All of the results established prior to 
theorem 3.8.4 hold for ANY predicate language.}\parm
\hs
{\bf Definitions 3.1 and 3.6.1} These are obtained in the exact same manner as 
is 2.2.2.\pars
\hm
{\bf Theorem on the existence and uniqueness of the process $\mod$ and 
$\not\mod$ described in 3.3.3 on Structure Valuation.} {\sl Given a structure 
$\str =\langle D,{P_i^n}\p \rangle,$ where the interpretation $I$ is a one-to-one  
correspondence from $N \subset C$ onto $D$ 
and every n-place predicate $P_i^n$ to an n-place relation $R_i^n.$ 
There exists a
$v\colon Pd \to \{\mod,\not\mod\}$ such that\pars
(a) for each $i \in \nat$ and $n \in \nat -\{0\},$ $v(P_i^n(c_1,\ldots,c_n) 
= \mod$ if and only if $(c_1\p,\ldots,c_n\p)\in R_i^n,$ where for any $c_i 
\in N,$ $I(c_i) = c_i^\prime.$  
 \pars
(b) If $A\to B \in Pd,$ then $v(A\to B)=\not\mod$ if and only if $v(A)=\mod$ and 
$v(B)=\not\mod.$ In all other cases$,$ $v(A\to B)=\mod.$\pars
(c) If $A \iff B \in Pd,$ then $v(A \iff B)=\mod$ if and only if $v(A)=\mod$ 
and $v(B)=\mod,$ or $v(A)=\not\mod$ and  $v(B)=\not\mod.$ \pars
(d) $v(A\lor B)=\mod$ if and only if $v(A)=\mod$ or $ 
v(B)=\mod.$\pars
(e) $v(A\land B)=\mod$ if and only if $v(A)=\mod$ and $ 
v(B)=\mod.$\pars
(f) $v(\neg A)=\mod$ if and only if $v(A)=\not\mod .$\parm
In what follows$,$ the constant $d$ is a general constant  and 
corresponds  to a general member $d\p$ of the set $D.$ 
Any constants that appear in the original 
predicates have been assigned FIXED members of $D$ by $I$ and never change their
corresponding elements throughout this valuation for a given structure.\parm
(g) For each sentence $C=\forall x A\in Pd,$ $v(\forall x A)=\mod$ 
if and only 
if for every $d\p \in D$ it follows that $v(S^x_d A])=\mod.$ Otherwise$,$
$v(\forall x A)=\not\mod.$ \pars
(h) For each sentence $C=\exists xA\in Pd,$ $v(\exists x A) =\mod$    
if and only 
if there is some $d\p \in D$ such that $v(S^x_d A])=\mod.$ Otherwise$,$ 
$v(\exists x A)=\not\mod .$ \pars
The function $v$ is unique in the sense that if any other function $f\colon Pd 
\to \{\mod,\not\mod\}$ such that $f =v$ for the statements in (a)$,$ then $f=v,$ in  
general.}\pars
Proof. The proof is by induction in the language levels $m.$ Let $m = 0.$ 
Define $v_0\colon Pd_0 \to \{\mod,\not\mod\}$ by condition (a) for all the 
predicates. Then $v_0$ satisfies (b) -- (h) vacuously.\pars
Suppose that $v_m\colon Pd_m\to \{\mod,\not\mod\}$ exists and satisfies (a) --- 
(h). Define $v_{m+1} \colon Pd_{m+1} \to \{\mod\not\mod\}$ as follows:
$v_{m+1}|Pd_m = v_m.$ For $F\in Pd_{m+1}-Pd_m,$ then $F =\neg A,\ A \land B,\ 
A \lor B,\ A \iff B,\forall x A,\ \exists x A,$ where $A,\ B \in Pd_m.$ Now 
define for the specific $F$ listed$,$ the function by the appropriate conditions 
listed in (b) --- (h). We note that $A,\ B$ are unique and this defines a 
function.\pars
We show that each $v_m$ is a unique function $f\colon Pd_m \to 
\{\mod,\not\mod\}$ satisfying (a) --- (h) by induction. If $f$ is another such 
function$,$ then letting $m = 0$ condition (a) implies that $f = v_0.$ Assume 
that $v_m$ is unique. Let $f\colon Pd_{m+1} \to \{\mod,\not\mod\}$ and $f$ 
satisfies (a) --- (h). Then $f|Pd_m =g$ satisfies 
(a) --- (h). Therefore$,$ $g = v_n.$ Now (b) --- (h)$,$ yields that $f = 
v_{m+1}.$\pars
The remainder of this proof follows in the exact same manner as the end of the 
proof of the {\bf Theorem on existence of a unique assignment dependent valuation 
function}.\qed\parm
\hs
NOTE: In the remaining portion of this appendix, it will be assumed that our language contains a non-empty countable set of constants $\cal C$.\parm\hm
We need for the proof of theorem 3.8.5 another conclusion. A 
set of sentences $\Gamma$ is {{\it universal}} for a language $Pd\p$ if 
 and only if $\forall x B \in \Gamma,$  whenever $S^x_dB\in \Gamma$ for all 
$d\in {\cal C}.$ \parm
{\bf Theorem 3.8.4} {\sl If consistent $\Gamma \subset \c S,$ then there 
is a language $Pd^{\prime\prime}$ that contains all of the symbols of 
$Pd\p,$ but with an 
additional set of constants and only constants adjoined$,$ and a set of 
sentences $\Gamma^{\prime\prime} \subset Pd^{\prime\prime}$ such that 
$\Gamma \subset \Gamma^{\prime\prime}$ 
and $\Gamma^{\prime\prime}$ is consistent$,$ negation complete$,$ and universal.}\pars
Proof. First$,$ we extend $\cal C$ by adjoining a new 
denumerable set of constants $\{b_0,\ldots \}$ to $\cal C$ giving a new language 
$Pd^{\prime\prime}.$ [For other languages$,$ the set of new constants may need 
to be a ``larger'' set than this.] 
This means that the set of sentences ${\c S}\p$ for  
$Pd^{\prime\prime}$ is denumerable and we can consider them as enumerated into 
an 
infinite sequence $S_1,S_2,\ldots $ and these are all of the members of 
${\c S}\p$. We now begin an inductive definition for an extension of 
$\Gamma_n$ 
by adjoining a finite set of sentences from  $S_1,S_2,\ldots $ which could 
mean that only finitely many members of $\{b_0,\ldots \}$ would appear in an
any $\Gamma_n.$ \pars
First$,$ for $n = 0,$ let $\Gamma_0 = \Gamma.$ Suppose that $\Gamma_n$ has been 
defined. We now define $\Gamma_{n+1}.$ \pars
\indent\indent (a) If $\Gamma_n \cup \{S_{n+1}\}$ is consistent$,$ then let
$\Gamma_{n+1} =  \Gamma_n \cup \{S_{n+1}\}$. \pars
\indent\indent (b) If $\Gamma_n \cup \{S_{n+1}\}$ is inconsistent and 
$S_{n+1}$ is NOT of the form $\forall xB,$ then let $\Gamma_{n+1} = \Gamma 
\cup \{\neg S_{n+1}\}.$ \pars
\indent\indent (c) If $\Gamma_n \cup \{S_{n+1}\}$ is inconsistent and 
$S_{n+1}$ is of the form $\forall xB,$ then let $\Gamma_{n+1}=
\Gamma \cup \{\neg S_{n+1}, \neg S^x_bB]\},$ where $b$ is the first constant 
in $\{b_0,\ldots \}$ that does not appear in $\Gamma_n.$\parm
We first show that for each $n$ the set of sentences $\Gamma_n$ is 
consistent.  Obviously$,$ for $n = 0,$ the result follows from the hypothesis 
that $\Gamma = \Gamma_0$ is consistent. Now assume that $\Gamma_n$ is 
consistent.\pars
\indent\indent (a)' Suppose that $\Gamma_{n+1}$ is obtained from case (a). 
Then $\Gamma_{n+1}$ is consistent.\pars
\indent\indent (b)' Suppose that $\Gamma_{n+1}$ is obtained from case (b). 
Then $\Gamma_n \cup \{S_{n+1}\}$ is inconsistent. Hence$,$
from corollary 3.8.3.1$,$ $\Gamma_n \vd \neg S_{n+1}.$ Now $\Gamma_{n+1} =
\Gamma_n \cup \{\neg S_{n+1}\}.$ Suppose that $\Gamma_n \cup \{\neg S_{n+1}\}$ is 
inconsistent. Then corollary 3.8.3.1 yields that $\Gamma_n \vd S_{n+1}.$ 
Hence$,$ $\Gamma_n$ is inconsistent. This contradiction yields that 
$\Gamma_{n+1}$ is consistent.\pars
\indent\indent (c)' Again $\Gamma_n \cup \{S_{n+1}\}$ is inconsistent yields 
that $\Gamma_n \vd \neg S_{n+1}.$ Now also suppose that  $\Gamma_{n+1}=
\Gamma \cup \{\neg S_{n+1}, \neg S^x_bB]\},$ where $b$ is the first member of 
new constants that does not appear in $\Gamma_n$ and assume that this is an
inconsistent collection of sentences. Then for some $C \in Pd^{\prime\prime}$ 
we have that  $\Gamma_n \cup \{\neg S_{n+1}, \neg S^x_bB]\}\vd C$ and $
\Gamma_n \cup \{\neg S_{n+1}, \neg S^x_bB]\}\vd \neg C.$  Then by the deduction 
theorem 3.7.4$,$ $\Gamma_n \cup \{ \neg S^x_bB]\}\vd (\neg S_{n+1})\to 
C$ and $
\Gamma_n \cup \{\neg S^x_bB]\}\vd (\neg S_{n+1})\to (\neg C).$  
Then by adjoining the proof that  $\Gamma_n \vd \neg S_{n+1}$ and two MP steps 
we have that  $\Gamma_n \cup \{\neg S^x_bB]\}$ is inconsistent. Thus by 
corollary 3.8.1$,$ $\Gamma_n \vd \neg S^x_bB].$ Now the constant 
$b$ does not occur anywhere in $\Gamma_n.$ Thus in the last demonstration we 
may substitute for $b$ some variable $y$ that does not appear anywhere in the 
demonstration for each occurrence of $b.$ This yields a demonstration that 
$\Gamma_n \vd S^x_yB.$ By Generalization$,$ $\Gamma_n \vd \forall 
y(S^x_yB]).$ By corollary 3.6.3.1$,$ we have that $\Gamma_n \vd \forall x B = 
S_{n+1}.$ This contradicts the consistency of $\Gamma_n.$ Since this is the 
last possible case$,$ then $\Gamma_{n+1}$ is consistent. Thus by induction for 
all $n,$ $\Gamma_n$ is consistent. \parm
We note that by definition $\Gamma=\Gamma_0 \subset \Gamma_1 \subset \cdots 
\subset \Gamma_n \subset \cdots.$ We now define $\Gamma^{\prime\prime} 
=\bigcup \{\Gamma_n\mid n\in \nat\}.$ We need to show that 
$\Gamma^{\prime\prime}$ is consistent$,$ negation complete and universal for 
$Pd^{\prime\prime}.$  \pars
(1)  Suppose that the set of 
sentences $\Gamma^{\prime\prime}$ is inconsistent. Then there is finite subset $F$ of 
$\Gamma^{\prime\prime}$  that is inconsistent. But$,$ from our remark above$,$ 
there is some $m \in \nat$ such that $F \subset \Gamma_m.$ This implies that 
$\Gamma_m$ is inconsistent. From this contradiction$,$ we have that 
$\Gamma^{\prime\prime}.$ [In fact$,$ its maximally consistent$,$ in that 
adjoining any sentence to  $\Gamma^{\prime\prime}$  that is not in 
$\Gamma^{\prime\prime},$ we get an inconsistent set of sentences.]\pars
(2) The set $\Gamma^{\prime\prime}$ is negation complete. Let $A$ be any 
sentence in $Pd^{\prime\prime}.$ Then $A$ is one of the $S_{n+1},$ where 
$n \in \nat.$ From definition (a)$,$ (b)$,$ (c)$,$ 
either $S_{n+1}\in \Gamma_{n+1}$ or $\neg S_{n+1} \in 
\Gamma_{n+1}.$  Thus $\Gamma^{\prime\prime}$ is 
negation complete.\pars
(3) We now show that $\Gamma^{\prime\prime}$ is universal. Let $\forall x B$ 
be a sentence in $Pd^{\prime\prime}$ such that $S^x_cB] \in 
\Gamma^{\prime\prime}$ for each $c \in {\cal C} \cup\{b_0,\ldots \}$. Suppose that 
$\forall x B \notin \Gamma^{\prime\prime}.$ We know that $\forall x B = 
S_{n+1}$ for some $n \in \nat.$ By case (a)$,$ $\Gamma_n \cup \{\forall x B\}$ 
is inconsistent$,$ by negation completeness$,$ $\neg \forall xB 
\in \Gamma^{\prime\prime}.$   
Now case (c) applies and $\Gamma_{n+1}=\Gamma_n \cup 
\{\neg S_{n+1}, \neg S^x_bB]\}\subset \Gamma^{\prime\prime}.$ This implies 
that $\neg S^x_bB]\in \Gamma^{\prime\prime}.$  But our assumption was that 
$S^x_cB]\in \Gamma^{\prime\prime}$ for all constants in $Pd^{\prime\prime}.$ 
Since $b$ is one of these constants$,$ we have a contradiction. Thus $\forall xB 
\in \Gamma^{\prime\prime}.$\qed
Prior to establishing our major theorem 3.8.5, we have the following Lemma 
and the method to construct the model  we seek. Let the domain of our structure $D= {\cal C}$ for a specific $Pd'$.   Let $\Gamma $ be any non-empty set 
of sentences from the language $Pd^{\prime}.$ 
For every n-place 
predicate $P_i^n$ in $Pd'$, we define the n-place relation $R_i^n$ by 
$(c_1,\ldots,c_n) \in R_i^n$ if and only if $P(c_1,\ldots,c_n) \in \Gamma.$ we 
denote the structure obtained from this definition by the notation 
$\str (\Gamma).$\parm
{\bf Lemma 3.8} {\sl Let $\Gamma$ be a consistent$,$ negation 
complete and universal set of sentences from $Pd^{\prime}.$ Then
for any sentence $A \in Pd^{\prime},$ $\str (\Gamma)\mod A$ if and only if $A \in \Gamma.$}\pars
Proof. This is established by induction on the size of a formula.\pars
(a) For $ n = 0.$  Let sentence $A\in Pd_0^{\prime}.$ The result follows 
from the definition of $\str (\Gamma) .$\pars
(b) Suppose that theorem holds for $n$. Suppose that $A\in 
Pd_{n+1}^{\prime}.$ Assume that $A =B \to C.$ Then $B,\ C \in 
Pd^{\prime}_n.$ We may assume by induction that the lemma holds for
$B$ and $C$. Suppose that $A \notin \Gamma^{\prime}.$ From negation 
completeness$,$ $\neg A \in \Gamma.$ But$,$ in general$,$ $\neg(B\to 
C) \vd B$ and $\neg (B\to C) \vd \neg C.$ Hence$,$ $\Gamma \vd B$ 
and $\Gamma \vd \neg C.$ If $\neg B \in \Gamma,$ 
$\Gamma$ is inconsistent. Then$,$ from negation completeness$,$
it must be that $B \in \Gamma.$ For the same reasons$,$ $\neg C 
\in \Gamma$ and$,$ thus$,$ $C \notin \Gamma.$ 
From the induction hypothesis$,$ $\str (\Gamma) \mod B$ and 
$\str (\Gamma)\not\mod C.$ Thus $\str (\Gamma)\not\mod A.$ Conversely$,$ assume that$\str (\Gamma) 
\not\mod A.$ Then $\str (\Gamma)\mod B$ and $\str (\Gamma)\not\mod C.$ By the induction 
hypothesis$,$ $B\in \Gamma$ and $C \notin \Gamma.$
By negation completeness$,$ $\neg C\in \Gamma.$ From$,$ 
$B,\ \neg C \vd \neg (B\to C)$$,$ we have that $\Gamma \vd \neg (B\to C).$ Again by negation completeness and consistency$,$ it follows that 
$\neg (B \to C) \in \Gamma.$ Hence$,$ $\neg A \in 
\Gamma$ and $A \notin \Gamma.$ \pars
(c) Let the sentence be $A = \forall x B.$ Suppose that $\str (\Gamma) \mod A.$ Then
for each $d \in D,$ $\str (\Gamma) \mod S^x_dB.$ By the induction hypothesis$,$ 
since $S^x_dB\in Pd_n^{\prime},$ then $S^x_dB \in 
\Gamma.$ Since $\Gamma$ is universal$,$
$\forall xB=A \in \Gamma.$  Conversely$,$ let $A \in 
\Gamma.$ Now$,$ in general$,$ $A \vd S^x_dB$ implies that 
$\Gamma \vd S^x_dB,$  which implies again by negation completeness 
and consistency that $S^x_dB \in \Gamma$ for each $d \in D.$ By 
the induction hypothesis$,$ we have that $\str (\Gamma) \mod S^x_dB$ for each $d \in 
D.$ From the definition of $\mod,$ we have that $\str (\Gamma)\mod A.$ \pars
This completes the proof.\qed
\vfil\eject
\hs\hm
{\bf Theorem 3.8.5} {\sl If consistent $\Gamma \subset \c S$$,$ then there 
exists a structure $\str$ for $\Gamma$ such that $\Gamma \vd A$ if and only if 
$\str \mod \forall A.$}\parm
 Proof. 
Let the consistent set of sentences $\Gamma \subset Pd^\prime.$ 
Then the language 
$Pd^\prime$ can be extended to a language $Pd^{\prime\prime}$  and the set of 
sentences $\Gamma$ extended to a set of sentences $\Gamma^{\prime\prime}$ such 
that $\Gamma^{\prime\prime}$ is consistent$,$ negation complete$,$ and universal.
Assume that $\Gamma \vd A.$ Then $\Gamma \vd \forall A$. Consequently$,$
$\Gamma^{\prime\prime} \vd \forall A.$ Thus $\forall A \in 
\Gamma^{\prime\prime}.$ Now we use the structure $\str (\Gamma^{\prime\prime}).$ Then lemma 3.8 
states that $\str (\Gamma^{\prime\prime})\mod \forall A,$ where $\str (\Gamma^{\prime\prime})$ is considered as restricted to 
$Pd\p$ since $\forall A \in Pd^\prime.$ \pars
Now suppose that $\Gamma \not\vd A.$ Then by repeated application of $P_5$ and 
MP$,$ we have that $\Gamma \not\vd \forall A.$ [See theorem 3.6.4.] Thus 
$\Gamma \cup \{\neg(\forall A)\}$ is consistent by theorem 3.8.2. 
Consequently$,$
$\neg (\forall A) \in \Gamma^{\prime\prime}.$ Therefore$,$ $\str (\Gamma^{\prime\prime}) \mod 
\neg(\forall A)$ implies that $\str (\Gamma^{\prime\prime}) \not\mod \forall A.$ Again$,$ by 
restriction to the language $Pd^\prime,$ the converse follows.\qed
\hs\hm
\vfil\eject

{\centerline{Answers to Some of the Exercise Problems}\parm
\centerline{\sl Exercise 2.2$,$  Section 2.2}
\medskip
\noindent [1] (a) $A \notin L.$ A required `` )'' is missing. (b) $A \in L.$
(c) $A \in L.$ (d) $A \notin L.$ The error is the symbol $(P)$. (e) $A \notin 
L.$ Another parenthesis error. (f) $A \notin L$. The  symbol string `` )P)'' 
 is in error. (g) $A \in L.$ (h) $A \notin L.$ Parenthesis error. 
(i) $A \notin L.$ Parenthesis error.
(j) $A \notin L.$ Parenthesis error. \pars

\noindent [3] (a) $(P \lor Q)$. (b) $(P\to (Q \lor R)).$ (c) $(P\land (Q \lor 
R)).$ (d) $(Q \lor R).$ (e) $((P\land Q)\lor R)$. (f) $\neg (P\lor Q).$ (g) 
$((P \iff Q)\land ((\neg Q)\to (R \land S))).$ (h) $((P\lor Q) \to R).$ (i)
$(((P\land Q) \to R) \land ((\neg P)\to (\neg R))).$\pars
\noindent [4] (a) If it is nice$,$ then it is not the case that it is hot and it 
is cold.\par
(b) It is small if and only if it is nice.\par
(c) It is small$,$ and it is nice or it is hot.\par
(d) If it is small$,$ then it is hot; or it is mice. \par
(e) It is nice if and only if; it is hot and it is not cold$,$ or it is 
small. [It may be difficult to express this thought nonambiguously in a single 
sentence unless this  ``strange'' punctuation is used. The ``;' indicates a 
degree of separation greater than a comma but less then a period.]  
\par
(f) If it is small$,$ then it is hot; or it is nice.\par 
\bigskip
\centerline{\sl Exercise 2.3$,$  Section 2.3} 
\medskip
\noindent [1] (a) The number of (not necessarily distinct connectives) =
the number of common pairs. \par
(b) The number of subformula = the number of common pairs.\pars
\noindent [2] (No problem$,$ sorry!)\pars
\noindent [3] (A) has common pairs (a$,$h)$,$ (b$,$e)$,$ (c$,$d)$,$ (f$,$g).\par
(B) has common pairs (a$,$m)$,$ (b$,$g)$,$ (c$,$f)$,$ (d$,$e)$,$ (h$,$k)$,$ (i$,$j). \par
(C) has common pairs (a$,$q)$,$ (b$,$p)$,$ (d$,$g)$,$ (e$,$f)$,$ (i$,$o)$,$ (j$,$k)$,$ (m$,$n)$,$ 
(c$,$h).\parm

\centerline{\sl Exercise 2.4$,$  Section 2.4}
\medskip
\noindent [1]. First$,$ we the assignment $\underline{a} = (T,F,F,T) \iff 
(P,Q,R,S).$ \pars
(a) $v((R\to(S\lor P)),\ass )= (F\to (T\lor T)) = (F \to T) = F.$ \par
(b) $v(((P\lor R)\iff (R \land(\neg S))),\ass )= ((T\lor F)\iff (F \land(\neg 
T))= (T\iff (F\land F))= F.$ \par
(c) $v((S \iff (P\to ((\neg P)\lor S))),\ass )=(T \iff (T\to ((\neg T)\lor 
T))= (T\iff (T \to T)) = (T \iff T) = T.$ \par
(d) $v((((\neg S)\lor Q)\to(P\iff S)),\ass )=(((\neg T)\lor F)\to(T\iff 
T)=(Fto T) = T.$\par
(e) $v((((P\lor (\neg Q))\lor R)\to ((\neg S)\land S)),\ass )=(((T\lor (\neg 
F))\lor F)\to ((\neg T)\land T)= (((T \lor T)\lor F)\to F) = (T \to F) = F.$
\pars

\noindent 3. 
(a) $(P \to Q)\to R,\ v(R) = T \Rightarrow ((P \to Q)\to T) = T$ always.\par
(b) $P\land (Q \to R),\ v(Q \to R) =F\Rightarrow (P\land F)=F$ always.\par
(c) $(P \to Q)\to ((\neg Q) \to (\neg P)),\ v(Q) = T\Rightarrow (P \to T)\to 
(F \to (\neg P))=T\to T=T$ always.\par
(d) $(R\to Q)\iff Q,\ v(R) = T\Rightarrow (T\to Q)\iff Q= T$  always.\par
(e) $(P \to Q) \to R,\ v(Q) = F\Rightarrow (P \to F) \to R= T \ {\rm or}\ 
F.$\par

(f) $(P\lor (\neg P))\to R,\ v(R) = F\Rightarrow (P\lor (\neg P))\to F= T\to 
F = F$ always.\parm
\medskip\vfil\eject
\centerline{\sl Assignment 4 Section 2.5}
\medskip
\noindent [1]. 
 

\newbox\medstrutbox
\setbox\medstrutbox=\hbox{\vrule height14.5pt depth9.5pt width0pt}
\def\medstrut{\relax\ifmmode\copy\medstrutbox\else\unhcopy\medstrutbox\fi}

\medskip

\newdimen\leftmargin
\leftmargin=0.0truein
\newdimen\widesize
\widesize=6truein
\advance\widesize by \leftmargin       
\hfil\vbox{\tabskip=0pt\offinterlineskip
\def\tablerule{\noalign{\hrule}}

\halign to \widesize{\medstrut\vrule#\tabskip=0pt plus2truein

&\hfil\quad#\quad\hfil&\vrule#
&\hfil\quad#\quad\hfil&\vrule#
&\hfil\quad#\quad\hfil&\vrule#
&\hfil\quad#\quad\hfil&\vrule#
&\hfil\quad#\quad\hfil&\vrule#
&\hfil\quad#\quad\hfil&\vrule#
&\hfil\quad#\quad\hfil&\vrule#
&\hfil\quad#\quad\hfil&\vrule#
&\hfil\quad#\quad\hfil&\vrule#

\tabskip=0pt\cr\tablerule

&P&&Q&&R&&$Q\to P$&&(1)&&$Q\to R$&&$P\to (Q\to R)$&&$P \to R$&&$P\to Q$&\cr\tablerule
&T&&T&&T&&T&&T&&T&&T&&T&&T&\cr\tablerule
&T&&T&&F&&T&&T&&F&&F&&F&&T&\cr\tablerule
&T&&F&&T&&T&&T&&T&&T&&T&&F&\cr\tablerule
&T&&F&&F&&T&&T&&T&&T&&F&&F&\cr\tablerule
&F&&T&&T&&F&&T&&T&&T&&T&&T&\cr\tablerule
&F&&T&&F&&F&&T&&F&&T&&T&&T&\cr\tablerule
&F&&F&&T&&T&&T&&T&&T&&T&&T&\cr\tablerule
&F&&F&&F&&T&&T&&T&&T&&T&&T&\cr\tablerule}}\pars

\newbox\medstrutbox
\setbox\medstrutbox=\hbox{\vrule height14.5pt depth9.5pt width0pt}
\def\medstrut{\relax\ifmmode\copy\medstrutbox\else\unhcopy\medstrutbox\fi}

\medskip

\newdimen\leftmargin
\leftmargin=0.0truein
\newdimen\widesize
\widesize=6truein
\advance\widesize by \leftmargin       
\hfil\vbox{\tabskip=0pt\offinterlineskip
\def\tablerule{\noalign{\hrule}}

\halign to \widesize{\medstrut\vrule#\tabskip=0pt plus2truein

&\hfil\quad#\quad\hfil&\vrule#
&\hfil\quad#\quad\hfil&\vrule#
&\hfil\quad#\quad\hfil&\vrule#
&\hfil\quad#\quad\hfil&\vrule#
&\hfil\quad#\quad\hfil&\vrule#
&\hfil\quad#\quad\hfil&\vrule#

\tabskip=0pt\cr\tablerule

&$(P\to 
Q)\to (P\to R)$&&(2)&&$\neg P$&&$\neg Q$&&$(\neg P)\to (\neg 
Q)$&&(21)&\cr\tablerule
&T&&T&&F&&F&&T&&T&\cr\tablerule
&F&&T&&F&&F&&T&&T&\cr\tablerule
&T&&T&&F&&T&&T&&T&\cr\tablerule
&T&&T&&F&&T&&T&&T&\cr\tablerule
&T&&T&&T&&F&&F&&T&\cr\tablerule
&T&&T&&T&&F&&F&&T&\cr\tablerule
&T&&T&&T&&T&&T&&T&\cr\tablerule
&T&&T&&T&&T&&T&&T&\cr\tablerule}}\pars

\newbox\medstrutbox
\setbox\medstrutbox=\hbox{\vrule height14.5pt depth9.5pt width0pt}
\def\medstrut{\relax\ifmmode\copy\medstrutbox\else\unhcopy\medstrutbox\fi}

\pars
\noindent Now use Theorem 2.5.2, with $P = A, \ Q = B,\ R = C.$\pars
\noindent [2] (a) A contradiction. (b) Not a contradiction. (c) A 
contradiction.
(d) Not a contradiction.\parm
\medskip

\centerline{\sl Exercise 2.6$,$ Section 2.6}
\medskip
\noindent [1] (A) Suppose that we assume that there is some $z \in [x]\cap 
[y].$ Then we have that $z \equiv x,\  z \equiv y.$ But symmetry yields that 
$x \equiv z.$ From the transitive property$,$ we have that $x \equiv y.$ Hence
$x \in [y].$ Now let $u \in [x].$ Then $u \equiv x\Rightarrow u \equiv y.$ 
Thus $u \in [y].$ Thus $[x] \subset [y].$ Since $y \equiv x,$ this last 
argument repeated for $y$ shows that $[y]\subset [x].$ Hence $[x]= [y].$ \par
(B) Well$,$ just note that for each $x\in X$ it follows that $x \equiv x$. 
Hence$,$ from the  definition$,$ $x \in [x].$\pars
\noindent [2] (A) Since $B$ is a binary relation on $X$$,$ it is defined for all 
members of $B$. We are given that $B$ is reflexive. Thus we need to show that 
it is transitive and symmetric. So$,$ let $x B y.$ Then $x \in (y).$ 
From reflexive$,$ we have that $x \in (x).$ From (A)$,$ we have that $(x) = (y).$  
Thus $y \in (x).$ Therefore$,$ $y B x.$ Thus yields that $B$ is symmetric. \par
To show that it is transitive$,$ assume that since $B$ in \underbar{on} the 
entire set $X$$,$ that $x \in {y},\ y \in (z).$ Then from the reflexive 
property$,$ $y \in (y)$. Hence $y \in (y) \cup (z).$ Thus $(y) = (z)$ from (A).
Thus $x B z$. Re-writing this in relation notation we$,$ have that if
$x B y$ and $y B z$$,$ then $x B z$. This $B$ is an equivalence relation.\par
(B) To show that reflexive is necessary$,$ consider the binary relation 
$R = \{(a,b),(b,a)\}$ on the two element set $\{a,b\}.$ Note that is can be 
re-written as $a R b,\ b R a.$ Then $(b) = \{a\}, \ 
(a)=\{b\}.$  Now this relation satisfies (A) vacuously (i.e (A) holds since 
that hypothesis never holds.) (B) is obvious. But$,$  
since neither $(a,a)$ nor $(b,b)$ are members of $R$$,$ then $R$ is 
not an equivalence relation.\pars
\noindent [3] (A) Let $D = (A \lor (A \lor B)),\ E = ((A \lor B)\lor C).$ Now 
from (3) Theorem 2.5.3$,$ we have that $C_D \equiv C_E.$\par
(B) Let $D = (A \lor B), \ E = ((\neg A) \to B).$ Then from part (46) of 
Theorem 2.5.3$,$ we have that $C_D \equiv C_E.$\par
(C) Let $D = (A \land B),\ E = (\neg (A\to (\neg B))).$ The part (43) of 
Theorem 2.5.3$,$ yields this result.\par
(D) First$,$ let $D = (A \iff B)$. Let $G = (A \to B) \land (B\to C).$ Then from 
part (47) of Theorem 2.5.3$,$ we have that $C_D \equiv C_G.$ We can now apply 
(C) twice$,$ and let $E = (\neg ((A \to B)\to (\neg (B\to C)))).$ This yields 
$C_D \equiv D_E.$ \par
(E) This actually requires a small induction proof starting at $n = 0.$ 
Let $n$ be the number of $\neg$ to the left of $A$. If $n = 0$ or $n = 1,$ 
then the property that there is either none or one such $\neg$ symbol to the 
left is established. Suppose that we have already established this for some 
$n >1.$ Now let there be $n+1$ such symbols to the left of $A$ and let this 
formula be denoted by $H$. First$,$ consider the subformula $G$ such that$\neg G 
= H. $ Then by the induction hypothesis$,$ there is a formula $K$ such that 
$K$ contains no $\neg$ to the left or at the most one $\neg$ to the left such 
that  $K \equiv G.$ If $K$ has no $\neg$ to the left$,$ then $C_H = C_{\neg 
G}\equiv C_{\neg K} \Rightarrow C_H\equiv C_{\neg K}. $ Now if $K$ has one 
$\neg$ to the left$,$ then$,$ in like manner$,$ $C_H \equiv C_{\neg(\neg K)} \equiv
C_K \Rightarrow C_H \equiv C_K,$ where $K$ has no $\neg$ symbols. The result 
follows by induction.\par
\medskip
\centerline{\sl Exercise 2.7$,$ Section 2.7}
\medskip
\noindent [1] (a) $P\iff (A \to (R \lor S)) \equiv ((\neg P) \lor (A \to 
(R \lor 
S))) \land ((\neg (A \to (R \lor S)))\lor P)\equiv ((\neg P) \lor ((\neg A) \lor(R \lor 
S))) \land ((\neg((\neg A)\lor (R \lor S)))\lor P)\equiv ((\neg P) \lor ((\neg A) \lor(R \lor 
S))) \land ( A\land ((\neg R) \land (\neg S))\lor P).$\par
(b) $((\neg P)\to Q) \iff R \equiv  (((\neg P)\to Q) \to R) \land (R \to  
((\neg P)\to Q)) \equiv  (\neg ((\neg P))\lor Q) \to R) \land ((\neg R) \lor  (\neg(\neg P)\lor Q))  \equiv((P\lor Q) \to R) \land ((\neg R) \lor  
(P\lor Q))\equiv ((\neg(P\lor Q)) \lor R) \land ((\neg R) \lor  
(P\lor Q))\equiv(((\neg P)\land(\neg Q)) \lor R) \land ((\neg R) \lor  
(P\lor Q)).$\par
(c) $(\neg ((\neg P)\lor (\neg Q))) \to R \equiv (P\land Q) \to R \equiv (\neg 
(P\land Q)) \lor R \equiv ((\neg P) \lor (\neg Q)) \lor R.$\par
(d) $((\neg P) \iff Q) \to R \equiv (\neg((\neg P) \iff Q))\lor R \equiv 
(\neg (((\neg P)\to  Q)\land (Q \to (\neg P))) \lor R\equiv ((\neg(\neg(\neg 
P ))\land(\neg Q))\lor Q)\lor((\neg Q)\lor (\neg P)))\lor R\equiv (((\neg 
P)\land (\neg Q))\lor Q) \lor ((\neg Q)\lor (\neg P)))\lor R.$\par
(e) $(S \lor Q)\to R \equiv (\neg (S \lor Q)) \lor R \equiv ((\neg S) \land 
(\neg Q)) \lor R.$\par
(f) $(P\lor(Q \land S)) \to R \equiv ((\neg P)\land((\neg Q)\lor (\neg S))) 
\lor R.$ \pars
\noindent [2]                                         
(a) $((\neg P)\lor Q)\land (((\neg Q)\lor P) \land R)\Rightarrow A_d =  
(P\land (\neg Q))\lor ((Q\land (\neg P)) \lor (\neg R))$\par
(b) $((P\lor (\neg Q))\lor R) \land(((\neg P)\lor Q)\land R)\Rightarrow A_d =
(((\neg P)\land Q)\land (\neg R)) \lor((P\land (\neg Q))\lor (\neg R)).$\par
(c) $((\neg R)\lor(\neg P)) \land (Q \land P)\Rightarrow A_d =
 (R\land P) \lor ((\neg Q) \lor (\neg P))$\par
(d) $(((Q\land (\neg R))\lor Q) \lor (\neg P))\land (Q \lor R)\Rightarrow A_d 
= ((((\neg Q)\lor R)\land (\neg Q)) \land  P)\lor ((\neg Q) \land (\neg R)).$ 
\pars
\noindent [3] (a) $(P \land Q \land R) \lor (P \land (\neg Q) \land R) \lor (P 
\land (\neg B) \land (\neg C)).$\par
(b) $(P \land Q \land R)\lor(P \land Q \land (\neg R))\lor(P \land (\neg Q)
 \land (\neg R))
\lor((\neg P) \land Q \land R)\lor((\neg P) \land Q \land (\neg R))\lor((\neg 
P) 
\land (\neg Q) \land R)
\lor((\neg P) \land (\neg Q) \land (\neg R) ).$\par
(c)  $(P \land Q \land (\neg R)) \lor(P \land (\neg Q) \land (\neg R))\lor
((\neg P) \land Q \land (\neg R))\lor((\neg P) \land (\neg Q) \land R).$\par
(d) $(P \land Q \land R)\lor(P \land Q \land (\neg R))\lor (P \land (\neg Q)
 \land R)\lor (P \land (\neg Q) \land (\neg R))\lor ((\neg P) \land Q \land R)
\lor ((\neg P) \land Q \land (\neg R))\lor ((\neg P) \land (\neg Q) 
\land R)\lor ((\neg P) \land (\neg Q) \land (\neg R)).$\pars

\noindent [5] (a) $C \land (A \lor(B \land ((\neg A)\lor B))).$\par
(b) $(C \lor ((A \lor B)\land (\neg(A \land B))))\land (\neg (C\land ((A \lor 
B) \land(\neg(A \land B))))).$\parm
\medskip
\centerline{\sl Exercise 2.8$,$  Section 2.8}
\medskip
\noindent [1]. 
\parm
\noindent (a) $P\to Q, \ (\neg P)\to Q \models Q.$\parm
\noindent (b) $P\to Q,\ Q\to R,\ P\models R.$\parm
\noindent (c) $(P\to Q)\to P,\ \neg P \models R.$ (Since premises are not 
satisfied$,$ $\models$ holds.)\parm
\noindent (d) $(\neg P)\to (\neg Q),\ P\not\models Q.$ (From row three.)\parm
\noindent (e) $(\neg P)\to (\neg Q),\ Q \models P.$\parm \medskip

\centerline{\sl Exercise 2.9$,$ Section 2.9}
\medskip
\noindent [1]. (a) $(\neg A) \lor B,\ C\to (\neg B) \models A \to C.$\pars
(1) $v(A \to C) = F \Rightarrow v(A) = T,\ v(C) = F \Rightarrow 
v(C\to (\neg B)) = T.$ Now (3)$,$ if $v(B) = T,$ then $v((\neg A) \lor B) = T.$ 
Consequently$,$ {\bf INVALID}. AND the letters represent atoms only. \pars
(b) $A\to (B\to C),\ (C\land D)\to E,\ (\neg G)\to (D \land (\neg E))\models A\to 
(B\to G).$ \pars
(1) Let $v(A\to (B\to G)) = F \Rightarrow v(A) = T, \ v(B \to G) = F 
\Rightarrow v(B) = T,\ v(G) = F.$ (2) Let $v(A\to (B\to C)) = T\Rightarrow 
v(C) = T.$ (3) Let $v((\neg G)\to (D \land (\neg E))) = T \Rightarrow v(D \land 
(\neg E)) = F \Rightarrow v(D) = F,\ v(E) = T\Rightarrow v((C\lor D) \to E) = 
F.$ Hence$,$ {\bf VALID.} \pars

(c) $(A\lor B)\to (C\land D),\ (D\lor E) \to G \models A \to G.$\pars
(1) Let $v(A \to G)= F \Rightarrow v(G) = F,\ v(A) = T.$ (2) Let
$v((D\lor E) \to G ) = T \Rightarrow v(D) = v(E) = F \Rightarrow  
v((A\lor B)\to (C\land D)) = F.$ Hence$,$ {\bf VALID.}\pars

(d) $A \to (B\land C),\ (\neg B) \lor D,\ (E \to (\neg G)) \to (\neg D),\ B\to 
(A \lor (\neg E)) \models B \to E.$\pars
(1) Let $v(B \to E)= F\Rightarrow v(B) =T, v(E) = F.$ (2) Let $v((\neg B)\lor 
D)= T \Rightarrow v(D) = T\Rightarrow v((E \to (\neg G)) \to (\neg D)) = F.$ 
Hence {\bf VALID}.\par\medskip\medskip
\noindent [2] (a) The argument is $H\lor S,\ (\neg H) \models S.$ Let $v(S) = 
F.$ Let $v(\neg H) = T\Rightarrow v(H) = F\Rightarrow v(H\lor S) = F.$ Hence$,$
{\bf VALID}.\pars

(b) The argument is $I \to C,\ (\neg I) \to D \models C\lor D.$ (1) Let 
$v( C \lor D) = F \Rightarrow v(C) = c(D) = F.$ (2) Let $v( (\neg I ) \to D) = 
T \Rightarrow v(I) = T\Rightarrow v(I \to C) = F.$ Hence$,$ {\bf VALID}.\pars
(c) The argument is $S \to I,\  I \to C,\  S \models C.$ (1) Let $v(C) = F.$ 
(2) Let $v(S) = T.$ (3) Let $v(S \to I) = T \Rightarrow v(I) = T\Rightarrow
v(I \to C) = F.$ Hence$,$ {\bf VALID.}\pars

(d) The argument is $P\to L,\  L\to N,\ N \models P.$ (1) Let $v(P)= 
F\Rightarrow v(P\to L) = T.$ (2) Let $v(N) = T\Rightarrow v(L \to N) = T.$ 
Hence$,$ {\bf INVALID}.\pars

(e) The argument is $W \lor C,\ W \to R,\ N \models W.$ (1) Let 
$v(W)=F\Rightarrow v(W \to R) = T.$ We have only one more premise remaining. 
This is NOT forced to be anything. We should be able to find values that will 
make it $T$. Let $v(C) = T\Rightarrow v(W \lor C) = T.$ Hence
{\bf INVALID.} [Note this is an over determined argument. The statement $N$ if 
removed will still lead to an invalid argument. However$,$ if we add the 
premises
$N\to W$$,$ then the argument is valid.\pars

(f) The argument is $C \to (M \to I),\  C \land (\neg M) \models \neg I.$
(1) Let $v(\neg I) = F \Rightarrow v(I) = T\Rightarrow v(C\to (M \to I))= T,$ 
independent of the values for $C,\ M.$ Thus for $v(C) = T,\ v(M) = 
F\Rightarrow v(C \land (\neg M) ) = T.$ Hence$,$ {\bf INVALID.}\pars

(g) The argument is $(L\lor C) \to (D \land S),\ D \to P,\ \neg P \models L.$
(1) Let $v(L) = F.$ (2) Let $V(\neg P) = T \Rightarrow v(P)= F.$ (3) Let $v(D \to 
P)) = T \Rightarrow v(D) = F.$ These values do not force $v((L\lor C) \to (D 
\land S))$ to be anything. This will immediately yield {\bf INVALID} since by 
special selection$,$ say $v(C) = F$$,$ we can always get a value 
$v((L\lor C) \to (D \land S))= T.$\parm
\medskip
\centerline{\sl Exercise 2.10$,$ Section 2.10}
\medskip
\noindent [1]. [Note: this is very important stuff. As will be shown in the 
next section$,$ if a set of premises is  inconsistent$,$ then there will always be 
a correct logical argument that will lead to any PRESELECTED conclusion.]\pars
(a) $A \to (\neg(B \land C)),\ (D \lor E) \to G,\ G \to (\neg (H \lor I)),\ 
(\neg C)\land E \land H.$\pars
(1) Let $v(A \to (\neg(B \land C)))= T \Rightarrow v(C) = F,\ v(E) = T,\ v(H) 
= T.$ (2) Let $v(G \to (\neg (H \lor I))) = T \Rightarrow v(G) = F \Rightarrow 
v((D \lor E) \to G) =F.$ Hence {\bf INCONSISTENT}.\pars
(b) $(A\lor B) \to (C \land D),\ (D\lor E) \to G,\ A \lor (\neg G).$\pars
(1) Let $v(A \lor (\neg G))=T \Rightarrow$ case studies (1i) $v(A) = T,\ v(G) 
= T,$ (1ii)$,$ $v(A) = T,\ v(G) = F.$ (1iii) $V(A) = F,\ v(G) = F.$ 
(2) Consider (1i). Then $v(G)= T \Rightarrow v(D\lor E) \to G) = T,$ 
independent of the values for $D,\ E,\ C.$ (3) So$,$ select $v(C) =T,\ v(D)= T,$ 
then $v((A\lor B) \to (C \land D)= T.$ Thus we have found an assignment that 
yields T for each member of the set of premises. Hence$,$ they are
 {\bf CONSISTENT} and the formula symbols represent ATOMS. \pars

(c) $(A\to B)\land (C \to D),\ (B \to D)\land((\neg C)\to A),\ (E \to G)\land 
(G \to (\neg D)),\ (\neg E)\to E.$\pars
(1) Let $v((\neg E)\to E) = T\Rightarrow V(E) = T.$ (2) Let $v((E \to G)\land 
(G \to (\neg D))) = T \Rightarrow v(G \to (\neg D))=T \Rightarrow v(G) = 
T\Rightarrow v(D) = F.$ (2) Let $v((A\to B)\land (C \to D))= T \Rightarrow
v(A \to B) = T,\ v(C \to D) = T \Rightarrow v(C) = F.$ (3) Now suppose that 
$v((B \to D)\land((\neg C)\to A)) = T.$ Then $v(B \to D) =T,\ v(((\neg C)\to 
A))= T \Rightarrow v(B) = F$ and from (2) $\rightarrow v(A) = F.$ From $
v(((\neg C)\to A))= T$ this yields that $v(\neg C) = F \Rightarrow v(C) = T.$ 
But this contradicts the statement just before (3). Hence$,$ 
$v((B \to D)\land((\neg C)\to A)) = F$ in all possible cases and the set of 
premises is {\bf INCONSISTENT}.\pars

(d) $(A\to (B\land C))\land(D\to (B\land E)), ((G \to (\neg A))\land H)\to I,
 (H \to I)\to (G\land D),  \neg((\neg C)\to E).$ \pars
(1) Let $v((A\to (B\land C))\land(D\to (B\land E))).= T \Rightarrow v(A\to 
(B\land C)) = T,\ v(D\to (B \land E))= T.$ (2) Let $v((\neg((\neg C)\to E)) = 
T\Rightarrow v((\neg C)\to E) = F\Rightarrow v(E) = F,\ v(C) =F$ and from (1) 
$\Rightarrow v(D) = F,\ v(A) = F.$ Also $\Rightarrow v(G\land D) = F.$ (3)
Let $v((H \to I)\to (G\land D)) = T \Rightarrow v(H\to I) = F \Rightarrow
v(H) = T,\ v(I) = F. $ (4) Let $V(((G \to (\neg A))\land H)\to I)= T 
\Rightarrow v(((G \to (\neg A))\land H))= F \Rightarrow v(G \to (\neg A))= 
F.$ But this implies that $v(\neg A)= F \Rightarrow v(A) = T$ which 
contradicts the result in (2). Thus {\bf INCONSISTENT}. \parm
\medskip
\centerline{\sl Exercise 2.11$,$ Section 2.11}
\medskip
\noindent [1] $\vdash(\neg (\neg A)) \to A$\parm
\line{\indent (1) $(\neg(\neg A))\to ((\neg(\neg(\neg(\neg A))))\to (\neg(\neg 
A)))$\leaderfill $P_1$}
\line{\indent (2) $( (\neg(\neg(\neg(\neg A))))\to (\neg (\neg A))) \to ((\neg 
A)\to 
(\neg (\neg (\neg A))))$ \leaderfill $P_3$}
\line{\indent (3) $(\neg (\neg A)) \to ((\neg A) \to (\neg (\neg (\neg 
A))))$\leaderfill $HS(1,2)$}
\line{\indent (4) $((\neg A) \to (\neg (\neg (\neg A))))\to ((\neg (\neg 
A))\to A)$\leaderfill $P_3$}
\line{\indent (5) $(\neg (\neg A)) \to ((\neg (\neg A)) \to A)$\leaderfill $HS(3,4)$}
\line{\indent (6) $((\neg (\neg A))\to ((\neg (\neg A)) \to A))\to$\hfil} 
\line{\indent\indent $(((\neg(\neg A)) \to (\neg (\neg A))) \to ((\neg (\neg 
A)) \to A))$\leaderfill $P_2$}
\line{\indent (7) $((\neg (\neg A))\to (\neg (\neg A))) \to ((\neg (\neg A)) \to
A)$\leaderfill $MP(5,6 )$}  
\line{\indent (8) $((\neg (\neg A))\to (((\neg(\neg A))\to (\neg(\neg A))) \to 
(\neg(\neg A)))) \to (((\neg(\neg A))\to$\hfil}
\line{\indent\indent $((\neg(\neg A)) \to (\neg (\neg A))))\to ((\neg(\neg A)) 
\to (\neg (\neg A))))$\leaderfill $P_2$}
\line{\indent (9) $(\neg(\neg A)) \to (((\neg (\neg A)) \to (\neg(\neg A))) 
\to (\neg (\neg A)))$ \leaderfill $P_1$}
\line{\indent (10) $((\neg(\neg A)) \to ((\neg(\neg A)) \to (\neg(\neg A)))) 
\to$\hfil} 
\line{\indent\indent $((\neg(\neg A)) \to (\neg(\neg A)))$\leaderfill $MP(8,9)$}
\line{\indent (11) $(\neg(\neg A)) \to ((\neg(\neg A)) \to (\neg(\neg 
A)))$\leaderfill $P_1$}
\line{\indent (12) $(\neg(\neg A)) \to (\neg(\neg A))$\leaderfill $MP(10,11)$}
\line{\indent (13) $(\neg(\neg A)) \to A$\leaderfill $MP(7,12)$} \parm

\noindent [2] (a) $\vdash A\to (\neg (\neg A))$\parm
\line{\indent (1) $((\neg(\neg(\neg A))) \to (\neg A)) \to (A \to (\neg (\neg 
A)))$\leaderfill $P_3$}
\line{\indent (2) $\vdash (\neg (\neg(\neg A)))\to (\neg A)$\leaderfill from [1].}
\line{\indent (3) $A \to (\neg (\neg A))$\leaderfill $MP(1,2)$}\parm

(b) $\vdash (\neg B) \to (B \to A).$ \parm
\line{\indent (1) $(\neg B) \to ((\neg A)\to (\neg B))$ \leaderfill $P_1$}
\line{\indent (2) $((\neg A) \to (\neg B)) \to (B \to A)$ \leaderfill $P_3$}
\line{\indent (3) $(\neg B) \to (B \to A)$ \leaderfill $HS(1,2)$}\parm
\medskip

\centerline{\sl Exercise 2.12$,$ Section 2.12}
\medskip
\indent [1] $(\neg(\neg A)) \vdash A.$\parm
\line{(1) $(\neg(\neg A))$\leaderfill Premise}
\line{(2) $(\neg(\neg A)) \to ((\neg (\neg (\neg (\neg A)))) \to (\neg(\neg 
A)))$\leaderfill $P_1$}
\line{(3) $ (\neg (\neg (\neg (\neg A)))) \to (\neg(\neg A))$\leaderfill 
$MP(1,2)$}
\line{(4) $((\neg(\neg(\neg(\neg A))))\to (\neg(\neg A))) \to ((\neg A) \to 
(\neg(\neg(\neg A))))$\leaderfill $P_3$}
\line{(5) $(\neg A) \to (\neg(\neg(\neg A)))$\leaderfill $MP(3,4)$}
\line{(6) $((\neg A) \to (\neg(\neg(\neg A))))\to ((\neg(\neg A)) \to 
A)$\leaderfill $P_3$}
\line{(7) $(\neg(\neg A)) \to A$\leaderfill$MP(5,6)$}
\line{(8) $A$\leaderfill $MP(1,7)$}
\parm
 \indent [2] (a) $A \to B,\ B \to C \vdash A\to C.$\parm
\line{(1) $(B \to C) \to (A \to (B \to C))$\leaderfill $P_1$}
\line{(2) $B \to C$\leaderfill Premise}
\line{(3) $A\to (B \to C)$\leaderfill $MP(1,2)$}
\line{(4) $(A\to (B \to C)) \to ((A \to B) \to (A \to C))$\leaderfill $P_2$}
\line{(5) $(A \to B) \to (A \to C)$\leaderfill $MP(3,4)$}
\line{(6) $(A \to B)$\leaderfill Premise}
\line{(7) $A \to C$\leaderfill $MP(5,6)$}\parm
\indent (b) $(\neg A) \vdash A \to B.$\parm
\line{(1) $\neg A$ \leaderfill Premise}
\line{(2) $\vdash (\neg A) \to (A \to B)$\leaderfill Ex. 2.11 (2b)}
\line{(3) $A \to B$ \leaderfill $MP(1,2)$}\parm
    
\indent (c) Given $\vdash (B\to A) \to ((\neg A) \to (\neg B)).$ Show that\par
\noindent  $\neg(A \to B) \vdash B \to A.$\parm
\line{(1) $\neg(A \to B)$\leaderfill Premise}
\line{(2) $\vdash (B\to A( \to B)) \to ((\neg (A \to B)) \to (\neg B))$\leaderfill Given}
\line{(3) $(B \to (A \to B)$ \leaderfill $P_1$}
\line{(4) $(\neg (A \to B))\to (\neg B)$\leaderfill $MP(2,3)$}
\line{(5) $(\neg B)$\leaderfill $MP(1,4)$}
\line{(6) $(\neg B) \to (B \to A)$\leaderfill Ex. 2.11 (2b)}
\line{(7) $B \to A$\leaderfill $MP(5,6)$}\parm
\medskip
\centerline{\sl Exercise 2.13$,$ Section 2.13}
\medskip
\noindent 1. Well$,$ we use the 
statement that $\vd A \to A$ in the proof of the deduction theorem.\parm
\noindent 2. 
\line{(A) Show that $\vd (B\to A) \to ((\neg A)\to (\neg B))$\hfil}\parm
\line{(1) $B\to A$ \leaderfill Premise and D. Thm.}\smallskip
\line{(2) $\vd (\neg(\neg B))\to B$\leaderfill Ex. 2.12.1 and D. Thm.}\smallskip
\line{(3) $(\neg (\neg B)) \to A$\leaderfill $HS(1,2)$}\smallskip
\line{(4) $A \to (\neg (\neg A))$\leaderfill Ex 2.11: 2a}\smallskip
\line{(5) $(\neg(\neg B)) \to (\neg(\neg A))$\leaderfill $HS(3,4)$}\smallskip
\line{(6) $((\neg(\neg B))\to (\neg A(\neg A)))\to ((\neg A) \to (\neg 
B))$\leaderfill $P_3$}\smallskip
\line{(7) $(\neg A) \to (\neg B)$\leaderfill $MP(5,6)$}\smallskip
\line{(8)  $\vd (B\to A) \to ((\neg A)\to (\neg B))$\leaderfill D. Thm.}\medskip

\line{(B) Show that $\vd ((A \to B)\to A) \to A$\hfil}\parm
\line{(1) $(A \to B) \to A$\leaderfill Premise D. Thm}\smallskip
\line{(2) $\vd (\neg A) \to (A \to B)$\leaderfill D.Thm. and Ex 2.12.2b}\smallskip
\line{(3) $(\neg A) \to A$\leaderfill $HS(1,2)$}\smallskip
\line{(4) $(\neg A)\to (\neg(\neg((\neg A) \to A)) \to (\neg A))$\leaderfill 
$P_1$}\smallskip
\line{(5) $(\neg(\neg((\neg A) \to A)) \to (\neg A)) \to (A \to (\neg 
((\neg A)\to A))$\leaderfill 
$P_3$}\smallskip
\line{(6) $(\neg A) \to (A \to (\neg ((\neg A)\to A)))$\leaderfill
$HS(4,5)$}\smallskip
\line{(7) $((\neg A)\to (A \to (\neg((\neg A) \to A))))\to(((\neg A)\to A)\to$\hfil}
\line{\hskip 1.5in $((\neg A) \to (\neg((\neg A) \to A))))$\leaderfill $P_2$}\s
\line{(8) $((\neg A)\to A)\to ((\neg A)\to (\neg((\neg A)\to A)))$\leaderfill 
$MP(6,7)$}\s
\line{(9) $(\neg A) \to (\neg ((\neg A) \to A))$\leaderfill $MP(3,8)$}\s
\line{(10) $((\neg A) \to (\neg((\neg A)\to A)))\to (((\neg A) \to A) \to 
A)$\leaderfill $P_3$}\s
\line{(11) $((\neg A)\to A) \to A$\leaderfill $MP(9,10)$}\s
\line{(12) $A$ \leaderfill $MP(3,11)$}\s
\line{(13) $\vd ((A \to B)\to B) \to A$ \leaderfill D. Thm.}\parm 
\noindent [3] 
\centerline{$A \to B,\ A \vd B$}\s
\line{(1) $\vd A\to A$\leaderfill Ex. 2.11.1}\s
\line{(2) $A \to B$ \leaderfill Premise}\s
\line{(3) $(A\to B) \to (A \to (A \to B))$\leaderfill $P_1$}\s
\line{(4) $A \to (A \to B)$\leaderfill $MP(2,3)$}\s
\line{(5) $(A \to (A\to B))\to ((A\to A) \to (A\to B)) $\leaderfill $P_2$}\s
\line{(6) $(A\to A) \to (A\to B)$\leaderfill $MP(4,5)$}\s
\line{(7) $A \to B$\leaderfill $MP(1,6)$}\s\parm
\medskip
\centerline{\sl Exercise 2.14$,$ Section 2.14}
\medskip
\noindent [1] (a) (i) $B,\ C \vd (\neg B) \to (\neg C)$. (ii) $B,\ (\neg C) \vd (\neg B) \to 
(\neg C).$
(iii) $(\neg B),\ C \vd \neg ((\neg B) \to (\neg C)).$  (iv) 
$(\neg B),\ (\neg C) \vd (\neg B) \to (\neg C).$ \pars   

(b) (i) $B,\ C \vd  B \to  C$. (ii) $B,\ (\neg C) \vd \neg (B \to 
 C).$
(iii) $(\neg B),\ C \vd  B \to  C.$  (iv) 
$(\neg B),\ (\neg C) \vd  B \to  C.$ \pars      

(c) (i) $B,\ C,\ D \vd B\to (C \to D).$ (ii) $B,\ C,\neg D \vd \neg (B\to (C 
\to D)).$ (iii) $B,\ \neg C,\ D \vd B\to (C \to D).$ (iv) $B,\ \neg C,\ \neg D 
\vd B\to (C \to D).$ (v) $\neg B,\ C,\ D \vd B\to (C \to D).$ (vi) $\neg B,\ C,\ 
\neg D \vd B\to (C \to D).$ (vii)$\neg B,\ \neg C,\ D \vd 
B\to (C \to D).$(viii) $\neg B,\ \neg C,\ \neg D \vd B\to (C \to D).$ \pars

(c) (i) $B,\ C,\ D \vd \neg (B\to (C \to D)).$ (ii) $B,\ C,\neg D \vd B\to (C 
\to D).$ (iii) $B,\ \neg C,\ D \vd \neg (B\to (C \to D)).$ (iv) $B,\ \neg C,\ \neg D 
\vd \neg (B\to (C \to D)).$ (v) $\neg B,\ C,\ D \vd B\to (C \to D).$ (vi) $\neg B,\ C,\ 
\neg D \vd B\to (C \to D).$ (vii)$\neg B,\ \neg C,\ D \vd 
B\to (C \to D).$(viii) $\neg B,\ \neg C,\ \neg D \vd B\to (C \to D).$ \pars
\medskip

\centerline{\sl Exercise 2.16$,$ Section 2.16}
\medskip
\noindent Almost all of the interesting stuff about consequence operators 
involves a much deeper use of set-theory than used in these questions.\parm
\noindent [1] Let ${\cal A} \subset {\cal B}.$  Suppose that$A \in C({\cal 
A}).$ Then there is a finite  ${\cal F}\subset {\cal A}$ such that $A \in 
C({\cal F}).$ But ${\cal F}\subset {\cal B}.$ Thus $A \in C({\cal B}).$ 
Therefore$,$ $C({\cal A}) \subset C({\cal B}).$ \parm
\noindent [2] Let ${\cal A} \subset C({\cal B}).$ Then $C({\cal A}) \subset 
C(C({\cal B})) = C({\cal B})$ from [1] and our axioms. \par
Conversely$,$ suppose that $C({\cal A}) \subset C({\cal B}).$ Since 
${\cal A} \subset C({\cal A}),$ then ${\cal A} \subset C({\cal B}).$\parm
\noindent [3] We know that ${\cal A} \cup {\cal B} \subset {\cal A} \cup 
C({\cal B}) \subset C({\cal A}) \cup C({\cal B}) \Rightarrow$
(a) $C({\cal A}\cup {\cal B}) \subset C({\cal A} \cup 
C({\cal B})) \subset C(C({\cal A}) \cup C({\cal B})).$ Now 
$C({\cal A}) \subset C({\cal A}\cup {\cal B}),\ C({\cal B}) \subset C({\cal A}\cup {\cal B})
\Rightarrow C({\cal A}) \cup C({\cal B})\subset C(C({\cal A}\cup {\cal B})) =
 C({\cal A}\cup {\cal B}).$  
Combining this with (a) implies that $C({\cal A}\cup {\cal B}) \subset C({\cal A} \cup 
C({\cal B})) \subset C(C({\cal A}) \cup C({\cal B}))\subset C({\cal A}\cup 
{\cal B}).$ Hence $C({\cal A}\cup {\cal B}) = C({\cal A} \cup 
C({\cal B})) = C(C({\cal A}) \cup C({\cal B})).$\par
Note: This is reasonable. You cannot get more deductions if you consider the 
premises broken up into subsets.\parm
\medskip
\centerline{\sl Exercise 2.17$,$ Section 2.17}
\medskip
\noindent [1] (ii) {\sl If $B \in L',$ then either $B \in \overline{\Gamma}$ 
or $\neg B \in \overline{\Gamma}$ not both.}\pars
Proof. Assume that $B\notin \overline{\Gamma}.$ Then $\overline{\Gamma}\cup 
\{B\}$ is inconsistent. Thus there is a finite $\{A_1,\ldots,A_n\} \subset  
\ol$ and $A_1,\ldots,A_n,B \vd C \land (\neg C),$ for some $C \in L'.$ We know 
that $C \land (\neg C)\mod D$ for any $D \in L'.$ Thus $C \land (\neg C)\mod 
\neg B.$ Hence$,$ $A_1,\ldots,A_n,B \mod \neg B \Rightarrow A_1,\ldots,A_n\mod B 
\to (\neg B) \mod \neg B \Rightarrow  A_1,\ldots,A_n \mod \neg B \Rightarrow
A_1,\ldots,A_n \vd \neg B \Rightarrow \neg B \in \ol.$ Suppose that 
$B, \  \neg B \in \ol.$ Then this contradicts Theorem 2.17.1 (e). \pars
 (iii) {\sl If $B \in \ol,$ then $A \to B \in \ol$ for each $A\in 
L'.$}\pars
Proof. Let $B\in \ol.$ The  $\ol \vd B.$ Thus there is a finite 
${\cal F}\subset \ol$ such that ${\cal F} \vd B.$ Hence $ {\cal F} \mod B.$
But $B \mod A \to B$ for any $A \in L'.$  Hence ${\cal F}\mod A \to B 
\Rightarrow {\cal F} \vd A \to B \Rightarrow \ol \vd A \to B \Rightarrow
A \to B \in \ol.$ \pars

(iv) {\sl If $A \notin \ol,$ then $A \to B \in \ol$ for each $B \in L'.$}\pars
Proof. If $A \notin \ol,$ then $\neg A \in \ol$ from (ii). But for any $B \in 
L',$ $\neg A \mod A \to B \Rightarrow \neg A \vd A \to B$ and since $\neg A 
\in \ol,$ $\ol \vd A\to B.$ Hence $A \to B \in \ol.$\pars

(v) {\sl If $A \in \ol$ and $B \notin \ol,$ then $A \to B \notin \ol.$}\pars
Proof. Let $A \in \ol,\ B \notin \ol.$ Suppose that $A \to B \in \ol.$ Then 
$\ol \vd A \to B.$ Since $\ol \vd A,$ by one MP step $\ol \vd B \Rightarrow 
B \in \ol$; a contradiction. The result follows.\parm
\noindent [2] (a) {\sl $\mod A$ iff $\{\neg A\}$ is not satisfiable.} \pars
Proof. Let $\mod A.$ Then for each $\ass$$,$ $v(A,\ass) = T.$ Hence for each 
$\ass,\ v(\neg A,\ass) = F.$ Hence $\{\neg A\}$ is not satisfiable.\par
Conversely$,$ let $\{\neg A\}$ not be satisfiable. Then for  each $\ass, v(\neg 
A,\ass) = F.$ This for each $\ass,\ v(A,\ass) = T.$ Thus $A$ is 
satisfiable.\pars

(b) {\sl $\{A\}$ is consistent iff $\not\vd \neg A.$}\pars
Proof. Let $\{A\}$ be consistent and $\vd \neg A.$ Hence $\mod \neg A.$ Hence 
for each $\ass,\ v(\neg A,\ass) = T\Rightarrow v(A,\ass) = F
.$ Hence $ A\mod B\land (\neg B)$ for some 
$B$. Thus $A\vd B\land (\neg B).$ Thus $\{A\}$ is inconsistent: a 
contradiction.\par
Conversely$,$ let $\not\vd \neg A.$ Then $\not\mod \neg A.$ Hence there is some
$\ass$ such that $v(\neg A,\ass) = F \Rightarrow v(A,\ass) = T \Rightarrow A 
\not\mod B\land (\neg B)$ for any $B$. $A\not\vd B\land (\neg B)$ for any $B$ 
implies that $\{A\}$ is consistent.\pars
(c) {\sl The Completeness Theorem is equivalent to the statement that every 
consistent formula is satisfiable.}\pars
Proof. Assume Completeness Theorem. Now let 
$\{A\}$ be consistent. Then from (b) $\not\vd \neg A.$ Hence from Completeness 
contrapositive$,$ we have that $\not\mod \neg A.$ Thus there is some $\ass$ such 
that $v(\neg A,\ass) = F.$ Hence $v(A,\ass) = T.$ Thus $\{A\}$ is 
satisfiable. \pars
Conversely$,$ assume that all consistent formula are satisfiable. [Note we 
cannot use the Completeness Theorem since this is what we want to establish.]
Let $\mod A.$ Then clearly $\neg A$ is not satisfiable. Hence$,$ from the 
contrapositive of our assumption$,$ $\{\neg A\}$ is inconsistent. Thus $\neg 
A\vd$anything. Consequently$,$ $\neg A \vd A \Rightarrow \ {\rm (a)}\ \vd (\neg A) \to A.$ 
But $\vd (A \to A) \to (((\neg A)\to A) \to A),$ by example 2.15.4. But $\vd A 
\to A.$ Hence by one MP step we have that $\vd ((\neg A )\to A) \to A.$ 
Adjoining this formal proof to (a)$,$ and we get $\vd A.$ This is but the 
Completeness Theorem.\parm
\noindent [3] (a) Let $\ass$ be any assignment to $A$. Then since $\mod A \to A.$
Now this could be established by formal induction$,$ but more simply$,$ note that 
if $\emptyset \not= {\cal F} \subset \Gamma,$ then one member of $\cal F$ contains 
that largest number of $A$s$,$ say $n \geq 2$ and the smallest sized subformula 
is $A \to A$ which has truth value $T$. All other formula$,$ if any$,$ have $A\to$ 
to the left etc. and$,$ hence$,$ have value $T$. It appears that $\Gamma$ is 
consistent.\pars
(b) $\{A_3\}$ is finite inconsistent subset of $\Gamma$. Hence $\Gamma$ is 
inconsistent.\pars
(c) $\{A_1,A_2\}$ is a finite inconsistent sunset of $\Gamma$. Hence$,$ $\Gamma$ 
is inconsistent.\pars\bigskip

\centerline{\sl Exercise 3.1$,$ Section 3.1}
\medskip
\noindent [1] (a) $A \in Pd$$,$ (b) $A \notin Pd$$,$ [Note the $\exists c.$] 
(c) $A \notin Pd$$,$ [Note the third (. (d) $A \notin Pd$$,$ [Note the symbols 
$\forall y).$] (e) $A \notin Pd$$,$ [Note the missing last two 
parentheses.]\pars
\noindent [2] (a) size (A) = 3$,$ (b) size (A) = 4 (c) size (A) = 4$,$ (d) 
size (A) = 4.\pars
\noindent [3] Please note the way I've translated these statements into 
first-order predicate form is not a unique translation. It is$,$ however$,$ in a 
mathematically standard form. Also note that from this moment on$,$ we will 
simply any formula that has a sequence of more than one 
$\land$s or $\lor$s by not putting parentheses about the inner subformula. 
Indeed$,$ the way we write or language tends not to include pauses between such 
combinations. [They are all equivalent no matter where we put the 
parentheses.] Finally notice that when we use the operators as abbreviations 
for the predicates we include ( and ) only to avoid confusion. If you make a 
complete substitution $,$ it will look like this $(x + 1 = y) \equiv R(x,1,y).$
\parm
\line{(a) $ \forall x((P(x) \land (x = 0)) \to (\exists y(Q(y,x) \land (y > 
x))))$\hfil}\pars                                              
\line{(b) $\forall x(R(x) \to (\exists y(R(x) \land (y>x))))$\hfil}\pars                                              
\line{(c) $\forall x(\exists y(\forall z((R(x) \land R(y) \land R(z)\land ( z+ 
1<y)) \to (x + 2 < 4)))) $\hfil}\pars                                              
\line{(d) $\forall x((W(x) \land L(x)) \to (\exists y(J(y) \land A(x,y))))$\hfil}\pars                                              
\line{(e) $ \exists x(L(x) \land S(x) \land A(x,j))$\hfil}\pars                                              
\line{(f) $\forall x(\forall y(\forall z(((P(x) \land P(y) \land P(z) \land 
R(x,z) \land R(y,z)) \to R(x,y))))$\hfil}\pars                                              
\line{(g) $\forall x(B(x) \to M(x))$\hfil}\pars                                              
\line{(h) $\forall x(P(x) \to (C(x) \land U(x)))$\hfil}\pars                                              
\line{(i) $\forall x(\forall y(\forall z((x=y)\land (y=z)\land (x=z)\land (x >0) \land (y > 0)\land (z > 0)\land(x + 
y + z >3)\land$\hfil}\line{\indent $(x + y+ z <9))\to
((1 < x < 3) \land (1 < y < 3) \land (1 < z < 
3)))))$\hfil}\pars                                            
\line{(j) $\forall x(\forall y((P(x) \land P(y) \land B(x,y)) \iff (M(x)\land 
M(y) \land (\neg (x = y))\land Q(x,y))))$\hfil}\pars
\noindent [4] [Note: In what follows$,$ you do not need to inquire whether of 
not the statement holds in order to translate into symbolic form.] \par
(a) Seven is a prime number and seven is an odd number.\pars
(b) (Two different ways). For each number$,$ if 2 divides the number$,$ then the 
number is even. [Note: A metalanguage variable symbol can also be used. 
For each X$,$ if 2 divides X$,$ then X is even. You can also write it more concisely 
as ``For each X$,$ if X is a number and 2 divides X$,$ then X is even.'' This is 
the very important  ``bounded form.'' The X is restricted to a particular 
set.]\pars
(c) Using the metasymbol method$,$ we have  ``There exists an X such that X is 
an even number and X is a prime number$,$ and there does not exist an X such 
that X is an even number and X is a prime$,$ and there exists some Y such that 
Y is not equal to X and Y is an even number and Y is a prime number.\pars
(d) For each X; if X is an even number$,$ then for each Y$,$ if X divides Y$,$ then 
Y is an even number.\pars
(e) For each X$,$ if X is an odd number$,$ then there exists a Y such that if Y is 
a prime number$,$ then Y divides X.\parm
\centerline{\sl Exercise 3.2$,$ Section 3.2}
\medskip
\noindent [1] (a) (1) scope = $\exists x Q(x,z)$. (2) scope = $Q(x,z)$. (3) 
scope = $Q(y,z).$ \pars
(b) (1) scope = $\forall y ( P(c) \land Q(y))$. (2) scope = $P(c) \land Q(y).$ 
(3) scope = $R(x).$\pars
(c) (1) scope = $(Q(y,z)\to (\forall x R(x))).$ (2) scope = $R(x).$\pars
(d) (1) scope = $(P(z) \land (\exists x Q(x,z)))\to (\forall z(Q(c) \lor 
P(z)))).$ (2) scope = $Q(x,z)$. (3) scope = $(Q(c) \lor P(z))).$\parm
 
\noindent [2]\pars 
(a) $\forall z_{\b 3}(\exists y_{\b 2}(P(z_{\b 3},y_{\b 2})\land (\forall z_{\b 1}Q(z_{\b 1},x)) \to 
M(z_{\b 3}))).$\par\medskip
(b) $\forall x_{\b 3}(\exists y_{\b 2}(P(x_{\b 3},y_{\b 2})\land (\forall y_{\b 1}Q(y_{\b 1},x_{\b 3})) \to M(x_{\b 3}))).$\par\medskip  
(c) $\forall z_{\b 3}(\exists x_{\b 2}(P(z_{\b 3},x_{\b 2})\land (\forall z_{\b 1}Q(z_{\b 1},y)) \to M(z_{\b 3}))).$\par\medskip 
(d) $\forall y_{\b 3}(\exists z_{\b 2}(P(y_{\b 3},z_{\b 2})\land (\forall z_{\b 1}Q(z_{\b 1},x)) \to M(y_{\b 3}))).$\par\medskip 
(e) $\forall y_{\b 3}(\exists z_{\b 2}(P(z_{\b 2},y_{\b 3})\land (\forall z_{\b 1}Q(z_{\b 1},x)) \to M(y_{\b 3}))).$\par\medskip
(f) $\exists x_{\b 3}(\forall z_{\b 2}(P(x_{\b 3},z_{\b 2}) \lor (\forall u_{\b 1} M(u_{\b 1},y,x_{\b 3})))).$\par\medskip
(g) $\exists y_{\b 3}(\forall x_{\b 2}(P(z,x_{\b 2}) \lor (\forall x_{\b 1} M(x_{\b 1},u,y_{\b 3})))).$\par\medskip
(h) $\exists y_{\b 3}(\forall x_{\b 2}(P(y_{\b 3},x_{\b 2}) \lor (\forall x_{\b 1} M(x_{\b 1},y_{\b 3},z)))).$\par\medskip
(i) $\exists z_{\b 3}(\forall x_{\b 2}(P(z_{\b 3},x_{\b 2}) \lor (\forall x_{\b 1} M(x_{\b 1},y,z_{\b 3})))).$\par\medskip 
(j) $\exists x_{\b 3}(\forall x_{\b 2}(P(z,x_{\b 2}) \lor (\forall z_{\b 1} 
M(x_{\b 2},y,z_{\b 1})))).$\par\medskip  
\noindent [3] $ (a) \cong (d);\ (f) \cong (i).$ \pars
\noindent [4] (a) Free$,$ $x,y,z$; Bound $x,z$.\pars
(b) Free$,$ $z$: Bound $x,z.$ \pars
(c) Free$,$ $x,y$; Bound $x$.\pars
(d) Free$,$ $z$; Bound $x,y,z.$\pars
(e) Free$,$ $x$; Bound $x,y$.\pars
(f) Free$,$ none; Bound $x,y.$\parm

\noindent [5] (a) They are (b)$,$ (d). (b) It is (b). (c) It is (f)\parm
\centerline{\sl Exercise 3.3$,$ Section 3.3}

\noindent [1] 
(a) $S^x_a\ (\exists x P(x)) \to R(x,y)]=(\exists x P(x)) \to R(a,y)$.
(b) $S^y_x\ (\exists y R(x,y))\iff (\forall x R(x,y))]= (\exists y R(x,y))\iff 
(\forall x R(x,x)).$
(c) $S^y_a\ (\forall x P(y,x)) \land (\exists y R(x,y))]=(\forall x P(a,x)) 
\land (\exists y R(x,y)),$
(d) $S^x_aS^x_b (\exists x P(x)) \to  R(x,y)]]=(\exists x P(x)) \to  R(b,y).$
(e) $S^x_a S^y_x (\exists y R(x,y))\iff (\forall x R(x,y))]]=(\exists y 
R(x,y))\iff (\forall x R(x,x)).$   
(f) $S^x_aS^y_b(\forall z P(y,x)) \land (\exists y R(x,y))]]=(\forall z 
P(b,a)) \land (\exists y R(a,y))$\parm

\noindent [2] (a) always true. (b) NOT true. Consider the example $\forall 
yP(x,y).$ (c) Always true. An argument is that$,$ for the left hand side$,$ 
we first substitute for all 
free occurrences of $z$$,$ if any$,$ a $w$. Whether or not the result gives no 
free occurrences of $w$ or not will not affect the substitution of $x$ for the 
free occurrences of $y$ since these are all distinct variables. The same 
argument goes for the right hand side. (d) Always true. These are the  ``do 
nothing'' operators. \pars
\noindent [3] (I have corrected the typo. by inserting the missing ( and ) in 
each formula.) \pars
 (a) $A = (\forall x(P(c) \lor Q(x,x))) \to  (P(c) \lor \forall 
xQ(x,x)).$
We need to determine what the value of $(\forall x(P(c) \lor Q(x,x)))$ is for 
this structure. Thus we need to determine the  ``value'' of the statement 
$P'(a')$ or $Q'(a',a')$ and $P'(a')$ or $Q'(b',b').$ Since $a'\in P'$ then  
${\cal M} \mod (\forall x(P(c) \lor Q(x,x))).$ For the statement $(P(c) \lor \forall 
xQ(x,x)),$ the same fact that $a' \in P'$ yields that ${\cal M}\mod (P(c) \lor \forall 
xQ(x,x)).$  Hence ${\cal M}\mod A.$\pars
 (b) $A = (\forall x(P(c) \lor Q(x,x))) \to  (P(c) \land \forall 
xQ(x,x)).$ We repeat the above for $(\forall x(P(c) \lor Q(x,x)))$ and get 
${\cal M}\mod (\forall x(P(c) \lor Q(x,x))).$ We now check $(P(c) \land \forall 
xQ(x,x)).$ We know that ${\cal M}\mod P(c).$ But $(b',b') \notin Q'.$ Hence
${\cal M} \not\mod (P(c) \land \forall 
xQ(x,x)) \Rightarrow  {\cal M}\not\mod A.$\pars
 (c) $A = (\forall x(P(c) \lor Q(x,x))) \to  (P(c) \land \exists xQ(x,x)).$
Again we know that ${\cal M}\mod (\forall x(P(c) \lor Q(x,x))).$ Also 
${\cal M} \mod P(c).$ Further$,$ $(a',a') \in Q'.$ Thus ${\cal M}\mod \exists xQ(x,x))
\Rightarrow {\cal M}\mod (P(c) \land \exists xQ(x,x))\Rightarrow {\cal M} \mod 
A.$\pars
 (d) $A=  (\forall x(P(c) \land Q(x,x)))\iff (P(c) \land \forall xQ(x,x)).$
First$,$ we know that ${\cal M}\mod P(c).$ But $(b',b') \notin Q'.$ Hence$,$
know that under our interpretation the mathematical statement corresponding to 
$(\forall x(P(c) \land 
Q(x,x)))$ is false. Thus 
${\cal M} \not\mod (\forall x(P(c) \land Q(x,x))).$ We now check $(P(c) \land 
\forall xQ(x,x)).$ Again since $(b',b') \notin Q'$ the mathematical statement 
is false. But then this implies that ${\cal M}\mod A.$\pars 
 (e) $A = (\forall x(P(c) \land Q(c,x)))\iff (P(c) \land \forall xQ(x,x)).$
Since $(a',a'),\ (a',b') \in Q'$ and $a' \in P'$ the mathematical statement 
$(\forall x(P(c) \land Q(c,x)))$ holds for this structure. Thus ${\cal M} \mod 
(\forall x(P(c) \land Q(c,x))).$ But as shown in (d) ${\cal M}\not\mod 
(P(c) \land \forall xQ(x,x)).$ Hence$,$ ${\cal M}\not\mod A.$\par
\medskip
\centerline{\sl Exercise 3.4$,$ Section 3.4}
\medskip
\noindent [1] (a) Note that it was not necessary to discuss special structures
in our definition for validity unless we wanted to find a countermodel. But in 
this case we need to also look at certain special structures. In particular$,$
the case for this problem that $P' = \emptyset.$ But in this case$,$
for any structure $\str \not\mod (\forall x(\exists yP(x,y))).$ Thus for this 
possibility$,$ $\str \mod (\forall x(\exists y P(x,y)))\to (\exists y(\forall x 
P(x,y))).$ \par
Now letting $D =\{a'\},$ then the only other possibility is that 
$P' = \{(a',a')\}.$ Mathematically$,$ it is true that there exists an $a'\in D,$ 
for all  $a' \in D$ we have that  $(a',a') \in P'.$ Hence$,$ in this case$,$
$\str \mod \exists y(\forall x 
P(x,y))).$ Consequently$,$ $\str \mod  (\forall x(\exists y P(x,y)))\to 
(\exists y(\forall x P(x,y))).$\pars
(b) Consider $D = \{a',b'\},\ P' = \{(a',a'),(b',b')\}.$ Then since
$(a',a'),(b',b') \in P',$ it follows that $\str \mod (\forall x(\exists y 
P(x,y))).$ However$,$ since $(b',a'), (b',b') \notin P'$  the mathematical 
statement  ``there exists some $d' \in D$ such that $(a',d'), (b',d')$''
does not hold. [Notice that the difference is that in the first case the 
second coordinate can be any member of $D$ while in the second case it must be 
a fixed member of $D$. Therefore$,$ $\str \not\mod (\forall x(\exists y P(x,y)))\to 
(\exists y(\forall x P(x,y)))$ and the formula is not 2-valid. [Hence$,$ not 
valid in general]\parm
\noindent [2]. In what follows I will make the substitution and see what 
happends in ecah case.
(a) $A = \forall w(P(x) \lor (\forall xP(x,y))\lor P(w,x));\ \lambda = 
y\Rightarrow \forall w(P(\underline{y}) \lor (\forall xP(x,y))\lor P(w,\underline{y}))$. 
Since $y \not= w$ then $y$ is free for $x$ in $A$.\pars                   
\noindent (b) $A = \forall w(P(x) \lor (\forall xP(x,y))\lor P(w,x));\ 
\lambda = w\Rightarrow \forall w(P(x) \lor (\forall xP(x,y))\lor 
P(w,\underline{w}))$. Now the variable has gone from a free occurrence in the 
underlined part to a bound occurrence. Hence $w$ is not free for $x$ in this 
$A.$\pars                   
\noindent (c) $A = (\forall x(P(x) \lor (\forall yP(x,y)))) \lor P(y,x);\ 
\lambda  = x.$ Yes$,$ any variable is always free for itself.\pars
\noindent (d) $A= (\forall x(P(x) \lor (\forall yP(x,y)))) \lor P(y,x);\  
\lambda = y\Rightarrow (\forall x(P(x) \lor (\forall yP(x,y)))) \lor 
P(y,\underline{y})$.  Since the only place that $x$ is free in this formula 
is in the underlined position$,$ the result of substitution still gives a free 
occurrence$,$ this time of $y$. Hence$,$ in this case$,$ $y$ is free for $x$ in 
$A$.\pars
\noindent (e) $A = (\forall x(\exists y P(x,y))\to 
(\exists yP(y,y,));\ \lambda = y.$  Since there are no free occurrences of $x$ 
in this formula$,$ then any variable is free for $x$ in this formula.
\pars
\noindent (f) $A= (\exists zP(x,z))\to (\exists zP(y,z));\ \lambda = 
z\Rightarrow (\exists zP(\underline{z},z))\to (\exists zP(y,z)).$  Since at 
the only position that $x$ was free$,$ the substitution now makes this position 
a bound occurrence$,$ then $z$ is not free for $x$ in $A.$\parm 
\noindent [3] See the above formula where I have made the substitutions in all 
cases.\parm
\noindent [4] 
Assume that $C$ does not have $x$ as a free variable and that $B$ may 
contain $x$ as a free variable. Further$,$ it's assumed that there are no other 
possible free variables. [This comes from$,$ our use of the special process 
(i).] 
$\exists x(C \land B)$ is a sentence. Let $\str$ be an 
arbitrary structure. Assume that $\str\mod \exists x(C \land B).$ 
Then there exists some $d\p \in D,\ \str\mod S^x_d(C \land B)] = C\land
 S^x_dB].$ Hence$,$
$\str\mod C$ and for some $d\p \in D,\ \str\mod S^x_d B].$  Hence$,$  
$\str\mod (C \land (\exists x B)).$ Then in like manner$,$ since $x$ is not 
free in 
$C$$,$ $\str\mod (C \land (\exists xB))\Rightarrow 
\str\mod \forall x(C \land B).$  [Note I have just copied the metaproof of 
(vii) and made appropriate changes.]\parm

\noindent [5] (a) $Q(x) \to (\forall x P(x)).$ Of course$,$ we first take the 
universal closure. This gives the sentence $\forall (Q(x) \to (\forall x 
P(x))).$ I have an intuitive feeling$,$ since $Q$ can be anything$,$ that this in 
invalid. So$,$ we must display a countermodel. To establish that $\str \mod Q(x) \to 
(\forall x P(x)).$ We consider $S^x_d(Q(x) \to (\forall x P(x))).$ Since $x$ 
is not free in $\forall xP(x),$ the valuation of $S^x_d(Q(x) \to (\forall x 
P(x)))$ is the same as the valuation for $(\forall xQ(x)) \to (\forall 
xP(x)).$  Now this makes sense. But this is not even 1-valid. For take 
$D = \{a'\},\  Q'= D,\  P' = \emptyset.$ Then $\str \mod (\forall xQ(x))$ but
$\str \not\mod (\forall xP(x)).$ Hence $\str \not\mod 
\forall (Q(x) \to (\forall x 
P(x))).$ Thus the formula is INVALID.\pars
(b) $(\exists xP(x))\to P(x).$ One again you consider the universal closure 
and we consider the formula $(\exists xP(x))\to (\forall xP(x)).$ Thus also 
does not seem ``logical'' in general. [Note that taking an empty relation will not 
do it.]  But take $D = \{a',b'\}, \ P' = \{a'\}.$ Now it follows that 
$\str \mod (\exists xP(x)).$ Since $b' \notin P',$ then $\str \not\mod 
(\forall xP(x)).$ Thus INVALID. \pars
(c) $(\forall x(P(x) \land Q(x)))\to ((\forall xP(x))\land (\forall 
xQ(x))).$ We don't need to do much work here. Simply consider Theorem 3.4.9 
part (vi). Then let $A= P(x),\ B = Q(x).$ Since that formula is valid$,$ then 
if $\str$ is any structure for (c)$,$ and $\str \mod (\forall x(P(x) \land 
Q(x))),$ then $\str \mod ((\forall xP(x))\land (\forall 
xQ(x))).$ Thus formula is VALID.\pars
(d) $(\exists x(\exists y P(x,y)))\to (\exists xP(x,x)).$ This seems to be   
be invalid since mathematically in the for $\str \mod(\exists x(\exists y 
P(x,y)))$ we do not need $x = y$ in the mathematical sense. 
This is that you can have some $a'$ and some $b' \not= a',$ which 
satisfy a binary relation but $(a', a')$ and $(b',b')$ do not and this is 
exactly how we construct a countermodel. Let $D = \{a',b'\},\  
P'=\{(a',b')\}.$ Then $\str \mod  (\exists x(\exists y P(x,y))),$ but 
$\str \not\mod (\exists xP(x,x)).$ Hence INVALID\pars
(e) $(\exists xQ(x)) \to (\forall xQ(x)).$ This also seems to be invalid.
Well$,$ take $D = \{a',b'\},\ Q' = \{a'\}.$ Then clearly$,$ $\str \mod (\exists 
xQ(x)).$ But since $b'\notin Q',$ then $\str \not\mod (\forall xQ(x)).$\parm

\noindent [6] [This is a very important process.]\pars 
 (a) $(\neg (\exists x P(x))) \lor (\forall xQ(x))\equiv 
(\neg (\exists x P(x))) \lor (\forall yQ(y))\equiv (\forall x (\neg 
P(x)))\lor (\forall yQ(y))\equiv \forall y(Q(y) \lor (\forall x(\neg 
P(x)))\equiv \forall y(\forall x(Q(y) \lor P(x))).$\pars
(b) $((\neg (\exists x P(x)))\lor (\forall xQ(x)))\land (S(c) \to 
(\forall xR(x)))\equiv ((\forall y (\neg P(y)))\lor (\forall xQ(x)))\land (S(c) \to 
(\forall xR(x)))\equiv (\forall y (\forall x(\neg P(y)\lor Q(x))))\land (S(c) 
\to(\forall xR(x)))\equiv (\forall y (\forall x(\neg P(y)\lor Q(x))))\land 
(S(c) \to(\forall zR(z)))\equiv (\forall y (\forall x(\neg P(y)\lor 
Q(x))))\land (S(c)\to(\forall zR(z)))\equiv (\forall y (\forall x(\neg P(y)\lor 
Q(x))))\land (\forall z(S(c) \to R(z)))\equiv \forall y(\forall x(\forall 
z((\neg P(x))\lor Q(x))\land (S(c) \to R(z)))).$\pars
(c) $\neg (((\neg (\exists xP(x)))\lor(\forall xQ(x))) \land 
(\forall xR(x)))\equiv \neg (\forall y(\forall x(\forall 
z((\neg P(x))\lor Q(y))\land R(z))))=  \exists y(\exists x(\exists   
z(P(x)\land (\neg Q(y)))\lor (\neg R(z)))).$\par
\medskip
\centerline{\sl Exercise 3.4$,$ Section 3.5}
\medskip
As suggested, I might try the deduction theorem for valid consequence 
determinations.
\noindent [1] (a) Consider $A_1 = \forall x(Q(x) \to R(x)),\ A_2 = \exists x 
Q(x), B = \exists x R(x).$\pars
Suppose that $\str$ is ANY structure, defined for $A_1,\ A_2,\ B$ such that $\str \not\mod \exists x R(x).$ 
Hence, for each $c'\in D,$ $c' \notin R'.$ Thus, we have that $R'$ is the empty set. Hence, under 
the hypothesis assume that $\str \mod  \forall x(Q(x) \to R(x)).$ This means 
that for each $c'\in D,\ \str\mod Q(c) \to R(c).$ However, we know that 
$\str\not\mod R(c)$ for each $c' \in D.$ Hence for each $c'\in D,$
$\str\not\mod \exists x Q(c).$ Therefore, $\str\not\mod Q(x).$  This implies that it is a VALID CONSEQUENCE.\par
Not using the deduction theorem, proceed as follows: suppose that
$\str $ is a structure defined for $A_1, \ A_2,\ B.$ Let $\str \mod A_1,\ 
\str\mod A_2.$ Then for every $c' \in D$, $\str \mod Q(c) \to R(c)$. Now since $\str \mod A_2$, then  
there is some $c_1' \in D$ such that $\str\mod Q(c_1)$ (i.e. $Q'$ is not empty). 
Since $\str \mod Q(c_1) \to R(c_1),$ then for this $c_1'$ we have that 
$\str \mod R(c_1).$  Hence $\str \mod B.$ Thus it is a VALID CONSEQUENCE.\par                                           
\noindent (b) $A_1 = \forall x(Q(x) \to R(x)),\ A_2 = \exists x (Q(x) \land 
Z(x)) \mod B = \exists x(R(x) \land Z(x)).$
(Without deduction theorem.) Suppose that there is any structure $\str$ 
defined for $A_1,\ A_2,\ B$ and $\str \mod A_1,\ \str \mod A_2.$ Hence, for 
each $c' \in D, \str\mod Q(c) \to R(c)$ and there exists some $c_1'\in D$ such 
that $\str \mod Q(c_1) \land Z(c_1).$ Thus $c_1' \in Q'$ and $c_1' \in Z'.$
Since $c_1' \in Q',$ then $c_1' \in R'.$   Thus using this $c_1'$ we have that
$\str\mod \exists x(R(x) \land Z(x)).$ Hence $\str\mod B$ implies that we have 
a VALID CONSEQUENCE. \parm
\noindent (c) $A_1 =\forall x(P(x) \to (\neg Q(x))),\  A_2 =\exists x (Q(x) 
\land R(x)) \mod B = \exists x ( R(x) \land (\neg Q(x))).$\pars We have a 
feeling that this might be invalid, so we need to construct a countermodel
Consider a one element domain $D = \{a'\}.$ Let $P' = \emptyset,\  Q' = R' = 
D.$ Let $\str = \langle D, P', Q', R' \rangle.$ Then $\str\mod A_1,$ and 
$\str\mod A_2.$ But, since there is no member of $c' \in D$ such that 
$\str \mod R(c)$ and $\str\mod \neg Q(x)$ (i.e. there is no member of $D$ that 
is a member of $D$ and not a member of $D$), it follows that $\str\not\mod B.$ 
Hence, argument if INVALID.\parm
\noindent (d) $A_1 = \forall x (P(x) \to Q(x)),\ A_2 = \exists x(Q(x) \land 
 R(x)) \mod B = \exists x (R(x) \land (\neg Q(x))).$ \pars
Take the same structure as defined in (c). The fact that $\str\mod A_1$ is 
not dependent upon the definition of $Q'.$ Hence, the argument is still 
INVALID.\par\medskip
\centerline{\sl Exercise 3.6$,$ Section 3.6}
\medskip
\centerline{(A) $\forall x(A \to B),\ \forall x(\neg B)\vd \forall x(\neg A).$}\s
\line{(1) $\forall x(A \to B)$\leaderfill Premise}\s
\line{(2) $\forall x(\neg B)$\leaderfill Premise }\s 
\line{(3) $\forall x(A \to B) \to (A \to B)$\leaderfill $P_5$}\s
\line{(4) $A \to B$\leaderfill MP(1$,$3)}\s
\line{(5) $(A \to B) \to ((\neg B)\to (\neg A))$\leaderfill Exer. 2.13$,$ 2A.}\s
\line{(6) $(\neg B) \to (\neg A)$\leaderfill MP(4$,$5)}\s
\line{(7) $(\forall x(\neg B)) \to (\neg B)$\leaderfill $P_5$}\s
\line{(8) $\neg B$\leaderfill MP(2$,$7)}\s
\line{(9) $\neg A$\leaderfill MP(6$,$8)}\s
\line{(10) $\forall x(\neg A)$\leaderfill G(9)}\medskip
\centerline{(B) $\forall x(\forall y A)\vd \forall y(\forall x A)$}\s
\line{(1) $\forall x(\forall yA)$\leaderfill Premise}\s
\line{(2) $(\forall x(\forall y A)) \to \forall y A$\leaderfill $P_5$}\s
\line{(3) $\forall y A$\leaderfill MP(1$,$2)}\s
\line{(4) $(\forall y A) \to A$\leaderfill $P_5$}\s
\line{(5) $A$\leaderfill MP(3$,$4)}\s
\line{(6) $\forall x A$\leaderfill G(5)}\s
\line{(7) $\forall y(\forall xA)$\leaderfill G(6)}\medskip
\centerline{(C) $A,\ (\forall xA) \to C \vd \forall x C$}\s
\line{(1) $A$\leaderfill Premise}\s
\line{(2) $\forall xA$\leaderfill G(1)}\s
\line{(3) $(\forall xA)\to C$\leaderfill Premise}\s
\line{(4) $C$ \leaderfill MP(2$,$3)}\s
\line{(5) $\forall xC$\leaderfill G(4)}\medskip
\centerline{(D) $\forall x(A \to B),\ \forall xA \vd \forall xB$}\s
\line{(1) $\forall x(A \to B)$\leaderfill Premise}\s
\line{(2) $\forall xA$\leaderfill Premise}\s
\line{(3) $(\forall x(A \to B)) \to (A \to B)$\leaderfill $P_5$}\s
\line{(4) $A \to B$\leaderfill MP(1$,$3)}\s
\line{(5) $(\forall xA) \to A$\leaderfill $P_5$}\s
\line{(6) $A$\leaderfill MP(2$,$5)}\s
\line{(7) $B$\leaderfill MP(4$,$6)}\s
\line{(8) $\forall xB$\leaderfill G(7)}\s
\par
\medskip
\centerline{\sl Exercise 3.7$,$ Section 3.7} \par
\medskip
1. Using the special process$,$ it may be assumed that $A$ has no free variables 
and the $B$ has at most one free variable $x$. We show that 
$\mod (\forall x(A \to B)) \to (A \to (\forall xB)).$ \pars
Proof. Let $\str$ be any structure defined for $A,\ B$ and we only need to 
suppose that $\str\not\mod  A \to (\forall x B).$ \pars
Consider the case that $x$ is not 
free in $B$. Then $B$ is a sentence. Now in this case $\str \mod \forall x(A 
\to B)$ iff $\str \mod A\to B.$ Since $\str\not\mod B$ and $\str \mod A$$,$ it 
follows that $\str \not\mod A \to B.$ Hence $\str \mod (A\to B) \to (A \to 
B).$\pars
Now assume that $x$ is free in $B.$ Then we know that for $\str \mod 
\forall x(A \to B)$ then $\str \mod A \to S^x_dB$ for each $c' \in D.$ 
But$,$ we have from our assumption that  $\str \mod A$ and there 
is some $c'\in D$ such that $\str \not\mod S^x_d B.$ Hence $\str \not\mod S^x_dB$
for each $c' \in D.$ Thus $\str\not\mod \forall x(A \to B).$ Therefore$,$
$\str \mod (\forall x(A \to B) ) \to (A \to B).$ \par
\medskip
\centerline{\sl Exercise 3.8$,$ Section 3.8} \par
\medskip
\noindent 1. Modify the argument given in example 3.8.1 as follows: let 
$L(x,y)$ correspond to the natural number binary relation of ``less than''
(i.e.  $<$). Give an argument that shows that there is a structure 
$\hyper{\str}$ that behaves like the natural numbers but in which there 
exists a member $b\p$ that is ``greater than'' any of the original 
natural numbers. \par
As in that example, let $\Gamma$ be the theory of natural numbers, described 
by a given $Pd$, 
and each 
member of $C$ denotes a member of the domain $D$ for a model $\str$ for 
$\Gamma,$ where $\str\mod$ models all of the theory definable predicts as well.  Again we let $b$ be a constant not in the original $C$ and adjoin 
this the $C$. Consider the set of sentences $\Phi = \{L(c,b)\mid c \in C\}.$ 
Now consider the set of sentences $\Gamma \cup \Phi$ and let $A$ be a finite 
subset of $\Gamma \cup \Phi.$ If
$\{a_1,\ldots,a_n\}\subset A$ and $\{a_1,\ldots,a_n\}\subset \Gamma$, then
$\str\mod $ is a model $\{a_1,\ldots,a_n\}$. Suppose that 
$\{a_{n+1},\ldots,a_m\}$ are the remaining members of $A$ that are not in 
$\Gamma.$ Now we investigate the actual members of $\Phi$. We know from the 
theory of natural numbers that for any finite set of natural numbers there is a 
natural number $b'$ greater than any member of that set. Now each member of $ 
\{a_{n+1},\ldots,a_m\}$ is but the sentence $L(c,b)$ where $c \in C$ and 
$b\not\in C.$ Thus there are at most finitely many different $c_i \in C$ 
contained in the formula in $\{a_{n+1},\ldots,a_m\}$. Each of these is 
interpreted as a name for a natural number. Hence let $b$ be interpreted as 
one of the $b'>'$ all of the $c_i'.$ Consequently, $\str \mod A.$ Thus from the compactness 
theorem there is a structure $\hyper \str$ that behaves in $Pd_b$ like the 
natural 
numbers but contains a type of natural number that is ``greater than'' all of 
the original natural numbers. \parm
\noindent 2. Let $\real$ denote the set of all real numbers. Let $C$ be a set 
of constants naming each member of $\real$ and suppose that $Pd$ is the 
language that describes the real numbers. Suppose that $b$ is a constant 
not a member of $C$. Let $\Gamma$ be the theory of real 
numbers. Let $Q(0,y,x)$ be the 3-place predicate that corresponds to the 
definable real number 3-place relation  $0\p <c\p <d\p$, where $0\p,\ c\p,\ 
d\p \in 
\real.$  Now in the real numbers there is a set of elements $G\p$ such that 
each member $c\p$ of $G\p$ has the property that $0\p< c\p$. Let $G$ be the set 
of constants that correspond to the members of the set $G\p.$ Consider the 
set of sentences $\Phi =\{Q(0,b,g)\mid g \in G\}$ in the language $Pd_b.$ 
Give an argument 
that  shows that there exists a structure $\hyper{\str}$ such that 
$\hyper{\str} \mod \Gamma \cup \Phi.$ That is there exists a mathematical 
domain $D$ that behaves like the real numbers, but $D$ contains a member
$b\p$ such that $b\p$ is ``greater than zero'' but $b\p$ is  ``less than''
every one of the original positive real numbers. \par
Well, simply consider any finite subset $A$ of $\Gamma \cup \Phi.$ Suppose 
that $\{a_1,\ldots,a_n\} \subset A$ and that $\{a_1,\ldots,a_n\}\subset 
\Gamma.$ Let $\str \mod \Gamma,$  where $\str\mod$ models all of the theory definable predicts as well. Then $\str \mod \{a_1,\ldots,a_n\}.$ Now 
suppose that $\{a_{n+1},\ldots,a_m\}$ are the remaining members of $A$ that are not in 
$\Gamma.$  Now each member of $\{a_{n+1},\ldots,a_m\}$ is of the form
$Q(0,b,g)$ where $g \in G$ is interpreted as a real number greater than $0'$. 
Since, in this case the $g$ represent a finite set of such real numbers, then 
one of these real numbers, say $g'_1,$ is the ``smallest one'' with respect 
to $<'$. Now consider the real number $g_1'/2$ and let $b$ be interpreted as 
this real number. Then $\str\mod \{a_{n+1},\ldots,a_m\}$. Hence we have that 
$\str \mod A.$ The compactness theorem states that there is a structure
$\hyper \str$ that behaves, in $Pd_b,$ like the real numbers. And, there 
exists in the structure a real like number $b'$ such that $0'<'b'<'g'$ for all 
of the original real numbers $g'$ such that $0' <' g'.$ \medskip\medskip
\centerline{\bf Discussion}
\centerline{In what follows prime notation has been charged to *}\parm

[1] (Note that we don't mention the ``constant'' in the structure notation since we have originally assigned all members of the domain and some members of the new domain constant names.) In the structure $\hyper \str =\ <\hypernat, \hyper +,\hyper =, \hyper <,\ldots >$ all the usual properties that can be expressed in a the appropriate first-order language hold for $\hyper \str.$ One extremely useful additional property that one would like $\hyper \str $ to possess is the embedding property. Is there a subset of $\hypernat$ that can be used in every way as the natural numbers themselves where the * operators restricted to this subset have the same properties as the original natural numbers? After some effort in model theory$,$ the answer is yes. Thus we can think of the natural numbers $\nat$ as a ``substructure'' of the $\hyper \str.$ The interesting part of all of this is that a simple comparison of properties can now be made. If there is one natural number $\Hyper b$ then you can consider the set of all such members of $\hypernat$.  This set is denoted by $\nat_\infty$ and is called the set of all ``infinite natural numbers.'' It has a algebra that behaves as Newton wished for such objects. For example$,$ if $\lambda,\ \beta \in \nat_\infty$ and nonzero $n \in \nat$$,$ then $n\lambda + \beta \in \nat_\infty.$ Further$,$ for any $n \in \nat$$,$ we have that $\lambda - n\in \nat_\infty,$ [can you show this?] where we now think of $\nat$ as the non-negative integers. But $\hypernat$ has a property that no set of natural numbers has and this property is why the set $\hypernat$ cannot be ``graphed'' in the usual manner. Every nonempty subset of $\nat$ has a first element. This means that if nonempty $A\subset \nat$$,$ then there is some
$a \in A$ such that $a\leq x$  for all $x \in A.$ But the set $\nat_\infty$ does not have a first element with respect to $\hyper <.$ Suppose that $\nat_\infty$ does have a first element $\lambda_1$. Then if you establish what has been written above$,$ since the same basic properties for $\hyper <$ hold as they do for $<$$,$ we have that $\lambda_1 - 1 \hyper < \lambda_1,$ which contradicts the concept of first element with respect to $\hyper <.$ Indeed$,$ a recent published paper that I am reviewing$,$ has forgotten this simple fact. 
\parm
[2] The structure here $\hyper \str =\  <\hyperreal,\hyper +,\hyper \cdot,\hyper <, \cdots>$ can have 
$\str =\ <\real,+,\cdot,<, \dots>$ embedded into it in such a manner that 
$\str$ is a substructure of $\hyper \str.$ And$,$ as in [1]$,$ this substructure behaves relative to the all the properties of the first-order theory of real numbers just like the real numbers. In this case$,$ we have that $\Hyper 0=0 \hyper < \Hyper b \hyper <\ r$ for each positive real number $r \in \real.$ Now we consider the entire set of all such $\hyper b$ and call this the set of positive infinitesimals
$\monad {+0}$. Since all the algebra holds one adjoins to this set the 0 and 
$\{-\eps \mid \eps \in \monad {+0}\} = \monad 0.$ Do these satisfy the theory of the ``infinitely small'' of Newton? Well$,$ here are a few properties. First$,$ the normal arithmetic of the real numbers holds for $\monad {0}$. This we have that if $\eps \in \monad {+0}$$,$ then $1/\eps  
\Hyper > \ r$ for any $r \in \real.$ Also  when Newton had an object he called infinitely small and squared it he claimed that this was ``more infinitely small'' in character. Well$,$ $0 \hyper < \ \eps^2 \hyper < \  \eps $ and they are not equal. Then Newton claimed that if he took any real number and multiplied it by an infinitely small number that the result was still infinitely small. One can establish that for each
$r \in \real,$ we have that $r\hyper \cdot \ \eps \in \monad 0.$ Indeed$,$ we can establish  that if $f$ is any function defined on $\real,$ continuous at $r = 0$ and $f(0)= 0,$ then all the properties of $f$ hold for $\eps$ and when they are applied to $\eps$ the result is always an infinitesimal. There is no where in the use of the Calculus that these infinitesimals contradict the intuitive procedures used by mathematicians through 1824. Their properties also correct the error discovered in 1824 that led to the introduction of the limit concept. \par}\vfil\eject
\newdimen\fullhsize
\fullhsize=6.5in \hsize=3.0in 
\def\fullline{\hbox to \fullhsize}

{\output={\if L\lr
         \global\setbox\leftcolumn=\columnbox \global\let\lr=R
\else \doubleformat \global\let\lr=L\fi
\ifnum\outputpenalty>-20000 \else\dosupereject\fi}
\def\doubleformat{\shipout\vbox{\makeheadline
\fullline{\box\leftcolumn\hfil\columnbox}
\makefootline}
\advancepageno}
\def\columnbox{\leftline{\pagebody}} 
\let\lr=L \newbox\leftcolumn
\def\sub#1{{\leftskip=0.2in \noindent #1 \par}\par}
\def\n{\noindent}
\indent\indent{\bf A}\par
{\quad}\par
\n and-gate 29. \par
\n assignment 15.\par
\n atomic formula$,$\par
\sub{propositional 13.}\par
\sub{predicate 64.}\par
\sub{first-order language 64.}\par
\n atoms$,$\par 
\sub{propositional 9, 10.}\par
\sub{predicate 63.}\par
\sub{first-order language 63.}\par
{\quad}\par
\indent\indent{\bf B}\par
{\quad}\par
\n binary 9$,$\par
\sub{relation 22.}\par
\n bound occurrence of the constant\par
 68.\par
{\quad}\par
\indent\indent{\bf C}\par
{\quad}\par
\n closed formula 69.\par
\n common pair 15, 13$,$\par
\sub{rule 13.}\par
\n complete set 59. \par
\n composite formula 64.\par
\n congruent formula 69.\par
\n connectives 9.\par
\n consequence$,$\par
\sub{of $\Gamma $ 57.}\par
\sub{consequence operator 55.}\par
\n consistent$,$\par
\sub{propositional 39.}\par
\sub{predicate$,$ first-order 91.}\par
\n constants$,$ insertion of 63, 57.\par
\n contradiction$,$\par
\sub{propositional 39.}\par
\sub{predicate$,$ first-order 91.}\par
\n countermodel 76.\par
{\quad}\par 
{\quad}\par 

\indent\indent{\bf D}\par
{\quad}\par
\n deduction$,$ \par
\sub{from a set of $\Gamma $ 45.}\par
\n deductive processes 98.\par
\n deductive system 59.\par
\n demonstration 45.\par
\n denial 26.\par
\n descriptions 98.\par
\n domain for a structure 72.\par
\n equality 22.\par
{\quad}\par
\indent\indent{\bf E}\par
{\quad}\par
\n equivalence relation 23.\par
\n existential quantifier 64.\par
\n extralogical symbols 10.\par
{\quad}\par
\indent\indent{\bf F}\par
{\quad}\par
\n formal demonstration$,$\par
\indent first-order$,$ 86.\par
\n formal proof$,$ propositional 42$,$\par
\indent 43.\par
\n formally consistent$,$\par
\sub{predicate 91.}\par
\sub{propositional 57.}\par
\n formula variables$,$\par
\sub{predicate 64.}\par
\sub{propositional 39.}\par
\n formulas$,$ \par
\sub{predicate 63.}\par
\sub{propositional 9.}\par
\n free for x in A 77. \par
\n free in the formula 68.\par
\n free substitution operator 72.\par
\n full disjunctive normal form 27.\par
\n fundamental conjunction 28.\par
{\quad}\par
\indent\indent{\bf G}\par
{\quad}\par
\n generalization$,$ first-order 86.\par
{\quad}\par
\indent\indent{\bf H}\par
{\quad}\par
\n hypotheses$,$ premises 32.\par
\n hypothetical syllogism 44.\par
{\quad}\par
\indent\indent{\bf I}\par
{\quad}\par
\n induction 102.\par
\n inconsistent$,$\par
\sub{predicate$,$ first-order 91.}\par
\sub{propositional 39.}\par
\n interpretation metasymbol 10.\par
\n interpretation$,$ for first-order\par
\indent structure 72.\par
\n inverter 29.\par
{\quad}\par
\indent\indent{\bf L}\par
{\quad}\par
\n language levels  9.\par
\n logic circuits 29.\par
\n logical connectives 10.\par
\n logical flow 29.\par
{\quad}\par\vfil\eject
\indent\indent{\bf M}\par
{\quad}\par
\n meta - 6.\par
\n model theory$,$ propositional 35.\par
\n modus ponens 43.\par
\n model$,$ first-order 73.\par
{\quad}\par
\indent\indent{\bf N}\par
{\quad}\par
\n names domain members$,$ first-order 94.\par
\n Natural event 98.\par
\n Natural system 97.\par
\n negation complete 91.\par
\n normally closed$,$ switch 14.\par
\n normally open$,$ switch 14.\par
\n not a model 73.\par
\n not satisfied 33.\par
{\quad}\par
\indent\indent{\bf O}\par
{\quad}\par
\n observer language 6.\par
\n or-gate 29.\par
{\quad}\par
\indent\indent{\bf P}\par
{\quad}\par
\n predicates 63$,$ \par
\sub{1-place predicates 63.} \par
\sub{2-place predicates 63.}\par
\sub{3-place predicates 63.}\par
\sub{n-place 63.}\par 
\n premises$,$ hypotheses 32.\par 
\n prenex normal form 82.\par 
{\quad}\par
\indent\indent{\bf R}\par
{\quad}\par
\n relative consistency 85.\par 
{\quad}\par
\indent\indent{\bf S}\par
{\quad}\par
\n satisfiable$,$\par
\sub{propositional 33$,$ 38$,$ 57}\par
\sub{predicate$,$ first-order 93.}\par
\n satisfy = satisfiable\par
\n scope$,$ 64\par
\sub{of that quantifier 67.}\par
\n semantical modus ponens 24.\par
\n semantical 19.\par
\n semantics 15.     \par
\n sentence 69.\par
\n set of all sentences 71.\par
\n simply consistent$,$ propositional\par
\indent 39.\par
\n size of a formula 11$,$ 64.\par
\n standard theory 100.\par
\n standard model 100.\par
\n stronger than for consequence\par
\indent operators 55.\par
\n structure$,$ first-order 72.\par
\n subformula 24.\par
\n substitution process$,$\par
\indent  propositional 20.\par
{\quad}\par
\indent\indent{\bf T}\par
{\quad}\par
\n theorem$,$ formal,\par
\sub{predicate$,$ first-order 86.}\par
\sub{propositional 42.}\par
\n truth-table 15 - 17.\par
{\quad}\par
\indent\indent{\bf U}\par
{\quad}\par
\n ultralogic 100.\par
\n ultraword 100.\par
\n unary connective 10.\par
\n unique equivalent form 27.\par
\n universal closure 76.\par
\n universal quantifier 64.\par
\n universal$,$ formula called 106.\par
{\quad}\par
\indent\indent{\bf V}\par
{\quad}\par
\n valid consequence of $\Gamma $$,$\par
\indent propositional 57.\par
\n valid consequence$,$ \par
\sub{predicate$,$ first-order 82.}\par
\sub{propositional 33, 57.}\par
\n valid$,$ \par
\sub{predicate$,$ first-order 82.}\par
\sub{propositional 19.}\par
\n valuation procedure$,$ 16.\par
\sub{in general for propositional 57.}\par
\n variable predicate forms 63$,$ 77\par
\n variable substitution$,$\par
\indent \indent propositional  40.\par
\n variables$,$ predicate 63 - 64.\par
{\quad}\par
\indent\indent{\bf W}\par
{\quad}\par
\n wffs = well formed formula 9.\par
\n well-ordered 102.\par
\vfil
\eject

\end